5/05/2015

# PSL(2; ℂ) connections on 3-manifolds with $L^2$ bounds on curvature [†]

Clifford Henry Taubes[††]

Department of Mathematics
Harvard University
Cambridge MA 02138

chtaubes@math.harvard.edu

Karen Uhlenbeck's compactness theorem for sequences of connections with $L^2$ bounds on curvature applies only to connections on principal bundles with compact structure group. This article gives an extension of Uhlenbeck's theorem that describes sequences of connections on principal PSL(2; ℂ) bundles over compact three dimensional manifolds.

[†]This version of the paper is a substantial revision of the published version in [T1]. This revision fixes some ultimately inconsequential errors and one serious error found by Thomas Walpuski.

[††] Supported in part by the National Science Foundation

Suppose that M is a compact Riemannian manifold and P is a principal bundle over M with fiber a Lie group to be denoted by G. Fix an integer no less than $\frac{1}{2}$ of the dimension of M, this denoted by p. In the case when G is compact, Karen Uhlenbeck's foundational paper *Connections with $L^p$ bounds on curvature* [U] explained the sense in which the space of connections on the given principal bundle with an a priori $L^p$ bound on the norm of the curvature is compact. This paper provides a generalization of Uhlenbeck's theorem in the case when M has dimension 2 or 3, the group G is $PSL(2;\mathbb{C})$ and p = 2. The generalization of Uhlenbeck's theorem for the case of dimension 3 is described by upcoming Theorems 1.1a and 1.1b. The dimension 2 case is subsumed by Theorems 1.1a,b by taking the manifold in Theorems 1.1a,b to be the product of the surface with the circle. The dimension 2 case is also stated separately as Theorem 1.2.

The space of automorphism equivalence classes of irreducible, flat $PSL(2;\mathbb{C})$ connections need not be compact. Morgan and Shalen [MS1]-[MS3] construct a compactification of the latter space using certain equivariant maps from the universal cover of M to certain sorts of $\mathbb{R}$-trees, an $\mathbb{R}$-tree being a metric space where by any two points are connected by a unique path. Daskoloupoulus, Dostoglu and Wentworth [DDW1-3] subsequently proved that the Morgan-Shalen maps can be taken to be harmonic. Some brief remarks are made near the end of Section 1 about what is said in Theorem 1.1a and what is said in [MS1]-[MS4] and [DDW1-3].

A generalization of Uhlenbeck's theorem to the case when G is $PSL(2;\mathbb{C})$ and M has dimension 3 is of specific, topical interest for reasons that are described at the end of Section 1 of this paper. Moreover, close kin to the techniques that are introduced to prove Theorem 1.1a,b will likely play crucial roles in proofs of Theorems 1.1a,b's analogs when G has rank greater than 2, when M has dimension greater than 3, and/or when p is not equal to 2.

**Acknowledgements**: A debt of gratitude is owed to Qing Han for explaining his work and those of his colleagues, Robert Hardt, and Fanghua Lin. The influence of their work and the earlier work of Almgren and others on the singular points of nodal curves can be seen in much of this article. A large debt of gratitude is also owed to John Morgan for his saintly patience with the author's questions about $\mathbb{R}$-trees and about his papers with Shalen that define their compactification of the space of automorphism classes of flat $PSL(2;\mathbb{C})$ connections. A substantial debt is owed to Ben Mares for his close reading of parts of an earlier version of this manuscript, and a similar debt is owed to Thomas Walpuski. Thomas Walpuski is also thanked for finding a serious error in the published version of this paper.

## 1. The PSL(2; $\mathbb{C}$) extension of Uhlenbeck's theorem

This section has eight subsections. Section 1a first presents Uhlenbeck's theorem and then states the latter's PSL(2;$\mathbb{C}$) generalization, Theorems 1.1a and Theorem 1.1b. These two theorems together revise and correct the assertions of Theorem 1.1 in the published version of this article [T1]. Most of the proof of Theorem 1.1b is in a companion to this article [T2].

Section 1b states and proves the dimension 2 analog of Theorems 1.1a and 1.1b. Section 1c first explains the relationship between what is said in Theorems 1.1a,b and what is said in [MS1]-[MS4] and [DDW1] about the flat connections. It then gives a very brief account how the data supplied by Theorems 1.1a,b leads to some central notions in 3-manifold topology. Section 1d provides a short outline of the proof of Theorems 1.1a and Section 1e makes the point that Theorem 1.1b follows as a special case of theorems in [T2]. Section 1f describes the salient differences between this article and the original. Section 1g describes some questions in 3 and 4 dimensional differential topology and geometry that may well require some sort of PSL(2;$\mathbb{C}$) extension of Uhlenbeck's theorem. This section also has a paragraph that says something about extensions of Uhlenbeck's theorem to PSL(n;$\mathbb{C}$) for $n > 2$.

Section 1h supplies a table of contents for the remaining sections of this article. Section 1h also states certain notational conventions that are subsequently invoked with no further comments.

### a) Theorems 1.1a and 1.1b

Theorem 1.1a is stated below after some necessary stage setting to define the notation and supply some needed background. The stage setting has eight parts.

*Part 1*: The group SL(2;$\mathbb{C}$) is viewed here as the group of $2 \times 2$ complex matrices with determinant 1. Viewed in this light, its Lie algebra is the vector space of trace zero, $2 \times 2$ complex matrices. The latter space can be written as the direct sum $\mathfrak{su}(2) \oplus i\,\mathfrak{su}(2)$ with $\mathfrak{su}(2)$ denoting the vector space of skew hermitian complex matrices and with i denoting the square root of -1. Keep in mind that the Lie algebra of SU(2) is the vector space $\mathfrak{su}(2)$. The linear form on the vector space of $2 \times 2$ complex matrices given by -1 times the trace is denoted in what follows by $\langle\ ,\ \rangle$. The square of the Hermitian norm on $\mathfrak{su}(2)$ can be written using this notation as the function $\mathfrak{u} \to \langle \mathfrak{u}\,\mathfrak{u} \rangle$.

*Part 2*: Assume henceforth that M is a compact 3-dimensional manifold with a given principal PSL(2;$\mathbb{C}$) bundle. As PSL(2;$\mathbb{C}$) bundles have reductions to principal SO(3) bundles, choose once and for all such a reduction so as to write the given principal PSL(2;$\mathbb{C}$) bundle as $P \times_{SO(3)} PSL(2;\mathbb{C})$ with $P \to M$ a principal SO(3) bundle.

Any given connection on $P \times_{SO(3)} PSL(2;\mathbb{C})$ can be written as $A + i\,\mathfrak{a}$ with A being a connection on P and with $\mathfrak{a}$ denoting a 1-form with values in the associated vector bundle $P\times_{SO(3)} \mathfrak{su}(2)$. The curvature of the connection $\mathbb{A} = A + i\,\mathfrak{a}$ is denoted by $F_{\mathbb{A}}$ and that of A by $F_A$. The former is a section of $(P\times_{SO(3)} \mathfrak{su}(2)) \otimes (\wedge^2 T^*M)$ and it can be written in terms of A's curvature as $F_{\mathbb{A}} = F_A - \mathfrak{a} \wedge \mathfrak{a} + i\, d_A\mathfrak{a}$ with $d_A$ denoting here the exterior covariant derivative that is defined by A.

*Part 3*: Fix once and for all a Riemannian metric on M. Unless instructed to the contrary, assume that all inner products on TM, T*M and their tensor products are defined by this metric. The metric Hodge star is denoted by $*$. Likewise, assume that all covariant derivatives on these tensor bundles are those induced by the associated Levi-Civita connection. These covariant derivatives are denoted by $\nabla$. The metric's Ricci curvature tensor defines a symmetric, bilinear form on T*M that is denoted by Ric. Integration on M is defined using the metric's volume element.

Norms of tensor bundle valued sections of $P\times_{SO(3)} \mathfrak{su}(2)$ are defined using the Riemannian metric and the norm on $\mathfrak{su}(2)$. Let A denote a given connection on P. The covariant derivative defined by A and the Levi-Civita connection on $P\times_{SO(3)} \mathfrak{su}(2)$ valued sections of tensor bundles is denoted by $\nabla_A$. The Hermitian adjoint of this operator is denoted by $\nabla_A^{\dagger}$.

*Part 4*: Sobolev spaces of connections are front and center in Uhlenbeck's theorem, and so they are front and center in Theorem 1.1a. This part of the subsection defines these spaces.

Fix a 'fiducial' connection on P to be denoted by $A_0$. The connection $A_0$ is used to define a given $k \in \{0, 1, 2, \ldots\}$ version of the Sobolev $L^2_k$ norm on tensors with values in $P\times_{SO(3)} \mathfrak{su}(2)$, this being the norm whose square assigns to a given tensor $\mathfrak{t}$ the integral over M of $\sum_{0\le m\le k} |(\nabla_{A_0})^{\otimes m}\mathfrak{t}|^2$. An *$L^2_k$ tensor* with values in $P\times_{SO(3)} \mathfrak{su}(2)$ is an almost everywhere defined section of the relevant vector bundle with finite Sobolev $L^2_k$ norm.

Let Conn(P) denote the space of smooth connections on P. Let k denote for the moment a given non-negative integer. The $L^2_k$ topology on Conn(P) is defined as follows: Write a given connection on P as $A_0 + \hat{a}$ with $\hat{a}$ being a section of $(P\times_{SO(3)} \mathfrak{su}(2)) \otimes T^*M$. Doing so identifies Conn(P) with the vector space of sections of $(P\times_{SO(3)} \mathfrak{su}(2)) \otimes T^*M$. The $L^2_k$ topology on Conn(P) is the topology that is induced by this identification from the $L^2_k$ metric metric topology on the space of sections of $(P\times_{SO(3)} \mathfrak{su}(2)) \otimes T^*M$. This topology does not depend on the chosen fiducial connection.

A connection on P is said to be an *$L^2_k$ connection* if it has the form $A_0 + \hat{a}$ with $\hat{a}$ denoting an $L^2_k$ section of $(P\times_{SO(3)} \mathfrak{su}(2)) \otimes T^*M$. This notion is likewise independent of the choice of $A_0$. A sequence $\{A_n\}_{n=1,2,..}$ of connections on P is said to converge *weakly* in the $L^2_k$ topology on P when it can be written as $\{A_n = A_0 + \hat{a}_n\}_{n=1,2,\ldots}$ with the sequence $\{\hat{a}_n\}_{n=1,2,..}$ having bounded $L^2_k$ norm and converging weakly with respect to the $L^2_k$ norm to an $L^2_k$ section of $(P\times_{SO(3)}\mathfrak{su}(2))\otimes T^*M$.

*Part 5*: Karen Uhlenbeck's [U] theorem applies to connections on P and in particular makes the following assertion:

**Uhlenbeck's Theorem**: *Suppose that* $\{A_n\}_{n=1,2,..}$ *is a sequence of connections on* P *with the corresponding sequence*

$$\{ \int_M |F_{A_n}|^2 \}_{n=1,2,...}$$

*being bounded. There is a subsequence of* $\{A_n\}_{n=1,2,..}$ *, hence renumbered consecutively from 1, and a corresponding sequence of automorphisms of* P*, this denoted by* $\{g_n\}_{n=1,2,..}$*, such that* $\{g_n{}^*A_n\}_{n=1,2,..}$ *converges weakly in the* $L^2_1$ *topology to an* $L^2_1$ *connection on* P.

By way of a reminder, the space of automorphisms of P acts on Conn(P) by pull-back. This space of automorphism can be identified in a canonical fashion with the space of sections of $P\times_{SO(3)}SO(3)$. Let g denote such an automorphism and let A denote a given connection. The pull-back g*A can be written as $A + \hat{g}^{-1}\, d_A\hat{g}$ where $\hat{g}$ is a (locally defined) section of the associated bundle to P with fiber the group of $2\times 2$ unitary matrices with determinant 1, this being SU(2). The group SO(3) acts on SU(2) via conjugation. The automorphism g need not lift over the whole of M to a section of the fiber bundle $P\times_{Ad(SO(3))}SU(2)$, but a lift does exists over any contractible subset of M. Even so, the section $\hat{g}^{-1}\, d_A\hat{g}$ of $(P\times_{SO(3)}\mathfrak{su}(2))\otimes T^*M$ is defined everywhere on M because any two lifts of g differ by the action of multiplication by 1 or -1.

*Part 6*: Uhlenbeck's theorem is a godsend because its assumptions are invariant under the action on Conn(P) of P's automorphism group; the reason being that the norm given by $\mathfrak{u} \to -\langle \mathfrak{u}\,\mathfrak{u}\rangle$ on $\mathfrak{su}(2)$ is ad-invariant. This implies that the norm of the curvature of a connection on P is pointwise identical to the norm of its pull-back via any automorphism of P. This being the case, the $L^2$ norm of the curvature of any given connection is is the same as that of its pull-back via an automorphism.

The Lie algebra of $SL(2;\mathbb{C})$ does not have a norm that is invariant under the adjoint action of $SL(2;\mathbb{C})$. Even so, there is a useful generalization of the $L^2$ norm of the curvature of an $SL(2;\mathbb{C})$ connection that is invariant under the action of the group of $PSL(2;\mathbb{C})$ automorphisms of $P\times_{SO(3)}PSL(2;\mathbb{C})$. The definition is given momentarily in (1.1). This generalization plays the role in Theorem 1.1a that is played by the curvature $L^2$ norm in Uhlenbeck's theorem.

The upcoming definition uses $\mathrm{Conn}(P\times_{SO(3)} PSL(2;\mathbb{C}))$ to denote the space of connections on $P\times_{SO(3)} PSL(2;\mathbb{C})$ The group of automorphisms of $P\times_{SO(3)} PSL(2;\mathbb{C})$ is denoted by $\mathcal{G}_{\mathbb{C}}$; and the $\mathcal{G}_{\mathbb{C}}$-orbit in $\mathrm{Conn}(P\times_{SO(3)} PSL(2;\mathbb{C}))$ of a given connection $\mathbb{A}$ is denoted by $\mathcal{G}_{\mathbb{C}}(\mathbb{A})$.

The generalization to $\mathrm{Conn}(P\times_{SO(3)} PSL(2;\mathbb{C}))$ of the square of the $L^2$ norm of the curvature is the function on $\mathrm{Conn}(P\times_{SO(3)} PSL(2;\mathbb{C}))$ that assigns to any given connection $\mathbb{A}$ the infimum over connections $A+i\mathfrak{a}\in\mathcal{G}_{\mathbb{C}}(A)$ of a function that is denoted by $\mathfrak{F}$ and defined by the rule

$$\mathbb{A}\to\mathfrak{F}(\mathbb{A}) = \int_M (\,|F_A - \mathfrak{a}\wedge\mathfrak{a}\,|^2 \;+\; |d_A\mathfrak{a}|^2 \;+\; |d_A{*}\mathfrak{a}|^2\,)\,. \tag{1.1}$$

*Part 7*: Fix $k\in\{0, 1, 2, \ldots\}$. The $L^2_k$ topology on $\mathrm{Conn}(P\times_{SO(3)} PSL(2;\mathbb{C}))$ is defined by first identifying the latter space with $\mathrm{Conn}(P)\times C^\infty(M;(P\times_{SO(3)}\mathfrak{su}(2))\otimes T^*M)$ with it understood that a given pair $(A, \mathfrak{a})$ in the latter space corresponds to the connection $A+i\,\mathfrak{a}$. Having done so, the topology is defined to be the product of the $L^2_k$ topologies on $\mathrm{Conn}(P)$ and $C^\infty(M;(P\times_{SO(3)}\mathfrak{su}(2))\otimes T^*M)$. An $L^2_k$ connection on the principal $PSL(2;\mathbb{C})$ bundle $P\times_{SO(3)} PSL(2;\mathbb{C})$ is defined by a pair $(A, \mathfrak{a})$ of $L^2_k$ connection on $P$ and $L^2_k$ section of $(P\times_{SO(3)}\mathfrak{su}(2))\otimes T^*M$.

A part of Theorem 1.1a describes a connection on an open set in M as being an $L^2_{k;loc}$ connection. Let U denote the given open set. A connection on $P|_U\times_{SO(3)} PSL(2;\mathbb{C})$ is said to be an $L^2_{k;loc}$ connection if it can be written as $A_0 + \hat{a} + i\,\mathfrak{a}$ with $\hat{a}$ and $\mathfrak{a}$ being almost everywhere defined sections over U of $(P\times_{SO(3)}\mathfrak{su}(2))\otimes T^*M$ with the following property: If $V\subset U$ is any given open set with compact closure, then the $L^2_k$ norm on V of $(\hat{a}, \mathfrak{a})$ is finite, the latter being the square root of the integral on V of $\sum_{0\le m\le k}|(\nabla_{A_0})^{\otimes m}\hat{a}|^2$ and $\sum_{0\le m\le k}|(\nabla_{A_0})^{\otimes m}\mathfrak{a}|^2$. This notion of an $L^2_{k;loc}$ connection does not depend on the fiducial connection $A_0$.

A sequence of connections on $P\times_{SO(3)} PSL(2;\mathbb{C})$ is said to converge weakly in the $L^2_k$ topology to an $L^2_k$ connection when the $\mathrm{Conn}(P)$ part of the corresponding sequence in $\mathrm{Conn}(P)\times C^\infty(M;(P\times_{SO(3)}\mathfrak{su}(2))\otimes T^*M)$ converges weakly in the $L^2_k$ topology to an $L^2_k$ connection on the principal SO(3) bundle P and the $C^\infty(M;(P\times_{SO(3)}\mathfrak{su}(2))\otimes T^*M)$ part has bounded $L^2_k$ norm and converges weakly to an $L^2_k$ section of $(P\times_{SO(3)}\mathfrak{su}(2))\otimes T^*M$. If $U\subset M$ is a given open set, a sequence of connections on $P\times_{SO(3)} PSL(2;\mathbb{C})$ is said to converge weakly to an $L^2_{k;loc}$ connection on U if both the $\mathrm{Conn}(P)$ part and the $C^\infty(M;(P\times_{SO(3)}\mathfrak{su}(2))\otimes T^*M)$ part of the connection converge weakly with respect to the $L^2_k$ topology on all open subsets in U with compact closure. This is to say that the Conn(P) part differ from $A_0$ on U by a sequence of sections over U of the vector bundle $(P \times_{SO(3)} \mathfrak{su}(2)) \otimes T^*M$ with bounded $L^2_k$ norm on open sets in U with compact closure; and that it converges weakly on such open sets to an $L^2_{k;loc}$ section of this vector bundle. Meanwhile, the $C^\infty(M; (P \times_{SO(3)} \mathfrak{su}(2)) \otimes T^*M)$ part of the sequence has bounded $L^2_k$ norm on open sets in U with compact closure and it also converges weakly on such open sets to an $L^2_{k;loc}$ section of $(P|_U \times_{SO(3)} \mathfrak{su}(2)) \otimes T^*M$.

*Part 8*: The term *real line bundle* is used here to describe the associated line bundle to a principal $\mathbb{Z}/2\mathbb{Z}$ bundle. With this term understood, let $U \subset M$ denote a given open set and suppose that $\mathcal{I} \to U$ is a real line bundle. An $\mathcal{I}$ valued tensor field on U is a section of the tensor product bundle with $\mathcal{I}$. The Riemannian metric defines the pointwise norm of an $\mathcal{I}$ valued p-form or any $\mathcal{I}$ valued tensor field. Meanwhile, the Levi-Civita covariant derivative defines a covariant derivative of such a tensor field.

Fix $p \in \{0, 1, 2, 3\}$and let $q$ denote for the moment an $\mathcal{I}$ valued p-form. The definition of the exterior derivative of $q$ is canonical and is denoted by $dq$. This $\mathcal{I}$ valued p-form is said to be closed if $dq = 0$ and it is said to be coclosed if $d{*}q = 0$. An $\mathcal{I}$ valued p-form that is closed and coclosed is said to be harmonic.

If A is a connection on P, then A defines a covariant derivative on sections of the bundle $(\mathcal{I} \otimes P) \times_{\mathbb{Z}/2\mathbb{Z} \times SO(3)} \mathfrak{su}(2)$. A given section is said to be A-covariantly constant if it is annihilated by A's covariant derivative.

Keep in mind that if $q$ is an $\mathcal{I}$-valued p-form and $\sigma$ is a section of the vector bundle $(\mathcal{I} \otimes P) \times_{\mathbb{Z}/2\mathbb{Z} \times SO(3)} \mathfrak{su}(2)$, then $q\sigma$ is a p-form on U with values in $P \times_{SO(3)} \mathfrak{su}(2)$.

**Theorem 1.1a**: *Suppose that* $\{\mathbb{A}_n = A_n + i\,\mathfrak{a}_n\}_{n=1,2,...}$ *is a sequence of connections on* $P \times_{SO(3)} PSL(2;\mathbb{C})$ *with the corresponding sequence* $\{\mathfrak{F}(\mathbb{A}_n)\}_{n=1,2,...}$ *being bounded. For each* $n \in \{1, 2, \ldots\}$, *use* $r_n$ *to denote the* $L^2$ *norm of* $\mathfrak{a}_n$.

- *If the sequence* $\{r_n\}_{i=1,2,...}$ *has a bounded subsequence, then there exists a subsequence of* $\{\mathbb{A}_n\}_{n=1,2,...}$, *hence renumbered consecutively from 1, and a corresponding sequence of automorphisms of* P, *this denoted by* $\{g_n\}_{n=1,2,...}$, *such that* $\{g_n{}^*\mathbb{A}_n\}_{n=1,2,...}$ *converges weakly in the* $L^2_1$ *topology to an* $L^2_1$ *connection on* $P \times_{SO(3)} PSL(2;\mathbb{C})$.
- *If the sequence* $\{r_n\}_{n=1,2,...}$ *has no bounded subsequence, then there exists a subsequence of* $\{\mathbb{A}_n\}_{n=1,2,...}$, *hence renumbered consecutively from* 1, *a corresponding sequence of automorphisms of* P, *this denoted by* $\{g_n\}_{n=1,2,...}$, *plus the following extra data: A closed set, nowhere dense set* $Z \subset M$, *a real line bundle* $\mathcal{I} \to M{-}Z$, *and a harmonic* $\mathcal{I}$ *valued 1-form on* $M{-}Z$. *The latter is denoted by* $\nu$. *These are such that*

1) *The norm* $|\nu|$ *of* $\nu$ *extends to the whole of* M *as a Hölder continuous,* $L^2_1$ *function with its zero locus being the set* Z.
2) *The sequence* $\{g_n{}^*A_n\}_{n=1,2,\ldots}$ *converges weakly in the* $L^2_{1;loc}$ *topology on* M–Z *to an* $L^2_{1;loc}$ *connection on* $P|_{M-Z}$, *this denoted by* A.
3) *The sequence* $\{r_n^{-1}\, g_n{}^*\mathfrak{a}_n\}_{n=1,2,\ldots}$ *converges weakly in the* $L^2_{1;loc}$ *topology on* M–Z *to* $\nu\,\sigma$ *with* $\sigma$ *being a unit length,* A-*covariantly constant homorphism over* M–Z *from* $\mathcal{I}$ *into* $P \times_{SO(3)} \mathfrak{su}(2)$. *Meanwhile,* $\{r_n^{-1}|\mathfrak{a}_n|\}_{n=1,2,\ldots}$ *converges to* $|\nu|$ *in the weak* $L^2_1$ *topology and the* $C^0$ *topology on the whole of* M.

The next theorem says more about the set Z from Theorem 1.1a. To set the notation for this upcoming Theorem 1.1b, suppose that (Z, $\mathcal{I}$, $\nu$) is a data set that comes from the second bullet of Theorem 1.1a. Theorem 1.1b introduces the notion of a *point of discontinuity* for the bundle $\mathcal{I}$. A point $p \in Z$ is a point of discontinuity for $\mathcal{I}$ if $\mathcal{I}$ is not isomorphic to the product $\mathbb{R}$ bundle on the complement of Z in any neighborhood of p. These are the interesting points in Z because if p is not a point of discontinuity for $\mathcal{I}$, then $\nu$ near p can be viewed as an honest $\mathbb{R}$-valued harmonic 1-form.

Theorem 1.1b also introduces the notion of a *geodesic arc* in a given ball in M. Supposing that B is the ball, a geodesic arc in B is a properly embedded geodesic segment in B through B's center point.

Theorem 1.1b introduces one final notion, this being the notion of a *1-dimensional, Lipshitz graph* in M. For the present such a set is characterized as follows: Let p denote any given point in the graph. Let $x = (x_1, x_2, x_3)$ denote Euclidean coordinates for $\mathbb{R}^3$. There is a coordinate chart for M centered at p that writes the graph near p as the small $|x|$ part of the map $t \to (x_1 = t,\ x_2 = \varphi_2(t),\ x_3 = \varphi_3(t))$ with $\varphi = (\varphi_1, \varphi_2)$ being a Lipshitz map from $\mathbb{R}$ to $\mathbb{R}^2$. By way of a reminder, a continuous map from an interval $I \subset \mathbb{R}$ into a Riemannian manifold is said to be Lipshitz under the following circumstances: Let $\gamma$ denote the map in question. Then $\gamma$ is Lipshitz when $\sup_{t,t'\in I} \mathrm{dist}(\gamma(t), \gamma(t')) \le c_\gamma |t - t'|$ with $c_\gamma$ being a constant.

**Theorem 1.1b**: *Let* (Z, $\mathcal{I}$, $\nu$) *denote a data set that comes via the second bullet of Theorem 1.1a from a sequence of PSL(2;* $\mathbb{C}$*) connections on* M.

- Z *has Hausdorff dimension at most* 1.
- *Given* $\varepsilon > 0$ *and* $\theta > 0$, *there are finite sets of balls* $\mathfrak{U}$ *and* $\mathfrak{V}$ *with the following properties: Their union contains* Z, *the balls in* $\mathfrak{U}$ *have pairwise disjoint closure, and* $\sum_{B\in\mathfrak{U}} (\mathrm{radius}(B))^\theta < \varepsilon$. *Meanwhile, if* $B \in \mathfrak{V}$ *and if* r *is the radius of* B, *then* $B \cap \mathfrak{V}$ *is contained in the radius* $r\varepsilon$ *tubular neighborhood of a geodesic arc in* B.
- *An open dense subset of* Z *is contained in a countable union of Lipshitz graphs.*

- *The points of discontinuity for* $\mathfrak{I}$ *are the closure of an open set in* Z *that is an embedded* $C^1$ *curve in* M.

As remarked at the outset, most of the proof of Theorem 1.1b is in [T2].

By way of a parenthetical remark, it is not out of the question that there is a 'generic metric' theorem to the effect that Z is necessarily the union of a finite set of points with a compact, embedded Lipshitz or even smooth curve if the metric on M is chosen from a suitable dense set. It is also possible that Z as just described if the sequence $\{\mathbb{A}_n\}_{n=1,2,\ldots}$ is chosen to have some additional properties; a case to consider is that where the sequence sits on an integral curve of $\mathfrak{F}$'s gradient vector field and the corresponding sequence of $\mathfrak{F}$ values converges to the infimum of $\mathfrak{F}$.

There is a certain topological significance in Z, $\mathfrak{I}$ and $\nu$ that are illustrated by two comparisons. The first comparison is that between what Theorem 1.1a says and what is said by their analogs for manifolds of dimension 2. The dimension 2 analog of these theorems is in the next subsection. This second comparison is that between what is said by Theorem 1.1a with no added assumptions and what can be said when the sequence $\{\mathbb{A}_n\}_{n=1,2,\ldots}$ is a sequence of flat $PSL(2;\mathbb{C})$ connections. Section 1c contains some very brief remarks about the latter case of Theorem 1.1a,b and then about the topological significance in Z, $\mathfrak{I}$ and $\nu$ in the general case.

**b) The case of dimension 2**

Let $\Sigma$ denote a compact, Riemann surface and let $P \to \Sigma$ denote a given principal SO(3) bundle. The function $\mathfrak{F}$ in (1.1) has its analog for connections on the principal $PSL(2;\mathbb{C})$ bundle $P \times_{SO(3)} PSL(2;\mathbb{C})$, the formula being identical but for the fact that the integration domain is $\Sigma$ and the Hodge dual is defined using the metric on $\Sigma$.

The upcoming Theorem 1.2 is the analog of Theorems 1.1a,b for oriented surfaces. The dimension 2 case of Theorems 1.1a,b for an non-orientable surface can be deduced from the theorem below by considering pull-backs to the double cover. The statement of the non-orientable version is left to the reader.

By way of background for Theorem 1.2, keep in mind first of all that the chosen orientation and metric for $\Sigma$ define a complex structure for $T\Sigma$ and thus a corresponding splitting of $T^*\Sigma \otimes_{\mathbb{R}} \mathbb{C}$ as the direct sum of two complex line bundles, $T^{1,0} \oplus T^{0,1}$. If $q$ is a given $P \times_{SO(3)} \mathfrak{su}(2)$ valued 1-form, then $q$ has a corresponding decomposition as $q^{1,0} + q^{0,1}$ with $q^{1,0}$ being the part of $q$ in $(P \times_{SO(3)} \mathfrak{su}(2)) \otimes_{\mathbb{R}} T^{1,0}$ and with $q^{0,1}$ denoting the part of $q$ in $(P \times_{SO(3)} \mathfrak{su}(2)) \otimes_{\mathbb{R}} T^{0,1}$.

Theorem 1.2 refers to what is known as a quadratic differential. The latter is a section of the complex line bundle $T^{2,0} = T^{1,0} \otimes_{\mathbb{C}} T^{1,0}$. This bundle has a canonical holomorphic structure and this is used to define the notion of a holomorphic quadratic differential. The Riemann-Roch theorem asserts that the space of holomorphic quadratic differentials is a complex vector space of dimension zero if $\Sigma$ is a sphere, dimension 1 over $\mathbb{C}$ if $\Sigma$ is a torus, and dimension $3G - 3$ over $\mathbb{C}$ with G denoting the genus of $\Sigma$ when this genus is greater than 1. Let $\mu$ denote a non-trivial, holomorphic quadratic differential. The zero locus of $\mu$ consists of $4G - 4$ points counted with multiplicity.

It may or may not be the case that $\mu$ is the square of a holomorphic section of $T^{1,0}$. This is the case if each of $\mu$'s zeros has even multiplicity. In any event, a square root of $\mu$ can be defined on the complement of its zero locus as a section of the tensor product of $T^{1,0}$ with a real line bundle on this complement. The zero locus of $\mu$ is denoted in Theorem 1.2 by $Z_\mu$, and the square root of $\mu$ on $\Sigma - Z_\mu$ is denoted by $\mu^{1/2}$. The corresponding real line bundle is denoted by $\mathcal{I}_\mu$.

**Theorem 1.2**: *Let* $\Sigma$ *denote a compact, oriented Riemann surface and let* $P \to \Sigma$ *denote a principal* SO(3) *bundle. Suppose that* $\{\mathbb{A}_n = A_n + i\,\mathfrak{a}_n\}_{n=1,2,...}$ *is a sequence of connections on* $P \times_{SO(3)} PSL(2;\mathbb{C})$ *with the corresponding sequence* $\{\mathfrak{F}(\mathbb{A}_n)\}_{n=1,2,...}$ *being bounded. For each* $n \in \{1, 2, \ldots\}$*, use* $r_n$ *to denote the* $L^2$ *norm of* $\mathfrak{a}_n$.

- *If the sequence* $\{r_n\}_{n=1,2,...}$ *has a bounded subsequence, then there exists a subsequence of* $\{\mathbb{A}_n\}_{n=1,2,...}$*, hence renumbered consecutively from 1, and a corresponding sequence of automorphisms of* P*, this denoted by* $\{g_n\}_{n=1,2,...}$*, such that* $\{g_n{}^*\mathbb{A}_n\}_{n=1,2,...}$ *converges weakly in the* $L^2_1$ *topology to an* $L^2_1$ *connection on* $P \times_{SO(3)} PSL(2;\mathbb{C})$.
- *If the sequence* $\{r_n\}_{n=1,2,...}$ *has no bounded subsequence, then there exists a subsequence of* $\{\mathbb{A}_n\}_{n=1,2,...}$*, hence renumbered consecutively from* 1*, a corresponding sequence of automorphisms of* P*, this denoted by* $\{g_n\}_{n=1,2,...}$*, and a non-trivial, holomorphic quadratic differential, this denoted by* $\mu$*. These are such that*
  1) *The sequence* $\{g_n{}^*A_n\}_{n=1,2,...}$ *converges weakly in the* $L^2_{1;loc}$ *topology on* $\Sigma - Z_\mu$ *to an* $L^2_{1;loc}$ *connection on* $P|_{M-Z_\mu}$ *this denoted by* A.
  2) *The sequence* $\{r_n^{-1} g_n{}^*\mathfrak{a}_n\}_{n=1,2,...}$ *converges weakly in the* $L^2_{1;loc}$ *topology on* $M - Z_\mu$ *to a section of* $(P \times_{SO(3)} \mathfrak{su}(2)) \otimes T^*\Sigma$ *whose* $T^{1,0}$ *part is* $\mu^{1/2}\,\sigma$ *with* $\sigma$ *being a unit length,* A*-covariantly homomorphism over* $M - Z_\mu$ *from* $\mathcal{I}_\mu$ to $P \times_{\mathbb{Z}/2\mathbb{Z} \times SO(3)} \mathfrak{su}(2)$.

By way of a parenthetical remark, a connection on a principal bundle over a Riemann surface defines a holomorphic structure on any associated bundle with complex fibers. This understood, the conclusions of Theorem 1.2 are foreshadowed by Simon Donaldson's paper [D1] about stable holomorphic bundles on surfaces.

***Proof of Theorem 1.2***: As explained momentarily, this theorem constitutes a special case to Theorem 1.1a. Even so, it can be proved independently of Theorem 1.1a and with much less effort. To obtain Theorem 1.2 from Theorem 1.1a, take M in Theorem 1.1a to be the product $S^1 \times \Sigma$ with the metric being the product metric. The pull-back of the principal SO(3) bundle P on $\Sigma$ to M via the projection map to $\Sigma$ defines a principal SO(3) bundle over M, the latter is denoted also by P. The bundle P is the SO(3) bundle to use for Theorem 1.1a. The corresponding pull-back of $\{\mathbb{A}_i\}_{i=1,2,\ldots}$ defines a sequence of connections on the incarnation of $P \times_{SO(3)} PSL(2;\mathbb{C})$ as a bundle on M. This pull-back sequence is also denoted by $\{\mathbb{A}_n\}_{n=1,2,\ldots}$. Keep in mind that the value of $\Sigma$'s version of $\mathfrak{F}$ on a connection is the same up to a multipicative factor as that of (1.1) on the pulled back connection. By the same token, the $\Sigma$ version of the sequence of $L^2$ norms $\{r_n\}_{n=1,2,\ldots}$ differs by an index independent multiplicative factor from the corresponding M version of $\{r_n\}_{n=1,2,\ldots}$.

Granted what was said in the preceding paragraph, the first bullet of Theorem 1.2 follows directly from the first bullet of Theorem 1.1a provided an argument can be made to the effect that the sequence $\{g_n\}_{n=1,2,\ldots}$ from Theorem 1.1a can be assumed to be the pull-backs of a corresponding sequence that is defined on $\Sigma$. The fact that such an assumption is valid ultimately follows from the fact that any given $L^2_1$ section of the bundle $(P \times_{SO(3)} \mathfrak{su}(2)) \otimes T^*(S^1 \times \Sigma)$ restricts to almost every constant $t \in S^1$ slice as an $L^2_1$ section along the slice, and it restricts to half of these slices with $L^2_1$ norm no greater than twice that of its $L^2_1$ norm on the whole of $S^1 \times \Sigma$.

A simliar sort of argument can be used to prove that Theorem 1.1a's set Z is the product of $S^1$ and a closed set $Z_\Sigma \subset \Sigma$ and that Theorem 1.1a's real line bundle $\mathcal{I}$ is isomorphic to the pull-back via the projection map of a real line bundle defined on the complement in $\Sigma$ of $Z_\Sigma$, this denoted for now by $\mathcal{I}_\Sigma$. Moreover, such an isomorphism identifies Theorem 1.1a's version of $\nu$ with the pull-back of a harmonic, $\mathcal{I}_\Sigma$ valued 1-form on $\Sigma - Z_\Sigma$ with $Z_\Sigma$ denoting the locus where its norm is zero. Let $\nu_\Sigma$ denote the latter. Use the splitting of $T^*\Sigma \otimes_{\mathbb{R}} \mathbb{C}$ as $T^{1,0} \oplus T^{0,1}$ to write $\nu_\Sigma$ as $e + \overline{e}$ with $e$ denoting the $T^{1,0} \otimes_{\mathbb{R}} \mathcal{I}_\mu$ part of $\nu_\Sigma$. The section $e$ is holomorphic, this being a consequence of the fact that $\nu_\Sigma$ is harmonic. Since $e$ is holomorphic, its square $e^2$ is a holomorphic section of $T^{2,0}$ over $\Sigma - Z_\Sigma$ whose norm vanishes on $\Sigma$. This last observation implies that $e^2$ extends over $Z_\Sigma$ so as to define a holomorphic quadratic differential on $\Sigma$. The latter is Theorem 1.2's quadratic differential $\mu$. Since the zero locus of $e$ is $Z_\Sigma$, this is likewise the zero locus of $\mu$. Thus $Z_\Sigma$ is Theorem 1.2's set $Z_\mu$ and Theorem 1.2's principal $\mathbb{Z}/2\mathbb{Z}$ bundle $\mathcal{I}_\mu$ is $\mathcal{I}_\Sigma$.

**c) The topological significance of Z, $\mathcal{I}$ and $\nu$**

There are two parts to what follows. The first part talks about Theorem 1.1a when the connections in the sequence $\{\mathbb{A}_n\}_{n=1,2,\ldots}$ are irreducible, flat connections. The second part gives a very brief account of how the data Z, $\mathcal{I}$ and $\nu$ in the general case lead to structures that are mainstays of 3-dimensional topology.

*Part 1*: Morgan and Shalen in [MS1]-[MS4] use $\mathbb{R}$-trees to define a compactification of the equivalence classes of irreducible representations of $\pi_1(M)$ in $PSL(2;\mathbb{C})$ with the equivalence relation defined by conjugation by $PSL(2;\mathbb{C})$. The latter space is a disjoint union of subspaces with each equivalent to the space of automorphism classes of irreducible, flat $PSL(2;\mathbb{C})$ connections on some principal $PSL(2;\mathbb{C})$ bundle over M. As noted in the introduction, the Morgan-Shalen compactification involves $\pi_1(M)$-equivariant maps from M's universal cover to $\mathbb{R}$-trees and Daskoloupoulus, Dostoglu and Wentworth [DDW1] proved that the Morgan-Shalen maps can be taken to be harmonic. The next paragraph gives a rough summary of [DDW1]'s construction of this map.

Suppose that $\{\mathbb{A}_n\}_{n=1,2,\ldots}$ is a sequence of flat connections on $P \times_{SO(3)} PSL(2;\mathbb{C})$ that is described by the second bullet in Theorem 1.1a. Donaldson [D2] and Corlette [Co] prove that there is a connection on the $\mathrm{Aut}(P\times_{SO(3)} PSL(2;\mathbb{C}))$ orbit of each member of this sequence with $\mathfrak{F}$ at most 1. This understood, Daskoloupoulus, Dostoglu and Wentworth start with a sequence $\{\mathbb{A}_n\}_{n=1,2,\ldots}$ with just this property. Daskoloupoulus, Dostoglu and Wentworth use what is known as the 'developing map' to obtain a $\pi_1$-equivariant harmonic map from M's universal cover to the hyperbolic 3-ball from each connection in the sequence. They then define a sequence of index n dependent rescalings of the hyperbolic metric on the hyperbolic 3-ball and then view the hyperbolic ball with this sequence of rescaled metrics as a sequence of pointed metric spaces. Having taken this view, a theorem of Korevaar and Schoen [KS] is invoked to conclude that the sequence of metric spaces converges in a suitable sense to an $\mathbb{R}$ tree with an action of $\pi_1(M)$ and that the sequence of harmonic maps to the hyperbolic ball has a subsequence that converges to a $\pi_1(M)$ equivariant harmonic map to this $\mathbb{R}$-tree. Daskoloupoulus, Dostoglu and Wentworth subsequently show that this limit $\mathbb{R}$-tree is of the sort that appears in the work of Morgan and Shalen.

To see the Daskoloupoulus, Dostoglu and Wentworth construction in context of Theorem 1.1a, note first that the hyperbolic 3-space can be viewed as the space of Hermitian, $2 \times 2$ complex matrices with determinant equal to 1, this denoted in what follows by $\mathbb{H}$. The hyperbolic metric is that defined by the norm on $T\mathbb{H}$ given by the trace of the square of a matrix. Let $\{\mathbb{A}_n\}_{n=1,2,\ldots}$ denote the sequence considered by Daskoloupoulus, Dostoglu and Wentworth. Fix an index $n \in \{1, 2, \ldots\}$. Let $\mathfrak{u}_n$ denote the map from M's universal cover $\mathbb{H}$ that is constructed by Daskoloupoulus, Dostoglu and Wentworth using the developing map with input being $\mathbb{A}_n$. This map appears in the context of Theorem 1.1a as follows: Write $\mathbb{A}_n$ as $A_n + i\,\mathfrak{a}_n$. The push-forward to M of the differential of $\mathfrak{u}_n$ is the 1-form $i\,\mathfrak{a}_n$.

The Daskoloupoulus, Dostoglu and Wentworth renormalization of the hyperbolic metric on $\mathbb{H}$ multiplies the latter by the inverse of Theorem 1.1a's constant $r_n$, this being the $L^2$ norm of the differential of $\mathfrak{u}_n$. Multiplying the metric on $\mathbb{H}$ by the factor $r_n^{-1}$ is accounted for in Theorem 1.1a by the appearance of $r_n^{-1}\mathfrak{a}_n$ in Item 3) of Theorem 1.1a's second bullet.

Theorem 1.1a's limit 1-form $\nu$ and its singular set $Z_S$ have the following interpretation in the context of [DDW1]: Use $\mathbb{T}$ to denote the [DDW1] limit $\mathbb{R}$ tree and $u$ their limit harmonic map from M's universal cover to $\mathbb{T}$. Gromov and Schoen [GS] (see also [S]) proved that $u$ can be viewed as an honest harmonic function on small balls in the complement of a set with Hausdorff dimension at most 1. The latter set appears in the context of Theorem 1.1a as the inverse image in M's universal cover of the set $Z_S$. Meanwhile, the differential of $u$ where it is an honest harmonic function is the pull-back of $\nu$ to M's universal cover.

Daskoloupoulus, Dostoglu and Wentworth [DDW2], [DDW3] tell a fascinating story about sequences of equivalence classes of flat $PSL(2;\mathbb{C})$ connections in the context of Theorem 1.2, these being sequences of connections on a principal $PSL(2;\mathbb{C})$ bundle over a Riemann surface.

*Part 2*: A tetrad of closely related notions in 3-manifold topology are singular measured foliations, measured laminations, weighted branched surfaces and maps to $\mathbb{R}$-trees. To paraphrase Hatcher and Oertel [HO], measured laminations are extremely useful generalizations of two central notions in 3-manifold topology, incompressible surfaces and foliations without Reeb components. Weighted branched surfaces and certain sorts of measured laminations are in some sense, two sides of the same coin. These notions were developed extensively by a number of people in concert and separately, among others Oertel [O], Hatcher [Hat], Morgan and Shalen [MS2], [MS3], and Gabai with Ortel [GO]. As explained by [HO], certain sorts of measured singular foliations give rise to measured laminations and weighted branched surfaces, and vice versa. Measured laminations are also closely related to maps from the universal cover to $\mathbb{R}$-trees, this being the central theme in [MS2]. Meanwhile, Bowditch [B] explains how to use measured, singular foliations to define maps to $\mathbb{R}$-trees. The reader should consult these references and the myriad of more recent articles to learn more about this tetrad of notions and their use in 3-manifold topology. As this author is a neophyte on this subject, no more will (or can) be said here except to point out that the notions from the tetrad serve as the dimension 3 generalization of notions that play central roles in research on the structure of Teichmuller space and the mapping class group for surfaces. An elegant account of the 2-dimensional story can found in the beautiful book by Calegari [Ca].

The singular foliation member of the tetrad appears in the context of Theorem 1.1a. To say more about how this comes about, suppose for the moment that $\Gamma \to M$ is a compact, embedded curve, that $\mathcal{I}_\Gamma$ is a real line bundle on $M-\Gamma$ and that $\nu_\Gamma$ is a smooth, closed $\mathcal{I}_\Gamma$-valued 1-form on $M-\Gamma$ with zero locus being the union of $\Gamma$ with a finite set of points in $M-\Gamma$. Data of this sort can be found if M has a suitable branched cover with branch locus being $\Gamma$. Let $Z_\Gamma$ denote the zero locus of $\nu_\Gamma$. The kernel of $\nu_\Gamma$ defines a 2-plane subbundle in M's tangent bundle over $M-Z_\Gamma$. This subbundle is integrable because $\nu_\Gamma$ is closed, and so it is everywhere tangent to the leaves of a foliation of $M-Z_\Gamma$. Moreover, the folation defined by $\nu_\Gamma$ is transversely measured with the measure given by integration of the pull-back of $|\nu_\Gamma|$ along curves that are transversal to the leaves. This foliation can be viewed as a singular, transversely measured foliation on M. More to the point, if the structure of $\nu_\Gamma$ near $Z_\Gamma$ is reasonable in a certain precise sense, then this singular foliation will have the local structure that is needed so as to invoke what is said in [HO].

Let Z, $\mathcal{I}$ and $\nu$ be as described in the second bullet of Theorem 1.1a. Being closed, the $\mathcal{I}$-valued 1-form $\nu$ defines a transversely measured foliation on $M-Z$ and so a singular, transversely measured foliation on M. It is possible that Z is a finite union of embedded curves when the metric on M is suitably generic, and in this case, it is likely that the $\mathcal{I}$-valued 1-form $\nu$ can be approximated by one that behaves in the desired manner near its zero set.

**d) An overview of the proof of Theorem 1.1a**

The first bullet of Theorem 1.1a is proved in Section 2a. A fundamental Bochner-Weitzenboch formula makes the first bullet little more than a corollary to Uhlenbeck's theorem for connections on principal SO(3) bundles.

The proof of the second bullet of Theorem 1.1a has three components. The first component obtains an $L^2_1$ limit of a subsequence of the rescaled sequence $\{r_n^{-1}|\mathfrak{a}_n|\}_{n=1,2,\ldots}$, this being the content of Lemma 2.1. The subsequence is hence renumbered from 1. As it turns out, the limit function is bounded, but the convergence is not $L^\infty$ convergence. For this and for other reasons, it proves necessary to modify the sequence $\{r_n^{-1}\mathfrak{a}_n\}_{n=1,2,\ldots}$ so as to obtain a new sequence, this denoted by $\{\hat{a}_n\}_{n=1,2,\ldots}$, with $\{|\hat{a}_n|\}_{n=1,2\ldots}$ being uniformly bounded, converging in $L^2_1$ to the same limit function, but with the following additional property: The limit $L^2_1$ function, now denoted by $|\hat{a}_\Diamond|$, is defined at each point in M by the rule $|\hat{a}_\Diamond| = \limsup_{n\to\infty}|\hat{a}_n|$. The properties of $\{\hat{a}_n\}_{n=1,2,\ldots}$ are described in Proposition 2.2.

This notion of pointwise convergence brings up a subtle but central issue, which is that pointwise convergence to an $L^\infty$ limit does not imply that the limit is $C^0$. The second component of the proof of Theorem 1.1a consists of a proof that $|\hat{a}_\diamond|$ is a continuous function. The assertion that $|\hat{a}_\diamond|$ is $C^0$ is made by Proposition 6.1. The intervening Sections 3, 4 and 5 develop the tools that are needed to prove Proposition 6.1. This proof that $|\hat{a}_\diamond|$ is continuous brings to bear, among other things, Uhlenbeck's theorem for SO(3) connections, the properties of a certain canonical, constant coefficient, first order elliptic operator on $\mathbb{R}^3$ that is defined by the non-linear structure of the curvature, and a gauge theoretic notion of the frequency function that was introduced by Almgren [Al] and used to great success by [HHL] and others to study the singularities of nodal sets of eigenfunctions of the Laplacian.

Uhlenbeck's theorem is brought to bear in Section 3 to study the behavior of the sequence $\{(A_n, \hat{a}_n)\}_{n=1,2,\ldots}$ on sequences of small radius balls about each point of M; the radii of the balls in any given such sequence depend on the chosen point and the index n. With a given point in M fixed, Section 4 uses the analysis in Section 3 to draw conclusions about the sequence whose n'th term is the curvature of the connection $A_n$ on the corresponding index n dependent radius ball about the point. This analysis brings to bear the aforementioned, constant coefficient first order operator. Section 5 defines the gauge theory analog of Almgren's frequency function and proves that it obeys an approximate monotonicity formula, the latter being the gauge theory analog of the monotonicity formula that is exploited by Almgren and others to study the nodal sets of eigenfunctions of the Laplacian. Section 6 uses the conclusions of Sections 3, 4 and 5 to prove that $|\hat{a}_\diamond|$ is continuous and that it is uniformly Hölder continuous.

The final component of the proof of Theorem 1.1a is in Section 7; it begins with the definition of Z and then the construction of Theorem 1.1a's real line bundle $\mathcal{I}$ over M–Z and the construction of Theorem 1.1a's $\mathcal{I}$ valued 1-form $\nu$. The set Z is defined to be the zero locus of $|\hat{a}_\diamond|$. The fact that $|\hat{a}_\diamond|$ is continuous implies that Z is a *closed* subset of M and thus that M–Z is an open set. The fact that Z is closed rules out all sorts of terrible pathologies, among them the possibility that M–Z has empty interior. The fact that Z is nowhere dense follows from Item a) of the fourth bullet of Proposition 7.1 and the fact that Z is closed.

Proposition 7.1 reasserts the existence of the real line bundle $\mathcal{I} \to$ M–Z; and it reasserts the existence of the harmonic section $\nu$ of $T^*M \otimes \mathcal{I}$ on M–Z with $|\nu| = |\hat{a}_\diamond|$. The proof of Proposition 7.1 is in Sections 7a and 7b. The proof that $|\nu|$ is Hölder continuous is in Section 7d. Items 2) and 3) of Theorem 1.1a follow directly from Lemma 7.2 which is in Section 7a.

**e) On the proof of Theorem 1.1b**

The first and second bullets of Theorem 1.1b constitute a special case of Theorem 1.3 in [T2], the third bullet of Theorem 1.1b is a special case of Theorem 1.4 in [T2] and the fourth bullet of Theorem 1.1b is a special case of Theorem 1.2 in [T2]. The input from this paper that is needed to invoke Theorems 1.2, 1.3 and 1.4 in [T2] consists of a data of the following sort:

- *A continuous non-negative function on* M *to be denoted by* $f$ *obeying*
  - i) $f > 0$ *somewhere*.
  - ii) *There exists* $\varepsilon > 0$ *such that if* $p \in M$ *and* $f(p) = 0$, *and if* r *is sufficiently small but positive, then* $\int_{\mathrm{dist}(p,\cdot)<r} f^2 \le r^{3+\varepsilon}$.
- *Let* Z *denote* $f^{-1}(0)$. *A real line bundle* $\mathcal{I} \to M - Z$.
- *A section* $\nu$ *of* $T^*M \otimes \mathcal{I}$ *over* $M-Z$ *that obeys*
  - i) $d\nu = 0$ *and* $d{*}\nu = 0$.
  - ii) $|\nu| = f$.
  - iii) *The function* $|\nabla\nu|^2$ *is integrable on* $M-Z$.

(1.2)

The data set $(Z, \mathcal{I}, \nu)$ from Proposition 7.1 with $f = |\nu|$ obeys the conditions in (1.2); the condition in Item ii) of the first bullet being a consequence of Items b) and c) of the fourth bullet of Proposition 7.1 and what is said by the second bullet of Lemma 7.7. The existence of the desired $\varepsilon$ also follows from Proposition 7.1's identification $|\nu| = |\hat{a}_\lozenge|$ and from what is said in Proposition 6.1 to the effect that $|\hat{a}_\lozenge|$ is uniformly Hölder continuous along its zero locus.

**f) On the differences between this paper and [T1]**

The statements in Theorem 1.1a corresponds to the parts of Theorem 1.1 in [T1] that were not affected by the error found by Thomas Walpuski. This error was in the proof of Lemma 8.8 in [T1]. Lemma 8.8 of [T1] makes a statement to the effect that Z has a unique tangent cone at each of its points and the error puts the assertion of that lemma in doubt. Thus, the statements in Theorem 1.1 of [T1] that were proved using [T1]'s Lemma 8.8 needed either a revision or another proof. These statements or their revisions are contained in Theorem 1.1b.

But for the doubtful Lemma 8.8, much of what is said Sections 8, 9 and 10 in [T1] now constitutes a special case of what is done in [T2]. This is why these sections do not have corresponding sections in this revision of [T1].

Other salient changes from the published version [T1] are listed below.

**SECTION 2:** Various relatively minor changes based on the comments of Ben Mares.

**SECTION 4**: Changes in Sections 4b and 4c correct some mystifying constructions in the corresponding parts of [T1] and they lead to a simpler proof of Lemma 4.3.

**SECTION 6**: The statement of Proposition 6.1 has been changed. The statement here makes an assertion to the effect that $|\hat{a}_\diamond|$ is uniformly Hölder continuous across its zero locus. This is stronger than the analogous statement in the [T1] version of Proposition 6.1. This stronger statement about Hölder continuity in Proposition 6.1 is restated and proved in a new version of Lemma 6.4. The proof given here of Lemma 6.4 also corrects some confusing statements in the proof of Lemma 6.4 in [T1].

**SECTION 7**: The proof of the second bullet of Lemma 7.7 is corrected. The proof of this bullet in the [T1] version of this lemma is circular. The corrected version uses the uniform Hölder continuity assertion in the new version Proposition 6.1.

### g) Extensions of Theorems 1.1a,b

This subsection briefly describes various contexts where the techniques and strategies that are used prove Theorem 1.1a,b could prove useful. See also [GU], [HW], [Tan] and [T3], for related papers that appeared after the writing of the original version of this paper.

GROUPS WITH RANK GREATER THAN ONE: As noted in the introduction, there is likely some sort of analog of Theorem 1.1a,b for Lie groups such as $PSL(n;\mathbb{C})$ for $n > 2$ The statement of a hypothetical $n > 2$ version of Theorem 1.1a,b will almost surely be more involved by virtue of the fact that the Lie algebra of such a group has non-Abelian subalgebras. In particular, the role that is played in Item 3 of Theorem 1.1a's second bullet by $\nu\sigma$ will likely involve a 1-form with values in a twisted vector bundle whose fiber is in a non-Abelian subalgebra of the group's Lie algebra. The singular locus Z may well extend beyond the zero locus of this $\nu\sigma$ analog to account for possible ranks of its stabilizer in the group.

MANIFOLDS OF DIMENSION GREATER THAN THREE: There is likely an $L^p$ version of Theorem 1.1a,b for any p greater than half the dimension of the ambient manifold. A case of special concern is the $p = 2$ case for a manifold of dimension 4. The reason being that this case is relevant to any attempt to define $PSL(2;\mathbb{C})$ analogs of Floer homology and $PSL(2;\mathbb{C})$ analogs of Donaldson's invariants. A bit more is said below about these analogs. The analog of Theorem 1.1a,b in the case where p is half the dimension of the manifold will be more complicated because this is so for the analogous version of Uhlenbeck's theorem. There may well be additional complications. In any event, the singular set in the case where the dimension is greater than 3 will likely involve some sort of union of rectifiable, codimension 2, Lipshitz submanifolds.

PSL(2;$\mathbb{C}$) FLOER HOMOLOGY: Of interest here are two sorts of equations for maps from $\mathbb{R}$ to the space of connections on $P \times_{SO(3)} PSL(2;\mathbb{C})$, these being

- $\frac{d}{dt} A = -*(F_A - \mathfrak{a} \wedge \mathfrak{a})$ *and* $\frac{d}{dt} \mathfrak{a} = *d_A\mathfrak{a}$ .
- $\frac{d}{dt} A = -*d_A\mathfrak{a}$ *and* $\frac{d}{dt} \mathfrak{a} = -*(F_A - \mathfrak{a} \wedge \mathfrak{a})$ .

(1.2)

The former is the gradient flow for the real part of the Chern-Simons functional

$$\mathbb{A} \to \mathrm{CS}(\mathbb{A}) = \frac{1}{2} \int_M \mathrm{tr}(\mathbb{A} \wedge (F_{\mathbb{A}} - \tfrac{1}{3} \mathbb{A} \wedge \mathbb{A})) \ . \tag{1.3}$$

The second equation in (1.2) is the gradient flow for the imaginary part of CS, this being the Hamiltonian flow for the real part as defined using a certain canonical symplectic form on the space of $P \times_{SO(3)} PSL(2;\mathbb{C})$ connections. The function CS is decreasing with respect to the gradient flow and constant along the Hamiltonian flow. The imaginary part of CS is constant under the gradient flow and is monotonic under the Hamiltonian flow. Witten [W1], [W2] conjectured that a certain $\mathbb{CP}^1$ parameterized family of linear combinations of the four equations in the first and second bullets of (1.2) can be used to give a gauge theoretic construction of Khovanov homology. See also [GW]. A similar suite of equations were introduced in [Hay]. Such linear combinations also enter in the work of Kapustin and Witten [KW] on the geometric Langlands program. Witten proposed in [W3], [W4] that the equations in (1.2) could be used to compute certain formal path integrals of the Chern Simons functional.

Solutions to the $\mathfrak{a} = 0$ version of the equation in the top bullet of (1.2) are used to define the differential for the SO(3) Floer homology on M; and the $L^2$ version of Uhlenbeck's theorem for the manifold $\mathbb{R} \times M$ plays a central role in the proof that this differential has square zero. See [Fl] and also [D3]. This being the case, it is almost sure bet that some sort of PSL(2;$\mathbb{C}$) extension of Uhlenbeck's theorem will be needed to define analogous algebraic structures using solutions to the equations in (1.2) and to the sorts of linear combinations that are introduced by Witten.

What follows is a likely relevant observation with regards to any such extension: The functional $\mathbb{A} \to \int_M |d_A * \mathfrak{a}|^2$ is constant along all of these flows.

PSL(2;$\mathbb{C}$) SELF DUALITY: Let X denote a smooth, compact and oriented 4 dimensional Riemannian manifold. As noted by Witten [W1], the equations in (1.2) have analogs on X, these being equations for a pair consisting of a connection on a principal SO(3) bundle over X and a 1-form with values in the associated vector bundle with fiber $\mathfrak{su}(2)$ given by the adjoint representation. Let (A, $\mathfrak{a}$) denote such a pair. Use the metric to define the respective bundles of self-dual and anti-self dual 2-forms. The orthogonal projection to these respective bundles are denoted by $\Pi^+$ and $\Pi^-$. The analog of the top equation in (1.2) reads

$$\Pi^+(F_A - \mathfrak{a}\wedge\mathfrak{a}) = 0 \quad \textit{and} \quad \Pi^- d_A\mathfrak{a} = 0 \quad \textit{and} \quad d_A{*}\mathfrak{a} = 0 \,. \tag{1.4}$$

There are corresponding analogs to the lower equation in (1.2) and to suitable linear combinations of the four equations in (1.2). See [DK] to read about the applications of the $\mathfrak{a} = 0$ version of (1.4)

Witten and Vafa [VW] proposed an alternate generalization of the equations in (1.2), this being a system of equations for a connection on a principal SO(3) bundle and a self-dual 2-form with values in the same associated bundle with fiber $\mathfrak{su}(2)$. Let (A, $\mathfrak{w}$) denote such a pair. The Witten-Vafa equations are written schematically as

$$\Pi^+(F_A - [\mathfrak{w};\mathfrak{w}]) = 0 \quad \textit{and} \quad d_A\mathfrak{w} = 0 \,, \tag{1.5}$$

where $[\cdot,\cdot]$ here denotes a certain canonical, bilinear, symmetric fiber preserving map that is defined by the metric's Hodge dual and the commutator on $\mathfrak{su}(2)$. The equations in (1.5) play a central role in recent work of Gaddes, Gukov and Putrov [GGP].

**h) Table of contents and conventions**

What follows is a table of contents for this article.

The paper employs two conventions throughout. The first convention has $c_0$ denoting a number that is greater than 100 whose value does not depend on any of the salient issues under consideration in a given assertion. The value of $c_0$ in any given appearance can depend on the particulars of M and its Riemannian metric; and also on the upper bound for Theorem 1.1a's sequence $\{\mathfrak{F}(\mathbb{A}_n)\}_{n=1,2,\ldots}$. However, under no circumstances does it depend on the index n. The value of $c_0$ in successive appearances can be assumed to increase.

The second convention concerns what is denoted by $\chi$. This is a fixed, smooth and nonincreasing function on $\mathbb{R}$ that equals 1 on $(-\infty, \frac{1}{4}]$ and equals 0 on $[\frac{3}{4}, \infty)$. A favorite version is chosen now and used throughout the paper.

## 2. $L^2_1$ and pointwise limits

The four subsections that follow comprise what Section 1c described as the Part 1 of the proof of Theorem 1.1a. The first Section 2a begins with a Bochner-Weitzenboch formula that plays a central role in the remaining subsections. It then uses this formula to prove the first bullet of Theorem 1.1a. The Section 2b begins the proof of Theorem 1.1a's second bullet with the construction an $L^2_1$ limit from the sequence $\{r_i^{-1}|\mathfrak{a}_i|\}_{i=1,2,\ldots}$. The third subsection modifies the sequence $\{\mathfrak{a}_i\}_{i=1,2,\ldots}$ to obtain a new sequence that allows greater control over the limit in Section 1b. The salient features of this new sequence are summarized by Proposition 2.2. The final subsection proves a central lemma that is used to prove Proposition 2.2.

### a) The Bochner-Weitzenboch formula

Let A denote a given connection on P and let $\mathfrak{a}$ denote a section of $P\times_{SO(3)}\mathfrak{su}(2)$. It is important to keep in mind that there is a Bochner-Weitzenbock formula that writes

$$*d_A*d_A\mathfrak{a} - d_A*d_A*\mathfrak{a} = \nabla_A{}^{\dagger}\nabla_A\mathfrak{a} + *(*F_A\wedge\mathfrak{a} + \mathfrak{a}\wedge *F_A) + \mathrm{Ric}(\mathfrak{a}) \,, \tag{2.1}$$

with $\mathrm{Ric}(\cdot)$ denoting here the Ricci curvature in its guise as a homomorphism of $T^*M$.

It proves useful to write (2.1) as an equality between integrals over M. To this end, let $f$ denote a chosen $C^2$ function on M and let $r\in[1,\infty)$ denote a chosen positive number. Take the inner product of both sides of (2.1) with the section $f\,\mathfrak{a}$ of $P\times_{SO(3)}\mathfrak{su}(2)$ and integrate over M. An integration by parts to obtain an expression with only first derivatives of $\mathfrak{a}$ leads from the resulting integral identity to the following one:

$$\begin{aligned}\int_M f(|d_A\mathfrak{a}|^2 + |d_A*\mathfrak{a}|^2 + |r^{-1}F_A - r\,\mathfrak{a}\wedge\mathfrak{a}|^2) = \\ \tfrac{1}{2}\int_M d^*df\,|\mathfrak{a}|^2 + \int_M f(|\nabla_A\mathfrak{a}|^2 + r^2|\mathfrak{a}\wedge\mathfrak{a}|^2 + r^{-2}|F_A|^2 + \mathrm{Ric}(\langle\mathfrak{a}\otimes\mathfrak{a}\rangle)) \\ - \int_M (df\wedge *\langle\mathfrak{a}(*d_A*\mathfrak{a})\rangle + df\wedge\langle\mathfrak{a}\wedge *d_A\mathfrak{a}\rangle)\,.\end{aligned} \tag{2.2}$$

The notation here has $\langle \mathfrak{a} \otimes \mathfrak{a} \rangle$ denoting the symmetric section of $T^*M \otimes T^*M$ that is defined as follows: Fix an orthonormal frame for $T^*M$ at any given point and use $\{\mathfrak{a}_\alpha\}_{\alpha \in \{1,2,3\}}$ to denote the coefficients of $\mathfrak{a}$ when written using the chosen frame. The corresponding coefficients of $\langle \mathfrak{a} \otimes \mathfrak{a} \rangle$ are $\{\langle \mathfrak{a}_\alpha \mathfrak{a}_\beta \rangle\}_{\alpha,\beta \in \{1,2,3\}}$. Meanwhile, the Ricci tensor in (2.2) is viewed using the metric as a linear functional on $T^*M \otimes T^*M$.

***Proof of the first bullet of Theorem 1.1a***: Take $f = 1$ and $r = 1$ in (2.2) to see that

$$\mathfrak{F}(A + i\,\mathfrak{a}) = \int_M (|\nabla_A \mathfrak{a}|^2 + |\mathfrak{a} \wedge \mathfrak{a}|^2 + |F_A|^2 + \mathrm{Ric}(\langle \mathfrak{a} \otimes \mathfrak{a} \rangle)) \ . \tag{2.3}$$

It follows as a consequence that

$$\int_M (|\nabla_A \mathfrak{a}|^2 + |\mathfrak{a} \wedge \mathfrak{a}|^2 + |F_A|^2) \le \mathfrak{F}(A + i\,\mathfrak{a}) + c_0 \int_M |\mathfrak{a}|^2 \ . \tag{2.4}$$

To exploit (2.4), suppose that $\{\mathbb{A}_n = A_n + i\,\mathfrak{a}_n\}_{n=1,2...}$ is a sequence with both $\{\mathfrak{F}(\mathbb{A}_n)\}_{n=1,2,..}$ and $\{\int_M |\mathfrak{a}_n|^2\}_{n=1,2,..}$ being bounded. Then the sequence $\{\int_M |F_{A_n}|^2\}_{n=1,2,...}$ is bounded and so Uhlenbeck's theorem finds a subsequence of $\{\mathbb{A}_n\}_{n=1,2,..}$ (hence renumbered consecutively from 1) and corresponding sequence of automorphisms of P, this denoted by $\{g_n\}_{n=1,2,...}$ such that $\{g_n{}^*A_n\}_{n=1,2,...}$ converges weakly in the $L^2_1$ topology to an $L^2_1$ connection on P. Denote the latter by A. Meanwhile, the sequence $\{\int_M |\nabla_{A_n} \mathfrak{a}_n|^2\}_{n=1.2,...}$ is also bounded, and this implies that $\{g_n{}^*\mathfrak{a}_n\}_{n=1,2,...}$ has a subsequence that converges weakly in the $L^2_1$ topology on the space of sections of $P \times_{SO(3)} \mathfrak{su}(2)$ to an $L^2_1$ section, this denoted by $\mathfrak{a}$. The pair $(A, \mathfrak{a})$ define the desired limit $L^2_1$ connection on the bundle $P \times_{SO(3)} PSL(2;\mathbb{C})$.

**b) Renormalization and $L^2_1$ convergence**

This subsection begins the proof of the second bullet of Theorem 1.1a. The input is a sequence $\{\mathbb{A}_n = A_n + i\,\mathfrak{a}_n\}_{n=1,2,...}$ with $\{\mathfrak{F}(\mathbb{A}_n)\}_{n=1,2,...}$ being bounded, but not so for the sequence whose n'th term is the $L^2$ norm of $\mathfrak{a}_n$. The latter sequence is assumed to be unbounded with no finite limits. Fix n and set $\hat{a}_n = r_n^{-1}\,\mathfrak{a}_n$ with $r_n$ denoting the $L^2$ norm of the section $\mathfrak{a}_n$. Fix a $C^2$ function, $f$, on M. Multiply the $(A, \mathfrak{a}) = (A_n, \mathfrak{a}_n)$ and $r = 1$ version of (2.2) by $r_n^{-2}$ to see that

$$\lim_{n\to\infty} \tfrac{1}{2} \int_M d^\dagger d f\, |\hat{a}_n|^2 + \lim_{n\to\infty} \int_M f(|\nabla_{A_n} \hat{a}_n|^2 + r_n^2\, |\hat{a}_n \wedge \hat{a}_n|^2 + r_n^{-2}\, |F_{A_n}|^2 + \mathrm{Ric}(\langle \hat{a}_n \otimes \hat{a}_n \rangle)) = 0. \tag{2.5}$$

This is so because the terms with $r_n^{-1} d_{A_n}\mathfrak{a}_n$ and $r_n^{-1} d_{A_n} {*\mathfrak{a}_n}$ and $r_n^{-2}|F_{A_n} - \mathfrak{a}_n\wedge\mathfrak{a}_n|$ have limit zero as $n \to \infty$.

**Lemma 2.1**: *There exists* $\kappa > 1$ *with the following significance: Let* $\{\mathbb{A}_n = A_n + i\,\mathfrak{a}_n\}_{n=1,2,\ldots}$ *denote a sequence of connections on* $P\times_{SO(3)}PSL(2;\mathbb{C})$ *with* $\{\mathfrak{F}(\mathbb{A}_n)\}_{n=1,2,\ldots}$ *being bounded, but with the sequence whose n'th term is the* $L^2$ *norm of* $\mathfrak{a}_n$ *being unbounded with no finite limits. For each* $n \in \{1, 2, \ldots\}$, *set* $r_n$ *to be the* $L^2$ *norm of* $\mathfrak{a}_n$ *and set* $\hat{a}_n = r_n^{-1}\mathfrak{a}_n$. *There is a subsequence in* $\{1, 2, \ldots\}$, *hence renumbered consecutively from 1 with the properties that are listed below.*

- *The sequence* $\{|\hat{a}_n|\}_{n=1,2,\ldots}$ *converges weakly in the* $L^2_1$ *topology and strongly in each* $p < 6$ *version of the* $L^p$ *topology to an* $L^2_1$ *function on M, this denoted by* $|\hat{a}|$. *The* $L^2$ *norm of* $|\hat{a}|$ *is* $1$ *and its* $L^2_1$ *norm is bounded by* $\kappa$. *Moreover,* $|\hat{a}|$ *defines an* $L^\infty$ *function with* $|\hat{a}| < \kappa$ *almost everywhere.*
- *The sequence* $\{\langle\hat{a}_n\otimes\hat{a}_n\rangle\}_{n=1,2,\ldots}$ *converges strongly in each* $q < 3$ *version of the* $L^q$ *topology on the space of sections of* $T^*M\otimes T^*M$ *and weakly in the* $L^3$ *topology. The limit section of* $T^*M\otimes T^*M$ *is denoted by* $\langle\hat{a}\otimes\hat{a}\rangle$*; it is in* $L^\infty$ *and its trace is the function* $|\hat{a}|^2$.
- *Let f denote any given* $C^2$ *function on* M. *The three sequences* $\{\int_M f\,|\nabla_{A_n}\hat{a}_n|^2\}_{n=1,2,\ldots}$, $\{\int_M r_n^{2} f\,|\hat{a}_n\wedge\hat{a}_n|^2\}_{n=1,2,\ldots}$ *and* $\{\int_M r_n^{-2} f\,|F_{A_n}|^2\}_{n=1,2,\ldots}$ *each converge. These limits are denoted by* $Q_{\nabla,f}$, $Q_{\wedge,f}$ *and* $Q_{F,f}$.
  a) *Each of* $Q_{\nabla,f}$, $Q_{\wedge,f}$ *and* $Q_{F,f}$ *is bounded by* $\kappa$ *times the sum of the supremum norm and the* $L^3_1$ *norm of f. Moreover,* $Q_{\wedge,f} = Q_{F,f}$*; and if* $f \ge 0$, *then* $Q_{\nabla,f} \ge \int_M f\,|d|\hat{a}||^2$.
  b) $\frac{1}{2}\int_M d^\dagger df\,|\hat{a}|^2 + Q_{\nabla,f} + 2Q_{\wedge,f} = -\int_M f\,\mathrm{Ric}(\langle\hat{a}\otimes\hat{a}\rangle)$.

***Proof of Lemma 2.1***: The assertion in the first bullet to the effect that $\{|\hat{a}_n|\}_{n=1,2,\ldots}$ has a subsequence which converges weakly in the $L^2_1$ topology follows if the sequence has uniformly bounded $L^2_1$ norm. That this is the case follows from the $f = 1$ version because

$$|d|\mathfrak{b}|| \le |\nabla_A\mathfrak{b}| \tag{2.6}$$

with A being any connection on P and $\mathfrak{b}$ being any tensor valued section of $P\times_{SO(3)}\mathfrak{su}(2)$. The $L^\infty$ assertions are proved momentarily. The proof of the remaining bullets do not require the $L^\infty$ assertion in the first bullet.

With regards to the second bullet, the inequality in (2.6) with (2.4) implies that the sequence $\{|\hat{a}_n|^{-1}\langle\hat{a}_n\otimes\hat{a}_n\rangle\}_{n=1,2,\ldots}$ is bounded in the $L^2_1$ topology. It follows as a consequence that a subsequence converges weakly in the $L^2_1$ topology and strongly in the

$L^p$ topology for $p < 6$. Since this is also the case for $\{|\hat{a}_n|\}_{n=1,2,\ldots}$, the product sequence $\{\langle \hat{a}_n \otimes \hat{a}_n \rangle\}_{n \in 1,2,\ldots}$ converges strongly in the $L^q$ topology for $q < 3$.

The existence of a subsequence that makes the third bullet true for any given $C^2$ function $f$ follows from (2.5). The third bullet's assertion is true for all $C^2$ functions $f$ simultaneously is due to the fact that the space of $C^2$ functions on M has a countable, dense subset. The equality $Q_{\wedge,f} = Q_{F,f}$ follows because $\lim_{n\to\infty} r_n^{-2} \int_M |F_{A_n} - r_n^2 \hat{a}_n \wedge \hat{a}_n|^2 = 0$. The $f \geq 0$ upper bound on $Q_{\nabla,f}$ follows from (2.6).

The four steps that follow momentarily prove that $|\hat{a}|$ defines an $L^\infty$ function with the asserted norm bound. By way of a reminder, a measurable function defines an $L^\infty$ function if it's norm is bounded on the complement of a measure zero set.

<u>Step 1</u>: Fix $p \in M$ and let $G_p$ denote the Green's function with pole at p for the operator $d^\dagger d + 1$, this being a smooth, non-negative function on M–p which extends to the whole of M as an $L^q$ function for any $q < 3$. More to the point, $G_p(\cdot) \leq c_0 \operatorname{dist}(p,\cdot)^{-1}$ which is a fact that is exploited below. Keep in mind as well that $|dG_p| \leq c_0 \operatorname{dist}(p,\cdot)^{-2}$.

The plan for what follows is to obtain an $L^\infty$ bound on $|\hat{a}|$ from a sequence of $C^2$ functions that converge in a suitable sense to the Green's function $G_p$. In particular, what follows directly describes a $(0, c_0^{-1}]$ parametrized sequence of $C^2$ approximations to $G_p$ that converge to $G_p$ as the parameter limits to 0 in the $L^q$ topology on $C^2(M)$ for $q < 3$ and in the $C^2$ topology on compact subsets of M–p. The $\varepsilon \in (0, c_0^{-1}]$ member of this sequence is denoted by $f_{p,\varepsilon}$. To define $f_{p,\varepsilon}$, let $\upsilon_{p,\varepsilon}$ denote the function that equals 1 on ball of radius $\varepsilon(1-\varepsilon^3)$ centered at p, is given by the rule $\varepsilon^{-3}(\varepsilon - \operatorname{dist}(p,\cdot))$ between on the spherical annulus between radius $\varepsilon(1-\varepsilon^3)$ centered at p and the radius $\varepsilon$ ball centered at p, and is zero on the complement of the radius $\varepsilon$ ball centered at p. This is a continuous, piecewise differentiable function that is very nearly the characteristic function of the radius $\varepsilon$ ball centered at p. Let $v_{p,\varepsilon}$ denote the integral of $\upsilon_{p,\varepsilon}$ and define $\delta_{p,\varepsilon}$ to be $v_{p,\varepsilon}^{-1} \upsilon_{p,\varepsilon}$. The function $f_{p,\varepsilon}$ is the solution on M to the equation $d^\dagger d f + f = \delta_{p,\varepsilon}$. The function $f_{p,\varepsilon}$ can be shown to be twice differentiable.

A depiction of $f_{p,\varepsilon}$ near p is given momentarily. This depiction introduces a constant, $z_p$, with norm bounded by $c_0$. This constant is defined by the depiction of $G_p$ in Gaussian coordinates near p, this having the form $x \to G_p(x) = \frac{1}{4\pi|x|} - z_p + \ldots$ where the unwritten terms have norm bounded by $c_0|x|$. Note in this regard that the $\mathbb{R}^3$ analog of $G_p$ with pole at the origin is the function $x \to \frac{1}{4\pi|x|} e^{-|x|}$ and so the analog on $\mathbb{R}^3$ of $z_p$ is $-\frac{1}{4\pi}$.

With $z_p$ understood, then the function $f_{p,\varepsilon}$ can be written in a Gaussian coordinate chart centered at p as $f_{p,\varepsilon} = g_\varepsilon + \mathfrak{e}_p$ where the function $x \to g_\varepsilon(x)$ is defined by the rule

$$g_\varepsilon(x) = \frac{1}{4\pi|x|} - z_p \quad \textit{for} \quad |x| > \varepsilon \quad \textit{and} \quad g_\varepsilon(x) = \frac{3}{8\pi\varepsilon}\left(1 - \frac{|x|^2}{3\varepsilon^2}\right) - z_p \quad \textit{for} \quad |x| \leq \varepsilon\,; \tag{2.7}$$

and where $\mathfrak{e}_p$ is a continuous function that is smooth on the complement of the origin and such that $|\mathfrak{e}_p| \le c_0(|x| + \varepsilon^2)$ and $|d\mathfrak{e}_p| \le c_0(1+\varepsilon)$.

<u>Step</u> <u>2</u>: The $f = f_{p,\varepsilon}$ version of the equality given by Item b) of the third bullet of Lemma 2.1 reads

$$\frac{1}{2}\int_M \delta_{p,\varepsilon}\, |\hat{a}|^2 + Q_{\nabla, f_{p,\varepsilon}} + 2Q_{\wedge, f_{p,\varepsilon}} = \int_M f_{p,\varepsilon}\, (\tfrac{1}{2}|\hat{a}|^2 - \mathrm{Ric}(\langle \hat{a} \otimes \hat{a}\rangle)) \ . \tag{2.8}$$

The right hand side of (2.8) converges as $\varepsilon \to 0$ because $|\hat{a}|^2$ is an $L^2$ function. It follows from the lemma's first two bullets that the various $\varepsilon \in (0, c_0^{-1}]$ versions converge as $\varepsilon \to 0$ with the $\varepsilon \to 0$ limit being the integral of $G_p(\frac{1}{2}|\hat{a}|^2 - \mathrm{Ric}(\langle \hat{a} \otimes \hat{a}\rangle))$

<u>Step</u> <u>3</u>: The functions in the family $\{g_\varepsilon\}_{\varepsilon\in(0,1/c_0)}$ obey $g_\varepsilon \ge g_{\varepsilon'}$ for $\varepsilon < \varepsilon'$. This has the following consequence for the left hand side of (2.8): Let $\bullet$ denote either $\nabla$ or $\wedge$. Then $Q_{\bullet, f_{p,\varepsilon}}$ can be written as $Q_{\bullet, g_\varepsilon} + P_{\bullet(p,\varepsilon)}$ where the sequence $\{P_{\bullet(p,\varepsilon)}\}_{\varepsilon\in(0,1/c_0)}$ converges as $\varepsilon$ limits to zero and where the sequence $\{Q_{\bullet(p,\varepsilon)}\}_{\varepsilon\in(0,1/c_0)}$ is bounded and monotonically increasing as $\varepsilon$ decreases to zero. It follows that this sequence also has a unique limit as $\varepsilon$ limits to zero. Thus, the sequence $\{Q_{\bullet(p,\varepsilon)}\}_{\varepsilon\in(0,1/c_0)}$ has a unique limit as $\varepsilon$ limits to zero. The limit is denoted by $Q_{\bullet, G_p}$.

<u>Step</u> <u>4</u>: Consider now the integral of $\delta_{p,\varepsilon}\, |\hat{a}|^2$ that appears on the left hand side of (2.8). This integral is positive, and it follows from (2.8) that the various $\varepsilon \in (0, c_0^{-1})$ versions are uniformly bounded as $\varepsilon \to 0$ with a $p \in M$ independent upper bound. This implies in particular that $|\hat{a}|$ is bounded by $c_0$ on the complement of a measure zero set. Meanwhile, if p is in the complement of a measure zero set in M, then the $\varepsilon \to 0$ limit of the integral of $\delta_{p,\varepsilon}|\hat{a}|^2$ converges to $|\hat{a}|^2$. The fact that this is so follows from Theorem 3.18 in [Fo] because any given $p \in M$ version of $\upsilon_{p,\varepsilon}$ is greater than the characteristic function of the radius $\varepsilon(1 - \varepsilon^3)$ ball centered at p and less than the characteristic function of radius $\varepsilon$ ball centered at p.

Given this $\varepsilon \to 0$ convergence of the integral of $\delta_{p,\varepsilon}$ to $|\hat{a}|^2$ for p in the complement of a measure zero set being the case, it then follows from what was said in Steps 2 and 3 that $|\hat{a}|$ can be modified on a measure zero set so that

$$\frac{1}{2}|\hat{a}|^2(p) + Q_{\nabla, G_p} + 2Q_{\wedge, G_p} = -\int_M G_p\, (\tfrac{1}{2}|\hat{a}|^2 - \langle \mathrm{Ric}(\hat{a}, \hat{a})\rangle) \tag{2.9}$$

for each $p \in M$. This formula implies the asserted norm bound.

**c) Second derivative bounds**

Let A denote for the moment a connection on P and let $\mathfrak{a}$ denote a section of $(P\times_{SO(3)}\mathfrak{su}(2))\otimes T^*M$. Define $q_A(\mathfrak{a})$ to be the section of $(P\times_{SO(3)}\mathfrak{su}(2))\otimes T^*M$ given by

$$q_A(\mathfrak{a}) = \nabla_A^{\dagger}\nabla_A\mathfrak{a} + *(*F_A\wedge\mathfrak{a} + \mathfrak{a}\wedge *F_A) + \mathrm{Ric}((\cdot)\otimes\mathfrak{a}) \ . \tag{2.10}$$

Note in particular that $q_A(\mathfrak{a})$ is the expression on the right hand side of (2.1).

Let $\{(A_n, \hat{a}_n)\}_{n\in 1,2,..}$ denote the renumbered subsequence given by Lemma 2.1. By way of a look ahead, the subsequent analysis of Lemma 2.1's limit function $|\hat{a}|$ requires a uniform bound for the sequence $\{\|q_{A_n}(\hat{a}_n)\|_2\}_{n=1,2,\ldots}$. The sad fact is that Lemma 2.1's assumptions are not strong enough to guarantee such a bound. The next proposition circumvents this problem. This proposition uses $\|\cdot\|_2$ to again denote the $L^2$ norm of an indicated tensor field valued section of $P\times_{SO(3)}\mathfrak{su}(2)$.

**Proposition 2.2**: *Suppose that* $\{\mathbb{A}_n = A_n + i\,\mathfrak{a}_n\}_{n=1,2,\ldots}$ *is a sequence of connections on* $P\times_{SO(3)}PSL(2;\mathbb{C})$ *with* $\{\mathfrak{F}(\mathbb{A}_n)\}_{n=1,2,\ldots}$ *being bounded. For each* $n\in\{1, 2, \ldots\}$, *let* $r_n$ *denote the* $L^2$ *norm of* $\mathfrak{a}_n$ *and assume that* $\{r_n\}_{n=1,2,\ldots}$ *is divergent with no finite limits. There exists a number* $\kappa > 1$ *that depends only on the upper bound for the sequence* $\{\mathfrak{F}(\mathbb{A}_n)\}_{n=1,2,\ldots}$, *and there exists a sequence* $\{\hat{a}_n\}_{n=1,2\ldots}$ *of sections of* $P\times_{SO(3)}\mathfrak{su}(2))\otimes T^*M$ *such that each* $n\in\{1, 2, \ldots\}$ *version of Items a)-e) below holds.*

a) $\|\nabla_{A_n}(\hat{a}_n - r_n^{-1}\mathfrak{a}_n)\|_2 + \kappa^2 r_n\|\hat{a}_n - r_n^{-1}\mathfrak{a}_n\|_2 < \kappa r_n^{-1}$ .

b) $\int_M (|\nabla_{A_n}\hat{a}_n|^2 + r_n^2|\hat{a}_n\wedge\hat{a}_n|^2 + r_n^{-2}|F_{A_n}|^2 + \mathrm{Ric}(\langle\hat{a}_n\otimes\hat{a}_n\rangle)) < \kappa r_n^{-2}$ .

c) $\|d_{A_n}\hat{a}_n\|_2^2 + \|d_{A_n}*\hat{a}_n\|_2^2 + r_n^{-2}\|F_{A_n} - r_n^2\hat{a}_n\wedge\hat{a}_n\|_2^2 < \kappa r_n^{-2}$ .

d) $\|q_{A_n}(\hat{a}_n)\|_2 < \kappa$ .

e) $\sup_M|\hat{a}_n| < \kappa$ .

*Moreover, there is a subsequence* $\Lambda\subset\{1, 2, \ldots\}$ *with the properties listed below.*

- *The sequence* $\{|\hat{a}_n|\}_{n\in\Lambda}$ *is bounded in* $L^2_1$, *it converges weakly in the* $L^2_1$ *topology and strongly in* $L^q$ *topology for* $q < 6$. *No member vanishes on an open set in* M.
- *The* $L^2$ *limit of* $\{|\hat{a}_n|\}_{n\in\Lambda}$ *is in* $L^\infty$. *The limit is denoted in what follows by* $|\hat{a}_\diamond|$. *This function is defined pointwise by the rule* $|\hat{a}_\diamond|(p) = \lim\sup_{n\in\Lambda}|\hat{a}_n|(p)$.
- *The sequence* $\{\langle\hat{a}_n\otimes\hat{a}_n\rangle\}_{n\in\Lambda}$ *converges strongly in any* $q < 3$ *version of the* $L^q$ *topology on the space of sections of* $(P\times_{SO(3)}\mathfrak{su}(2))\otimes T^*M$, *and it converges weakly in the* $L^3$ *topology. The limit section is* $\langle\hat{a}_\diamond\otimes\hat{a}_\diamond\rangle$.

- *Let $f$ denote a given $C^2$ function. The three sequences* $\{\int_M f\,|\nabla_{A_n}\hat{a}_n|^2\}_{n\in\Lambda}$, $\{\int_M r_n^{\,2} f\,|\hat{a}_n\wedge\hat{a}_n|^2\}_{n\in\Lambda}$ *and* $\{\int_M r_n^{\,-2} f\,|F_{A_n}|^2\}_{n\in\Lambda}$ *converge with respective limits that are denoted in what follows by* $Q_{\nabla,f}$, $Q_{\wedge,f}$ *and* $Q_{F,f}$. *These are such that* $Q_{\wedge,f} = Q_{F,f}$ *and any* $f\geq 0$ *version of* $Q_{\nabla,f}$ *is no less than* $\int_M f\,|d|\hat{\mathrm{a}}_\Diamond||^2$ . *Moreover,*

$$\tfrac{1}{2}\int_M d^*df\,|\hat{\mathrm{a}}_\Diamond|^2 + Q_{\nabla,f} + 2Q_{\wedge,f} + \int_M f\,\mathrm{Ric}(\langle\hat{\mathrm{a}}_\Diamond\otimes\hat{\mathrm{a}}_\Diamond\rangle) = 0.$$

- *The sequence* $\{q_{A_n}(\hat{a}_n)\}_{n\in\Lambda}$ *has a weak limit in the $L^2$ topology. Moreover, if $f$ is any $L^2$ function, then* $\lim_{n\in\Lambda}\int_M f\langle\hat{a}_n\wedge *q_{A_n}(\hat{a}_n)\rangle = 0$.
- *Fix* $p\in M$. *The two sequences indexed by* $\Lambda$ *with respective n'th terms given by the integral of* $G_p|\nabla_{A_n}\hat{a}_n|^2$ *and the integral of* $G_p r_n^{\,2}|\hat{a}_n\wedge\hat{a}_n|^2$ *are less than* $\kappa$. *Introduce by way of notation* $Q_{\Diamond,p}$ *to denote the lim-inf of the sequence with n'th term the integral of* $G_p(|\nabla_{A_n}\hat{a}_n|^2 + 2r_n^{\,2}|\hat{a}_n\wedge\hat{a}_n|^2)$. *The function* $|\hat{\mathrm{a}}_\Diamond|^2$ *obeys the equation*

$$\tfrac{1}{2}|\hat{\mathrm{a}}_\Diamond|^2(p) + Q_{\Diamond,p} = -\int_M G_p\,(\tfrac{1}{2}|\hat{\mathrm{a}}_\Diamond|^2 - \mathrm{Ric}(\langle\hat{\mathrm{a}}_\Diamond\otimes\hat{\mathrm{a}}_\Diamond\rangle))\ .$$

The proof of Proposition 2.2 requires a preliminary lemma that directly asserts a part of the proposition.

**Lemma 2.3**: *Let* $\{(A_n, \hat{\mathrm{a}}_n)\}_{n\in 1,2,..}$ *denote the renumbered subsequence from Lemma 2.1. There exists* $\kappa > 1$ *that depends only on the upper bound for the sequence* $\{\mathfrak{F}(\mathbb{A}_n)\}_{n=1,2,\ldots}$ *and, given* $z > \kappa$, *there exists a sequence* $\{\hat{a}_n\}_{n=1,2\ldots}$ *of sections of* $(P\times_{SO(3)}\mathfrak{su}(2))\otimes T^*M$ *such that for each* $n\in\{1, 2, \ldots\}$,

- $\|\hat{a}_n - \hat{\mathrm{a}}_n\|_2 \leq z^{-1/2} r_n^{\,-2}$ .
- $\|d_{A_n}\hat{a}_n\|_2^{\,2} + \|d_{A_n}*\hat{a}_n\|_2^{\,2} \leq \kappa\, r_n^{\,-2}$ .
- $\|q_{A_n}(\hat{a}_n)\|_2 \leq \kappa z$ .

This lemma is proved in the next subsection. Accept as gospel truth for the moment.

***Proof of Proposition 2.2***: The proof has three parts. Fix $z \geq c_0$ so as to invoke Lemma 2.3. As is evident in the proof, a large choice for $z$, but in any event less than $c_0$ suffices to prove the assertions of the proposition. This said, view $z$ for now as a chosen parameter. By way of notation, the proof introduces $\|\cdot\|_q$ to denote a given $q\in[1,\infty]$ version of the $L^q$ norm on tensor valued sections of $P\times_{SU(2)}\mathfrak{su}(2)$. The proof also denotes any given $n\in\{1, 2, \ldots\}$ version of $r_n^{-1}\mathfrak{a}_n$ by $\hat{\mathrm{a}}_n$.

*Part 1*: This first part proves Item e) of Proposition 2.2 and the assertion in the proposition's sixth bullet that the integrals of each $n\in\Lambda$ version of $G_p|\nabla_{A_n}\hat{a}_n|^2$ and $G_p r_n^2|\hat{a}_n\wedge\hat{a}_n|^2$ has a $p\in M$ and index n independent upper bound. There are four steps.

<u>Step 1</u>: This step proves that the $L^2$ norms of $\nabla_{A_n}\hat{a}_n$ and $r_n\hat{a}_n\wedge\hat{a}_n$ are both bounded by $c_0$. To start, integrate by parts to rewrite the second bullet of Lemma 2.3 as

$$\int_M (|\nabla_{A_n}\hat{a}_n|^2 + 2\langle *F_{A_n}\wedge\hat{a}_n\wedge\hat{a}_n\rangle + \mathrm{Ric}(\langle\hat{a}_n\otimes\hat{a}_n\rangle)) \le \kappa r_n^{-2}. \qquad (2.11)$$

To exploit (2.11), use the first bullets of Lemmas 2.1 and 2.3 to see that $\|\hat{a}_n\|_2 = 1+\mathfrak{e}$ with $|\mathfrak{e}|\le c_0 z^{-1/2} r_n^{-2}$. Next, use the $f = 1$ version of Lemma 2.1's third bullet to bound the integral of $|\langle *F_{A_n}\wedge\hat{\mathrm{a}}_n\wedge\hat{\mathrm{a}}_n\rangle|$ by $c_0$. Then write $\hat{a}_n$ as $\hat{\mathrm{a}}_n+(\hat{a}_n-\hat{\mathrm{a}}_n)$ and use the first bullet of Lemma 2.3 to go from the preceding bound to a bound on the integral of $|\langle *F_{A_n}\wedge\hat{a}_n\wedge\hat{a}_n\rangle|$ by $c_0 + c_0\|F_{A_n}\|_2\|\hat{a}_n-\hat{\mathrm{a}}_n\|_6^{1/2}(\|\hat{a}_n\|_6+\|\hat{\mathrm{a}}_n\|_6)^{3/2}$. Note that $\|F_{A_n}\|_2\|\hat{a}_n-\hat{\mathrm{a}}_n\|_6^{1/2}\le c_0$ by Lemmas 2.1 and 2.3; and $\|\hat{\mathrm{a}}_n\|_6\le c_0$ due to Lemma 2.1. The $L^6$ norm of $\hat{a}_n$ is bounded by $c_0(\|\nabla_{A_n}\hat{a}_n\|_2+\|\hat{a}_n\|_2)$ using (2.6) and the Sobolev inequality that bounds a function's $L^6$ norm by $c_0$ times its $L^2_1$ norm. These bounds give a $c_0(1+\|\nabla_{A_n}\hat{a}_n\|_2^{3/2})$ bound for the integral of $|\langle *F_{A_n}\wedge\hat{a}_n\wedge\hat{a}_n\rangle|$ and thus (2.11) gives a $c_0$ bound for $\|\nabla_{A_n}\hat{a}_n\|_2$. Note for use below that this leads back to a $c_0$ bound $\|\hat{a}_n\|_6$.

Use the $c_0$ bounds for the $L^6$ norm of $\hat{a}_n$ and the $L^2$ norm of $F_{A_n} - r_n^2\hat{\mathrm{a}}_n\wedge\hat{\mathrm{a}}_n$ in (2.11) to see that

$$\int_M (|\nabla_{A_n}\hat{a}_n|^2 + 2r_n^{2}\langle *(\hat{\mathrm{a}}_n\wedge\hat{\mathrm{a}}_n)\wedge\hat{a}_n\wedge\hat{a}_n\rangle) \le c_0 . \qquad (2.12)$$

Write $\hat{\mathrm{a}}_n = \hat{\mathrm{a}}_n - \hat{a}_n + \hat{a}_n$ and use this decomposition in (2.12) to see that

$$\int_M (|\nabla_{A_n}\hat{a}_n|^2 + 2r_n^{2}|\hat{a}_n\wedge\hat{a}_n|^2) \le c_0(1 + r_n^2\|\hat{a}_n-\hat{\mathrm{a}}_n\|_2(\|\hat{\mathrm{a}}_n\|_6+\|\hat{a}_n\|_6)\|\hat{a}_n\|_6^{2}) . \qquad (2.13)$$

Lemma 2.3's first bullet, the $c_0$ bounds for the $L^6$ norms of $\hat{\mathrm{a}}_n$ and $\hat{a}_n$ and (2.13) imply that

$$\int_M (|\nabla_{A_n}\hat{a}_n|^2 + 2\,r_n^{\,2}\,|\hat{a}_n\wedge\hat{a}_n|^2) \le c_0 \ .$$

(2.14)

This last inequality implies that $\|\nabla_{A_n}\hat{a}_n\|_2^{\,2}$ and $r_n^{\,2}\|\hat{a}_n\wedge\hat{a}_n\|_2^{\,2}$ are both bounded by $c_0$.

<u>Step 2</u>: Let $(A, \mathfrak{a})$ denote for the moment a given pair of connection on P and section of $(P\times_{SO(3)}\mathfrak{su}(2))\otimes T^*M$. Let $f$ denote a given $C^2$ function. Integrate by parts in the right most integral on its right hand side of (2.2) and then use (2.1) with the definition of $q_A(\mathfrak{a})$ to obtain the equality

$$\tfrac{1}{2}\int_M d^\dagger df\,|\mathfrak{a}|^2 + \int_M f(|\nabla_A\mathfrak{a}|^2 + 2\langle *F_A\wedge\mathfrak{a}\wedge\mathfrak{a}\rangle + \mathrm{Ric}(\langle\mathfrak{a}\otimes\mathfrak{a}\rangle)) - \int_M f\,\langle\mathfrak{a}\wedge *q_A(\mathfrak{a})\rangle = 0.$$

(2.15)

Since A and $\mathfrak{a}$ are smooth, this equality also holds when $f$ is such that $d^\dagger df$ is a distribution. In particular, it holds with $f$ being the Green's function of $d^\dagger d + 1$ with pole at any given point. This understood, fix $p \in M$ and let $G_p$ again denote the Green's function for $d^\dagger d + 1$ with pole at p. The corresponding version of (2.11) reads

$$\tfrac{1}{2}|\mathfrak{a}|^2(p) + \int_M G_p(|\nabla_A\mathfrak{a}|^2 + 2\langle *F_A\wedge\mathfrak{a}\wedge\mathfrak{a}\rangle) = \int_M G_p(\tfrac{1}{2}|\mathfrak{a}|^2 - \mathrm{Ric}(\langle\mathfrak{a}\otimes\mathfrak{a}\rangle) + \langle\mathfrak{a}\wedge *q_A(\mathfrak{a})\rangle)\,.$$

(2.16)

A bound for the right hand side of this inequality can be had by using the fact that $G_p(\cdot)$ is bounded by $c_0\,\mathrm{dist}(p,\cdot)^{-1}$ and thus has bounded $L^q$ norm for any $q < 3$. This being the case, the left most two terms on the right hand side of (2.16) contribute at most $c_0\|\mathfrak{a}\|_4^{\,2}$ to the absolute value of the right hand side. What with (2.6), a dimension 3 Sobolev inequality bounds this by $c_0(\|\nabla_A\mathfrak{a}\|_2^{\,2} + \|\mathfrak{a}\|_2^{\,2})$.

Meanwhile, the term with $q_A(\mathfrak{a})$ in (2.16) contributes at most

$$c_0\,(\sup_{p\in M}\|\mathrm{dist}(p,\cdot)^{-1}\,\mathfrak{a}\|_2)\,\|q_A(\mathfrak{a})\|_2$$

(2.17)

to the absolute value of the right hand side of (2.16). As explained in the next paragraph, $\sup_{p\in M}\|\mathrm{dist}(p,\cdot)^{-1}\,\mathfrak{a}\|_2$ is no greater than $c_0(\|\nabla_A\mathfrak{a}\|_2 + \|\mathfrak{a}\|_2)$. Granted such a bound, then the absolute value of the right hand side of (2.16) is no greater than

$$c_0(\|\nabla_A\mathfrak{a}\|_2^{\,2} + \|\mathfrak{a}\|_2^{\,2} + \|q_A(\mathfrak{a})\|_2^{\,2}).$$

(2.18)

The assertion in the preceding paragraph about the supremum in (2.17) invokes Hardy's inequality: Let $f$ denote any given $L^2_1$ function on M. Then

$$\sup_{p\in M}\int_M \frac{1}{\mathrm{dist}(p,\cdot)^2} f^2 \ \le c_0\,(\|df\|_2^{\,2} + \|f\|_2^{\,2})$$

(2.19)

The latter and (2.6) imply the asserted bound for $\sup_{p\in M}\|\mathrm{dist}(p,\cdot)^{-1}\,\mathfrak{a}\|_2$.

Step 3: Fix $n \in \{1, 2, \ldots\}$. Use the first and third bullets of Lemma 2.3 and the bound for $\|\nabla_{A_n}\hat{a}_n\|_2$ from Step 1 to bound the $(A = A_n, \mathfrak{a} = \hat{a}_n)$ version of (2.18) by $c_0 z^2$ when $z > c_0$. Given this bound, it then follows that the absolute value of the right hand side of the $(A = A_n, \mathfrak{a} = \hat{a}_n)$ version of (2.16) is also bounded by $c_0 z^2$.

The integral of $G_p|\nabla_{A_n}\hat{a}_n|^2$ that appears on the left hand side of the $(A = A_n, \mathfrak{a} = \hat{a}_n)$ version of (2.16) is nonnegative because $G_p$ is nonnegative; and this is all that need be said about this integral for the time being. The other integral on the left hand side of this same version of (2.16) is that of $2G_p\langle *F_{A_n} \wedge \hat{a}_n \wedge \hat{a}_n\rangle$. This integral is not manifestly nonnegative and so more needs to be said about it. To start the story on this integral, write $*F_{A_n}$ as the sum of two terms, these being $*F_{A_n} - r_n^{\,2}\hat{\mathrm{a}}_n \wedge \hat{\mathrm{a}}_n$ and $r_n^{\,2}\hat{\mathrm{a}}_n \wedge \hat{\mathrm{a}}_n$. Use this decomposition to write

$$\int_M G_p\left\langle *F_{A_n} \wedge \hat{a}_n \wedge \hat{a}_n\right\rangle = r_n^{\,2}\int_M G_p\left\langle *(\hat{\mathrm{a}}_n \wedge \hat{\mathrm{a}}_n) \wedge (\hat{a}_n \wedge \hat{a}_n)\right\rangle + \mathfrak{e}\ ,$$

(2.20)

with $\mathfrak{e}$ being the contribution from $*F_{A_n} - r_n^{\,2}\hat{\mathrm{a}}_n \wedge \hat{\mathrm{a}}_n$. The bound by $c_0$ on the latter's $L^2$ norm leads to the bound $|\mathfrak{e}| \le c_0\|\mathrm{dist}(p,\cdot)^{-1}\hat{a}_n\|_2\|\hat{a}_n\|_\infty$. This last bound with (2.19), (2.16) and Step 1's bound for $\|\nabla_{A_n}\hat{a}_n\|_2$ implies that $|\mathfrak{e}| \le c_0\|\hat{a}_n\|_\infty$.

To see about the integral on the right hand side of (2.20), write $\hat{\mathrm{a}}_n \wedge \hat{\mathrm{a}}_n$ as a sum of two terms, these being $\hat{a}_n \wedge \hat{a}_n$ and $\hat{\mathrm{a}}_n \wedge \hat{\mathrm{a}}_n - \hat{a}_n \wedge \hat{a}_n$. Use this splitting to write (2.20) as

$$\int_M G_p\left\langle *F_{A_n} \wedge \hat{a}_n \wedge \hat{a}_n\right\rangle = r_n^{\,2}\int_M G_p\,|\hat{a}_n \wedge \hat{a}_n|^2 \ + \mathfrak{e}' + \mathfrak{e}\ ,$$

(2.21)

with $\mathfrak{e}'$ being the contribution from $\hat{\mathrm{a}}_n \wedge \hat{\mathrm{a}}_n - \hat{a}_n \wedge \hat{a}_n$. Of particular note is that

$$|\mathfrak{e}'| \le c_0\, r_n^{\,2}\, \|\hat{a}_n - \hat{\mathrm{a}}_n\|_2\, (\|\mathrm{dist}(p,\cdot)^{-1}\hat{\mathrm{a}}_n\|_2 + \|\mathrm{dist}(p,\cdot)^{-1} a_n\|_2)\,\|\hat{a}_n\|_\infty^{\,2}\ ,$$

(2.22)

and thus $|\mathfrak{e}'|$ is no greater that $c_0 z^{-1/2}\|\hat{a}_n\|_\infty^{\,2}$. This bound follows from the first bullet of Lemma 2.3 using what was said already about the integrals that involve $\mathrm{dist}(p,\cdot)^{-1}$.

Step 4: Use the bounds on $|\mathfrak{e}|$ and $|\mathfrak{e}'|$ in (2.21) with the bound in Step 2 for the right hand side of (2.16) to see that the latter equation implies the bound

$$(1 - c_0 z^{-1/2}) \| \hat{a}_n \|_\infty^2 + \sup_{p \in M} \int_M G_p (|\nabla_{A_n} \hat{a}_n|^2 + r^2 |\hat{a}_n \wedge \hat{a}_n|^2) \leq c_0 z^2 . \tag{2.23}$$

Any $z > c_0$ version of (2.23) supplies an index n independent bound for the $\| \hat{a}_n \|_\infty$ and an index n and $p \in M$ independent bound for the integral of $G_p |\nabla_{A_n} \hat{a}_n|^2$ and $G_p r_n^2 |\hat{a}_n \wedge \hat{a}_n|^2$.

*Part 2*: This part of the subsection proves Items a), b) and c) of Proposition 2.2. This is done in three steps.

Step 1: This step explains why Item b) follows from Item c) and the first bullet of Lemma 2.3. To start, fix $n \in \{1, 2, \ldots\}$ for the moment and define $\mathcal{A}_n$ to be the connection $\mathcal{A}_n = A_n + i r_n \hat{a}_n$, this being a connection on $(P \times_{SO(3)} PSL(2; \mathbb{C})$. The assertion in Item c) is equivalent to the assertion that $\mathfrak{F}(\mathcal{A}_n) \leq c_0$. Meanwhile, the first bullet of Lemma 2.3 implies that

$$\| \hat{a}_n \|_2 = 1 + \mathfrak{e} \quad \textit{with} \quad |\mathfrak{e}| \leq c_0 r_n^{-2}. \tag{2.24}$$

Assuming the bounds in the preceding paragraph, then the sequence $\{\mathcal{A}_n\}_{n=1,2,\ldots}$ can be used in lieu of $\{\mathbb{A}_n\}_{n=1,2,\ldots}$ as input for Lemma 2.1. Item b) of Proposition 2.2 follow directly from the second bullet of the $\{\mathcal{A}_n\}_{n=1,2,\ldots}$ version of Lemma 2.1.

Step 2: This step and Step 3 prove Items a) and c) in tandem. Note in this regard that it suffices to prove that both $\| \nabla_{A_n} (\hat{a}_n - \hat{\mathrm{a}}_n) \|_2^2$ and $r_n^{-2} \| F_{A_n} - r_n^2 \hat{a}_n \wedge \hat{a}_n \|_2^2$ are at most $c_0 r_n^{-2}$. The fact that Items a) and c) follow from bounds of this sort follows from the $z \geq c_0$ versions of the first and second bullets of Lemma 2.3.

To start the proof, use the second bullet of Lemma 2.3 with the triangle inequality to see that

$$\| d_{A_n} (\hat{a}_n - \hat{\mathrm{a}}_n) \|_2^2 + \| d_{A_n} * (\hat{a}_n - \hat{\mathrm{a}}_n) \|_2^2 \leq c_0 r_n^{-2} . \tag{2.25}$$

Integration by parts leads from (2.25) to the inequality

$$\| \nabla_{A_n} (\hat{a}_n - \hat{\mathrm{a}}_n) \|_2^2 + 2 \int_M \langle *F_{A_n} \wedge (\hat{a}_n - \hat{\mathrm{a}}_n) \wedge (\hat{a}_n - \hat{\mathrm{a}}_n) \rangle + \int_M \mathrm{Ric}(\langle (\hat{a}_n - \hat{\mathrm{a}}_n) \otimes (\hat{a}_n - \hat{\mathrm{a}}_n) \rangle) \leq c_0 r_n^{-2} . \tag{2.26}$$

The absolute value of the integral with Ricci curvature tensor is bounded by $c_0 r_n^{-4}$, this being a consequence of the first bullet of Lemma 2.3. To see about the middle integral on the left hand side of (2.25), write $F_{A_n}$ as $(F_{A_n} - r_n^2 \hat{a}_n \wedge \hat{a}_n) + r_n^2 \hat{a}_n \wedge \hat{a}_n$. The contribution to the middle integral on the left hand side of (2.25) from $(F_{A_n} - r_n^2 \hat{a}_n \wedge \hat{a}_n)$ is no greater than $2\|F_{A_n} - r_n^2 \hat{a}_n \wedge \hat{a}_n\|_2 \|\hat{a}_n - \hat{\mathrm{a}}_n\|_4^2$ . This, in turn, is no greater than

$$c_0\|F_{A_n} - r_n^2 \hat{a}_n \wedge \hat{a}_n\|_2 \ \|\hat{a}_n - \hat{\mathrm{a}}_n\|_2^{1/2} \ \|\hat{a}_n - \hat{\mathrm{a}}_n\|_6^{3/2} . \tag{2.27}$$

To say more about the expression in (2.27), use the first bullet of Lemma 2.3 to bound it by $c_0$ times the product $z^{-1/2}(r_n^{-1}\|F_{A_n} - r_n^2 \hat{a}_n \wedge \hat{a}_n\|_2)\|\hat{a}_n - \hat{\mathrm{a}}_n\|_6^{3/2}$. Meanwhile,

$$\|\hat{a}_n - \hat{\mathrm{a}}_n\|_6^{3/2} \le c_0 (\|\nabla_{A_n} (\hat{a}_n - \hat{\mathrm{a}}_n)\|_2 + r_n^{-1}) . \tag{2.28}$$

To derive this, note that $\|\hat{a}_n - \hat{\mathrm{a}}_n\|_6^{3/2} \le \|\hat{a}_n - \hat{\mathrm{a}}_n\|_6 (\|\hat{a}_n\|_6 + \|\hat{\mathrm{a}}_n\|_6)^{1/2}$. This understood, then (2.28) follows from (2.6) and a standard Sobolev inequality given that the $L^2$ norms of $\hat{a}_n$ and $\hat{\mathrm{a}}_n$ and their respective $A_n$-covariant derivatives are bounded by $c_0$. Such a bound for the $A_n$-covariant derivative of $\hat{a}_n$ is supplied by Step 1 of Part 1; and that of $\hat{\mathrm{a}}_n$ is supplied by Lemma 2.1.

Use what is said in the preceding two paragaphs with (2.26) to see that

$$\|\nabla_{A_n} (\hat{a}_n - \hat{\mathrm{a}}_n)\|_2^2 \le c_0 z^{-1/2}\|F_{A_n} - r_n^2 \hat{a}_n \wedge \hat{a}_n\|_2 (\|\nabla_{A_n} (\hat{a}_n - \hat{\mathrm{a}}_n)\|_2 + r_n^{-1}) + c_0 r_n^{-2} . \tag{2.29}$$

This last inequality is invoked in at the very end of Step 3.

<u>Step</u> <u>3</u>: The $L^2$ norm of $r_n^{-1}(F_{A_n} - r_n^2 \hat{a}_n \wedge \hat{a}_n)$ is no greater than the sum of the $L^2$ norms of $r_n^{-1}(F_{A_n} - r_n^2 \hat{\mathrm{a}}_n \wedge \hat{\mathrm{a}}_n)$ and $r_n(\hat{\mathrm{a}}_n \wedge \hat{\mathrm{a}}_n - \hat{a}_n \wedge \hat{a}_n)$. The $L^2$ norm of the former is bounded by $c_0 r_n^{-1}$. Meanwhile, that of the latter is no greater than the sum of the $L^2$ norms of $r_n(\hat{a}_n - \hat{\mathrm{a}}_n) \wedge (\hat{a}_n - \hat{\mathrm{a}}_n)$ and twice that of $r_n(\hat{\mathrm{a}}_n - \hat{a}_n) \wedge \hat{a}_n$.

The preceding observations imply directly that

$$r_n^{-2}\|F_{A_n} - r_n^2 \hat{a}_n \wedge \hat{a}_n\|_2^2 \le c_0 (r_n^2\|\hat{a}_n - \hat{\mathrm{a}}_n\|_4^4 + r_n^2 \|\hat{a}_n\|_\infty^2 \|\hat{a}_n - \hat{\mathrm{a}}_n\|_2^2 + r_n^{-2}). \tag{2.30}$$

Part 1 bounds $\|\hat{a}_n\|_\infty$ by $c_0$ and Lemma 2.3 bounds $\|\hat{a}_n - \hat{\mathrm{a}}_n\|_2$ by $c_0 r_n^{-2}$. This leads to a bound on the right hand side of (2.29) by $c_0 r_n^2 \|\hat{a}_n - \hat{\mathrm{a}}_n\|_4^4 + c_0 r_n^{-2}$. The latter is no greater than $c_0(r_n^2 \|\hat{a}_n - \hat{a}_n\|_2 \|\hat{a}_n - \hat{a}_n\|_6^3 + r_n^{-2})$. This understood, use the first bullet of Lemma 2.3 and the bound on $\|\hat{a}_n - \hat{a}_n\|_6$ given by (2.28) to see that

$$r_n^{-2} \| F_{A_n} - r_n^2 \hat{a}_n \wedge \hat{a}_n \|^2 \le c_0 z^{-1} \| \nabla_{A_n} (\hat{a}_n - \hat{a}_n) \|_2^2 + c_0 r_n^{-2} . \tag{2.31}$$

Taken together, the inequalities in (2.29) and (2.31) imply that

$$(1 - c_0 z^{-1})( \| \nabla_{A_n} (\hat{a}_n - \hat{a}_n) \|_2^2 + r_n^{-2} \| F_{A_n} - r_n^2 \hat{a}_n \wedge \hat{a}_n \|_2^2 \le c_0 r_n^{-2} . \tag{2.32}$$

This inequality supplies the desired bounds if $z > c_0$.

*Part 3*: The nine steps that follow in this part of the subsection prove the bulleted assertions of Proposition 2.2.

Step 1: Reintroduce the sequence $\{\mathcal{A}_n = A_n + i\, r_n \hat{a}_n\}_{n=1,2,\ldots}$ from Step 1 of Part 2, this being a sequence of connections on $P \times_{SO(3)} PSL(2;\mathbb{C})$. As noted therein, the corresponding sequence $\{\mathfrak{F}(\mathcal{A}_n)\}_{n=1,2,\ldots}$ is bounded and so it can be used in lieu of $\{\mathbb{A}_n\}_{n=1,2,\ldots}$ as input for Lemma 2.1.

Except for the assertion about vanishing on an open set, what is said by the first bullet of Proposition 2.2 constitutes a part of the first bullet the $\{\mathcal{A}_n\}_{n=1,2,\ldots}$ version of Lemma 2.1. Note that the limit $L^2_1$ function has $L^2$ norm equal to 1. This is because the top bullet of Lemma 2.3 implies that the $L^2$ norm of each $n \in \{1, 2, \ldots\}$ version of $\hat{a}_n$ differs from 1 by at most $z^{-1/2} r_n^{-2}$. The nonvanishing condition can readily be satisfied by making a very small perturbation of a sequence that has all of the other properties that are required by the proposition. This understood, no more will be said about the nonvanising on an open set requirement.

The first bullet of the $\{\mathcal{A}_n\}_{n=1,2,\ldots}$ version of Lemma 2.1 asserts what is said in the second bullet of Proposition 2.2 to the effect that the limit $L^2_1$ function of the sequence $\{|\hat{a}_n|\}_{n=1,2,\ldots}$ is an $L^\infty$ function. In fact, it follows from Item a) of Proposition 2.2 that the limit functions of $\{|\hat{a}_n|\}_{n=1,2,\ldots}$ and $\{|\hat{a}_n|\}_{n=1,2,\ldots}$ have the same weak $L^2_1$ limits with it understood that the subsequences that are chosen for the respective $\{\mathcal{A}_n\}_{n=1,2,\ldots}$ and $\{\mathbb{A}_n\}_{n=1,2,\ldots}$ versions of Lemma 2.1 are identical.

Step 2. The definition of the function $|\hat{a}_\Diamond|$ by the second bullet of Proposition 2.2 raises a subtle point, this being the distinction between elements in $L^2_1 \cap L^\infty$ and functions that are defined pointwise. An element in the Banach space of $L^2_1 \cap L^\infty$ is an *equivalence class* of functions that are defined almost everywhere with two functions being equivalent if they agree on the complement of a set of measure zero. The distinction between a pointwise defined function and an equivalence class of functions that differ on sets with measure zero is at issue with regards to the definition in the second bullet of Proposition 2.2 of $|\hat{a}_\Diamond|$. In particular, this bullet of Proposition 2.2 *defines* an honest function, $|\hat{a}_\Diamond|$; and in so doing, this bullet makes the following implicit assertion:

*The function defined by the rule* $p \to \limsup_{n\to\infty} |\hat{a}_n|(p)$ *is in the equivalence class of the* $L^2_1$ *limit of the sequence* $\{|\hat{a}_n|\}_{n=1,2,\ldots}$.

(2.33)

The proof of the second bullet of Proposition 2.2 requires a proof of (2.33).

The distinction between a pointwise defined function and an equivalence class of functions that agree on the complement of a measure zero set is also at issue with regards to the proof of the sixth bullet of Proposition 2.2. The assertion in (2.33) and the fifth bullet of Proposition 2.2 are proved simultaneously in Steps 4-8. In the meantime, let $|\hat{a}_\Diamond|$ denote a chosen function from the $L^2_1$ equivalence class of the weak limit of $\{|\hat{a}_n|\}_{n=1,2,\ldots}$.

<u>Step 3</u>: The assertions made by the third bullet of Proposition 2.2 follow because the sequence $\{\langle \hat{a}_n \otimes \hat{a}_n\rangle\}_{n\in\{1,2\ldots\}}$ is bounded in the $L^2_1$ topology; a priori bounds on the $L^2_1$ norms of its elements come via Items b) and e) of Proposition 2.2.

The proof of the assertions made by the fourth bullet of Proposition 2.2 are almost verbatim identical to those made by the third bullet of Lemma 2.1.

To prove the assertion made by the fifth bullet, let $f$ denote for the moment a smooth function. Fix $n \in \{1, 2, \ldots\}$. The right hand side (2.1) is $q_A(\mathfrak{a})$. This being the case, an integration by parts writes

$$\int_M f\left\langle \hat{a}_n \wedge * q_{A_n}(\hat{a}_n)\right\rangle = \int_M \left(f\,(|d_{A_n}\hat{a}_n|^2 + |d_{A_n} * \hat{a}_n|^2) + df\wedge *\left\langle \hat{a}_n(*d_{A_n} * \hat{a}_n) + \hat{a}_n \wedge * d_{A_n}\hat{a}_n\right\rangle\right).$$

(2.34)

Granted (2.34), then the second bullet of Lemma 2.3 and Item e) of Proposition 2.2 finds

$$\left|\int_M f\left\langle \hat{a}_n \wedge * q_{A_n}(\hat{a}_n)\right\rangle\right| \le c_0\left(\|f\|_\infty r_n^{-2} + \|df\|_2\, r_n^{-1}\right).$$

(2.35)

Let g denote a given $L^2$ function. Fix $\varepsilon > 0$ and let $f$ denote a smooth function such that $\|f - g\|_2 \le \varepsilon$. Then

$$|\int_M g\left\langle \hat{a}_n \wedge * q_{A_n}(\hat{a}_n)\right\rangle| \le c_0\,\varepsilon\,\|\hat{a}_n\|_\infty \|q_{A_n}(\hat{a}_n)\|_2 + |\int_M f\left\langle \hat{a}_n \wedge * q_{A_n}(\hat{a}_n)\right\rangle| \,. \tag{2.36}$$

Now use (2.35) with (2.36) and Item e) of Proposition 2.2 to conclude that

$$\int_M g\left\langle \hat{a}_n \wedge * q_{A_n}(\hat{a}_n)\right\rangle \le c_0\,(\varepsilon + r_n^{-2}\,\|f\|_\infty + r_n^{-1}\|df\|_2) \,. \tag{2.37}$$

Taking n sufficiently large bounds the right hand side of (2.37) by $\varepsilon$ and so (2.37) implies that

$$\lim_{n\in\Lambda} \int_M g\left\langle \hat{a}_\Diamond \wedge * q_{A_n}(\hat{a}_n)\right\rangle = 0. \tag{2.38}$$

<u>Step</u> <u>4</u>: This step with Steps 5-9 prove the assertion in (2.33) and the assertion made by the sixth bullet of Proposition 2.2. To start, fix $n \in \{1, 2, \ldots\}$ and a point $p \in M$ so as to consider the ($A = A_n$, $\mathfrak{a} = \hat{a}_n$) version of (2.16). Of particular concern in this step is the integral of $G_p\langle \hat{a}_n \wedge * q_{A_n}(\hat{a}_n)\rangle$ that appears on the right hand side of the $(A_n, \hat{a}_n)$ version. As explained directly, the various $p \in M$ versions are such that

$$\sup_{p\in M} |\int_M G_p\left\langle \hat{a}_n \wedge * q_{A_n}(\hat{a}_n)\right\rangle| < c_0\, r_n^{-1/5} \,. \tag{2.39}$$

To see that (2.39) holds, fix for the moment $p \in M$ and $\rho \in (0, 1)$. Having done so, let $\chi_{p,\rho}$ denote the function on M given by $\chi(2 - \rho^{-1}\,\mathrm{dist}(p,\cdot))$. This function equals 1 where the distance to p is greater than $2\rho$ and it equals 0 where the distance is less than $\rho$. Write $G_p$ as $(1 - \chi_{p,\rho})G_p + \chi_{p,\rho}G_p$ so as to split a given $n \in \{1, 2, \ldots\}$ version of the integral in (2.39) into two integrals.

The absolute value of the integral of $(1 - \chi_{p,\rho})\,G_p\langle \hat{a}_n \wedge * q_{A_n}(\hat{a}_n)\rangle$ is no greater than $c_0\,\rho^{1/2}\,\|\hat{a}_n\|_\infty \,\|q_{A_n}(\hat{a}_n)\|_2$, this because $G_p \le c_0\,\mathrm{dist}(p,\cdot)^{-1}$. Meanwhile, the absolute value of the integral of $\chi_{p,\rho}G_p\langle \hat{a}_n \wedge * q_{A_n}(\hat{a}_n)\rangle$ is no greater than $c_0\,\rho^{-2}\, r_n^{-1}$, this being a consequence of (2.35). Granted these bounds, take $\rho = r_n^{-2/5}$ and invoke Items d) and e) of Proposition 2.2 to conclude that the absolute value in (2.39) is no greater than $c_0\, r_n^{-1/5}$ as claimed.

Step 5: Fix $n \in \{1, 2, \ldots\}$ and $p \in M$ again so as to return to the $(A = A_n, \mathfrak{a} = \hat{a}_n)$ version of (2.16). Of particular concern here is the integral of $G_p\langle *F_{A_n} \wedge \hat{a}_n \wedge \hat{a}_n\rangle$ that appears on the left hand side of this version. As explained directly,

$$\int_M G_p \left\langle *F_{A_n} \wedge \hat{a}_n \wedge \hat{a}_n \right\rangle = \int_M G_p \, r_n^{\,2} \, |\hat{a}_n \wedge \hat{a}_n|^2 \; + \mathfrak{e} \;, \tag{2.40}$$

with $\mathfrak{e}$ having absolute value no greater than $c_0 \, r_n^{-2/3}$. To see why this is, decompose $F_{A_n}$ as the sum $r_n^2 \hat{a}_n \wedge \hat{a}_n + (F_{A_n} - r_n^2 \hat{a}_n \wedge \hat{a}_n)$ so as to decompose the integral on the left hand side of (2.40) as a sum of two integrals. The term designated as $\mathfrak{e}$ in (2.40) is the second of these two integrals. To bound $|\mathfrak{e}|$, fix for the moment $\rho > 0$ and again write $G_p$ as the sum $(1 - \chi_{p,\rho}) G_p + \chi_{p,\rho} G_p$. The absolute value of the contribution of $(1 - \chi_{p,\rho}) G_p$ to $\mathfrak{e}$ is no greater than $c_0 \rho^{-1} \| F_{A_n} - r_n^2 \hat{a}_n \wedge \hat{a}_n \|_2 \| \hat{a}_n \wedge \hat{a}_n \|_2$ because $G_p \le c_0 \, \mathrm{dist}(p, \cdot)^{-1}$. This understood, Items b) and c) of Proposition 2.2 bound this contribution by $c_0 \rho^{-1} r_n^{-2}$. Meanwhile, the absolute value of the contribution of $\chi_{p,\rho} G_p$ to $\mathfrak{e}$ is no greater than $c_0 \rho^{1/2} \| F_{A_n} - r_n^2 \hat{a}_n \wedge \hat{a}_n \|_2 \| \hat{a}_n \|_\infty^{\,2}$. Items c) and e) of Proposition 2.2 imply that this is no greater than $c_0 \rho^{1/2}$. Granted these bounds, take $\rho = r_n^{-4/3}$ to see that $|\mathfrak{e}| \le c_0 \, r_n^{-2/3}$.

Step 6: Fix once again $n \in \{1, 2, \ldots\}$ and $p \in M$. Use the $(A = A_n, \mathfrak{a} = \hat{a}_n)$ version of (2.16) with what is said in Steps 4 and 5 to see that

$$\tfrac{1}{2} |\hat{a}_n|^2(p) + \int_M G_p (|\nabla_{A_n} \hat{a}_n|^2 + 2 \, r_n^{\,2} \, |\hat{a}_n \wedge \hat{a}_n|^2) = - \int_M G_p (\tfrac{1}{2} |\hat{a}_n|^2 - \mathrm{Ric}(\langle \hat{a}_n \otimes \hat{a}_n \rangle)) + \mathfrak{e}_n \;, \tag{2.41}$$

with $\mathfrak{e}_n$ having absolute value no greater than $c_0 \, r_n^{-1/5}$. Let $I_n(p)$ denote the integral that appears on the right hand side of (2.41). Introduce by way of notation

$$I_\Diamond(p) = - \int_M G_p \, (\tfrac{1}{2} |\hat{a}_\Diamond|^2 - \mathrm{Ric}(\langle \hat{a}_\Diamond \otimes \hat{a}_\Diamond \rangle)) \;. \tag{2.42}$$

As explained directly, $\lim_{n\to\infty} \sup_{p \in M} |I_\Diamond(p) - I_n(p)| = 0$. That this is so follows from the first and third bullets of Proposition 2.2 with the fact that $G_p$ is square integrable.

Step 7: Introduce for the time being $\mathrm{v}$ to denote the function on M that is defined by the rule $p \to \mathrm{v}(p) = \limsup_{n\to\infty} |\hat{a}_n|(p)$. Let $Q_\Diamond$ denote the function on M that is defined by the sixth bullet of Proposition 2.2. As explained momentarily, the identity in (2.41) and the fact that $\lim_{n\to\infty} \sup_{p\in M} |I_\Diamond(p) - I_n(p)| = 0$ imply that

$$\tfrac{1}{2}\,\nu(\mathrm{p})^2 + Q_\lozenge(\mathrm{p}) = I_\lozenge(\mathrm{p})\ . \tag{2.43}$$

The final assertion of sixth bullet of Proposition 2.2 follows immediately from (2.43) if it is the case that $\nu = |\hat{a}_\lozenge|$ on the complement of a measure zero set.

To see about (2.43), fix $\Lambda \in \{1, 2, \ldots\}$ such that $\lim_{n\in\Lambda} |\hat{a}_n|(p) = \nu(p)$. Suppose that there exists $\delta > 0$ and a subsequence $\Lambda' \subset \Lambda$ such that

$$\int_M G_p\big(|\nabla_{A_n} \hat{a}_n|^2 + 2\, r_n^{\,2}\, |\hat{a}_n \wedge \hat{a}_n|^2\big) > Q_\lozenge(p) + 2\,\delta \tag{2.44}$$

when $n \in \Lambda'$. If this is the case, then the left hand side of (2.41) for $n \in \Lambda'$ will be greater than $\frac{1}{2}\nu(p)^2 + Q_\lozenge(p) + 2\delta$ when n is large, and this implies that $I_\lozenge(p)$ is no less than $\frac{1}{2}\nu(p)^2 + Q_\lozenge(p) + 2\delta$. Now fix a second subsequence, $\Theta \subset \{1, 2, \ldots\}$ such that $n \in \Theta$ versions of the integral on the left hand side of (2.45) converges to $Q_\lozenge(p)$. The $n \in \Theta$ versions of $|\hat{a}_n|(p)$ must in any event be less than $\nu(p) + \delta$ when n is large, and this implies that $I_\lozenge(p)$ is no greater than $\frac{1}{2}\nu(p)^2 + Q_\lozenge(p) + \delta$. This last conclusion is incompatible with the lower bound just stated for $I_\lozenge$.

<u>Step 8</u>: This step proves that $|\hat{a}_\lozenge| \ge \nu$ on a set of full measure. To start the proof, fix $\varepsilon > 0$ and reintroduce the functions $\delta_{(\cdot),\varepsilon}$ and $f_{(\cdot),\varepsilon}$ that are defined in Step 1 of the proof of Lemma 2.1. The function on M given by the rule

$$p \to \int_M \delta_{p,\varepsilon}\, |\hat{a}_\lozenge|^2 \tag{2.45}$$

is smooth. As noted in Step 3 of Lemma 2.1's proof, the function depicted in (2.45) converges as $\varepsilon \to 0$ to an $L^\infty$ function that equals $|\hat{a}_\lozenge|^2$ on the complement of a measure zero set. Of particular note is that the function defined by (2.45) obeys (2.8).

Write $f_{p,\varepsilon}$ as $g_\varepsilon + \mathfrak{e}_p$ with $g_\varepsilon$ defined using Gaussian coordinates centered at p by the formulas in (2.7). The functions $g_\varepsilon$ and $G_p$ are such that

$$g_\varepsilon \le G_p + \mathfrak{r}_p\ , \tag{2.46}$$

with $\mathfrak{r}_p$ such that $|\mathfrak{r}_p| \le c_0 |x|$. The inequality in (2.46) and the fact that $\{f_{p,\varepsilon}\}_{\varepsilon\in(0,1]}$ converges to $G_p$ as $\varepsilon \to 0$ in the $C^2$ topology on compact subsets of $M{-}p$ has the following consequence: Fix $\rho \in (0, c_0^{-1}]$ and there exists $c_\rho > 1$ that is independent of the point p and such that there is a continuous function on M, to be denoted by $\mathfrak{r}_{\mathfrak{p},\varepsilon,\rho}$, with norm less than $\rho$ and such that

$$f_{p,\varepsilon} \le G_p + \mathfrak{r}_{p,\varepsilon,\rho} \quad \textit{when } \varepsilon < c_\rho^{-1} .$$

(2.47)

Fix $\rho \in (0, c_0^{-1}]$ and $\varepsilon < c_\rho^{-1}$ Given $n \in \{1, 2, \ldots\}$ and $p \in M$, multiply both sides of (2.41) by $\delta_{p,\varepsilon}$, and integrate the result over M and use (2.47) to see that

$$\tfrac{1}{2} \int_M \delta_{p,\varepsilon} |\hat{a}_n|^2 + \int_M G_p (|\nabla_{A_n} \hat{a}_n|^2 + 2\, r_n^{\,2}\, |\hat{a}_n \wedge \hat{a}_n|^2) \ge I_\Diamond(p) - e_{n,\varepsilon,\rho} ,$$

(2.48)

where $e_{n,\varepsilon,\rho}$ is such that $\lim_{n\to\infty} \sup_{p\in M} |e_{n,\varepsilon,\rho}| \le c_0\, \rho$. Meanwhile, what is said in the third bullet of Proposition 2.2 has the following implication: Given $\varepsilon < c_\rho^{-1}$ and $\varepsilon' \in (0, 1]$, there exists $n_{\varepsilon,\varepsilon'} \ge 1$ such that

$$\sup_{p\in M} \Big| \int_M \delta_{p,\varepsilon} |\hat{a}_\Diamond|^2 - \int_M \delta_{p,\varepsilon} |\hat{a}_n|^2 \Big| < \varepsilon' \quad \textit{when } n > n_{\varepsilon,\varepsilon'} .$$

(2.49)

This follows from the fact that $\{|\hat{a}_n|\}_{n\in\{1,2,\ldots\}}$ converges strongly in the $L^2$ topology to $|\hat{a}_\Diamond|$.

Fix $p \in X$, take $\varepsilon' = \rho$ and then $n > n_{\varepsilon,\varepsilon'=\rho}$ with two additional constraints, both giving lower bounds: The first is that n's version of the $G_p$ integral on the right hand side of (2.48) differs from $Q_\Diamond(p)$ by at most $\rho$. The second is that $\sup_{p\in M} |e_{n,\varepsilon,\rho}| \le c_0\, \rho$. With n so chosen, invoke the $\varepsilon' = \rho$ version of (2.49) to see that

$$\tfrac{1}{2} \int_M \delta_{p,\varepsilon} |\hat{a}_\Diamond|^2 + Q_\Diamond(p) \ge I_\Diamond(p) - c_0 \rho \ .$$

(2.50)

Since $\rho$ and then $\varepsilon$ can be chosen as small as desired, and since this inequality holds for each $p \in X$ a comparison between (2.50) and (2.43) leads to the conclusion that $|\hat{a}_\Diamond| \ge \nu$ on the complement of a measure zero set.

<u>Step 9</u>: To see that $|\hat{a}_\Diamond| \le \nu$ on a set of full measure, fix $\varepsilon > 0$, $p \in M$ and a positive integer n. Having done so, multiply n's version of (2.41) by $\delta_{p,\varepsilon}$ and integrate the result over M. Do this for choices of $n > n_{\varepsilon,\varepsilon'=\varepsilon}$ that obey two additional lower bound constraints. The first constraint asks that n's version of the $G_p$ integral on the right hand side of (2.48) is greater than $Q_\Diamond(p) - \varepsilon$.; and the second asks that $\sup_{p\in M} |I_\Diamond(p) - I_n(p)| < \varepsilon$. With n so chosen, invoke the $\varepsilon' = \varepsilon$ version of (2.49) to see that

$$\frac{1}{2}\int_M \delta_{p,\varepsilon}|\hat{a}_\Diamond|^2 + Q_\Diamond(p) \le I_\Diamond(p) + c_0\varepsilon\ .$$
(2.51)

Since $\varepsilon$ can be as small as desired, a comparison between (2.51) and (2.43) leads to the conclusion that $|\hat{a}_\Diamond| \le v$ on the complement of a measure zero set.

**d) Proof of Lemma 2.3**

The proof has four parts. By way of a look ahead, any given $\hat{a}_n$ for $n \in \{1, 2, \ldots\}$ is a particular $t \in [0, \infty)$ value of the solution to a certain heat equation whose initial value is $\hat{a}_n$.

*Part 1*: Fix a pair $(A, \mathfrak{a})$ with A being a connection on P and $\mathfrak{a}$ being a section of $(P\times_{SO(3)}\mathfrak{su}(2))\otimes T^*M$. Standard existence and uniqueness theorems for parabolic differential equations prove that there is a unique section of $(P\times_{SO(3)}\mathfrak{su}(2))\otimes T^*M$ over $[0,\infty)\times M$, this denoted by $a$, that obeys the linear heat equation

$$\frac{\partial}{\partial t} a = -(\nabla_A{}^\dagger\nabla_A a + *(*F_A\wedge a + a\wedge *F_A) + \mathrm{Ric}(a) \quad \textit{with}\ a|_{t=0} = \mathfrak{a}\ .$$
(2.52)

This equation can be written equivalently in two ways,

$$\frac{\partial}{\partial t} a = -\, q_A(a) \quad \textit{and} \quad \frac{\partial}{\partial t} a = -(*d_A*d_A a - d_A*d_A*a)\ .$$
(2.53)

Use $\|\cdot\|_2$ as before to denote the $L^2$ norm on M. It follows from the right most identity in (2.53) that the function $[0,\infty)$ given by $\|a\|_2$ is non-increasing, and in particular, that

$$\frac{d}{dt}\|a\|_2{}^2 = -2\,(\|d_A a\|_2{}^2 + \|d_A*a\|_2{}^2).$$
(2.54)

Let $\mathfrak{E}$ denote the function on $[0, \infty)$ given by $\frac{1}{2}(\|d_A a\|_2{}^2 + \|d_A*a\|_2{}^2)$ It follows from the left most identity in (2.53) that

$$\frac{d}{dt}\mathfrak{E} = -\,\|q_A(a)\|_2{}^2\ .$$
(2.55)

The equations in (2.54) and (2.55) play central roles in Parts 2 and 3 of the proof.

*Part 2*: Define the function $\mathfrak{n}$ on $[0, \infty)$ by the rule

$$t \to \mathfrak{n}(t) = t^{-1} \int_0^t \|q_A(a)\|_2^2 \ .$$

(2.56)

Having done so, integrate (2.56) to see that

$$\mathfrak{E}(t) = \mathfrak{E}(0) - t\,\mathfrak{n}(t) \ .$$

(2.57)

As $\mathfrak{E}(t) \ge 0$ in any event, the preceding identity requires that

$$\mathfrak{n}(t) \le \mathfrak{E}(0)\, t^{-1} \ .$$

(2.58)

Given the definition of $\mathfrak{n}$, this implies in particular that there exists $s \in (0, t]$ such that

$$\| q_A(a|_s) \|_2^2 \le \mathfrak{E}(0)\, t^{-1} \ .$$

(2.59)

Meanwhile, for any $s \in [0, \infty)$, use of the left most identity in (2.53) gives the bound

$$\| a|_s - \mathfrak{a} \|_2^2 \le s^2\, \mathfrak{n}(s) \ .$$

(2.60)

Indeed, (2.60) follows by noting that $\| a|_s - \mathfrak{a} \|_2 \le \int_0^s \|\frac{\partial}{\partial s} a\|_2 = \int_0^s \|q_A(a)\|_2$ .

*Part 3*: Let $\mathrm{E}$ denote an upper bound for the sequence $\{\mathfrak{F}(A_n + i\,\mathfrak{a}_n)\}_{n=1,2,\ldots}$. Now fix $n \in \{1, 2, \ldots\}$ and take $(A, \mathfrak{a})$ in Parts 1 and 2 to be $(A_n, \hat{\mathfrak{a}}_n)$. If $\mathrm{E} = 0$, take $\hat{a}_n = \hat{\mathfrak{a}}_n$.

Assume now that $\mathrm{E}$ is positive. The relevant version in this case of the function $\mathfrak{E}$ has $\mathfrak{E}(0) \le \mathrm{E}\, r_n^{-2}$. Granted that such is the case, take t in (2.59) to be $z^{-1}\, \mathrm{E}\, r_n^{-2}$. The resulting version of (2.59) asserts that there exists $s \in (0, t]$ such that $\| q_A(a|_s) \|_2 \le z$; and for such s, the resulting version of (2.60) asserts that $\| a|_s - \mathfrak{a} \|_2^2 \le z^{-1}\, r_n^{-4}$. Meanwhile, (2.57) implies that $\| d_A a|_s \|_2^2 + \| d_A{*}a|_s \|_2^2 \le \mathrm{E}\, r_n^{-2}$ for this same choice of s.

What is said in the preceding paragraph implies directly that all requirements of Lemma 2.3 are met by taking $\hat{a}_n$ to equal $\hat{a}|_s$.

**3. Scaling limits**

Theorem 1.3 and Lemma 2.5 from [U] play a central role in what is done in this section. The required parts are stated momentarily in (3.1). The notation uses $B \subset \mathbb{R}^3$ to denote the radius 1 ball centered on the origin and $\theta_0$ to denote the product connection on the product principal bundle $B \times SO(3)$. The assertion that follows restates the relevant parts Theorem 1.3 and Lemma 2.5 in [U].

*There exists* $\kappa_U > 1$ *with the following significance: Suppose that* A *is a connection on the product* SO(3) *bundle over* B *whose curvature has* $L^2$ *norm at most* $\kappa_U^{-1}$*. There is an automorphism of this bundle that pulls* A *back as* $\theta_0 + \hat{a}_A$ *with* $\hat{a}_A$ *an* $\mathfrak{su}(2)$ *valued 1-form on* B *with the properties listed below.*

- *The 1-form* $\hat{a}_A$ *is coclosed, thus* $d{*}\hat{a}_A = 0$.
- *The 2-form* ${*}\hat{a}_A$ *pulls back as zero to the boundary of the closure of* B.
- *The Sobolev* $L^2_1$ *norm of* $\hat{a}_A$ *is bounded by* $\kappa_U$ *times the* $L^2$ *norm of* $F_A$.

*Conversely, if* $A = \theta_0 + \hat{a}_A$ *is a connection on* $B \times SO(3)$ *such that* $\hat{a}_A$ *obeys the first two bullets and has* $L^2_1$ *norm bounded by* $\frac{1}{2}\kappa_U^{-1}$*, then the third bullet is also obeyed.*

(3.1)

By way of a guide for those unfamiliar with [U], the proof of (3.1) uses an open/closed argument that exploits four facts, the first being that the $L^2_1$ norm dominates the $L^4$ norm. The second fact is that the curvature of $\theta_0 + \hat{a}$ differs from $d\hat{a}$ by a term that is quadratic in $\hat{a}$. The third fact is as follows: If $\hat{a}$ is an $\mathfrak{su}(2)$ valued 1-form on B with ${*}\hat{a}$ pulling back as zero to the boundary of B, then the sum of the $L^2$ norms $d\hat{a}$ and $d{*}\hat{a}$ is greater than $c_0^{-1}$ times the $L^2_1$ norm of $\hat{a}$. The fourth fact asserts that the differential of an $\mathfrak{su}(2)$ valued function can be added to any given $\mathfrak{su}(2)$ valued 1-form so that the result obeys the first and second bullets in (3.1).

This section uses Uhlenbeck's theorem to analyze the behavior of suitably constrained connections on P and sections of the bundle $(P \times_{SO(3)} \mathfrak{su}(2)) \otimes T^*M$ on balls with radius chosen to guarantee an apriori $L^2$ bound for the curvature of the connection.

**a) The constraints**

Fix $r \geq 1$ and $E \geq 1$, and let $(A, \hat{a}) \in Conn(P) \times C^\infty(M; (P \times_{SO(3)} \mathfrak{su}(2)) \otimes T^*M)$ denote a pair that obeys the following constraints:

a) $\int_M |\hat{a}|^2 = 1.$

b) $\int_M (|d_A\hat{a}|^2 + |d_A{*}\hat{a}|^2) < E r^{-2}$

c) $\int_M |F_A - r^2 \hat{a} \wedge \hat{a}|^2 < E$ .

d) $\int_M |q_A(\hat{a})|^2 < E$ .

(3.2)

The subsections that follow and Sections 4 and 5 use pairs that obey Items a)-d) in (3.2). The rest of this subsection derives some direct implications, these summarized by

- $\int_M (|\nabla_A \hat{a}|^2 + r^2|\hat{a}\wedge\hat{a}|^2 + r^{-2}\,|F_A|^2) \le c_0 E$ ,
- $\sup_{p\in M} \int_{dist(p,\cdot)\le r} G_p(|\nabla_A \hat{a}|^2 + 2r^2|\hat{a}\wedge\hat{a}|^2) \le c_0 E^2$ ,
- $\sup_{p\in M}|\hat{a}|(p) \le c_0 E$ .

(3.3)

Items a)-c) with the $f = 1$ and $(A, \hat{a})$ version of (2.2) lead to the bounds in the first bullet. The bounds in the second bullet follow with the addition of Item d) using arguments that are very much like those in Steps 2-4 of the proof of Proposition 2.2. The next paragraph gives the details.

Since $(A,\hat{a})$ are smooth, the function $f$ in the $(A,\hat{a})$ version of (2.15) can be replaced by the Green's function $G_p$ so as to obtain the corresponding version of (2.16). Given the top bullet of (3.3), then what is said in Step 2 of the proof of Proposition 2.2 bounds the right hand side of the $(A,\hat{a})$ version of (2.6) by $c_0 E$. The term on the integral on the left hand side with the integrand $G_p\langle *F_A\wedge\hat{a}\wedge\hat{a}\rangle$ is written as a sum of two integrals, the first with integrand $G_p\langle *(F_A - r^2\hat{a}\wedge\hat{a})\wedge\hat{a}\wedge\hat{a}\rangle$, and the second with integrand $G_p r^2|\hat{a}\wedge\hat{a}|^2$. The latter is non-negative and it contributes in any event to the left hand side of the second bullet in (3.3). The absolute value of the former is at most

$$c_0\|F_A - r^2\hat{a}\wedge\hat{a}\|_2\,\|\hat{a}\|_\infty\,\|dist(\cdot,p)^{-1}\,\hat{a}\|_2 .$$

(3.4)

Use the first bullet of (3.3) with (2.19) to bound $\|dist(\cdot,p)^{-1}\,\hat{a}\|_2$ by $c_0 E^{1/2}$. This bound and Item c) in (3.2) bound (3.4) by $c_0 E\|\hat{a}\|_\infty$. Given what was said about the right hand side of the $(A,\hat{a})$ version of (2.16), the latter bound leads directly to the bound

$$|\hat{a}|^2(p) + \int_{dist(p,\cdot)\le r} G_p(|\nabla_A \hat{a}|^2 + 2r^2|\hat{a}\wedge\hat{a}|^2) \le c_0 E^2$$

(3.5)

with $p \in M$ being any given point.

**b) The parameter $r_\Diamond$**

Fix a point $p \in M$. With $r \in (0, c_0^{-1})$ specified, introduce by way of notation $B_r$ to denote the ball of radius r centered at p. Denote by $r_\Diamond$ the largest value of r such that

$$\int_{B_r} |F_A|^2 \le \tfrac{1}{100}\,\kappa_U^{-2}\, r^{-1}$$

(3.6)

Note that $r_\Diamond \geq c_0^{-1} r^{-1}$, this being a consequence of Item c) in (3.2) and the third bullet of (3.3). The upcoming Propositions 3.1 and 3.2 give an indication of the significance of $r_\Diamond$.

To set the stage for Proposition 3.1 and for the discussion in the subsequent subsections, fix $p \in M$ and Gaussian coordinates centered on p. Use these coordinates to identify a radius $c_0^{-1}$ ball centered at p with the radius $c_0^{-1}$ ball centered at the origin in $\mathbb{R}^3$. Let $\phi$ denote the map from the radius $c_0^{-1} r_\Diamond^{-1}$ ball about the origin in $\mathbb{R}^3$ to the radius $c_0^{-1}$ ball about p that is obtained by composing first the map $x \to r_\Diamond x$ from $\mathbb{R}^3$ to itself, and then the map that is defined by the Gaussian coordinates.

Pull the pair $(A, r_\Diamond^{-1}\hat{a})$ back to the radius $c_0^{-1}r_\Diamond^{-1}$ ball about the origin in $\mathbb{R}^3$ using $\phi$ to define a pair consising of a connection on $\phi^*P$ and a $\phi^*(P\times_{SO(3)} \mathfrak{su}(2))$ valued 1-form. Denote this pair by $(A_\Diamond, \hat{a}_\Diamond)$. The definition of $\hat{a}_\Diamond$ as $r_\Diamond^{-1}\phi^*\hat{a}$ implies that $|\hat{a}_\Diamond| \leq \|\hat{a}\|_\infty$ and that $|\hat{a}_\Diamond|(0) = |\hat{a}|(p)$ with it understood that the norm of $\hat{a}_\Diamond$ is defined by using the metric that is given by $r_\Diamond^{-2}$ times $\phi$'s pull-back of M's metric. The latter metric is denoted by $\mathfrak{m}_\phi$. Note in any event that the metric $\mathfrak{m}_\phi$ and the Euclidean metric differ by a term whose norm and first two derivatives are bounded by $c_0 r_\Diamond^2$, and whose derivatives to each order $k > 2$ are bounded by a k-dependent constant time $r_\Diamond^k$.

The number $r_\Diamond$ is defined so that the curvature of $A_\Diamond$ obeys

$$\int_{|x|\leq 1} |F_{A_\Diamond}|^2 \leq \tfrac{1}{100}\kappa_U^{-2} \ \textit{with equality if}\ r_\Diamond < c_0^{-1}, \tag{3.7}$$

with it again understood that $\mathfrak{m}_\phi$ is used to define the norm and volume form.

The proposition also introduces $\theta_0$ to denote the product connection on principal SO(3) bundle $\mathbb{R}^3 \times SO(3)$. Proposition 3.1's Hodge star operator is the Euclidean metric's Hodge star, not $\mathfrak{m}_\phi$'s Hodge star. Meanwhile, the norm and volume form used in Proposition 3.1 can be either those defined by $\mathfrak{m}_\phi$ of those defined by the Euclidean metric.

**Proposition 3.1**: *There exists* $\kappa > 1$ *with the following significance: Suppose that* $E \geq 1$, $r > 1$ *and* $(A, \hat{a})$ *is a pair of connection on* P *and section of* $(P\times_{SO(3)} \mathfrak{su}(2)) \otimes T^*M$ *that obey (3.2). Fix* $p \in M$. *There is an isomorphism, to be denoted by* g, *from the product* SO(3) *bundle over the* $|x| < \kappa^{-1} r_\Diamond^{-1}$ *ball in* $\mathbb{R}^3$ *to* $\phi^*P$ *such that*

- *The connection* $g^*A_\Diamond$ *can be written as* $\theta_0 + \hat{a}_{A_\Diamond}$ *where* $\hat{a}_{A_\Diamond}$ *is an* $\mathfrak{su}(2)$*-valued 1-form on the* $|x| \leq 1$ *ball in* $\mathbb{R}^3$ *with* $L^2_1$ *norm bounded by* $\kappa \int_{|x|\leq 1} |F_{A_\Diamond}|^2$. *Moreover,* $\hat{a}_{A_\Diamond}$ *obeys* $d_{A_\Diamond} * \hat{a}_{A_\Diamond} = 0$ *and* $* \hat{a}_{A_\Diamond}$ *pulls back as zero to the* $|x| = 1$ *sphere.*

- *The* $L^2_1$ *norm of* $g^*\hat{a}_\Diamond$ *on the* $|x| < 1$ *ball is bounded by* $\kappa E^2$. *Moreover, given* $r \in [\frac{1}{2}, 1)$*, there exists* $\kappa_{E,r} > \kappa$ *which is independent of* $(A, \hat{a})$ *and is such that the* $L^2_2$ *norm of* $g^*\hat{a}_\Diamond$ *on the* $|x| \le r$ *ball in* $\mathbb{R}^3$ *is bounded by* $\kappa_{E,r}$.

Granted (3.7) and some standard Sobolev inequalities, Proposition 3.1 amounts to little more than a corollary to Uhlenbecks theorem in (3.1). In any event, its proof is in the next subsection.

Proposition 3.1 can be said to see only the SO(3) subgroup in PSL(2;$\mathbb{C}$), this being the maximal compact subgroup. The upcoming Proposition 3.2 can be said to see the whole of PSL(2;$\mathbb{C}$). Proposition 3.2 is the key to Theorem 1.1a's extension of Uhlenbeck's theorem to PSL(2; $\mathbb{C}$).

**Proposition 3.2**: *Given* $E \ge 1$, $\mu \in (0, \frac{1}{2}]$ *and* $\varepsilon \in (0, 1]$*, there exists* $\kappa_{E,\mu,\varepsilon} > 1$ *with the following significance: Suppose that* $r > 1$ *and* $(A, \hat{a})$ *is a pair of connection on* P *and section of* $(P \times_{SO(3)} \mathfrak{su}(2)) \otimes T^*M$ *that obey (3.2). Fix* $p \in M$. *If both* $r_\Diamond < \kappa_{E,\mu,\varepsilon}^{-1}$ *and* $\int_{B_{r_\Diamond}} |\nabla_A \hat{a}|^2 \le \kappa_{E,\mu,\varepsilon}^{-1} r_\Diamond^{-2} \int_{B_{r_\Diamond}} |\hat{a}|^2$, *then the* $r = (1 - \mu) r_\Diamond$ *version of* $\int_{B_r} |F_A|^2$ *is less than* $\varepsilon\, r^{-1}$.

The proof of Proposition 3.2 will invoke Proposition 3.1. This being the case, it proves convenient to restate the proposition in terms of $A_\Diamond$ and $\hat{a}_\Diamond$. To this end, introduce $z$ to denote the $L^2$ norm of $\hat{a}_\Diamond$ on the $|x| < 1$ ball; then define $\hat{a}_*$ to equal $z^{-1}\hat{a}_\Diamond$. Proposition 3.2's assertion is equivalent to the following:

*Given* $E \ge 1$, $\varepsilon \in (0, 1]$ *and* $\mu \in (0, \frac{1}{2}]$ *and there exists* $\kappa_{E,\mu,\varepsilon} > 1$ *with the following significance: Suppose that* $(A, \hat{a}) \in \mathrm{Conn}(P) \times C^\infty(M; P \times_{SO(3)} \mathfrak{su}(2)) \otimes T^*M$ *obeys (3.2) Fix* $p \in M$. *If* $r_\Diamond < \kappa_{E,\mu,\varepsilon}^{-1}$ *and* $\int_{|x| \le 1} |\nabla_{A_\Diamond} \hat{a}_*|^2 \le \kappa_{E,\mu,\varepsilon}^{-1}$*, then* $\int_{|x| \le 1-\mu} |F_{A_\Diamond}|^2 < \varepsilon$.

(3.8)

The equivalence between the assertion in (3.8) and Proposition 3.2 follows directly from the scaling identities given in the next subsection. Note that the norms and volume form used in (3.8) can be either those defined by $\mathfrak{m}_\phi$ or those defined by the Euclidean metric. The proof of (3.8) will use those defined by the Euclidean metric.

The proof of Proposition 3.2 is in Section 4. Section 3c states some rescaling identities and, as noted previously, it proves Proposition 3.1. Sections 3d and 3e state and prove a pair of lemmas that are used in Section 4 for the proof of Proposition 3.2.

**c) Scaling identities**

Proposition 3.1 is proved at the end of this subsection. The observations that follow directly in (3.9) are used in the proof and in the subsequent subsections. The upcoming (3.9) lists some rescaling identities. What is denoted by $\mathfrak{e}_G$ in the fourth bullet has norm bounded by $c_0 \mathrm{r}_\lozenge$. All bullets refer to a chosen $R \in (0, c_0^{-1}\mathrm{r}_\lozenge^{-1})$. The norms, inner products, volume form and Hodge dual in (3.9) are all defined using $\mathfrak{m}_\phi$.

- $\int_{|x|\le R} (|d_{A_\lozenge} \hat{a}_\lozenge|^2 + |d_{A_\lozenge} {*\hat{a}_\lozenge}|^2) \le \mathrm{E}\,\mathrm{r}_\lozenge^{-1} r^{-2}$ .
- $\int_{|x|\le R} |F_{A_\lozenge} - \mathrm{r}_\lozenge^2 r^2 \hat{a}_\lozenge \wedge \hat{a}_\lozenge|^2 \le \mathrm{E}\,\mathrm{r}_\lozenge$ .
- $\int_{|x|\le R} |\nabla_{A_\lozenge} \hat{a}_\lozenge|^2 = \mathrm{r}_\lozenge^{-1} \int_{B_{R\mathrm{r}_\lozenge}} |\nabla_A \hat{a}|^2$ .
- $\int_{|x|\le R} |F_{A_\lozenge}|^2 = \mathrm{r}_\lozenge \int_{B_{R\mathrm{r}_\lozenge}} |F_A|^2$ .
- $\int_{|x|\le R} |\mathcal{q}_{A_\lozenge}(\hat{a}_\lozenge)|^2 = \mathrm{r}_\lozenge \int_{B_{R\mathrm{r}_\lozenge}} |\mathcal{q}_A(\hat{a})|^2$ .

(3.9)

Note that $|\hat{a}_\lozenge| \le c_0 \mathrm{E}$, this being a consequence of the third bullet in (3.3). Moreover,

- $\int_{|x|\le R} |\nabla_{A_\lozenge} \hat{a}_\lozenge|^2 \le c_0 \mathrm{E}^2 R$
- $\mathrm{r}_\lozenge^2 r^2 \int_{|x|\le R} |\hat{a}_\lozenge \wedge \hat{a}_\lozenge|^2 \le c_0 \mathrm{E}^2 R$,
- $\int_{|x|\le R} |\mathcal{q}_{A_\lozenge}(\hat{a}_\lozenge)|^2 \le c_0 \mathrm{E}\,\mathrm{r}_\lozenge$ .

(3.10)

The first and second bullets are consequences of the second bullet in (3.3); and the third bullet is a consequence of Item d) of (3.2). These inequalities hold with the norms, volume form and covariant derivative defined by either $\mathfrak{m}_\phi$ or the Euclidean metric.

A standard, dimension 3 Sobolev inequality is used in the proofs of both Proposition 3.1 and Proposition 3.2; and a two other inequalities are used only in the proof of Proposition 3.2. The first inequality is a version of the assertion that the $L^2_1$ norm dominates the $L^4$ norm, the second inequality is a local version of Hardy's inequality in (2.19) and the third is a version of the assertion that $L^4_1$ norm dominates the Holder exponent $\frac{1}{4}$ norm. These inequalities are given directly for future reference. To set the stage, fix $p \in M$ and $r \in (0, c_0^{-1})$. Let $f$ denote a Lipshitz function on $B_r$. Then

- $\int_{B_r} |f|^4 \le c_0 \left( \int_{B_r} |f|^2 \right)^{1/2} \left( \int_{B_r} (|df|^2 + r^{-2} |f|^2) \right)^{3/2}$ .
- $\int_{B_r} \frac{1}{\mathrm{dist}(p,\cdot)^2} f^2 \le c_0 \int_{B_r} (|df|^2 + r^{-2} |f|^2)$ .
- *If* $q_1$ *and* $q_2$ *are any two points in* $B_r$, *then* $|f(q_2) - f(q_1)| \le c_0 \, \mathrm{dist}(q_2, q_1)^{1/4} \left( \int_{B_r} |df|^4 \right)^{1/4}$ .

(3.11)

These same inequalities hold when $f$ is a function on a ball about the origin in $\mathbb{R}^3$ with the only change being the use of the Euclidean metric to define the norms, volume form, and distance function.

It proves convenient at this point to use henceforth the following notational conventions. Unless stated explicitly to the contrary, the Euclidean metric is used to define the norms, Hodge star, covariant derivatives and volume form on the $|x| < c_0^{-1} r_\lozenge^{-1}$ part of $\mathbb{R}^3$. Covariant derivatives on $\mathfrak{su}(2)$ valued tensors on $\mathbb{R}^3$ that are defined using the connection $\theta_0$ are denoted by $\nabla$. What is denoted subsequently as $c_{\mathrm{E}}$ is a number that is greater than 16 and depends only on $\mathrm{E}$. It can be assumed to increase between successive appearances. If $\mu \in (0, 1)$ has also been specified, $c_{\mathrm{E},\mu}$ is used to denote a number that is greater than 16 and depends only on $\mathrm{E}$ and $\mu$. It can also be assumed to increase between successive appearances

***Proof of Proposition 3.1***: Granted (3.3), then Uhlenbeck's theorem in the guise of (3.1) can be invoked. This theorem supplies an isomorphism from the product SO(3) bundle over the $|x| < 1$ ball to $\phi^*P$ that pulls $A_\lozenge$ back as $\theta_0 + \hat{a}_{A_\lozenge}$ with $\hat{a}_{A_\lozenge}$ as described by Proposition 3.1. Use g to denote this isomorphism.

The $L^2$ norm of $\nabla_{A_\lozenge} \hat{a}_\lozenge$ on any radius $R < c_0 r_\lozenge^{-1}$ ball in $\mathbb{R}^3$ centered at the origin is a priori bounded by $c_{\mathrm{E}} R$ this being a consequence of the first bullet in (3.10). This $L^2$ norm bound on $\nabla_{A_\lozenge} \hat{a}_\lozenge$ with what was said about $g^*A_\lozenge$ implies that the $L^2_1$ norm of $g^*\hat{a}_\lozenge$ on the $|x| < 1$ ball is bounded by $c_{\mathrm{E}}$.

Fix $\mu \in (0, \frac{1}{2}]$ and introduce by way of notation $\chi_\mu$ to denote the function on $\mathbb{R}^3$ given by the rule $x \to \chi(\frac{1}{\mu}(|x| - 1 + \mu))$. This function equals 1 where $|x| \le 1 - \frac{3}{4}\mu$ and it equals 0 where $|x| \ge 1 - \frac{1}{4}\mu$. Integrate by parts and use the $L^\infty$ bound on $|\hat{a}_\lozenge|$ to see that

$$\int_{|x|\le 1} |\nabla_{A_\lozenge}(\chi_\mu \nabla_{A_\lozenge} \hat{a}_\lozenge)|^2 \le$$

$$\int_{|x|\le 1} \chi_\mu |q_\lozenge|^2 + c_0\Big(\mu^{-2} \|\nabla_{A_\lozenge} \hat{a}_\lozenge\|_2^2 + \int_{|x|\le 1} |F_{A_\lozenge}| \, \chi_\mu^2 |\nabla_{A_\lozenge} \hat{a}_\lozenge|^2 + r_\lozenge^2\Big) + c_0 \mathrm{E}^2 \int_{|x|\le 1} |F_{A_\lozenge}|^2 \chi_\mu^2$$

(3.12)

To exploit this inequality, let $x$ denote for the moment $\chi_\mu\, \nabla_{A_\Diamond} \hat{a}_\Diamond$. Then (3.7), the first and third bullets of (3.10) and (3.12) lead to the bound

$$\| \nabla_{A_\Diamond} x \|_2^2 \le c_E(1+\mu^{-2} + \|x\|_4^2) , \tag{3.13}$$

Given the Sobolev inequality asserted by the top bullet of (3.11), the inequality in (3.13) implies that $\chi_\mu|\nabla_{A_\Diamond} \hat{a}_\Diamond|$ is an $L^2_1$ function on the $|x| \le 1 - \mu$ ball and a second appeal to (3.11) and (3.12) bounds $\|\nabla_{A_\Diamond} x\|_2$ by $c_{E,\mu}$.

Granted the latter $L^2$ norm bound, write $g^*A_\Diamond = \theta_0 + \hat{a}_{A_\Diamond}$ to see that

$$\|\nabla(g^*x)\|_2 \le c_0 \| \hat{a}_{A_\Diamond} \|_4 \|x\|_4 . \tag{3.14}$$

Use this last bound, the a priori $L^2_1$ bound for $\hat{a}_{A_\Diamond}$ and the $f = |\hat{a}_{A_\Diamond}|$ version of the top bullet of (3.11) to see that the $L^2_1$ norm of $g^*x$ is also bounded by $c_{E,\mu}$. Write $g^*(\nabla_{A_\Diamond} \hat{a}_\Diamond)$ as $\nabla_{g^*A_\Diamond}(g^*\hat{a}_\Diamond)$ to see that

$$|\nabla(\nabla(g^*\hat{a}_\Diamond))| \le c_0(|\nabla(g^*(\nabla_{A_\Diamond}\hat{a}_\Diamond))| + (|\nabla \hat{a}_{A_\Diamond}| + |\hat{a}_{A_\Diamond}|^2)|\hat{a}_\Diamond| + |\hat{a}_{A_\Diamond}||\nabla_{A_\Diamond}\hat{a}_\Diamond|) . \tag{3.15}$$

Given this last inequality, then the apriori $L^2_1$ bounds for $g^*x$ and $\hat{a}_{A_\Diamond}$ with the apriori $L^\infty$ bound for $|\hat{a}_\Diamond|$ lead directly to the desired $c_{E,\mu}$ bound for the $L^2_2$ norm of $g^*\hat{a}_\Diamond$.

**d) When $r_z = z r_\Diamond r$ is small**

Let $z$ again denote the $L^2$ norm of $\hat{a}_\Diamond$ over the $|x| \le 1$ ball in $\mathbb{R}^3$. Introduce by way of notation $r_z$ to denote the combination $z r_\Diamond r$. This combination appears when writing the integral in the second bullet of (3.9) in terms of $\hat{a}_*$. Do so and it asserts the following:

$$\int_{|x|\le R} |F_{A_\Diamond} - r_z^2\, \hat{a}_* \wedge \hat{a}_*|^2 \le E r_\Diamond . \tag{3.16}$$

The lemma that follows makes a formal statement to the effect that (3.8) is true if there is an a priori bound on $r_z$.

**Lemma 3.3**: *Given* $E \ge 1$, $K \ge 1$, $\mu \in (0, \frac{1}{2}]$ *and* $\varepsilon \in (0, 1]$, *there exists* $\kappa > 1$ *with the following significance: Fix* $r > 1$ *and a pair* $(A, \hat{a})$ *of connection on* P *and section of* $(P\times_{SO(3)}\mathfrak{su}(2))\otimes T^*M$ *obeying (3.2). Suppose that* $p \in M$ *is a point where* $r_z \leq \kappa$. *If both* $r_\Diamond < \kappa^{-1}$ *and* $\int_{|x|\leq 1} |\nabla_{A_\Diamond}\hat{a}_*|^2 < \kappa^{-1}$, *then* $\int_{|x|\leq 1-\mu} |F_{A_\Diamond}|^2 < \varepsilon$.

***Proof of Lemma 3.3***: Define $\chi_\mu$ to be the function on $\mathbb{R}^3$ given by $x \to \chi(\frac{1}{\mu}(|x| - 1+\mu))$. This function equals 1 where $|x| \leq 1 - \frac{3}{4}\mu$ and it equals 0 where $|x| \geq 1 - \frac{1}{4}\mu$. Fix constant orthonormal vectors $e_1$ and $e_2$ on $\mathbb{R}^3$ and integrate $\chi_\mu{}^2\langle\hat{a}_*(e_1)\,[F_{A_\Delta}(e_2,e_1), \hat{a}_*(e_2)]\rangle$ over the radius 1 ball about the origin. The resulting integral is the same as the integral of $\chi_\mu{}^2\langle F_{A_\Delta}(e_2,e_1)\,[\hat{a}_*(e_2), \hat{a}_*(e_1)]\rangle$. It follows from (3.16) that $c_0$ times the latter integral plus $c_E\, r_z{}^{-2} r_\Diamond$ bounds $r_z{}^{-2}\|\chi_\mu\, F_{A_\Delta}(e_1,e_2)\|_2{}^2$.

Meanwhile, $[F_{A_\Delta}(e_2,e_1), \hat{a}_*(e_2)]$ can be written as the commutator of $A_\Diamond$ covariant derivatives of $\hat{a}_*(e_2)$. Do so and use this depiction with an integration by parts to see that

$$r_z{}^{-2}\|\chi_\mu\, F_{A_\Delta}(e_1,e_2)\|_2{}^2 \leq c_0(\|\nabla_{A_\Diamond}\hat{a}_*\|_2{}^2 + \mu^{-1}\|\nabla_{A_\Diamond}\hat{a}_*\|_2 + r_z{}^{-2} r_\Diamond c_E)\,, \tag{3.17}$$

Suppose that $r_z < \kappa$. Then (3.16) asserts the desired bound $\|\chi_\mu\, F_{A_\Delta}\|_2{}^2 < \varepsilon$ when $r_\Diamond < c_E{}^{-1}\varepsilon$ and when $\|\nabla_{A_\Diamond}\hat{a}_*\| \leq c_0{}^{-1}\kappa^{-2}\mu\,\varepsilon$.

The remaining subsections assume unless stated to the contrary that $r_z \geq 1$,

**e) First order equations**

The 4 parts of this subsection explain how a certain almost tautological system of inhomogeneous, semi-linear first order differential equations can be used to prove Proposition 3.2. The linear terms in the homogeneous version of these equations define the first order, constant coefficient system of elliptic equations on $\mathbb{R}^3$ that is described in Part 1 of the subsection. The fully non-linear, inhomogeneous equation is described in Part 2. Part 3 explains how a suitable a priori estimate for a solution to these equations can be used to prove Proposition 3.2. The final part of the subsection states and then proves a lemma that is subsequently used to invoke the arguments in Part 3.

*Part 1*: To set the stage for what is to come, let $\tau$ denote a given, unit length element in $\mathfrak{su}(2)$. With $\tau$ chosen, introduce by way of notation $\mathbb{V} \subset (T^*\mathbb{R}^3\oplus\mathbb{R})\otimes\mathfrak{su}(2)$ to denote the subbundle that is annihilated by the homomorphism to $(T^*\mathbb{R}^3\oplus\mathbb{R})$ that is defined by the rule $\mathfrak{f} \to \langle\tau\mathfrak{f}\rangle$. Let $e$ denote a chosen, constant 1-form on $\mathbb{R}^3$ with norm 1.

Fix $m > 1$. Define $\mathcal{L}_m$: $C^\infty(\mathbb{R}^3; \mathbb{V}\oplus\mathbb{V}) \to C^\infty(\mathbb{R}^3; \mathbb{V}\oplus\mathbb{V})$ as follows: Write an element in $\mathbb{V}$ as $(\mathfrak{a}, \mathfrak{a}_0)$ with $\mathfrak{a}$ being an $\mathfrak{su}(2)$ valued 1-form and $\mathfrak{a}_0$ an element in $\mathfrak{su}(2)$. The operator $\mathcal{L}_m$ sends a given element $((\mathfrak{a}, \mathfrak{a}_0), (\mathfrak{b}, \mathfrak{b}_0))$ to one whose components in the left and right most factor of $\mathbb{V}\oplus\mathbb{V}$ are the respective pairs of $\mathfrak{su}(2)$ valued 1-form and valued function given by

- $*(d\mathfrak{a} - m e\wedge[\tau, \mathfrak{b}]) + d\mathfrak{a}_0 + m e[\tau, \mathfrak{b}_0]$ *and* $-*(d{*}\mathfrak{a} + m\, e \wedge [\tau, *\mathfrak{b}])$ ,
- $-*(d\mathfrak{b} + m e\wedge[\tau, \mathfrak{a}]) - d\mathfrak{b}_0 + m e[\tau, \mathfrak{a}_0]$ *and* $*(d{*}\mathfrak{b} - m\, e\wedge[\tau, *\mathfrak{a}])$ ,

$$(3.18)$$

with $*$ denoting the Euclidean metric's Hodge star operator. The operator $\mathcal{L}_m$ is symmetric with respect to the $L^2$ inner product defined by the Euclidean metric and such that $\mathcal{L}_m{}^2 = (\nabla^\dagger\nabla + m^2)\,\mathbb{I}$ with $\mathbb{I}$ denoting the identity endomorphism of $\mathbb{V}$, with $\nabla$ denoting the covariant derivative on $C^\infty(\mathbb{R}^3; \mathbb{V}\oplus\mathbb{V})$ as defined by the Euclidean metric with $\nabla^\dagger$ being its formal, $L^2$ adjoint.

This depiction of $\mathcal{L}_m{}^2$ has the following useful corollary: Let $\|\cdot\|_2$ denote the $L^2$ inner product on the space of $L^2_1$ map from $\mathbb{R}^3$ to $\mathbb{V}\oplus\mathbb{V}$. If $\mathfrak{k}$ is any such map, then

$$\|\mathcal{L}_m \mathfrak{k}\|_2{}^2 = \|\nabla\mathfrak{k}\|_2{}^2 + m^2\|\mathfrak{k}\|_2{}^2 \,. \tag{3.19}$$

The positivity of the square of $\mathcal{L}_m$ implies that the equation $\mathcal{L}_m \mathfrak{f} = \mathfrak{h}$ has a unique, $L^2_1$ solution in $C^\infty(\mathbb{R}^3; \mathbb{V}\oplus\mathbb{V})$ when $\mathfrak{h} \in C^\infty(\mathbb{R}^3; \mathbb{V}\oplus\mathbb{V})$ is square integrable. This solution $\mathfrak{h}$ can be written explicitly using the Green's function for $-d^\dagger d + m^2$. The version of the latter with pole at a given point $y \in \mathbb{R}^3$ is denoted in what follows by $G_y$; it is the function on $\mathbb{R}^3 - \{y\}$ given by

$$G_y(\cdot) = \frac{1}{4\pi}\frac{1}{|(\cdot)-y|}\, e^{-m|(\cdot)-y|} \,. \tag{3.20}$$

As can be seen from (3.20), the function $G_y$ obeys

$$|G_y| + |dG_y| + |\nabla(dG_y)| \le c_0\, e^{-m\delta/2} \tag{3.21}$$

at points with distance greater than $\delta$ from y. The Green's function for $\mathcal{L}_m$ with pole at y is the $\mathrm{End}(\mathbb{V}\oplus\mathbb{V})$-valued function on $\mathbb{R}^3 - y$ that is defined by $x \to (\mathcal{L}_m G_{(\cdot)})|_y(x)$. Note in particular that the bounds given above for $G_y$ and its derivatives lead directly to the following observation: The Green's function for $\mathcal{L}_m$ with pole at y and those of its first derivatives are also bounded by $c_0\,\delta^{-2} e^{-m\delta/2}$ at points with distance greater than $\delta$ from y.

*Part 2*: Let $(A_\Delta, a_\Delta)$ denote a pair consisting of a connection on the product SO(3) bundle over the $|x| \le 1$ ball in $\mathbb{R}^3$ and an $\mathfrak{su}(2)$ valued 1-form on this ball. The connection $A_\Delta$ is written as $\theta_0 + \hat{a}_{A_\Delta}$. Fix $m \ge 1$ and, with $\tau \in \mathfrak{su}(2)$ as in (3.18), define

$$\mathfrak{a} = m^{-1}(\hat{a}_{A_\Delta} - \tau\langle\tau\,\hat{a}_{A_\Delta}\rangle) \quad \textit{and} \quad \mathfrak{b} = \sqrt{\tfrac{4\pi}{3}}\,(a_\Delta - \tau\langle\tau\, a_\Delta\rangle)\,.$$

(3.22)

With $\mathfrak{b}$ understood, define the 1-form $e_\Delta$ by writing $a_\Delta$ as $a_\Delta = \sqrt{\tfrac{3}{4\pi}}\,(\tau e_\Delta + \mathfrak{b})$.

Let $\mathfrak{su}_\perp \subset \mathfrak{su}(2)$ denote the kernel of the homorphism $\mathfrak{f} \to \langle\tau\mathfrak{f}\rangle$. Use $\mathfrak{s}$ to denote the $\mathfrak{su}_\perp$ part of the $\mathfrak{su}(2)$ valued 1-form $m^{-1}{*}(F_{A_\Delta} - \tfrac{4\pi}{3} m^2 a_\Delta\wedge a_\Delta)$. By the same token, use $\mathfrak{r}$ to denote the $\mathfrak{su}_\perp$ parts of $-\sqrt{\tfrac{4\pi}{3}}\, d_{A_\Delta} a_\Delta$. Let $*_\phi$ denote the Hodge dual that is defined by the metric $\mathfrak{m}_\phi$ and introduce $\mathfrak{r}_0$ to denote the $\mathfrak{su}_\perp$ part of $\sqrt{\tfrac{4\pi}{3}}\, *(d_{A_\Delta} *_\phi a_\Delta)$. Granted this notation, define an $\mathfrak{su}_\perp$ valued 1-form $\mathrm{S}$ by writing $\mathfrak{s}$ as $*(d\mathfrak{a} - me\wedge[\tau,\mathfrak{b}]) - \mathrm{S}$. Introduce a second $\mathfrak{su}_\perp$ valued 1-form $\mathrm{R}$ by writing $\mathfrak{r}$ as $-*(d\mathfrak{b} - me\wedge[\tau,\mathfrak{a}]) - \mathrm{R}$, and introduce an $\mathfrak{su}_\perp$ valued function $\mathrm{R}_0$ by writing $\mathfrak{r}_0$ as $*(d{*}\mathfrak{b} - me\wedge[\tau,*\mathfrak{a}]) - \mathrm{R}_0$. By way of a summary the troika $(\mathrm{S}, \mathrm{R}, \mathrm{R}_0)$ are written below in the first three bullets of (3.23). Write the $\mathfrak{su}_\perp$ valued function $-*(d{*}\mathfrak{a} + me\wedge[\tau,*\mathfrak{b}])$ as $\mathrm{S}_0$ and use $\mathrm{S}_0$ to define the $\mathfrak{su}_\perp$ valued function $\mathfrak{s}_0$ using the fourth bullet of (3.23).

- $*\mathrm{S} = m(e_\Delta - e)\wedge[\tau,\mathfrak{b}] - \langle\tau\hat{a}_{A_\Delta}\rangle\wedge[\tau,\mathfrak{a}] + \mathfrak{s}$ ,
- $*\mathrm{R} = -m(e_\Delta - e)\wedge[\tau,\mathfrak{a}] + \langle\tau\,\hat{a}_{A_\Delta}\rangle\wedge[\tau,\mathfrak{b}] + \mathfrak{r}$ ,
- $*\mathrm{R}_0 = *(d(\pi_\phi\mathfrak{b}) - me\wedge[\tau,\pi_\phi\mathfrak{a}]) + m(e_\Delta - e)\wedge[\tau, *_\phi\,\mathfrak{a}] + \langle\tau\,\hat{a}_{A_\Delta}\rangle\wedge[\tau, *_\phi\,\mathfrak{b}] + \mathfrak{r}_0$ ,
- $*\mathrm{S}_0 = *(-d(\pi_\phi\mathfrak{a}) - me\wedge[\tau,\pi_\phi\mathfrak{b}]) + m(e_\Delta - e)\wedge[\tau, *_\phi\,\mathfrak{b}] - \langle\tau\,\hat{a}_{A_\Delta}\rangle\wedge[\tau,*\mathfrak{a}] + \mathfrak{s}_0$ ,

(3.23)

where the notation has $\pi_\phi$ denoting $(*_\phi - *)$.

Granted all of these definitions, view $\mathfrak{f} = ((\mathfrak{a}, 0), (\mathfrak{b}, 0))$ as mapping the $|x| < 1$ ball to $\mathbb{V}\oplus\mathbb{V}$. The map $\mathfrak{f}$ obeys the tautological equation $\mathcal{L}_m\mathfrak{f} = \mathfrak{h}$ with $\mathfrak{h} = ((\mathrm{S}, \mathrm{S}_0), (\mathrm{R}, \mathrm{R}_0))$. The respective left and write $\mathbb{V}$ summands of this equation are

- $*(d\mathfrak{a} - me\wedge[\tau,\mathfrak{b}]) = \mathrm{S}$ *and* $-*(d{*}\mathfrak{a} + me\wedge[\tau,*\mathfrak{b}]) = \mathrm{S}_0$ .
- $-*(d\mathfrak{b} + me\wedge[\tau,\mathfrak{a}]) = \mathrm{R}$ *and* $*(d{*}\mathfrak{b} - me\wedge[\tau,*\mathfrak{a}]) = \mathrm{R}_0$ .

(3.24)

The equation $\mathcal{L}_m\mathfrak{f} = \mathfrak{h}$ is introduced so that the properties of the Green's function for $\mathcal{L}_m$ can be used to obtain a priori bounds on $\mathfrak{f}$. This is done in Part 4 of the subsection.

*Part 3*: With the proof of Proposition 3.2 in mind, what follows constitutes a digression that explains how bounds on the $L^2$ norm of $\mathfrak{b}$ and $L^\infty$ norm of $a_\Delta$ lead to a bound for the $L^2$ norm of $F_{A_\Delta}$. This comes about by writing $a_\Delta$ as $\sqrt{\frac{3}{4\pi}}\,(\tau e_\Delta + \mathfrak{b})$ so as to write $F_{A_\Delta}$ as

$$F_{A_\Delta} = m^2 e_\Delta \wedge [\tau, \mathfrak{b}] + m^2\, \mathfrak{b} \wedge \mathfrak{b} + (\,F_{A_\Delta} - \tfrac{4\pi}{3} m^2 a_\Delta \wedge a_\Delta)\ . \tag{3.25}$$

Fix $\mu \in (0, \frac{1}{2}]$. It follows directly from this depiction of $F_{A_\Delta}$ that its $L^2$ norm on the $|x| < 1 - \mu$ ball in $\mathbb{R}^3$ obeys

$$\int_{|x|\le 1-\mu} |F_{A_\Delta}|^2 \le 2\,(\sup_{|x|\le 1-\mu} |a_\Delta|^2)\,(m^4 \int_{|x|\le 1-\mu} |\mathfrak{b}|^2\,) \ + \int_{|x|\le 1-\mu} |F_{A_\Delta} - \tfrac{4\pi}{3} m^2 a_\Delta \wedge a_\Delta|^2 \tag{3.26}$$

This inequality has the following consequence: Fix $\varepsilon \in (0, 1]$. The $L^2$ norm of $F_{A_\Delta}$ on the $|x| < 1 - \mu$ ball will be less than $\varepsilon$ if $\sup_{|x|\le 1-\mu} |a_\Delta| < D$, and if the $L^2$ norm of $|\mathfrak{b}|$ on the $|x| \le 1 - \mu$ ball is less than $\frac{1}{4} m^{-2} D^{-1} \varepsilon$, and if that of $(F_{A_\Delta} - \frac{4\pi}{3} m^2 a_\Delta \wedge a_\Delta)$ is less than $\frac{1}{2}\varepsilon$.

*Part 4*: The tautological equation $\mathcal{L}_m \mathfrak{f} = \mathfrak{h}$ leads to a priori bound for $\mathfrak{f}$ when $\mathfrak{h}$ is small in a suitable sense. A precise statement as to what is meant by 'small' requires the introduction of parameters $M \ge 1$, $\mu \in (0, 2^{-20}]$ and $\rho \in (0, 1]$. Having chosen their values, make the following assumptions:

- $\sup_{|x|\le 1-\mu} (|\mathfrak{b}| + |e_\Delta - e|) < \rho$ .
- *The* $L^2$ *norm of* $\mathfrak{a}$ *on the* $|x| \le 1 - \mu$ *ball is less than* $\rho$.
- *The* $L^4$ *norm of* $\langle \tau\, \hat{a}_{A_\Delta} \rangle$ *on the* $|x| \le 1 - \mu$ *ball is less than* $M$ .
- *The* $L^2$ *norms of* $\mathfrak{s}, \mathfrak{r}, \mathfrak{r}_0$ *and* $\mathfrak{s}_0$ *on the* $|x| \le 1 - \mu$ *ball are less than* $m^{-1}\rho$.
- *The endomorphism* $\pi_\phi = *_\phi - *$ *and its covariant derivative obey* $|\nabla \pi_\phi| + |\pi_\phi| < \rho$.

(3.27)

The equation $\mathcal{L}_m \mathfrak{f} = \mathfrak{h}$ is used to prove the next lemma.

**Lemma 3.4**: *There exists* $\kappa > 1$ *and given* $\mu \in (0, 2^{-20}]$ *and* $\varepsilon \in (0, 1)$*, there exists* $\kappa_* > 1$, *these with the following significance*: *Fix* $M$ *and suppose that* $m > \kappa M + \kappa_*^{-1}$ *and* $\rho < \kappa_*^{-1}$. *Define* $(\mathfrak{a}, \mathfrak{b})$ *as in Part 2 and suppose that the bounds in (3.27) are satisfied*. *Then*

$$\int_{|x|\le 1-2048\mu} (\,|\nabla \mathfrak{a}|^2 + |\nabla \mathfrak{b}|^2 + m^2 |\mathfrak{a}|^2 + m^2 |\mathfrak{b}|^2\,) < \varepsilon\, m^{-2}\ .$$

***Proof of Lemma 3.4***: The proof has seven steps.

Step 1: Define $\mathfrak{f}_\mu$: $\mathbb{R}^3 \to \mathbb{V} \oplus \mathbb{V}$ to be the function $\chi_{4\mu}\,\mathfrak{f}$. The latter obeys an equation that has the schematic form $\mathcal{L}_m \mathfrak{f}_\mu = \mathfrak{h}_\mu$ with $\mathfrak{h}_\mu = \mathfrak{S}(d\chi_{4\mu})\mathfrak{f} + \chi_{4\mu}\,\mathfrak{h}$ where $\mathfrak{S}$ denotes the principal symbol of the operator $\mathcal{L}_m$. In this case, $\mathfrak{S}$ is constant and so an element in Hom($\mathbb{R}^3$; End($\mathbb{V} \oplus \mathbb{V}$)). The introduction of $\mathfrak{f}_\mu$ facilitates the use of the Green's function for $\mathcal{L}_m$ that is described in Part 1. The parameter $4\mu$ is used with $\chi_{(\cdot)}$ because each point where $\chi_{4\mu} > 0$ has distance at least $\frac{1}{4}\mu$ from the $|x| > 1 - \mu$ part of $\mathbb{R}^3$ and in particular from the support of $1 - \chi_\mu$.

Step 2: Let $\mathfrak{g}$ denote the square integrable map from $\mathbb{R}^3$ to $\mathbb{V} \oplus \mathbb{V}$ that solves the equation $\mathcal{L}_m \mathfrak{g} = \mathfrak{S}(d\chi_\mu)\mathfrak{f}$. To bound the size of $|\mathfrak{g}|$, first use Part 1's Green's function for $\mathcal{L}_m$ to see that

- $|\mathfrak{g}|(y) \le c_0 \int_{\mathbb{R}^3} (\frac{1}{|(\cdot)-y|^2} + m\frac{1}{|(\cdot)-y|})e^{-m|(\cdot)-y|}\ |d\chi_{4\mu}|\,|\mathfrak{f}|$ *at any given* $y \in \mathbb{R}^3$.
- $|\nabla\mathfrak{g}|(y) \le c_0 \int_{\mathbb{R}^3} (\frac{1}{|(\cdot)-y|^3} + m\frac{1}{|(\cdot)-y|^2} + m^2\frac{1}{|(\cdot)-y|})e^{-m|(\cdot)-y|}\ |d\chi_{4\mu}|\,|\mathfrak{f}|$ *if* $d\chi_{4\mu} = 0$ *at* $y$.

(3.28)

Let $D_\mu$: $\mathbb{R}^3 \to [0, \infty)$ denote the distance to the support of $d\chi_{4\mu}$. Since $|\mathfrak{f}| \le |\mathfrak{a}| + |\mathfrak{b}|$, the bounds in (3.27) and (3.28) imply that

- $|\mathfrak{g}|(y) \le c_0 D_\mu(y)^{-2}\, e^{-mD_\mu(y)/2}\mu^{-1}\rho$ *if* $d\chi_{4\mu} = 0$ *at* $y$.
- $|\nabla\mathfrak{g}|(y) \le c_0 D_\mu(y)^{-3}\, e^{-mD_\mu(y)/2}\mu^{-1}\rho$ *if* $d\chi_{4\mu} = 0$ *at* $y$.

(3.29)

Save these bounds for the moment.

Step 3: Introduce from (3.23) and (3.24) the terms that are denoted by $\mathfrak{s}$, $\mathfrak{r}$, $\mathfrak{r}_0$ and $\mathfrak{s}_0$. Let $\mathfrak{t}$ denote the solution in $L^2_1(\mathbb{R}^3; \mathbb{V} \oplus \mathbb{V})$ to the equation

$$\mathcal{L}_m \mathfrak{t} = \chi_{4\mu}((\mathfrak{s}, \mathfrak{s}_0), (\mathfrak{r}, \mathfrak{r}_0))\ . \tag{3.30}$$

Take the $L^2$ norm of both sides of (3.30) and (3.29) with (3.27)'s fourth bullet find

$$\|\nabla\mathfrak{t}\|_2^{\,2} + m^2\|\mathfrak{t}\|_2^{\,2} \le c_0\, m^{-2}\rho^2\ . \tag{3.31}$$

This bound should also be saved.

Step 4: Let $\mathfrak{q} = \mathfrak{f}_\mu - \mathfrak{g} - \mathfrak{t}$. The latter obeys an inhomogeneous differential equation that can be written schematically as

$$\mathcal{L}_m \mathfrak{q} - \chi_{4\mu} \mathrm{P}(\mathfrak{q}) = \mathfrak{d} \ , \tag{3.32}$$

where $\mathrm{P}$ is defined momentarily and $\mathfrak{d} = \chi_{4\mu} \mathrm{P}(\mathfrak{g} + \mathfrak{t})$. The homomorphism $\mathrm{P}$ of $\mathbb{V} \oplus \mathbb{V}$ sends a given element $\mathfrak{s} = ((\mathfrak{s}_1, \mathfrak{s}_{01}), (\mathfrak{s}_2, \mathfrak{s}_{02}))$ to $\mathrm{P}(\mathfrak{s}) = ((\mathrm{P}_1(\mathfrak{s}), \mathrm{P}_{01}(\mathfrak{s})), (\mathrm{P}_2(\mathfrak{s}), \mathrm{P}_{02}(\mathfrak{s}))$ where

- $*\mathrm{P}_1(\mathfrak{s}) = m(e_\Delta - e) \wedge [\tau, \mathfrak{s}_2] - \langle \tau\, \hat{\mathrm{a}}_{\mathrm{A}_\Delta} \rangle \wedge [\tau, \mathfrak{s}_1]$ .
- $*\mathrm{P}_2(\mathfrak{s}) = -m(e_\Delta - e) \wedge [\tau, \mathfrak{s}_1] + \langle \tau\, \hat{\mathrm{a}}_{\mathrm{A}_\Delta} \rangle \wedge [\tau, \mathfrak{s}_2]$ .
- $*\mathrm{P}_{01}(\mathfrak{s}) = *(\mathrm{d}(\pi_\phi \mathfrak{s}_1) - m e \wedge [\tau, \pi_\phi \mathfrak{s}_2]) + m(e_\Delta - e) \wedge [\tau, *_\phi \mathfrak{s}_2] + \langle \tau\, \hat{\mathrm{a}}_{\mathrm{A}_\Delta} \rangle \wedge [\tau, *_\phi \mathfrak{s}_1]$ .
- $*\mathrm{P}_{02}(\mathfrak{s}) = *(-\mathrm{d}(\pi_\phi \mathfrak{s}_2) - m e \wedge [\tau, \pi_\phi \mathfrak{s}_1]) + m(e_\Delta - e) \wedge [\tau, *_\phi \mathfrak{s}_1] - \langle \tau\, \hat{\mathrm{a}}_{\mathrm{A}_\Delta} \rangle \wedge [\tau, *\mathfrak{s}_2]$ .

(3.33)

As explained momentarily, the operator $\mathcal{L}_m - \chi_{4\mu}\mathrm{P}$ is invertible if $\rho < c_0^{-1}$ and $m > c_0 \mathrm{M}^2$. This implies in particular that (3.32) has a unique solution.

The asserted invertibility of $\mathcal{L}_m - \chi_{4\mu}\mathrm{P}$ is proved using the bounds stated below in (3.34). The notation in (3.34) has $\mathfrak{v}$ denoting a given $\mathrm{L}^2_1$ section of $T^*\mathbb{R}^3 \otimes \mathfrak{su}(2)$.

- $m \| \chi_{4\mu} |e_\Delta - e|\, \mathfrak{v} \|_2 \le c_0 \rho\, m \| \mathfrak{v} \|_2$.
- $\| \chi_{4\mu} |\langle \tau\, \hat{\mathrm{a}}_{\mathrm{A}_\Delta} \rangle|\, \mathfrak{v} \|_2 \le c_0 \mathrm{M} (m^{-1} \|\nabla \mathfrak{v}\|_2 + m \| \mathfrak{v} \|_2)$.
- $\| \chi_{4\mu} \mathrm{d}(\pi_\phi \mathfrak{v}) \|_2 \le c_0 \rho (\| \nabla \mathfrak{v} \|_2 + \| \mathfrak{v} \|_2)$ *and* $m \| \chi_{4\mu} \pi_\phi \mathfrak{v} \|_2 \le m \rho \| \mathfrak{v} \|_2$.

(3.34)

These bounds are derived in the next paragraph. The assertion that $\mathcal{L}_m - \chi_{4\mu}\mathrm{P}$ is invertible follows directly from the following assertion: If $\rho < c_0^{-1}$ and $m > c_0 \mathrm{M}^2$, then

$$\| (\mathcal{L}_m - \chi_{4\mu} \mathrm{P}) \mathfrak{k} \|_2^2 \ge \tfrac{3}{4} (\| \nabla \mathfrak{k} \|_2^2 + m^2 \| \mathfrak{k} \|_2^2) \tag{3.35}$$

when $\mathfrak{k}$ is any given $\mathrm{L}^2_1$ map from $\mathbb{R}^3$ to $\mathbb{V} \oplus \mathbb{V}$. This last assertion follows from (3.34) because the (3.34) implies the following: If $\mathfrak{k}$ is an $\mathrm{L}^2_1$ map from $\mathbb{R}^3$ to $\mathbb{V} \oplus \mathbb{V}$, then

$$\| \chi_{4\mu} \mathrm{P}(\mathfrak{k}) \|_2^2 \le c_0 (m^{-1} \mathrm{M} + \rho^2) (\| \nabla \mathfrak{k} \|_2^2 + m^2 \| \mathfrak{k} \|_2^2) \ . \tag{3.36}$$

If $\rho < c_0^{-1}$ and $m > c_0 \mathrm{M}$, then the right hand side is no greater than $\frac{1}{100} \| \mathcal{L}_m k \|_2^2$. The latter bound with (3.29) lead directly to the bound in (3.35).

The assertions of the first and third bullets in (3.34) follow directly from the bounds in the respective first and fifth bullets of (3.27). To see about the middle bullet, keep in mind that the $L^2$ norm of $\chi_{4\mu}|\langle \tau\,\hat{a}_{A_\Delta}\rangle|\,|\mathfrak{v}|$ is bounded by the product of their $L^4$ norms. Use the $f = |\mathfrak{v}|$ version of the top bullet in (3.11) to bound the $L^4$ norm of $\mathfrak{v}$. Having done so, then the inequality in the middle bullet of (3.34) follows directly from the bound in the third bullet of (3.27).

<u>Step 5</u>: Assume henceforth that $\rho < c_0^{-1}$ and that $m > c_0 M$ so as to use (3.29) to bound the right hand side of (3.36) by $\frac{1}{100}\|\mathcal{L}_m \mathfrak{k}\|_2^2$. This being the case, then (3.35) holds. Use the $\mathfrak{k} = \mathfrak{q}$ version of (3.35) with (3.32) to conclude that

$$\|\nabla \mathfrak{q}\|_2^2 + m^2\|\mathfrak{q}\|_2^2 = \tfrac{4}{3}\|\chi_{4\mu}\,P(\mathfrak{g}+\mathfrak{t})\|_2^2\ . \tag{3.37}$$

To exploit this last bound, first apply the bound for $\|P(\cdot)\|_2 \le \frac{1}{10}\|\mathcal{L}_m(\cdot)\|_2$ to conclude that

$$\|\nabla \mathfrak{q}\|_2 + m\|\mathfrak{q}\|_2 \le \tfrac{1}{10}\left(\|\mathcal{L}_m \mathfrak{g}\|_2 + \|\mathcal{L}_m \mathfrak{t}\|_2\right)\ . \tag{3.38}$$

Use (3.44) to bound $\|\mathcal{L}_m \mathfrak{t}\|_2$ by $c_0\, m^{-1}\rho$. Use the identity $\mathcal{L}_m \mathfrak{g} = \mathfrak{S}(d\chi_{4\mu})\,\mathfrak{f}$ to bound $\|\mathcal{L}_m \mathfrak{g}\|_2$ by $c_0\,\mu^{-1}\|\chi_\mu \mathfrak{f}\|_2$. Substituting these bounds on the right hand side of (3.38) finds

$$\|\nabla \mathfrak{q}\|_2 + m\|\mathfrak{q}\|_2 \le c_0\left(m^{-1}\rho + \mu^{-1}\|\chi_\mu \mathfrak{f}\|_2\right) \tag{3.39}$$

when $\rho < c_0^{-1}$ and $m > c_0 M^2$.

<u>Step 6</u>: Write $\mathfrak{f}$ on the $|x| \le 1 - 4\mu$ part of $\mathbb{R}^3$ as $\mathfrak{f} = \mathfrak{g} + \mathfrak{t} + \mathfrak{q}$ to conclude that

$$\|\chi_{16\mu}\nabla \mathfrak{f}\|_2 + m\|\chi_{16\mu}\mathfrak{f}\|_2 \le \|\chi_{16\mu}\nabla \mathfrak{g}\|_2 + m\|\chi_{16\mu}\mathfrak{g}\|_2 + \|\nabla \mathfrak{t}\|_2 + m\|\mathfrak{t}\|_2 + \|\nabla \mathfrak{q}\|_2 + m\|\mathfrak{q}\|_2\ . \tag{3.40}$$

Use (3.31) to bound the terms in (3.40) with $\mathfrak{t}$ and thus by $m^{-1}\rho$. Use (3.39) to bound those with $\mathfrak{q}$ by $c_0(m^{-1}\rho + \mu^{-1}\|\chi_\mu \mathfrak{f}\|_2)$. The pointwise norms in (3.29) lead to a bound on the terms with $\mathfrak{g}$ by $c_0\,\mu^{-4}\,e^{-m\mu/4}\rho$. These substitutions imply the bound

$$\|\chi_{16\mu}\nabla \mathfrak{f}\|_2 + m\|\chi_{16\mu}\mathfrak{f}\|_2 \le c_0\,\mu^{-1}\left(\mu^{-4}\,e^{-m\mu/4}\,m + m^{-1}\right)\rho + c_0\,\mu^{-1}\|\chi_\mu \mathfrak{f}\|_2\ . \tag{3.41}$$

Note in particular that $\mu^{-4}\,e^{-m\mu/4}\,m \le m^{-1}$ when $m > c_0\,\mu^{-2}$. Granted that such is the case, then (3.41) implies that

$$\|\chi_{16\mu}\nabla\mathfrak{f}\|_2 + m\|\chi_{16\mu}\mathfrak{f}\|_2 \le c_0\,\mu^{-1}(m^{-1}\rho + \|\chi_\mu\mathfrak{f}\|_2).$$

(3.42)

The first and second bullets of (3.27) bound $\|\chi_\mu\mathfrak{f}\|_2$ by $c_0\rho$ and so the right hand side of (3.42) is no larger than $c_0\mu^{-1}\rho$.

<u>Step</u> <u>7</u>: Assume henceforth that $m > c_0\mu^{-2}$ so as to conclude from (3.42) that

$$\|\chi_{16\mu}\mathfrak{f}\|_2 \le c_0\mu^{-1}m^{-1}\rho.$$

(3.43)

With (3.43) in hand, repeat Steps 1-6 but with $\mu$ replaced by $\mu' = 128\mu$. The salient difference is the replacement of the factor $\|\chi_\mu\mathfrak{f}\|_2$ in (3.39), (3.41) and (3.42) by $\|\chi_{16\mu}\mathfrak{f}\|_2$. Where as the former can be bounded only by $c_0\rho$, the latter is bounded courtesy of (3.43) by $c_0 m^{-1}\rho$, this being an improvement by a factor of $m^{-1}$. Granted this replacement, then the $\mu' = 128\mu$ version of (3.42) reads

$$\|\chi_{2048\mu}\nabla\mathfrak{f}\|_2 + m\|\chi_{2048\mu}\mathfrak{f}\|_2 \le c_0\mu^{-1}m^{-1}\rho.$$

(3.44)

The assertion made by Lemma 3.4 follows directly from (3.44).

## 4. Unexpectedly small curvature

Lemma 3.3 and a version of Lemma 3.4 are brought to bear in Section 4d to prove Proposition 3.2. The intervening subsections supply data that can be used by Lemma 3.4. The latter requires as input a parameter $m$, a connection on the product bundle over the $|x| < 1$ ball, and an $\mathfrak{su}(2)$ valued 1-form over this ball, these being $A_\Delta$ and $a_\Delta$. The parameter $m$ is taken to be $m = \sqrt{\frac{3}{4\pi}}\, r_z$ and the connection $A_\Delta$ is chosen to have the form $h^*A_\lozenge$ with h being a suitable isomorphism from the product bundle on the $|x| < 1$ ball to the bundle $\phi^*P$ over this ball. The $\mathfrak{su}(2)$ valued 1-form $a_\Delta$ will be the pull back via h of a suitable perturbation of $\hat{a}_*$. Section 4a constructs this perturbation and Section 4c constructs h. The intervening subsection says more about the perturbation.

### a) The heat equation on the |x| < 1 ball

The three parts of this subsection modify $\hat{a}_*$ over the $|x| < 1$ ball in $\mathbb{R}^3$ so as to obtain a $\phi^*(P\times_{SO(3)}\mathfrak{su}(2))$ valued 1-form whose pointwise norm and $L^2_2$ norm on concentric balls with radius less than 1 obey a priori bounds that can not be assumed to hold for $\hat{a}_*$. The modified version is denoted by $\tilde{a}_*$. The following proposition summarizes the salient features of $\tilde{a}_*$. The proposition uses the metric $\mathfrak{m}_\phi$ to define the norms, volume form, Hodge star and covariant derivatives on tensors.

**Proposition 4.1**: *Given* $\mathrm{E} \geq 1$, *there exists* $\kappa_{\mathrm{E}} > \kappa$, *and given also* $\mu \in (0, \frac{1}{2}]$, *there exists* $\kappa_{\mathrm{E},\mu} > 1$*; these having the following significance: Suppose that* $r > 1$ *and* $(\mathrm{A}, \hat{a})$ *is a pair of connection on* P *and section of* $(\mathrm{P} \times_{\mathrm{SO}(3)} \mathfrak{su}(2)) \otimes T^*M$ *that obey (3.2). Fix* $\mathrm{p} \in \mathrm{M}$ *such that* $r_z = z\, \mathrm{r}_\Diamond\, r \geq 1$. *There exists a* $\phi^*(\mathrm{P} \times_{\mathrm{SO}(3)} \mathfrak{su}(2))$ *valued 1-form on the* $|\mathrm{x}| \leq 1$ *ball in* $\mathbb{R}^3$ *to be denoted by* $\tilde{a}_*$ *with the properties in the list below*.

- $\int_{|\mathrm{x}|\leq 1} (|d_{\mathrm{A}_\Diamond} \tilde{a}_*|^2 + |d_{\mathrm{A}_\Diamond} {*\tilde{a}_*}|^2) < \kappa_{\mathrm{E}}\, \mathrm{r}_\Diamond\, r_z^{-2}$ .
- $\int_{|\mathrm{x}|\leq 1} |\hat{a}_* - \tilde{a}_*|^2 < \kappa_{\mathrm{E}}\, \mathrm{r}_\Diamond\, r_z^{-4}$.
- $\int_{|\mathrm{x}|\leq 1} |\nabla_{\mathrm{A}_\Diamond} (\hat{a}_* - \tilde{a}_*)|^2 < \kappa_{\mathrm{E}}\, \mathrm{r}_\Diamond r_z^{-2}$.
- $\int_{|\mathrm{x}|\leq 1} |q_{\mathrm{A}_\Diamond}(\tilde{a}_*)|^2 < \kappa_{\mathrm{E}}\, \mathrm{r}_\Diamond$ .
- *If* $\mu \in (0, \frac{1}{2}]$, *then* $\int_{|\mathrm{x}|\leq 1-\mu} |\nabla_{\mathrm{A}_\Diamond} (\nabla_{\mathrm{A}_\Diamond} \tilde{a}_*)|^2 < \kappa_{\mathrm{E},\mu}$ .
- *If* $\mu \in (0, \frac{1}{2}]$, *then* $\sup_{|\mathrm{x}|\leq 1-\mu} |\tilde{a}_*| < \kappa_{\mathrm{E},\mu}$ .
- *If* $\mu \in (0, \frac{1}{2}]$, *then* $\int_{|\mathrm{x}|\leq 1-\mu} |F_{\mathrm{A}_\Diamond} - r_z^{\,2}\, \tilde{a}_* \wedge \tilde{a}_*|^2 \leq \kappa_{\mathrm{E},\mu}\, \mathrm{r}_\Diamond$ .

***Proof of Proposition 4.1***: The desired $\tilde{a}_*$ is constructed from $\hat{a}_*$ by mimicking what is done in Sections 2c and 2d. The construction has six steps. As in the statement of the proposition, these steps implicitly use the metric $\mathfrak{m}_\phi$ to define the norms, the volume form, the Hodge star and the covariant derivatives on tensors. This metric is also used to define the formal, $L^2$ adjoint of the covariant derivative.

Step 1: The section $\tilde{a}_*$ is obtained from the solution to an analog of the heat equation in (2.52). The heat equation in this case specifies a $\phi^*(\mathrm{P} \times_{\mathrm{SO}(3)} \mathfrak{su}(2))$ valued 1-form over the product of $[0, \infty)$ with the $|\mathrm{x}| \leq 1$ ball in $\mathbb{R}^3$. The 1-form in question is denoted by $a$ and obeys the following:

- $\frac{\partial}{\partial t} a = -(\nabla_{\mathrm{A}_\Diamond}^{\dagger} \nabla_{\mathrm{A}_\Diamond} a + {*({*F_{\mathrm{A}_\Diamond}}} \wedge a + a \wedge {*F_{\mathrm{A}_\Diamond}})) + \mathrm{Ric}_{\mathfrak{m}_\phi}((\cdot) \otimes a)$ *where* $|\mathrm{x}| < 1$.
- $a|_{t=0} = \hat{a}_*$ *for all* $|\mathrm{x}| \leq 1$.
- $a|_{|\mathrm{x}|=1} = \hat{a}_*|_{|\mathrm{x}|=1}$ *for all* $t \geq 0$.

(4.1)

What is denoted by $\mathrm{Ric}_{\mathfrak{m}_\phi}$ in the top bullet of (4.1) is the Ricci curvature tensor of $\mathfrak{m}_\phi$. Note in particular that its pointwise norm is bounded by $c_0\, \mathrm{r}_\Diamond^{\,2}$. The third bullet in (4.1) specifies both tangential and normal components of $a$ on the boundary $|x| = 1$ sphere. Standard results about parabolic equations prove that there is a unique solution to (4.1).

The desired $\tilde{a}_*$ is given by $a|_s$ for an appropriate stopping time $s \in [0, \infty)$.

Step 2: The analog of the function that is denoted by $\mathfrak{E}$ in Part 1 of Section 2d is the function on $[0, \infty)$ given by

$$\mathfrak{E} = \int_{|x|\le 1} (|d_{A_\Diamond} a|^2 + |d_{A_\Diamond} {*a}|^2) \ . \tag{4.2}$$

Since $a$ is constant on the $|x| = 1$ sphere, integration by parts writes

$$\tfrac{d}{dt}\mathfrak{E} = - 2 \, \|q_{A_\Diamond}(a)\|_2^2 \ , \tag{4.3}$$

where $\|\cdot\|_2$ denotes the $L^2$ norm on the $|x| < 1$ sphere and where $q_{A_\Diamond}(a)$ is defined by writing the top bullet of (4.1) schematically as $\frac{\partial}{\partial t} a = -q_{A_\Diamond}(a)$. Note that the integration by parts has no boundary contribution because $a = \hat{a}_*$ for all t where $|x| = 1$.

Step 3: The first bullet in (3.9) implies that $\mathfrak{E}(\hat{a}_*) \le c_E r_\Diamond r_z^{-2}$, and so (4.3) finds

$$\mathfrak{E}(t) \le c_0 r_\Diamond r_z^{-2} - t\,\mathfrak{n}(t) \ , \tag{4.4}$$

where $\mathfrak{n}(t)$ is defined by analogy with (2.54) as

$$\mathfrak{n}(t) = t^{-1} \int_0^t \|q_{A_\Diamond}(a)\|_2^2 \ . \tag{4.5}$$

The inequality in (4.5) implies the assertion in the first bullet of Proposition 4.1 if $\tilde{a}_*$ is defined to be any $s \ge 0$ version of $a|_s$. In any event, (4.4) requires that $\mathfrak{n}(t) \le c_E r_\Diamond r_z^{-2} t^{-1}$ if t is positive. Meanwhile, the formula in the first bullet of (4.1) for $\frac{\partial}{\partial t} a$ implies that

$$\|\hat{a}_* - a|_t\|_2^2 \le t^2 \mathfrak{n}(t) \quad \textit{for any } t \in [0, \infty). \tag{4.6}$$

Step 4: The definition in (4.5) with the bound $\mathfrak{n}(t) \le c_E r_\Diamond r_z^{-2} t^{-1}$ have the following consequence: There exists $s \in [\frac{1}{2} t, t]$ such that $a|_s$ obeys

$$\| q_{A_\Diamond}(a|_s)\|_2^2 \le c_E r_\Diamond r_z^{-2} t^{-1} \,.$$

(4.7)

This understood, take $t = r_z^{-2}$ and fix $s \in [\frac{1}{2}t, t]$ so that (4.7) holds. Set $\tilde{a}_*$ to equal $a|_s$. It follows from (4.7) that $\| q_{A_\Diamond}(\tilde{a}_*)\|_2 \le c_E r_\Diamond^{1/2}$ and (4.6) finds $\| \hat{a}_* - \tilde{a}_{*n}\|_2 \le c_E r_\Diamond^{1/2} r_z^{-2}$. The former bound is the assertion in the fourth bullet and the latter bound is the assertion in the second bullet of the proposition.

To see about the third bullet of the proposition, use the fact that both $\mathfrak{E}(\hat{a}_*)$ both $\mathfrak{E}(\tilde{a}_*)$ are bounded by $c_E r r_z^{-2}$ to see that

$$\int_{|x|\le 1} (|d_{A_\Diamond}(\hat{a}_* - \tilde{a}_*)|^2 + |d_{A_\Diamond} {*}(\hat{a}_* - \tilde{a}_*)|^2) \le c_E r_\Diamond r_z^{-2}$$

(4.8)

also. This understood, an integration by parts with the fact that $\hat{a}_* = \tilde{a}_*$ where $|x| = 1$ leads from (4.8) to the bound

$$\int_{|x|\le 1} (|\nabla_{A_\Diamond}(\hat{a}_* - \tilde{a}_*)|^2 + 2\left\langle {*}F_{A_\Diamond} \wedge (\hat{a}_* - \tilde{a}_*) \wedge (\hat{a}_{*n} - \tilde{a}_{*n})\right\rangle) \le c_E r_\Diamond^2 \|\hat{a}_* - \tilde{a}_*\|_2^2 + c_E r_\Diamond r_z^{-2} \,.$$

(4.9)

In turn, this last inequality, the bound on the $L^2$ norm of $F_{A_\Diamond}$ and the Sobolev inequality given in the top bullet in (3.11) lead to the bound

$$\int_{|x|\le 1} |\nabla_{A_\Diamond}(\hat{a}_* - \tilde{a}_*)|^2 \le c_E r_\Diamond r_z^{-2} \,.$$

(4.10)

This is the bound that is asserted by the third bullet of the proposition

<u>Step 5</u>: To see about the fifth and sixth bullets, fix $\mu \in (0, \frac{1}{2}]$. Set $\chi_\mu$ to denote again denote the function on $\mathbb{R}^3$ that is given by $x \to \chi(\frac{1}{\mu}(|x| - 1 + \mu))$. The bound on the $L^2$ norm of $q_{A_\Diamond}(\tilde{a}_*)$, the bound on $|\mathrm{Ric}_{m_\phi}|$ by $c_0 r_\Diamond^2$ and the fact that (4.10) implies that $\|\nabla_{A_\Diamond} \tilde{a}_*\|_2 \le c_E$ implies via an integration by parts that

$$\|\nabla_{A_\Diamond}(\nabla_{A_\Diamond}(\chi_\mu \tilde{a}_*))\|_2^2 \le c_{E,\mu} + c_0 \|F_{A_\Diamond}\|_2^2 \|\chi_\mu \tilde{a}_*\|_\infty^2 \,,$$

(4.11)

Note that (4.11) also uses the fact that $A_\Diamond = \theta + \hat{a}_{A_\Diamond}$ with the bound given by Uhlenbeck's theorem for the $L^2_1$ norm of $\hat{a}_{A_\Diamond}$. Since $\|F_{A_\Diamond}\|_2 \le 1$, the bound in (4.11) implies that

$$\|d|\nabla_{A_\Diamond}(\chi_\mu\tilde{a}_*)|\|_2 \le c_{E,\mu} + c_0\|\chi_\mu\tilde{a}_*\|_\infty\ ,$$

(4.12)

and so $|\nabla_{A_\Diamond}(\chi_\mu\tilde{a}_*)|$ is an $L^2_1$ function with $L^2_1$ norm bounded by $(c_{E,\mu} + c_0\|\chi_\mu\tilde{a}_*\|_\infty)$.

Invoke the top bullet in (3.11) yet again to see that the $L^4$ norm of $|\nabla_{A_\Diamond}(\chi_\mu\tilde{a}_*)|$ is also bounded by $(c_{E,\mu} + c_0\|\chi_\mu\tilde{a}_*\|_\infty)$, and so this is also the case for the $L^4$ norm of $d|\chi_\mu\tilde{a}_*|$. With the preceding understood, invoke the the third bullet in (3.11) using $|\chi_\mu\tilde{a}_*|$ for $f$ and then invoke the first bullet of (3.11) using $|\nabla_{A_\Diamond}(\chi_\mu\tilde{a}_*)|$ conclude the following: Fix $\delta \in (0, 1)$. Then

$$\|\chi_\mu\tilde{a}_*\|_\infty \le \delta^{-1}(c_{E,\mu} + \|\nabla_{A_\Diamond}\tilde{a}_*\|_2) + \delta\|\nabla_{A_\Diamond}(\nabla_{A_\Diamond}(\chi_\mu\tilde{a}_*))\|_2\ .$$

(4.13)

A $\delta = c_0^{-1}$ version of (4.13), the fact that $\|F_{A_\Diamond}\|_2 \le 1$ with (3.27) and (3.26) lead directly to the assertion in the fifth bullet of Proposition 4.1. The assertion in the sixth bullet follows from the $\delta = 1$ version of (4.13) and the assertion in the fifth bullet.

<u>Step 6</u>: This step proves the assertion in the seventh bullet of the proposition. Given what is said by (3.16), it is only necessary to prove that

$$r_z^4\|\chi_\mu((\hat{a}_* \wedge \hat{a}_*) - (\tilde{a}_* \wedge \tilde{a}_*))\|_2^2 \le c_r r_\Diamond\ .$$

(4.14)

The left hand side of (4.14) is no greater than

$$c_0\, r_z^4\|\,|\hat{a}_* - \tilde{a}_*|^2\|_2^4 + c_0\, r_z^4\|\hat{a}_* - \tilde{a}_*\|_2^2\|\chi_\mu\tilde{a}_*\|_\infty^2\ .$$

(4.15)

Use the first bullet in (3.11) with the first two bullets of the proposition to see that the left most term in (4.15) is no greater than $c_E r_\Diamond^{5/2} r_\Diamond^{-1}$. Invoke the first bullet of the proposition to bound the right most term by $c_E r_\Diamond\|\chi_\mu\tilde{a}_*\|_\infty^2$; then use the proposition's sixth bullet to bound $c_0 r_\Diamond\|\chi_\mu\tilde{a}_*\|_\infty^2$ by $c_{E,\mu} r_\Diamond$.

**b) The decomposition of $\tilde{a}_*$**

Lemma 3.4 refers to parameters $m$, $\rho$ and $\mathrm{M}$. As noted at the outset, $m$ will be taken to be $\sqrt{\frac{3}{4\pi}}\, r_z$. The parameters $\rho$ and $\mathrm{M}$ are such that (3.27) holds. An appeal to a given $\varepsilon > 0$ version of Lemma 3.4 requires an upper bound for $\rho$. The upcoming Lemma 4.2 is the first of two lemmas that are used to obtain the required upper bound. It proves useful for notational reasons to fix an isomorphism between $\phi^*P$ and the product principal SO(3) bundle over the |x| < 1 ball. Having done so, view $A_\Diamond$ henceforth as a connection on this product bundle over the |x| < 1 ball and to view $\tilde{a}_*$ as a $\mathfrak{su}(2)$ valued 1-form on this same ball.

**Lemma 4.2**: *There exists* $M_0 > 1$, *and given* $E \geq 1$, $\mu \in (0, \frac{1}{2}]$ *and* $\rho \in (0, 1]$, *there exists* $\kappa_\rho > 1$ *with the following significance: Suppose that* $r > 1$ *and* $(A, \hat{a})$ *is a pair of connection on* P *and section of* $(P \times_{SO(3)} \mathfrak{su}(2)) \otimes T^*M$ *that obey (3.2). Fix* $p \in M$ *where* $r_z \geq \kappa_\rho$. *Suppose in addition that* $r_\Diamond < \kappa_\rho^{-1}$ *and that* $\| \nabla_{A_\Diamond} \tilde{a}_* \|_2 < \kappa_\rho^{-1}$. *There is a constant, unit length 1-form* $\mathfrak{e}$ *on* $\mathbb{R}^3$ *and a map* $\tau_\Diamond$ *from the* $|x| \leq 1 - \mu$ *ball in* $\mathbb{R}^3$ *to the unit sphere in* $\mathfrak{su}(2)$ *with the properties that are listed below.*

- *Write* $\tilde{a}_*$ *as* $\sqrt{\frac{3}{4\pi}}\,(\tau_\Diamond e_\Delta + \mathfrak{b}_*)$ *with* $\langle \tau_\Diamond \mathfrak{b}_* \rangle = 0$. *The* $\mathfrak{m}_\phi$ *metric pairing between* $\mathfrak{e}_\Delta$ *and* $\mathfrak{b}_*$ *is zero. Moreover,* $\sup_{|x| \leq 1-\mu} |e_\Delta - e| < \rho$ *and* $\sup_{|x| < 1-\mu} |\mathfrak{b}_*| < \kappa_\rho r_z^{-2}$.
- *Neither the norm of the covariant derivative* $\mathfrak{e}_\Delta$ *nor that of either the* $A_\Diamond$ *covariant derivative of* $\tau_\Diamond$ *or the* $A_\Diamond$ *covariant derivative of* $\mathfrak{b}$ *are larger than* $10\, |\nabla_{A_\Diamond} \tilde{a}_*|$.

***Proof of Lemma 4.2***: The proof has six steps. The first five steps prove the top bullet and the sixth step proves the lower bullet. These steps use $\delta$ to denote $\| \nabla_{A_\Diamond} \tilde{a}_* \|_2$.

<u>Step 1</u>: Use the top bullet in (3.11) with the bounds in Proposition 4.1's third and fifth bullets to bound the $L^4$ norm of $|\nabla_{A_\Diamond} \tilde{a}_*|$ by $c_{E,\mu}\, \delta^{1/4}$ on the $|x| \leq 1 - \mu$ ball. Let $e$ denote for the moment any given constant, unit length 1-form on $\mathbb{R}^d$ and let $\tilde{a}_*(e)$ denote the pairing between $\tilde{a}_*(e)$ and $e$. The $c_{E,\mu}\, \delta^{1/4}$ bound on the $L^4$ norm of $|\nabla_{A_\Diamond} \tilde{a}_*|$ leads to a corresponding $c_{E,\mu}\, \delta^{1/4}$ bound on the $L^4$ norm of $|d|\tilde{a}_*(e)||$. This in turn implies via a Sobolev inequality that $|\tilde{a}_*(e)|$ is Holder continuous for an exponent $\upsilon \in (c_0^{-1}, \frac{1}{4})$ and that

$$\big|\, |\tilde{a}_*(e)|(x) - |\tilde{a}_*(e)|(y) \big| \leq c_{E,\mu}\, \delta^{1/c_{E,\mu}}\, |x - y|^{\upsilon} \tag{4.16}$$

if both |x| and |y| are less than $1 - \mu$.

<u>Step 2</u>: This step explains why $\big|\, |\tilde{a}_*| - \sqrt{\frac{3}{4\pi}}\ \big| \leq c_{E,\mu}\, (\delta^{1/4} + r_\Diamond r_z^{-2})$ where $|x| < 1 - \mu$. To prove this, suppose for the moment that $c > 1$ is such that the $|\tilde{a}_*|(0) < \sqrt{\frac{3}{4\pi}} - c^{-1}$. If so, then (4.16) implies that $|\tilde{a}_*|^2 < (\sqrt{\frac{3}{4\pi}} - c^{-1} + c_{E,\mu}\, \delta^{1/4})^2$ on the whole $|x| < 1 - \mu$ ball. It follows that the integral of $|\tilde{a}_*|^2$ over this ball is no greater than $(1 - \sqrt{\frac{4\pi}{3}}\, c^{-1} + c_{E,\mu} \delta^{1/4})^2 (1 - \mu)^3$. Moreover, it follows that the integral of $|\tilde{a}_*|^2$ on the $|x| = 1 - \mu$ sphere is no greater than

$3(1-\sqrt{\frac{4\pi}{3}}\, c^{-1} + c_{E,\mu}\delta^{1/4})^2(1-\mu)^2$. Granted this last bound, it then follows using the fundamental theorem of calculus that the integral of $|\tilde{a}_*|^2$ over the spherical annulus where $1-\mu \le |x| \le 1$ is no greater than $3(1-\sqrt{\frac{4\pi}{3}}\, c^{-1} + c_{E,\mu}\delta^{1/4})^2\mu + c_0\,\delta\mu^{1/2}$. Add these various bounds to conclude that $\|\tilde{a}_*\|_2^2 \le (1 - c_0 c^{-1} + c_{E,\mu}\delta^{1/4})^2 + c_0\,\delta\mu^{1/2}$. This is nonsense if $c^{-1}$ is less than $c_{E,\mu}(\delta^{1/4} + r_\Diamond r_z^{-2})$ because Proposition 4.1's second bullet implies that the $L^2$ norm of $\tilde{a}_*$ is no less than $1 - c_0 r_\Diamond r_z^{-2}$. Very much the same argument derives nonsense if $|\tilde{a}_*|(0) > \sqrt{\frac{3}{4\pi}} + c_{E,\mu}(\delta^{1/4} + r_\Diamond r_z^{-2})$.

Step 3: It follows from what is said in Steps 1 and 2 that the $L^4$ norm of the function $|\nabla_{A_\Diamond}(\tilde{a}_* \wedge \tilde{a}_*)|$ on the $|x| \le 1-\mu$ ball is also bounded by $c_{E,\mu}\,\delta^{1/4}$. This implies the same for the $L^4$ norm of $d|\tilde{a}_* \wedge \tilde{a}_*|$; and so the function $|\tilde{a}_* \wedge \tilde{a}_*|$ is also Holder continuous with exponent $\upsilon \in (c_0^{-1}, \frac{1}{4})$. In addition,

$$||\tilde{a}_* \wedge \tilde{a}_*|(x) - |\tilde{a}_* \wedge \tilde{a}_*|(y)| \le c_{E,\mu}\,\delta^{1/c_{E,\mu}}\,|x-y|^\upsilon \tag{4.17}$$

if $|x|$ and $|y|$ are less than $1-\mu$. Since the integral of $|\tilde{a}_* \wedge \tilde{a}_*|^2$ over the $|x| < 1-\mu$ ball is less than $c_0 r_z^{-4}$, it follows from (4.17) that $|\tilde{a}_* \wedge \tilde{a}_*| \le c_{E,\mu} r_z^{-2}$ at each point in the $|x| < 1-\mu$ ball.

Step 4: Fix a point $x \in \mathbb{R}^3$ with $|x| < 1 - \mu$ and fix an $\mathfrak{m}_\phi$ orthonormal frame for $T^*\mathbb{R}^3$ at x. Use this frame to write the components of $\tilde{a}_*|_x$ as $\{\tilde{a}_{*a}\}_{a=1,2,3}$ and use these components to define a symmetric, non-negative endomorphism of $\mathfrak{su}(2)$ by the rule

$$\sigma \to \mathbb{T}(\sigma) = \tilde{a}_{*a}\langle \tilde{a}_{*a}\sigma\rangle\,, \tag{4.18}$$

with it understood that repeated indices are to be summed. This endomorphism does not depend on the chosen orthonormal frame. As explained in the next paragraph, if $\delta$ is less than $c_{E,\mu}^{-1}$ and $r_\Diamond \le c_{E,\mu}^{-1}$ and $r_z \ge c_{E,\mu}$, then $\mathbb{T}$ has only one eigenvalue greater than the minimum of $\frac{1}{1000}$ and $c_0 r_z^{-4}$, and the corresponding eigenspace has dimension one. Note in this regard that it follows from what is said in Step 1 that $\mathbb{T}$ has at least one eigenvalue greater than $\frac{1}{4\pi}(1 - c_{E,\mu}\,\delta^{1/4} - r_\Diamond r_z^{-2})$; and this is greater than $\frac{1}{8\pi}$ when $\delta < c_{E,\mu}^{-1}$ and $r_z > c_{E,\mu}$.

Suppose that $c > 1$ and that $\mathbb{T}$ is greater than $c\, r_z^{-4}$ on a two dimensional vector space. Use $\tau_\Diamond$ to denote a length 1 eigenvector of $\mathbb{T}$ with eigenvalue greater than $\frac{1}{8\pi}$ and use $\sigma_\Diamond$ to denote an orthogonal, length 1 eigenvector of $\mathbb{T}$ with eigenvalue greater than $c\, r_z^{-4}$. Given that $\mathbb{T}$ is symmetric, the eigenvalue equations imply that the 1-form $\langle \tau_\Diamond \tilde{a}_*\rangle$ is $\mathfrak{m}_\phi$-orthogonal to $\langle \sigma_\Diamond \tilde{a}_* \rangle$. This being the case, it then follows that $|\tilde{a}_* \wedge \tilde{a}_*|(x) \ge c_0^{-1} c\, r_z^{-2}$. What with (4.17), such a lower bound with (4.17) runs afoul of (3.6) if $c > c_0^{-1}$.

Step 5: The map $\tau_\Diamond$ for Lemma 4.2 is defined so that at any given point, it is a norm 1 eigenvector of $\mathbb{T}$ with eigenvalue greater than $\frac{1}{8\pi}$. Define $e_\Delta$ to be $\sqrt{\frac{4\pi}{3}}\,\langle \tau_\Diamond \tilde{a}_* \rangle$. It follows from the definition of $\tau_\Diamond$ as an eigenvector of $\mathbb{T}$ and what is said in Step 4 about $\mathbb{T}$ having only one large eigenvalue that $|e_\Delta| > 1 - c_{E,\mu}(\delta^{1/4} + r_\Diamond r_z^2)$ and thus greater than $\frac{1}{2}$ when $\delta < c_{E,\mu}^{-1}$ and $r_z > c_{E,\mu}$. Let $\mathfrak{b}_* = \tilde{a}_* - \sqrt{\frac{3}{4\pi}}\, \tau_\Diamond e_\Delta$. The definition of $\tau_\Diamond$ as an eigenvector of $\mathbb{T}$ implies that $\mathfrak{b}_*$ is orthogonal to $\tau_\Diamond$ in $\mathfrak{su}(2)$ and that its $\mathfrak{m}_\phi$ metric pairing with $e_\Delta$ is zero. Moreover, the bound $|\mathfrak{b}_*| \le c_{E,\mu}\, r_z^{-2}$it follows from what Step 4 says about $\mathbb{T}$ having but one large eigenvalue. Let $e$ denote the constant 1-form given that is equal to the 1-form $|e_\Delta|^{-1} e_\Delta$ at the origin. The bound $|e - e_*| \le c_{E,\mu}\, \delta^{1/c_{E,\mu}}$ follows from (4.16).

Step 6: The assertion made by the second bullet of the lemma follows by writing the covariant derivative of $g^*\tilde{a}_*$ as

$$\nabla_{A_\Diamond} \tilde{a}_* = \sqrt{\tfrac{3}{4\pi}}\,(\nabla e_\Delta \tau_\Diamond + e_\Delta \nabla_{g^*A_\Diamond} \tau_\Diamond + \nabla_{g^*A_\Diamond} \mathfrak{b}_*)\ . \tag{4.19}$$

The fact that $\langle \tau_\Diamond \mathfrak{b}_* \rangle = 0$ and the fact that $\mathfrak{b}_*$ has zero metric pairing with $e_\Delta$ can be used to see that

$$|\nabla_{A_\Diamond} \tilde{a}_*|^2 \ \ge\ \tfrac{1}{6}\,(1 - c_0|\mathfrak{b}|^2)\,(|\nabla e_\Delta|^2 + |e_\Delta|^2\, ||\nabla_{g^*A_\Diamond} \tau_\Diamond|^2\,) \ +\ \tfrac{1}{6}\,|\nabla_{g^*A_\Diamond} \mathfrak{b}_*|^2\ . \tag{4.20}$$

Given the bounds on $|\mathfrak{b}|$ in Lemma 4.2's first bullet, the inequality in (4.20) implies the claim that is made by Lemma 4.2's second bullet.

**c) The definition of h**

The first four parts of this subsection construct the automorphism h. The results of this construction are summarized in Part 5 by Lemma 4.3. Fix $\rho > 0$ and suppose in what follows that $r_\Diamond$ and $\|\nabla_{A_\Diamond} \tilde{a}_*\|_2$ are small enough and $r_z$ is large enough to invoke Lemma 4.2 with the given, small value of $\rho$.

*Part 1*: Define a connection $\hat{A}_*$ on the $|x| < 1 - \mu$ ball by setting it to be

$$\hat{A}_* = g^*A_\Diamond - \tfrac{1}{4}[\tau_\Diamond, \nabla_{g^*A_\Diamond}\tau_\Diamond] \,.$$

(4.21)

This connection is such that $\nabla_{\hat{A}_*}\tau_\Diamond = 0$. Note also that the curvature of $\hat{A}_*$ obeys

$$\|F_{\hat{A}_*}\|_2^2 \le c_0\|F_{A_\Diamond}\|_2^2 + c_0\|\nabla_{A_\Diamond}\tilde{a}_*\|_4^4$$

(4.22)

with it understood that the integrals that define the $L^2$ and $L^4$ norms on both sides of this equation are restricted to the $|x| \le 1-\mu$ ball.

*Part 2*: Since the ball where $|x| < 1-\mu$ is contractible, there exists an automorphism of the product SU(2) bundle over this ball that pulls back $\tau_\Diamond$ so as to be the constant, unit length element $\tau \in \mathfrak{su}(2)$. Choose such an automorphism for the moment and denote it by $u$. A point to keep in mind is that $u$ is not uniquely determined because there are automorphisms that fix $\tau$. These can be written as fiber multiplication by $e^{x\tau}$ with $x$ being an $\mathbb{R}$ valued function on the $|x| < 1-\mu$ ball. This flexibility with the choice of $u$ is used momentarily.

*Part 3*: Since $\nabla_{\hat{A}_*}\tau_\Diamond = 0$, it follows that $\tau$ is $u^*\hat{A}_*$-covariantly constant; and this implies that $u^*\hat{A}_*$ can be written as $u^*\hat{A}_* = \theta_0 + \hat{A}_*\tau$ with $\hat{A}_*$ being an $\mathbb{R}$ valued 1-form on the $|x| < 1-\mu$ ball. The $L^2$ norm of $d\hat{A}_*$ on the $|x| < 1-\mu$ ball is bounded by $c_0(1+c_{E,\mu}\delta^{1/4})$. What with (4.22) and (3.11), this follows from (3.6) and Proposition 4.1's fifth bullet.

Granted this $c_0(1+c_{E,\mu}\delta^{1/4})$ bound on the $L^2$ norm bound of $d\hat{A}_*$, Hodge theory (the Abelian version of Uhlenbeck's theorem) can be used to find a function on the $|x| < 1-\mu$ ball, this denoted by $x$, such that $\hat{A} = \hat{A}_* + dx$ is coclosed on the $|x| < 1-\mu$ ball and such that its $L^2_1$ norm is bounded by $c_0(1+c_{E,\mu}\delta^{1/4})$. This bound with the top bullet in (3.11) leads to a $c_0(1+c_{E,\mu}\delta^{1/4})$ bound for the $L^4$ norm of $\hat{A}$.

*Part 4*: Define $\hbar$ to be the composition of first $u$ and then the automorphism on the product SO(3) bundle that defined by $e^{x\tau}$. Writing $\hbar^*A_\Diamond$ as $\theta + \hat{a}_{A_\Delta}$, it follows that $\langle\tau\hat{a}_{A_\Delta}\rangle = \hat{A}$. Note that this definition of h implies that the $\mathfrak{a}$ and $\mathfrak{b}$ are given by

$$\mathfrak{a} = \tfrac{1}{4}\, m^{-1}\,[\tau, \hbar^*(\nabla_{A_\Diamond}\tau_\Diamond)] \;\; \textit{and} \;\; \mathfrak{b} = \hbar^*\mathfrak{b}_* \;\; ;$$

(4.23)

the first identity being a consequence of (4.21).

Define an isomorphism, this denoted by h, from the product principal SO(3) bundle over the $|x| < 1-\mu$ ball to $\phi^*P$ by composing first the isomorphism chosen at the start of this Section 4b and then the automorphism h.

*Part 5*: The following lemma summarizes the salient features of the preceding definitions.

**Lemma 4.3**: *Given* $\mathrm{E} \geq 1, \mu \in (0, \frac{1}{2}]$ *and* $\rho \in (0, 1]$, *there exists* $c > 1$ *with the following significance: Suppose that* $r > 1$ *and* $(\mathrm{A}, \hat{a})$ *is a pair of connection on* P *and section of* $(\mathrm{P} \times_{\mathrm{SO}(3)} \mathfrak{su}(2)) \otimes \mathrm{T^*M}$ *that obey (3.2). Fix* $\mathrm{p} \in \mathrm{M}$ *where* $r_z \geq c$. *Suppose in addition that* $\mathrm{r}_\Diamond < c^{-1}$ *and that* $\|\nabla_{\mathrm{A}_\Diamond} \tilde{a}_*\|_2 < c^{-1}$. *Use Parts 1-4 above to define the isomorphism* h *from the product principal* SO(3) *bundle over the* $|x| \leq 1 - \mu$ *ball to* $\phi^*\mathrm{P}$ *over this ball such that the pair* $(\mathrm{A}_\Delta = \mathrm{h}^*\mathrm{A}_\Diamond, a_\Delta = \mathrm{h}^*\tilde{a}_*)$ *satisfies the conditions in the* $(\rho, \mathrm{M} = \rho\, r_z)$ *version (3.27).*

***Proof of Lemma 4.3***: The proof has seven steps.

<u>Step 1</u>: Keeping in mind that $\mathfrak{b}$ is given by (4.23), then bound in the top bullet of (3.27) follows from the top bullet of Lemma 4.2. With $\mathfrak{a}$ given by (4.23), then the bound in the second bullet of (3.27) follows from what is said by Lemma 4.2's second bullet. The lemma's second bullet leads to a $c_{\mathrm{E},\mu}\, r_z^{-1} \|\nabla_{\mathrm{A}_\Diamond} \tilde{a}_*\|_2$ bound for the $\mathrm{L}^2$ norm of $\mathfrak{a}$.

What the third bullet of (3.27) denotes by $\langle \tau\, \hat{\mathrm{a}}_{\mathrm{A}_\Delta} \rangle$ is the 1-form $\hat{\mathrm{A}}$ from Part 2 in the previous subsection. It follows from what is said there that the third bullet of (3.27) holds with $\mathrm{M} \leq c_{\mathrm{E},\mu}$.

The fifth bullet in (3.27) is obeyed when $\mathrm{r}_\Diamond < c_0^{-1} \rho^{1/2}$ because $\pi_\phi$ depends only on the Riemannian metric and the radius $\mathrm{r}_\Diamond$. It has no dependence on either $(\mathrm{A}_\Diamond, \tilde{a}_\Diamond)$ or h.

<u>Step 2</u>: To see about the fourth bullet in (3.27), note first that $\mathfrak{r}$, $\mathfrak{s}$ and $\mathfrak{r}_0$ are described in Part 2 of Section 3e. In particular, it follows from what is said there and from the first and last bullets of Proposition 4.1 that their $\mathrm{L}^2$ norms on the $|x| < 1 - \mu$ ball are no greater than $c_{\mathrm{E},\mu}\, m^{-1} \mathrm{r}_\Diamond$. It follows as a consequence that $\mathfrak{r}$, $\mathfrak{s}$ and $\mathfrak{r}_0$ obey the bounds in (3.27) if $\mathrm{r}_\Diamond$ is less than $c_{\mathrm{E},\mu}^{-1} \rho$.

<u>Step 3</u>: The derivation of a suitable bound for $\mathfrak{s}_0$ starts with the reminder that $\mathrm{s}_0$ is defined as follows: $\mathrm{s}_0 = -*(\mathrm{d}*\mathfrak{a} + m\, e \wedge [\tau, *\mathfrak{b}])$, and that $\mathfrak{s}_0$ is defined by $\mathrm{s}_0$ using the fourth bullet of (3.23). By way of notation, write $\mathrm{s}_0 - \mathfrak{s}_0$ as $\mathfrak{q}$. Thus, $\mathfrak{q}$ is given by setting $\mathfrak{s}_0$ equal to zero on the right hand side of the fourth bullet in (3.23). The term $*(m\, e \wedge [\tau, *\mathfrak{b}])$ in $\mathrm{s}_0$ contributes only to $\mathfrak{q}$ because the Euclidean metric pairing between $\mathfrak{b}$ and $e_\Delta$ is zero. Indeed, this term would be zero were it the case that $e = e_\Delta$ and $* = *_\phi$. The difference between $e$ and $e_\Delta$ and between $*$ and $*_\phi$ are accounted for the terms in the fourth bullet of (3.23) that comprise $\mathfrak{q}$.

Step 4: Write the $d*\mathfrak{a}$ term from the expression $-*(d*\mathfrak{a} + m e\wedge[\mathfrak{t}, *\mathfrak{b}])$ (which is $\mathfrak{s}_0$) as $d(*_\phi \mathfrak{a}) + d(\pi_\phi \mathfrak{a})$. The $d(\pi_\phi \mathfrak{a})$ contribution is directly accounted for on the right hand side of the fourth bullet of (3.23) by a term in $\mathfrak{q}$. Write $A_\Delta$ as $\theta_0 + \hat{a}_{A_\Delta}$ to write

$$d(*_\phi \mathfrak{a}) = d_{A_\Delta}(*_\phi \mathfrak{a}) - \hat{a}_{A_\Delta} \wedge *_\phi \mathfrak{a} + \mathfrak{a} \wedge *_\phi \hat{a}_{A_\Delta} \tag{4.24}$$

Of interest here is the $\mathfrak{su}_\perp$ part of this expression. Write $\hat{a}_{A_\Delta}$ as $\hat{a}_{A_\Delta} = \mathfrak{t}\langle \mathfrak{t}\, \hat{a}_{A_\Delta} \rangle + m\mathfrak{a}$. Only the $\mathfrak{t}\langle \mathfrak{t}\, \hat{a}_{A_\Delta} \rangle$ part of $\hat{a}_{A_\Delta}$ contributes to the $\mathfrak{su}_\perp$ part of the right most terms in (4.24); and this contribution is accounted for in $\mathfrak{q}$. The $\mathfrak{su}_\perp$ part of the left most term on the right hand side of (4.24) is $\mathfrak{s}_0$.

Step 5: The $\mathfrak{su}_\perp$ part of the left most term on the right hand side of (4.24) is the pull-back via the automorphism $\hbar^*$ of the $\mathfrak{m}_\phi$ Hodge dual of $-\frac{1}{4} m^{-1}[\mathfrak{t}_\Diamond, \nabla_{A_\Diamond}^\dagger \nabla_{A_\Diamond} \mathfrak{t}_\Diamond]$, the formal $L^2$ adjoint here being defined by the metric $\mathfrak{m}_\phi$. To exploit this observation, use Lemma 4.2 to write $\tilde{a}_*$ as $\sqrt{\frac{3}{4\pi}}\,(e_\Delta \mathfrak{t}_\Diamond + \mathfrak{b}_*)$. Use this rewriting to write $\nabla_{A_\Diamond}^\dagger \nabla_{A_\Diamond} \tilde{a}_*$ in terms of $e_\Delta$, $\mathfrak{t}_\Diamond$ and $\mathfrak{b}_*$. The resulting expression can be used with the fact that $\mathfrak{t}_\Diamond$ is orthogonal to $\mathfrak{b}_*$ in $\mathfrak{su}(2)$ and the fact that $e_\Delta$ is $\mathfrak{m}_\phi$ orthogonal to $\mathfrak{b}_*$ to the formula for $\nabla_{A_\Diamond}^\dagger \nabla_{A_\Diamond} \mathfrak{t}_\Diamond$ that is given below in (4.25). The notation in this formula uses $\{e_{\Delta a}\}_{a=1,2,3}$ and $\{\mathfrak{b}_{*a}\}_{a=1,2,3}$ and $\{(\nabla_{A_\Diamond}^\dagger \nabla_{A_\Diamond} \tilde{a}_*)_a\}_{a=1,2,3}$ to denote the respective components of $e_\Delta$, $\mathfrak{b}_*$ and $\nabla_{A_\Diamond}^\dagger \nabla_{A_\Diamond} \tilde{a}_*$ with respect to a chosen $\mathfrak{m}_\phi$ orthonormal frame for $T^*\mathbb{R}^3$ at any given $|x| < 1 - \mu$ point. The formula in (4.25) uses the convention that repeated indices from $\{1, 2, 3\}$ are summed.

$$\nabla_{A_\Diamond}^\dagger \nabla_{A_\Diamond} \mathfrak{t}_\Diamond - \mathfrak{b}_{*a}\langle \mathfrak{b}_{*a} \nabla_{A_\Diamond}^\dagger \nabla_{A_\Diamond} \mathfrak{t}_\Diamond \rangle = e_{\Delta a}(\nabla_{A_\Diamond}^\dagger \nabla_{A_\Diamond} \tilde{a}_*)_a + \mathfrak{b}_{*a}\langle \mathfrak{t}_\Diamond (\nabla_{A_\Diamond}^\dagger \nabla_{A_\Diamond} \tilde{a}_*)_a \rangle + \mathfrak{R}\ , \tag{4.25}$$

with $\mathfrak{R}$ being a sum of terms that are quadratic functions of the covariant derivatives of $e_\Delta$, $\mathfrak{t}_\Diamond$ and $\mathfrak{b}_*$. The essential point being that the norm of $\mathfrak{R}$ is bounded curtesy of the second bullet of Lemma 4.2 by $c_0\, |\nabla_{A_\Diamond} \tilde{a}_*|^2$.

Step 6: If $|\mathfrak{b}_*| < \frac{1}{2}$, then (4.25) leads to the bound

$$|\nabla_{A_\Diamond}^\dagger \nabla_{A_\Diamond} \mathfrak{t}_\Diamond| \le c_0 |e_{\Delta a}(\nabla_{A_\Diamond}^\dagger \nabla_{A_\Diamond} \tilde{a}_*)_a| + c_0 |\mathfrak{b}_*|\, |\nabla_{A_\Diamond}^\dagger \nabla_{A_\Diamond} \tilde{a}_*| + c_0 |\nabla_{A_\Diamond} \tilde{a}_*|^2\ . \tag{4.26}$$

Step 7 derives a $c_{E,\mu}(r_\lozenge + r_z^{-2} + \|\nabla_{A_\lozenge}\tilde{a}_*\|_2^{1/2})$ bound for the $L^2$ norm on the $|x| < 1-\mu$ ball of the right hand side in (4.26). This bound leads directly to the desired $m^{-1}\rho$ bound for for the $L^2$ norm of $\mathfrak{s}$ on the $|x| < 1-\mu$ ball if $r_\lozenge < c_{E,\mu}^{-1}\rho$ and $\|\nabla_{A_\lozenge}\tilde{a}_*\|_2 \le c_{E,\mu}^{-1}\rho^2$ and $r_z^2 > c_{E,\mu}\rho^{-1}$.

<u>Step 7</u>: The fifth bullet of Proposition 4.1 with the top bullet of (3.11) give a $c_{E,\mu}\|\nabla_{A_\lozenge}\tilde{a}_*\|_2^{1/2}$ bound for the $L^2$ norm of $|\nabla_{A_\lozenge}\tilde{a}_*|^2$ over the $|x| < 1-\mu$ ball. Meanwhile, this same fifth bullet of Proposition 4.1 and the $c_{E,\mu} r_z^{-2}$ bound for $|\mathfrak{b}_*|$ from Lemma 4.2 lead to a $c_{E,\mu} r_z^{-2}$ bound for the $L^2$ norm on the $|x| < 1-\mu$ ball of the $|\mathfrak{b}_*||\nabla_{A_\lozenge}^\dagger\nabla_{A_\lozenge}\tilde{a}_*|$ term from the right hand side of (4.26).

To say something about the $|e_{\Delta a}(\nabla_{A_\lozenge}^\dagger\nabla_{A_\lozenge}\tilde{a}_*)_a|$ term, first invoke the definition in (2.10) to bound this term by

$$|q_{A_\lozenge}(\tilde{a}_*)| + |e_{\Delta a}\wedge(*F_{A_\lozenge}\wedge\tilde{a}_* + \tilde{a}_*\wedge *F_{A_\lozenge})| + c_0 r_\lozenge^2|e_\Delta||\tilde{a}_*| \tag{4.27}$$

The fourth bullet of Proposition 4.1 supplies a $c_{E,\mu} r_\lozenge$ bound for the $L^2$ norm of $q_{A_\lozenge}(\tilde{a}_*)$ on the $|x|<1-\mu$ ball. The middle term in (4.27) is no greater than $c_0|\mathfrak{b}_*||F_{A_\lozenge}|$ and so its $L^2$ norm is at most $c_{E,\mu} r_z^{-2}$. The $L^2$ norm of the right most term in (4.27) is at most $c_{E,\mu} r_\lozenge^2$.

**d) The proof of Proposition 3.2**

Fix small $\mu > 0$ and $\varepsilon' > 0$. Lemma 4.3 with the first and second bullets of Proposition 4.1 and Lemma 3.4 lead directly to the following: There exists $c > 1$ such that if $r_z > c$, $r_\lozenge < c^{-1}$ and $\|\nabla_{A_\lozenge}\hat{a}_*\|_2 < c^{-1}$, then $|u(g^*\tilde{a}_* g^{-1})u^{-1}| \le c_0$ on the $|x|<1-\mu$ ball and the $(A_\Delta = u^*g^*A_\lozenge,\ a_\Delta = u(g^*\tilde{a}_* g^{-1})u^{-1})$ version of what is denoted by $\mathfrak{b}$ in (3.22) obeys

$$\int_{|x|\le 1-\mu} |\mathfrak{b}|^2 < \varepsilon' r_z^{-4}\,. \tag{4.28}$$

Given what is said in the second bullet of (3.9), these bounds with the $m = \sqrt{\frac{3}{4\pi}}\ r_z$ version of (3.26) lead directly to the bound

$$\int_{|x|\le 1-\mu} |F_{A_\lozenge}|^2 \le c_0(\varepsilon' + r_\lozenge)\,. \tag{4.29}$$

If $\varepsilon'$ is taken equal to $c_0^{-1}\varepsilon$ and if $r_\lozenge < c_0^{-1}\varepsilon$, then (4.29) is the assertion of Proposition 3.2. Therefore, if $r_z > c$ and $r_\lozenge < \min(c^{-1}, c_0^{-1}\varepsilon)$ and $\|\nabla_{A_\lozenge}\hat{a}_*\|_2 < c^{-1}$, then Proposition 3.2 is true.

Let $\kappa_c$ denote the ($\kappa = c, \mu, \varepsilon$) version of Lemma 3.3's constant $\kappa$. If $r_z \le c$, then Lemma 3.3's conclusion holds if both $r_\lozenge < \kappa_c^{-1}$ and $\|\nabla_{A_\lozenge} \hat{a}_*\|_2 < \kappa_c^{-1}$. As Lemma 3.3's conclusion is the assertion of Proposition 3.2, it follows that Proposition 3.2 is true no matter the value of $r_z$ when both $r_\lozenge$ and $\|\nabla_{A_\lozenge} \hat{a}_*\|_2$ are less than $\min(c^{-1}, c_0^{-1}\varepsilon, \kappa_c^{-1})$.

## 5. Integral identities and monotonicity

This section has two purposes, the first is to state and then prove a proposition that asserts a monotonicity property of two functions on $(0, c_0^{-1})$ that are associated to any given point in M and a pair $(A, \hat{a})$ of connection on P and section of $P\times_{SO(3)}\mathfrak{su}(2)$ that obeys the constraints in (3.2) for a given $\mathrm{E}$ and $r$. This monotonicity assertion constitutes the upcoming Proposition 5.1. The second purpose of the section is to establish various integral inequalities, some of which are used in the proof of Proposition 5.1 and some are used in Section 6.

To set the stage for Proposition 5.1, fix $c_0 > 1$ so that the ball of radius $c_0^{-1}$ centered about any given point in M has compact closure in a Gaussian coordinate chart about the chosen point. Fix for the moment a point $p \in M$. Given $r \in (0, c_0^{-1})$, let $B_r \subset M$ denote radius r ball centered at p and let $\partial B_r$ denote the boundary of the closure of $B_r$.

Reintroduce Uhlenbeck's constant $\kappa_U$ from (3.1). Having specified a connection, A, on P and $p \in M$, define $r_\lozenge \in (0, c_0^{-1}]$ as in (3.6) to be the largest value of r such that

$$\int_{B_r} |F_A|^2 \le \tfrac{1}{100}\, \kappa_U^{-2} r^{-1}, \tag{5.1}$$

To continue the stage setting, suppose that $p \in M$ and $r \in (0, c_0^{-1})$ are given and that $r$ is a chosen, positive number. Fix a pair $(A, \hat{a})$ with A being a connection on P's restriction to the closure of $B_r$ and with $\hat{a}$ being a section over this closure of the bundle $(P\times_{SO(3)}\mathfrak{su}(2))\otimes T^*M$. The parameter $r$ with A and $\hat{a}$ define functions h and H on $(0, c_0^{-1})$ by the rules

$$r \to h(r) = \int_{\partial B_r} |\hat{a}|^2 \quad \textit{and} \quad r \to H(r) = \int_{B_r} (|\nabla_A \hat{a}|^2 + 2r^2 |\hat{a}\wedge\hat{a}|^2)\ . \tag{5.2}$$

Use these two to define the function $\mathrm{N}$ on the $h \ne 0$ part of $(0, c_0^{-1})$ by the rule

$$r \to \mathrm{N}(r) = \frac{r\,H(r)}{h(r)}\ . \tag{5.3}$$

The function $\mathrm{N}$ is called the *frequency* function because it plays the role here that is played by an eponymous function that was introduced by Almgren [Al]. The latter incarnation of $\mathrm{N}$ plays a central role in proving various theorems about the nodal sets of eigenfunctions of the Laplacian on Riemannian manifolds, and on solutions to certain nonlinear equations, second order differential equations with Laplace symbol. See for example [Han], [DF] and [HHL].

**Proposition 5.1**: *Fix* $r_0 > 0$ *such that each* $p \in M$ *version of the functions* $h, H$ *and* $\mathrm{N}$ *are defined on the radius* $r_0$ *ball centered at* $p$. *Given* $r_0$, *there exists* $\kappa > 100$ *whose significance is explained by what follows. Suppose that* $\mathrm{E} \geq 1$, *that* $r > \kappa \mathrm{E}^{\kappa}$ *and that* $(A, \hat{a})$ *is a pair in* $\mathrm{Conn}(P) \times C^{\infty}(M; (P \times_{SO(3)} \mathfrak{su}(2)) \otimes T^*M)$ *obeying Items a)-d) in (3.2). Fix a point* $p \in M$ *so as to define the function* $h$, $H$ *and* $\mathrm{N}$ *using* $r$ *and* $(A, \hat{a})$ *on* $(0, r_0]$.

- *If* $r \in [\frac{1}{2} r_{\Diamond}, \kappa^{-1} r_0]$, *then* $\frac{d}{dr} h = 2 r^{-1}(1 + \mathrm{N}) h + \mathfrak{e}$ *with* $\mathfrak{e}$ *being a function of* $r$ *whose absolute value is less than* $\kappa(\mathrm{E} r^{-2} + \mathrm{E} r^{-1/2} h^{1/2} + \mathrm{E}^{1/2} r^{-1} H^{1/2})$.
- *Let* $u$ *denote the function on* $[\frac{1}{2} r_{\Diamond}, \kappa^{-1} r_0]$ *that is defined by the rule* $r \to u(r) = \kappa \mathrm{E}^{\kappa} r^{1/\kappa}$. *If* $r_{\dagger} \in [\frac{1}{2} r_{\Diamond}, \kappa^{-1} r_0]$ *is such that* $h(r) \geq r^{3-1/\kappa}$ *when* $r \in [r_{\dagger}, \kappa^{-1} r_0]$. *Then*

$$\mathrm{N}(r_2) \geq e^{-u(r_2)+u(r_1)} \mathrm{N}(r_1) - u(r_2) + u(r_1) - r^{-1/8} .$$

  *when* $r_2$ *and* $r_1$ *are such that* $r_{\dagger} \leq r_1 \leq r_2 \leq \kappa^{-1} r_0$.
- *If* $r_* \in [\frac{1}{2} r_{\Diamond}, \kappa^{-1})$ *is such that* $h(r_*) = r_*^{2+1/16}$ *and* $h(r) > r^{2+1/16}$ *when* $r$ *is greater than* $r_*$ *but sufficiently close to* $r_*$, *then* $h(r) \leq \kappa r^{2+1/16}$ *when* $r \in [r_*, \kappa^{-1} r_0]$.
- *If* $r_* \in [\frac{1}{2} r_{\Diamond}, \kappa^{-1})$ *and* $h(r_*) \leq r_*^{2+1/16}$, *then* $h(r) \leq \kappa r^{2+1/16}$ *for all* $r \in [r_*, \kappa^{-1} r_0]$.

Note that the version of $\kappa$ in this proposition depends on $r_0$. The case when $r$ or $r_{\dagger}$ or $r_*$ are less than $r_{\Diamond}$ is of no interest to the subsequent applications. Keep in mind in what follows that $r$ is less than $r_{\Diamond}$ if $r > c_0 \mathrm{E}$ and $r \leq c_0^{-1} \mathrm{E}^{-1} 10^{-1/2} \kappa_U^{-1/2} r^{-1}$, this being a consequence Item c) in (3.2) and the third bullet of (3.3).

By way of an outline for the rest of this section, the proof of the first bullet of the proposition is in Section 5b and that of the second is in Section 5d. Section 5e contains the proofs Proposition 5.1's third and fourth bullets. The intervening subsections establish some facts and observations that are used in the proof.

### a) Integration by parts identities

The upcoming Lemma 5.2 states the two fundamental integration by parts identities. The second refers to the parameter $r$ that is used to define the function H. By way of notation, the lemma introduces functions $\mathscr{h}$ and $\mathfrak{f}$ on $(0, c_0^{-1})$, these defined by

$$\mathfrak{h}(r) = \int_{B_r} (|\nabla_A \hat{a}|^2 + 2\langle *F_A \wedge \hat{a} \wedge \hat{a}\rangle) \quad and \quad \mathfrak{f}(r) = \int_{B_r} (r^{-2}|F_A - r^2\hat{a}\wedge\hat{a}|^2 + |d_A\hat{a}|^2 + |d_A *\hat{a}|^2) \ . \tag{5.4}$$

The lemma also uses $\|\cdot\|_{2,r}$ to denote the $L^2$ norm on $B_r$, and it $\|\cdot\|_{\infty,r}$ to denote the supremum norm on $B_r$.

To set the rest of the notation, keep in mind that if $p \in M$ and if $r \in (0, c_0^{-1})$, then the ball $B_r$ has compact closure in a Gaussian coordinate chart centered at p. This being the case, the outward pointing, unit length tangent vectors to the geodesic rays that start at p define a smooth vector field on B−p. This vector field is denoted by $\partial_r$. When $\hat{a}$ denotes a section over $B_r$ of $(P\times_{SO(3)}\mathfrak{su}(2))\otimes T^*M$, then $\hat{a}_r$ is used to denote the section of $(P\times_{SO(3)}\mathfrak{su}(2))$ overy $B_r$−p that is given by pairing $\hat{a}$ with $\partial_r$. When A denotes a connection on P's restriction to $B_r$, then $\partial_{A,r}$ is used to denote the A's directional covariant derivative along the vector field $\partial_r$.

**Lemma 5.2**: *There exists* $\kappa > 1$ *with the following significance: Fix* $p \in M$ *and fix* $r \in (0, \kappa^{-1})$. *Suppose that* A *is a connection on* P*'s restriction to the closure of* $B_r$ *and that and* $\hat{a}$ *is a section over this closure of* $(P\times_{SO(3)}\mathfrak{su}(2))\otimes T^*M$. *Then*

- $\int_{\partial B_r} \langle \partial_{A,r}\hat{a} \wedge *\hat{a}\rangle = \mathfrak{h} + \int_{B_r} \mathrm{Ric}(\langle \hat{a}\otimes\hat{a}\rangle) - \int_{B_r}\langle \hat{a}\wedge *q_A(\hat{a})\rangle$ .
- $\int_{\partial B_r} (|\nabla_A\hat{a}|^2 + 2r^2|\hat{a}\wedge\hat{a}|^2) = 2\int_{\partial B_r} (|\partial_{A,r}\hat{a}|^2 + r^2|[\hat{a},\hat{a}_r]|^2) + \frac{1}{r}H + \mathfrak{R}$ ,

*with* $\mathfrak{R}$ *obeying* $|\mathfrak{R}| \le c_0(\|q_A(\hat{a})\|_{2,r} + r\,\mathfrak{f}^{1/2}\|\hat{a}\|_{\infty,r})H^{1/2} + c_0(h + rH) + c_0(\mathfrak{f}^{1/2}\|\hat{a}\|_{2,r} + \|\hat{a}\|_{2,r}^{2})$.

***Proof of Lemma 5.2***: Given the definitions of $\mathfrak{h}$ and $q_A$, the identity in the first bullet of the lemma follows using Stoke's theorem. The proof of the assertion in the second bullet has two steps.

Step 1: Fix $p \in M$ and an orthonormal frame $\{e^i\}_{i=1,2,3}$ for $T^*M$ on the radius $c_0^{-1}$ ball centered at p. Write $\hat{a} = \hat{a}_i e^i$ and $\nabla_A\hat{a} = (\nabla_{A,i}\hat{a})_k e^k\otimes e^i$ with it understood that repeated indices are to be summed. View the Riemann curvature tensor as a section of $\otimes^4 T^*M$ and write it as $R_{ikmn}\, e^i\otimes e^k\otimes e^m\otimes e^n$. Likewise, write the Ricci tensor as $\mathrm{Ric}_{ik} e^i\otimes e^k$. The metric tensor is $\delta_{ik} e^i\otimes e^k$ in this notation with $\delta_{ik} = 0$ if $i \ne k$ and $\delta_{11} = \delta_{22} = \delta_{33} = 1$. Granted this notation, a symmetric section of $T^*M\otimes T^*M$ to be denoted by $T = T_{ij} e^i\otimes e^j$ is defined by setting any given $i,j \in \{1, 2, 3\}$ version of $T_{ij}$ to be

$$T_{ij} = \langle(\nabla_{A,i}\hat{a})_k(\nabla_{A,j}\hat{a})_k\rangle + r^2\langle[\hat{a}_i,\hat{a}_k][\hat{a}_j,\hat{a}_k]\rangle + R_{ikjm}\langle\hat{a}_k\hat{a}_m\rangle$$
$$- \tfrac{1}{2}\,\delta_{ij}\,(\langle(\nabla_{A,m}\hat{a})_k(\nabla_{A,m}\hat{a})_k\rangle + r^2\langle[\hat{a}_k,\hat{a}_m][\hat{a}_k,\hat{a}_m]\rangle + \mathrm{Ric}_{mk}\langle\hat{a}_k\hat{a}_m\rangle) \ , \tag{5.5}$$

with it again understood that repeated indices are summed.

Let $d^\dagger$ denote the map from $C^\infty(M; T^*M \otimes T^*M)$ to $C^\infty(M; T^*M)$ given by the formal $L^2$ adjoint of the Levi-Civita covariant derivate. The 1-form $d^\dagger T$ can be written schematically as

$$(d^\dagger T)_i = -\langle (\nabla_{A,i}\hat{a})_k q_A(\hat{a})_k \rangle + \mathcal{Q}_{\wedge i}(r\mathfrak{q} \otimes \hat{a} \otimes r(\hat{a} \wedge \hat{a})) + \mathcal{Q}_{\nabla i}(r\mathfrak{q} \otimes \hat{a} \otimes \nabla_A \hat{a}) + \mathcal{R}_i(\hat{a} \otimes \mathfrak{q}) + \mathcal{S}_i(\hat{a} \otimes \hat{a}), \tag{5.6}$$

where the notation has $\mathfrak{q}$ denoting the triple

$$\mathfrak{q} = (r^{-1} * (F_A - r^2 \hat{a} \wedge \hat{a}),\ * d_A \hat{a},\ d_A * \hat{a}), \tag{5.7}$$

this being a section of $(P \times_{SO(3)} \mathfrak{su}(2)) \otimes (T^*M \oplus T^*M \oplus (\wedge^3 T^*M))$. What are denoted by $\mathcal{Q}_\wedge$, $\mathcal{Q}_\nabla$, $\mathcal{R}$ and $\mathcal{S}$ in (5.6) are homomorphisms that involve only the metric, the Riemann curvature tensor and the latter's covariant derivative. In any event, each has norm bounded by $c_0$.

<u>Step</u> <u>2</u>: Fix $r \in (0, c_0^{-1})$. Let $\Delta$ denote for the moment the function $\text{dist}(\cdot, p)^2$ on the ball $B_r$. Integrate $d^\dagger T \wedge * d\Delta$ on $B_r$ and use Stoke's theorem to write

$$\int_{B_r} d^\dagger T \wedge * d\Delta = -2r \int_{\partial B_r} T(\partial_r \otimes \partial_r) + \int_{B_r} T_{ij} \nabla_i \nabla_j \Delta . \tag{5.8}$$

The identity in the second bullet of Lemma 5.2 follows directly from (5.8) given what is said in the subsequent paragraph about the integral on the left hand side and the two integrals on the right hand side.

The norm of $|d\Delta|$ on $B_r$ is bounded by r and so it follows from (5.5) that the absolute value of the integral on the left hand side of (5.6) is no greater than r times

$$c_0 (\| q_A(\hat{a}) \|_{2,r} + r \mathfrak{f}^{1/2} \| \hat{a} \|_{\infty,r}) H^{1/2} + c_0 (\mathfrak{f}^{1/2} \| \hat{a} \|_{2,r} + \| \hat{a} \|_{2,r}{}^2) . \tag{5.9}$$

Meanwhile, the boundary integral on the right hand side of (5.8) differs from r times

$$\int_{\partial B_r} (|\nabla_A \hat{a}|^2 + 2r^2 |\hat{a} \wedge \hat{a}|^2) - 2 \int_{\partial B_r} (|\partial_{A,r} \hat{a}|^2 + r^2 |[\hat{a}, \hat{a}_r]|^2) \tag{5.10}$$

by at most $c_0 r h$, this being a direct consequence of (5.5)'s formula for T. The right most integral in (5.8) differs from the integral of the trace of T by at most $c_0 (r^2 H + \| \hat{a} \|_{2,r}{}^2)$, this because $\nabla_i \nabla_k \rho$ differs from $2\delta_{ik}$ by at most $c_0 r^2$. A look at (5.5) finds that the trace of T differs from $-\frac{1}{2}(|\nabla_A \hat{a}|^2 + 2r^2 |\hat{a} \wedge \hat{a}|^2)$ by at most $c_0 |\hat{a}|^2$.

**b)** ***Proof of the first bullet of Proposition 5.1***

The proof has five parts.

*Part 1*: Invoke the first bullet of Lemma 5.2 to write

$$\frac{d}{dr} h = 2r^{-1}(1+\mathfrak{r})h + 2\hbar + \int_{B_r} \mathrm{Ric}(\langle \hat{a} \otimes \hat{a} \rangle) - \int_{B_r} \langle \hat{a} \wedge * q_A(\hat{a}) \rangle \ , \tag{5.11}$$

with $\mathfrak{r}$ denoting a term whose absolute value is bounded by $c_0 r^2$. This term accounts for the fact that the second fundamental form of $\partial B_r$ may differ from $r^{-1}$ times the induced metric. The bound on $|\mathfrak{r}|$ is due to the fact that these two tensors can differ by at most $c_0 r$.

*Part 2*: Assume that A is a connection on P's restriction to the closure of $B_r$ and that $\hat{a}$ is a section over this closure of $(P \times_{SO(3)} \mathfrak{su}(2)) \otimes T^*M$. Assume in addition that a positive number $\mathrm{E}$ has been specified and that $\| F_A - r^2 \hat{a} \wedge \hat{a} \|_{2,r} \le \mathrm{E}^{1/2}$. Write $\langle * F_A \wedge \hat{a} \wedge \hat{a} \rangle$ as the sum of $r^{-1} \langle *(F_A - r^2 \hat{a} \wedge \hat{a}) \wedge (r \hat{a} \wedge \hat{a}) \rangle$ and $r^2 |\hat{a} \wedge \hat{a}|^2$ to derive the bound given below for the difference between $\hbar$ and H:

$$|\hbar - H| \le c_0 \mathrm{E}^{1/2} r^{-1} H^{1/2} \ . \tag{5.12}$$

Use this bound with the definition of $\mathrm{N}$ to write (5.11) as

$$\frac{d}{dr} h = 2r^{-1}(1+\mathrm{N}+\mathfrak{r})h + \mathfrak{r}_1 + \int_{B_r} \mathrm{Ric}(\langle \hat{a} \otimes \hat{a} \rangle) - \int_{B_r} \langle \hat{a} \wedge * q_A(\hat{a}) \rangle \tag{5.13}$$

with $\mathfrak{r}_1$ being a term with norm at most $c_0 \mathrm{E}^{1/2} r^{-1} H^{1/2}$.

*Part 3*: Assume as in Part 2 that $(A, \mathfrak{a})$ is a pair of connection on P's restriction to the closure of $B_r$ and that $\hat{a}$ is a section over this closure of $(P \times_{SO(3)} \mathfrak{su}(2)) \otimes T^*M$. Then

$$\int_{B_r} |\hat{a}|^2 \le 2r\, h(r) + 8 r^2 H(r) \ . \tag{5.14}$$

This inequality is a direct consequence of the following integration by parts formula for the radius 1 ball in $\mathbb{R}^3$: Let $f$ denote a smooth on the $|x| \le 1$ ball. Fix $\varepsilon \in (0, 1)$. Then

$$(1-\varepsilon)\int_{|x|\le 1}\frac{1}{|x|^2}f^2 \le \int_{|x|=1} f^2 \;+\; \varepsilon^{-1}\int_{|x|\le 1}|df|^2 \;.$$

(5.15)

The inequality in (5.15) leads to (5.14) by writing the integrals in (5.14) using a Gaussian coordinate chart centered at p. Having done so, appeal to (2.6) when invoking the $f = |\hat{a}|$ and $\varepsilon = \frac{1}{4}$ version of (5.15).

It follows directly from (5.14) that the contribution to (5.13) from the Ricci curvature integral on (5.13)'s right hand side has absolute value at most $c_0\mathrm{r}(1+\mathrm{N})\mathrm{h}$.

*Part 4*: Assume that A is a connection on P's restriction to the closure of $B_\mathrm{r}$ and that $\hat{a}$ is a section over this closure of $(P\times_{SO(3)}\mathfrak{su}(2))\otimes T^*M$. Assume in addition that a number $\mathrm{E} > 1$ has been specified and that each of $\|q_A(\hat{a})\|_{2,\mathrm{r}}$, $r\|d_A\hat{a}\|_{2,\mathrm{r}}$ and $r\|d_A{*}\hat{a}\|_{2,\mathrm{r}}$ is less than $\mathrm{E}^{1/2}$. As explained momentarily, $r > c_0\mathrm{E}$ and if $\mathrm{r} > c_0^{-1}\mathrm{E}^{-1}10^{-1/2}\kappa_U^{-1/2}r^{-1}$, then

$$\Big|\int_{B_\mathrm{r}}\big\langle \hat{a}\wedge{*}q_A(\hat{a})\big\rangle\Big| \le c_0(\mathrm{E}\,r^{-2} + \mathrm{E}\,r^{-1/2}\mathrm{h}^{1/2} + \mathrm{E}^{1/2}r^{-1}\mathrm{H}^{1/2}) \;.$$

(5.16)

By way of a parenthetical remark, (5.16) holds with $\mathrm{r} > m^{-1}r$ for any $m > 0$ with $c_0$ replaced by an $m$-dependent constant. The given lower bound on r is used because r is necessarily less than $\mathrm{r}_\Diamond$ if $r > c_0\mathrm{E}$ and $\mathrm{r} \le c_0^{-1}10^{-1/2}\kappa_U^{-1/2}r^{-1}$. The proof of (5.16) has two steps.

<u>Step 1</u>: Fix $\varepsilon \in (0, \frac{1}{2})$ for the moment and let $\chi_\varepsilon$ denote the function on $B_\mathrm{r}$ given by $\chi(\varepsilon^{-1}(\mathrm{r}^{-1}\mathrm{dist}(\cdot,p) - 1 + \varepsilon))$. This function equals 1 where $\mathrm{dist}(\cdot,p) \le (1-\varepsilon)\mathrm{r}$ and it equals 0 where $\mathrm{dist}(\cdot,p) \ge \mathrm{r}$. Write the integral over $B_\mathrm{r}$ of $\langle\hat{a}\wedge{*}q_A(\hat{a})\rangle$ as the sum of two integrals, the first being the integral over $B_\mathrm{r}$ of $\chi_\varepsilon\langle\hat{a}\wedge{*}q_A(\hat{a})\rangle$ and the second being the integral of $(1-\chi_\varepsilon)\langle\hat{a}\wedge{*}q_A(\hat{a})\rangle$. The absolute values of these two integrals are bounded respectively by

$$c_0\mathrm{E}\,r^{-2} + c_0\mathrm{E}^{1/2}\varepsilon^{-1/2}\mathrm{r}^{-1/2}r^{-1}(\sup\nolimits_{(1-\varepsilon)\mathrm{r}\le s\le \mathrm{r}}\mathrm{h}(s))^{1/2} \quad \textit{and} \quad c_0\mathrm{E}^{1/2}\varepsilon^{1/2}\mathrm{r}^{1/2}(\sup\nolimits_{(1-\varepsilon)\mathrm{r}\le s\le \mathrm{r}}\mathrm{h}(s))^{1/2} \;.$$

(5.17)

To prove the bound for the first integral, write $q_A(\hat{a})$ as ${*}d_A{*}d_A\hat{a} - d_A{*}d_A{*}\hat{a}$, integrate by parts and invoke the assumption that $\|d_A\hat{a}\|_{2,\mathrm{r}} + \|d_A{*}\hat{a}\|_{2,\mathrm{r}}$ is less than $r^{-1}\mathrm{E}^{1/2}$. The bound for the second integral follows directly from the assumption that $\|q_A(\hat{a})\|_{2,\mathrm{r}} \le \mathrm{E}^{1/2}$.

<u>Step 2</u>: Use the fundamental theorem of calculus to see that

$$\sup_{(1-\varepsilon)\mathrm{r} \le s \le \mathrm{r}} \mathrm{h}(s) \le c_0 \mathrm{h}(\mathrm{r}) + c_0 \varepsilon \mathrm{r} \int_{B_{\mathrm{r}}} |\nabla_A \hat{a}|^2 \quad . \tag{5.18}$$

This inequality with (5.17) leads directly to the following bound:

$$\Big| \int_{B_{\mathrm{r}}} \big\langle \hat{a} \wedge * q_A(\hat{a}) \big\rangle \Big| \le c_0 \mathrm{E}\, r^{-2} + c_0 \mathrm{E}^{1/2} (\varepsilon^{-1/2} \mathrm{r}^{-1/2} r^{-1} + \varepsilon^{1/2} \mathrm{r}^{1/2}) \mathrm{h}^{1/2} + c_0 (r^{-1} + \varepsilon \mathrm{r}) (\mathrm{H}^{1/2} + \mathrm{E}^{1/2} r^{-1}) . \tag{5.19}$$

Assume henceforth that $\mathrm{r} > c_0^{-1} \mathrm{E}^{-1} 10^{-1} \kappa_U^{-1/2} r^{-1}$. Use $\varepsilon = c_0^{-1} \kappa_U^{-1/2} \mathrm{E}^{-1} \mathrm{r}^{-1} r^{-1}$ in the right hand side of (5.19) to go from the latter bound to the bound in (5.16).

*Part 5*: Use what is said in Parts 2-4 to conclude that the right hand side of (5.11) differs from $2\frac{1}{\mathrm{r}}(1+\mathrm{N})\mathrm{h}$ by at most $c_0 \mathrm{r}(1+\mathrm{N})\mathrm{h} + c_0(\mathrm{E}\, r^{-2} + \mathrm{E}^{1/2} r^{-1} \mathrm{H}^{1/2}) + c_0 \mathrm{E}\, r^{-1/2} \mathrm{h}^{1/2}$. This is what is claimed by the first bullet of Proposition 5.1.

**c) A differential equation for N**

The next lemma asserts a formula for $\frac{d\mathrm{N}}{d\mathrm{r}}$. This lemma also refers to Uhlenbeck's constant $\kappa_U$ from (3.1).

**Lemma 5.3**: *There exists* $\kappa \ge 100$ *with the following significance: Suppose that* $\mathrm{E} \ge 1$*, that* $r > \kappa \mathrm{E}$ *and that* $(A, \hat{a}) \in \mathrm{Conn}(P) \times C^\infty(M; (P \times_{SO(3)} \mathfrak{su}(2)) \otimes T^*M)$ *obeys the conditions in (3.2). Fix* $p \in M$ *so as to define the functions* h, H *and* N *using* $r$ *and* $(A, \hat{a})$*. If* $\mathrm{r} \ge \kappa^{-1} \kappa_U^{-1/2} r^{-1}$*, then*

$$\frac{d\mathrm{N}}{d\mathrm{r}} \ge 2 \frac{\mathrm{r}}{\mathrm{h}} \int_{\partial B_{\mathrm{r}}} r^2 \|[\hat{a}, \hat{a}_{\mathrm{r}}]\|^2 - \mathfrak{x} \mathrm{N} - \kappa \mathrm{E} \left(\frac{\mathrm{r}}{\mathrm{h}}\right)^{1/2} \mathrm{N}^{1/2} - \kappa \left(\mathrm{r} + \kappa \mathrm{E}\, r^{-1/2} \mathrm{r} \left(\frac{\mathrm{r}}{\mathrm{h}}\right)^{1/2}\right) \quad ,$$

*with* $\mathfrak{x}$ *obeying* $|\mathfrak{x}| \le \kappa \left( \mathrm{E}\, r^{-1} \left(\frac{1}{\mathrm{r}\mathrm{h}}\right)^{1/2} \mathrm{N}^{1/2} + \mathrm{r}\mathrm{N} + \mathrm{E} \frac{1}{\mathrm{h}} r^{-2} + c_0 \mathrm{E}\, r^{-1/2} \left(\frac{1}{\mathrm{h}}\right)^{1/2} \right)$

***Proof of Lemma 5.3***: The proof does its best to mimic what is done in [Al], [HHL] and [Han] to prove an analoguous monotonicity assertion for the version of N that arises when studying eigenfunctions of $d^\dagger d$. In any event, the starting point is the identity

$$\frac{d\mathrm{N}}{d\mathrm{r}} = \frac{\mathrm{H}}{\mathrm{h}} + \frac{\mathrm{r}}{\mathrm{h}} \frac{d\mathrm{H}}{d\mathrm{r}} - \frac{\mathrm{r}\mathrm{H}}{\mathrm{h}^2} \frac{d\mathrm{h}}{d\mathrm{r}} . \tag{5.20}$$

The derivation of the lemma's lower bound from (5.20) has four parts.

*Part 1*: The derivative of H is

$$\frac{dH}{dr} = \int_{\partial B_r} (|\nabla_A \hat{a}|^2 + 2r^2 |\hat{a}\wedge\hat{a}|^2) \ , \tag{5.21}$$

and the second bullet of Lemma 5.2 is used to write this as

$$\frac{dH}{dr} = 2 \int_{\partial B_r} (|\partial_{A,r}\hat{a}|^2 + r^2 |[\hat{a}, \hat{a}_r]|^2) + \tfrac{1}{r} H + \mathfrak{R} \ . \tag{5.22}$$

Use this last identity to write N's derivative as

$$\frac{dN}{dr} = 2\frac{r}{h} \int_{\partial B_r} (|\partial_{A,r}\hat{a}|^2 + r^2 |[\hat{a}, \hat{a}_r]|^2) + \frac{r\mathfrak{R}}{h} + 2\frac{H}{h} - \frac{rH}{h^2}\frac{dh}{dr} \ . \tag{5.23}$$

Save this identity for Part 3.

*Part 2*: Use the definition of h to write

$$\frac{dh}{dr} = 2(1+\mathfrak{r})\frac{h}{r} + 2\int_{\partial B_r} \langle \partial_{A,r}\hat{a} \wedge *\hat{a}\rangle \ , \tag{5.24}$$

with $\mathfrak{r}$ being the same as its namesake in (5.11). Thus, $|\mathfrak{r}| \le c_0 r^2$. This then leads to

$$2\frac{H}{h} - \frac{rH}{h^2}\frac{dh}{dr} = -2\frac{rH}{h^2}\int_{\partial B_r} \langle \partial_{A,r}\hat{a} \wedge *\hat{a}\rangle - 2\tfrac{1}{r}\mathfrak{r} N \ . \tag{5.25}$$

To continue, use (5.12) to write $H = \hslash - \mathfrak{z}_1$ with $\mathfrak{z}_1$ obeying $|\mathfrak{z}_1| \le c_0 E^{1/2} r^{-1} H^{1/2}$ and then use this to write the right hand side of (5.25) as twice

$$-\frac{r\hslash}{h^2}\int_{\partial B_r} \langle \partial_{A,r}\hat{a} \wedge *\hat{a}\rangle + \frac{r}{h^2}\mathfrak{z}_1 \int_{\partial B_r} \langle \partial_{A,r}\hat{a} \wedge *\hat{a}\rangle - \tfrac{1}{r}\mathfrak{r} N \ , \tag{5.26}$$

Now use the first bullet of Lemma 5.2 to replace $\hslash$ in (5.26) and in doing so, replace (5.26) with the expression

$$-\frac{r}{h^2}\Big(\int_{\partial B_r} \langle \partial_{A,r}\hat{a} \wedge *\hat{a}\rangle\Big)^2 + \frac{r}{h^2}\Big(\mathfrak{z}_1 + \int_{B_r} \mathrm{Ric}(\langle \hat{a}\otimes\hat{a}\rangle) - \int_{B_r} \langle \hat{a}\wedge * q_A(\hat{a})\rangle\Big) \int_{\partial B_r} \langle \partial_{A,r}\hat{a} \wedge *\hat{a}\rangle - \tfrac{1}{r}\mathfrak{r} N \ . \tag{5.27}$$

One last replacement is needed, this the replacement in (5.27) of the right most integral over $\partial B_r$ of $\partial_{A,r}\hat{a}\wedge *\hat{a}$ using the first bullet of Lemma 5.2 to equate (5.27) with

$$-\frac{\mathrm{r}}{h^2}\Big(\int_{\partial B_r}\big\langle \partial_{A,r}\hat{a}\wedge *\hat{a}\big\rangle\Big)^2 + \frac{\mathrm{r}}{h^2}\,\mathfrak{z}\,(H+\mathfrak{z}) - \frac{1}{\mathrm{r}}\,\mathfrak{t}\,\mathrm{N}\,, \tag{5.28}$$

with $\mathfrak{z}$ denoting $\mathfrak{z} = \mathfrak{z}_1 + \int_{B_r}\mathrm{Ric}(\langle\hat{a}\otimes\hat{a}\rangle) - \int_{B_r}\big\langle \hat{a}\wedge *q_A(\hat{a})\big\rangle$.

*Part 3*: Twice the expression in (5.28) is the same as the right hand side of (5.25) and thus the same as the left hand side of (5.25). The latter appears in (5.23). Replace its incarnation in (5.23) with (5.28) to write the derivative of N as

$$\frac{d\mathrm{N}}{d\mathrm{r}} = 2\,\frac{\mathrm{r}}{h^2}\,\mathfrak{X} + \frac{\mathrm{r}\,\mathfrak{N}}{h} + 2\,\frac{\mathrm{r}}{h^2}\,\mathfrak{z}\,(H+\mathfrak{z}) - 2\,\frac{1}{\mathrm{r}}\,\mathfrak{t}\,\mathrm{N}\,, \tag{5.29}$$

with $\mathfrak{X} = \Big(\int_{\partial B_r}|\hat{a}|^2\Big)\Big(\int_{\partial B_r}(|\partial_{A,r}\hat{a}|^2 + r^2|[\hat{a},\hat{a}_r]|^2)\Big) - \Big(\int_{\partial B_r}\big\langle \partial_{A,r}\hat{a}\wedge *\hat{a}\big\rangle\Big)^2$.

The key observation with regards to (5.29) is that $\mathfrak{X}$ is non-negative. More to the point, writing the derivative of N as in (5.29) leads to the lower bound

$$\frac{d\mathrm{N}}{d\mathrm{r}} \ge 2\,\frac{\mathrm{r}}{h}\int_{\partial B_r} r^2|[\hat{a},\hat{a}_r]|^2 + \frac{\mathrm{r}\,\mathfrak{N}}{h} + 2\,\frac{\mathrm{r}}{h^2}\,\mathfrak{z}\,(H+\mathfrak{z}) - 2\,\frac{1}{\mathrm{r}}\,\mathfrak{t}\,\mathrm{N}\,. \tag{5.30}$$

The left most term on the right hand side of (5.30) is nonnegative. Part 4 gives bounds for the absolute values of the remaining terms.

*Part 4*: With regards to (5.30), keep in mind that $|\mathfrak{t}| \le c_0\mathrm{r}^2$. Keep in mind also that Lemma 5.2 with (5.14), the assumptions in Items a)-d) of (3.2) and the first plus third bullet of (3.3) give the bound

$$|\mathfrak{N}| \le c_0\,\mathrm{E}\,H^{1/2} + c_0\,h + c_0\,\mathrm{r}\,H + c_0\,\mathrm{E}\,r^{-1}\,\mathrm{r}^{1/2}\,h^{1/2} \tag{5.31}$$

Meanwhile, what is said above about $|\mathfrak{z}_1|$ and what is said in (5.14) and (5.16) imply that

$$|\mathfrak{z}| \le c_0\,\mathrm{E}^{1/2}\,r^{-1}\,H^{1/2} + c_0(\mathrm{r}\,h + \mathrm{r}^2 H) + c_0\,(\mathrm{E}\,r^{-2} + \mathrm{E}\,r^{-1/2}\,h^{1/2})\,. \tag{5.32}$$

Use these bounds for $|\mathfrak{t}|$, $|\mathfrak{N}|$ and $|\mathfrak{z}|$ to bound the derivative of N from below by

$$\frac{dN}{dr} \geq 2\frac{r}{h}\int_{\partial B_r} r^2 |[\hat{a}, \hat{a}_r]|^2 - c_0 E r H^{1/2} h^{-1} - c_0 r - c_0 E r^{-1/2} r^{3/2} h^{-1/2}$$

$$- c_0(r + E r^{-1} H^{1/2} h^{-1} + c_0 r^2 H h^{-1} + c_0 E r^{-2} h^{-1} + c_0 E r^{-1/2} h^{-1/2}) N \,. \tag{5.33}$$

Write each occurrence of H in (5.33) as $r^{-1} N h$ to obtain the bound asserted by Lemma 5.3.

**d) Proof of the second bullet of Proposition 5.1**

The proof of the second bullet has two parts.

*Part 1*: This part of the subsection states and then proves a lemma that asserts the conclusions of the second bullet on any interval in $[r_\ddagger, r_0]$ where N is not too large.

**Lemma 5.4**: *Given* $m > 1$, *there exists* $\kappa \geq 100$ *with the following significance: Suppose that* $E \geq 1$, $r > \kappa E^\kappa$ *and that* $(A, \hat{a}) \in \mathrm{Conn}(P) \times C^\infty(M; (P \times_{SO(3)} \mathfrak{su}(2)) \otimes T^*M)$ *obeys the conditions in Items a)-e) in the statement of Proposition 5.1. Suppose that* $r_2 > r_1$ *are from the interval* $[r_\Diamond, r_0]$ *and such that both* $h(r) \geq r^{3-1/m}$ *and* $N(r) \leq r^{1/2m}$ *when* $r \in [r_1, r_2]$. *Let* u *denote the function on* $[r_1, r_2]$ *that is defined by the rule* $r \to u(r) = \kappa E^\kappa r^{1/\kappa}$. *Then* $N(r_2) \geq e^{-u(r_2)+u(r_1)} N(r_1) - u(r_2) + u(r_1)$.

***Proof of Lemma 5.4***: The proof has two steps.

Step 1: Keep in mind that the function H is uniformly bounded by $c_0 E$, this being a consequence of the first bullet of (3.3). Write each explicit occurrence of $H^{1/2}$ on the right hand side of (5.33) as $(r^{-1} N h)^{1/2}$ to obtain the bound

$$\frac{dN}{dr} \geq - c_0 E r^{1/2} h^{-1/2} N^{1/2} - c_0 r - c_0 E r^{-1/2} r^{3/2} h^{-1/2}$$

$$- c_0(r + E r^{-1} (rh)^{-1/2} N^{1/2} + c_0 r^2 H h^{-1} + c_0 E r^{-2} h^{-1} + c_0 E r^{-1/2} h^{-1/2}) N \,. \tag{5.34}$$

Replace the left most occurrence of $N^{1/2}$ on the right hand side of (5.34) by $(1 + N)$, replace the explicit occurrence of H on the right hand side by $c_0 E$ to obtain the lower bound

$$\frac{dN}{dr} \geq - c_0 E r^{1/2} h^{-1/2} - c_0 r - c_0 E r^{-1/2} r^{3/2} h^{-1/2}$$

$$- c_0(r + E r^{1/2} h^{-1/2} + E r^{-1} (rh)^{-1/2} N^{1/2} + c_0 E r^2 h^{-1} + c_0 E r^{-2} h^{-1} + c_0 E r^{-1/2} h^{-1/2}) N \,. \tag{5.35}$$

Now suppose $h(r) \geq r^{3-1/m}$ when $r \in [r_1, r_2]$. Use this assumption to replace all occurrences of h on the right hand side of (5.35) with $r^{3-1/m}$. Meanwhile, replace the occurrence of $N^{1/2}$ on the right hand side by $r^{1/4m}$. These replacements lead to the lower bound

$$\frac{dN}{dr} \geq - c_0 E\, \mathrm{r}^{-1+1/2m} - c_0 \mathrm{r} - c_0 E\, r^{-1/2} \mathrm{r}^{1/2m}$$

$$- c_0(\mathrm{r} + E\mathrm{r}^{-1+1/2m} + E\, r^{-1+1/4m} \mathrm{r}^{-2+1/2m} + c_0 E\, r^{-2} \mathrm{r}^{-3+1/m} + c_0 E\, r^{-1/2} \mathrm{r}^{-3/2+1/2m})\, N \tag{5.36}$$

This inequality plays the role of the monotonicity inequalities for N's namesakes in [Al], [HHL] and [Han].

<u>Step 2</u>: Integrate (5.36) to obtain the inequality

$$N(\mathrm{r}_2) \geq e^{-u(\mathrm{r}_2)+u(\mathrm{r}_1)} (N(\mathrm{r}_1) - c_0 (m + E\, r^{-1/2}) (\mathrm{r}_2^{1/2m} - \mathrm{r}_1^{1/2m}) \;, \tag{5.37}$$

where u is defined by the rule

$$\mathrm{r} \to u(\mathrm{r}) = c_0 E\, (m\, \mathrm{r}^{1/2m} - r^{-1+1/4m} \mathrm{r}^{-1+1/2m} - r^{-2} \mathrm{r}^{-2+1/m} + r^{-1/2} \mathrm{r}^{-1/2+1/2m}) \;. \tag{5.38}$$

The right hand side of (5.38) is an increasing function of r, but it is less than $c_0 E\, m\, \mathrm{r}^{1/2m}$ unless $\mathrm{r} \leq c_m E^{c_m} r^{-1-1/8m}$ with $c_m$ depending only on m. This with the third bullet in (3.3) implies that the integral of $|F_A|^2$ over the radius r ball centered at p is no greater than $c_0$ $(c_m E^{c_m})^3 r^{-1-3/8m}$. Note in particular that the latter $L^2$ norm bound violates the bound in (5.1) when $r$ is large. If the lower bound $r \geq c_0 (c_m E^{c_m})^{2m}$ is obeyed, then Lemma 5.4's assertion follows directly from (5.37) and (5.38).

*Part 2*: The second bullet of Proposition 5.1 is an immediate corollary of Lemma 5.4 and the subsequent lemma.

**Lemma 5.5**: *Given* $m > 1$, *there exists* $\kappa \geq 100$ *with the following significance: Suppose that* $E \geq 1$, $r > \kappa E^{\kappa}$ *and that* $(A, \hat{a}) \in \mathrm{Conn}(P) \times C^{\infty}(M; (P \times_{SO(3)} \mathfrak{su}(2)) \otimes T^*M)$ *obeys the conditions in (3.2). Fix* $p \in M$ *and suppose that* $\mathrm{r}_\dagger \in [\mathrm{r}_\Diamond, \frac{1}{2}\mathrm{r}_0]$ *is such that* $h(\mathrm{r}) \geq \mathrm{r}^{3-1/m}$ *when* $\mathrm{r} \in [\mathrm{r}_\dagger, \mathrm{r}_0]$. *Then* $N(\cdot) \leq r^{1/2m}$ *on the interval* $[\mathrm{r}_\dagger, \frac{1}{2}\mathrm{r}_0]$.

***Proof of Lemma 5.5***: Let $\kappa_m$ denote the version of $\kappa$ that is given by Lemma 5.4, and suppose that $r \geq \kappa_m E^{\kappa_m}$ so as to invoke Lemma 5.4. Suppose that $\mathrm{r}_1 \in [\mathrm{r}_\dagger, \frac{1}{2}\mathrm{r}_0]$ is a point where $N(\mathrm{r}_1) \geq r^{1/2m}$. It follows as a consequence of Lemma 5.4 that $N(\cdot) \geq \frac{1}{2}\kappa_m^{-1} r^{1/2m}$ on the interval $[\mathrm{r}_1, \mathrm{r}_0]$. Use this fact with the first bullet of Proposition 5.1 to see that

$$\tfrac{d}{dr}\,\mathrm{h} \geq c_*^{-1} r^{1/2m}\, \mathrm{r}^{-1}\mathrm{h} - c_0 E\,(r^{-2} + r^{-1-1/2m}) \quad \textit{when } \mathrm{r} \in [\mathrm{r}_1, \mathrm{r}_0]\ , \tag{5.39}$$

where $c_*$ denotes a number that is greater than 1 and depends only on m. Introduce by way of shorthand $\textsc{z}$ denote $c_*^{-1} r^{1/2m}$. The differential inequality in (5.39) implies that

$$\mathrm{h}(\mathrm{r}_0) \geq 2^{\textsc{z}}(\mathrm{h}(\mathrm{r}) - c_m E\,(r^{-2-1/2m} + r^{-1-1/m})) \ \textit{ for any } \mathrm{r} \in [\mathrm{r}_1, \tfrac{1}{2}\mathrm{r}_0]\ , \tag{5.40}$$

where $c_m$ denotes here and in what follows a number that is greater than 1 and depends only on m. The assigned value for $c_m$ can be assumed to increase between successive appearances.

The fact that $|\hat{a}| \leq c_0 E$ means that $\mathrm{h}(\mathrm{r}_0) \leq c_0 E\, \mathrm{r}_0^2$. This runs afoul of the $\mathrm{r} = \frac{1}{2}\mathrm{r}_0$ version of (5.40) when $r \geq c_m$ unless $\mathrm{h}(\frac{1}{2}\mathrm{r}_0) < c_m E\, r^{-1-1/m}$. But the latter bound is not allowed if $r > c_m E^{c_m}$ because $\mathrm{h}(\frac{1}{2}\mathrm{r}_0)$ is no smaller than $\mathrm{r}_0^{3-1/m}$.

**e) Proof of the third and fourth bullets of Proposition 5.1**

The fourth bullet follows directly from the third. The proof of the proposition's third bullet has three parts.

*Part 1*: Fix $\delta \in (0, 1)$ for the moment. With $\delta$ chosen define the function $x$ on $[\mathrm{r}_\Diamond, \frac{1}{2}\mathrm{r}_0]$ by the rule $\mathrm{r} \to x(\mathrm{r}) = \mathrm{h}(\mathrm{r}) - \mathrm{r}^{2+2\delta}$. The first bullet of Proposition 5.1 implies that $x$ obeys the differential inequality

$$\tfrac{d}{dr}\, x \leq 2\,\mathrm{r}^{-1}(1 + \textsc{n})\, x + 2\,(\textsc{n} - \delta)\,\mathrm{r}^{1+2\delta} + c_0 E\,(r^{-2} + r^{-1/2}\mathrm{h}^{1/2} + r^{-1}\mathrm{H}^{1/2})\ . \tag{5.41}$$

Write $\mathrm{H}^{1/2}$ as $(\mathrm{r}^{-1}\textsc{n}\,\mathrm{h})^{1/2}$ in (5.41) to see that (5.41) leads to the bound

$$\tfrac{d}{dr}\, x \leq 2\mathrm{r}^{-1}(1 + \textsc{n})\, x + 2(\textsc{n} - \delta + r^{-1/8})\,\mathrm{r}^{1+2\delta} + c_0 E\,(1 + \delta^{-1})\, r^{-7/8} + c_0 E\, r^{-1/2}\mathrm{h}^{1/2}\ . \tag{5.42}$$

Let $\mathrm{r}_\delta \in [\mathrm{r}_\Diamond, \frac{1}{2}\mathrm{r}_0]$ denote a zero of $x$, thus a value of r where the function h obeys $\mathrm{h}(\mathrm{r}) = \mathrm{r}^{2+2\delta}$. It follows from (5.42) that $x$ is decreasing at $\mathrm{r}_\delta$ unless

$$\textsc{n} \geq \delta - r^{-1/8} - c_0 E\,(1 + \delta^{-1})\, r^{-2}\,\mathrm{r}^{-1-2\delta} - c_0 E\, r^{-1/2}\,\mathrm{r}^{-\delta}\ . \tag{5.43}$$

In particular, if $\delta \le \frac{1}{4}$, then $x$ is decreasing at $\mathrm{r}_\delta$ unless $\mathrm{N}(\mathrm{r}_\delta) > \frac{1}{4}\delta$ or else $\mathrm{r} \le c_0^{-1}\delta^4 \mathrm{E}^{-4} r^{-4/3}$. The latter option is of no concern if $r > c_0\, \mathrm{E}^{c_0}$ because it runs afoul of the $\mathrm{r} \ge \mathrm{r}_\Diamond$ assumption. This understood, assume henceforth this lower bound for $r$ so as to conclude that $x$ is decreasing at any zero in $[\mathrm{r}_\Diamond, \frac{1}{2}\mathrm{r}_0]$ where $\mathrm{N} \le \frac{1}{4}\delta$.

*Part 2*: Take $\delta = \frac{1}{16}$. Suppose that $\mathrm{r}_* \in [\mathrm{r}_\Diamond, \frac{1}{2}\mathrm{r}_0]$ is a zero of $x$ with the property that $x(\mathrm{r}_*+\mathrm{t}) > 0$ if t is positive and sufficiently small. Let $\mathrm{r}_1 \in (\mathrm{r}_*, c_0^{-1}\mathrm{r}_0]$ denote the largest value of r such that $\mathrm{h}(\mathrm{s}) \ge \mathrm{s}^{2+3/4}$ when $\mathrm{s} \in [\mathrm{r}_*, \mathrm{r}]$. Note in particular that $\mathrm{r}_1$ is strictly larger than $\mathrm{r}_*$ because $\mathrm{h}(\mathrm{r}_*) = \mathrm{r}_*^{2+1/16}$. It follows from the second bullet of Proposition 5.1 that

$$\mathrm{N}(\mathrm{r}) \ge \tfrac{1}{16} - r^{-1/8} - c_{\mathrm{E}}(\mathrm{s}^{1/c_0} - \mathrm{r}_*^{1/c_0}) \quad \textit{for all } \mathrm{s} \in [\mathrm{r}_*, \mathrm{r}_1]. \tag{5.44}$$

Keeping this in mind, use the first bullet of Proposition 5.1 with to see that

$$\tfrac{d}{dr}\mathrm{h} \ge 2\mathrm{r}^{-1}(1 + \tfrac{1}{16} - r^{-1/8} - c_{\mathrm{E}}(\mathrm{r}^{1/c_0} - \mathrm{r}_*^{1/c_0}))\mathrm{h} - c_{\mathrm{E}}(r^{-2} + r^{-7/8}\mathrm{r}) \tag{5.45}$$

on $[\mathrm{r}_*, \mathrm{r}_1]$. Since $\mathrm{r}_* \ge c_0^{-1} r^{-1}$, this inequality implies that

$$\tfrac{d}{dr}\mathrm{h} \ge 2\mathrm{r}^{-1}(1 + \tfrac{1}{16} - r^{-1/8} - c_{\mathrm{E}}(\mathrm{r}^{1/c_0} - \mathrm{r}_*^{1/c_0}))\mathrm{h} - c_{\mathrm{E}}\, r^{-3/4}\mathrm{r}^{9/8} \tag{5.46}$$

on $[\mathrm{r}_*, \mathrm{r}_1]$. Integrating (5.46) finds that

$$\mathrm{h}(\mathrm{r}) \ge e^{-v(\mathrm{r}) + v(\mathrm{r}_*)}\, \mathrm{r}^{2+1/16} \quad \textit{when } \mathrm{r} \in [\mathrm{r}_*, \mathrm{r}_1], \tag{5.47}$$

where v is a non-negative function on $[0, c_0^{-1}\mathrm{r}_0]$ that obeys $|v|(\mathrm{r}) \le c_0\, \mathrm{r}^{1/c_0}$. The lower bound in (5.47) implies that $\mathrm{h}(\mathrm{r}_1) \ge \mathrm{r}_1^{2+3/4}$ if $\mathrm{r}_1 \le c_0^{-1}\mathrm{r}_0$. This being the case, it follows that $\mathrm{r}_1$ must equal $c_0^{-1}\mathrm{r}_0$ and so both (5.44) and (5.46) hold on $[\mathrm{r}_*, c_0^{-1}\mathrm{r}_0]$.

*Part 3*: Fix $\mathrm{r} \in [\mathrm{r}_*, c_0^{-1}\mathrm{r}_0]$ and use $\mathrm{r}_{0*}$ for the moment to denote the upper end point of this interval. Integrate (5.46) from r to $\mathrm{r}_{0*}$ to see that

$$\mathrm{h}(\mathrm{r}_{0*}) \ge c_0^{-1}\left(\tfrac{\mathrm{r}_{0*}}{\mathrm{r}}\right)^{2+1/16}\mathrm{h}(\mathrm{r}). \tag{5.48}$$

As noted previously, $\mathrm{h}(\mathrm{r}_0)$ is at most $c_0\mathrm{E}^2\mathrm{r}_{0*}^2$. Use the latter bound in (5.47) to see that $\mathrm{h}(\mathrm{r}) \le c_0\mathrm{r}^{2+1/16}$ when $\mathrm{r} \in [\mathrm{r}_*, c_0^{-1}\mathrm{r}_0]$. This last conclusion is the assertion of the third bullet of Proposition 5.1.

## 6. Continuity of the limit

This section uses the results from Sections 3-5 to prove that Proposition 2.2's limit function $|\hat{a}_\diamond|$ is continuous. The proposition below makes a formal statement to this effect.

**Proposition 6.1**: *Fix a subsequence* $\Lambda \subset \{1, 2, \ldots\}$ *so that* $\{(A_n, \hat{a}_n)\}_{n\in\Lambda}$ *is described by Proposition 2.2. The limit function* $|\hat{a}_\diamond|$ *given by the second bullet of Proposition 2.2 is continuous. This function is also Hölder continuous with exponent* $\frac{1}{4}$ *on compact where it is bounded away from zero. Meanwhile, there exists* $\kappa > 1$ *such that if* $p \in M$ *and if* $|\hat{a}_\diamond|(p) = 0$, *then* $|\hat{a}_\diamond| \le \mathrm{dist}(p,\cdot)^{1/\kappa}$ *on a sufficiently small radius ball centered at* $p$.

The proof of Proposition 6.1 occupies the remainder of this section. Section 6a proves that $|\hat{a}_\diamond|$ is continuous where it is positive and it proves the Hölder continuity assertion for the points in compact sets where $|\hat{a}_\diamond|$ is positive. The arguments in Section 6a assume Lemma 6.2, this being the crucial input. Section 6b contains the proof of Lemma 6.2. Section 6c proves that $|\hat{a}_\diamond|$ is continuous near its zeros, and it proves that $|\hat{a}_\diamond|$ has the asserted Hölder continuity property at each zero.

### a) Continuity where $|\hat{a}_\diamond| > 0$

The assertion in Proposition 6.1 to the effect that $|\hat{a}_\diamond|$ where positive is Hölder continuous with exponent $\frac{1}{4}$ is seen momentarily to be a consequence of the upcoming Lemma 6.2. The notation for Lemma 6.2 and for the subsequent subsections refers back to notation that was introduced in Section 3b. To say more, fix for the moment $p \in M$. Each $n \in \Lambda$ version of $(A_n, \hat{a}_n)$ has a corresponding version of what is denoted in Section 3b as $r_\diamond$. The $(A_n, \hat{a}_n)$ version is denoted in what follows by $r_{\diamond n}$. Section 3b describes a map that it denotes by $\phi$ from the radius $c_0^{-1} r_\diamond^{-1}$ ball in $\mathbb{R}^3$ to the radius $c_0^{-1}$ ball in M centered at p. The $(A_n, \hat{a}_n)$ version of this map is denoted by $\phi_n$. The corresponding pair $(\phi_n{}^*A_n, r_{\diamond n}{}^{-1}\hat{a}_n)$ is denoted by $(A_{\diamond n}, \hat{a}_{\diamond n})$, this being a pair of connection on $\phi_n{}^*P$ and $\phi_n{}^*P \times_{SO(3)} \mathfrak{su}(2)$ valued 1-form.

**Lemma 6.2**: *Let* $p$ *denote a point in M where* $|\hat{a}_\diamond|(p) > 0$. *Choose a subsequence in* $\Lambda$ *such that the corresponding sequence with n'th term* $|\hat{a}_n|(p)$ *converges to* $|\hat{a}_\diamond|(p)$. *This chosen subsequence has a subsequence, to be denoted by* $\Lambda_p$, *such that* $\liminf_{n\in\Lambda_p} r_{\diamond n} > 0$. *Moreover, there exists data consisting of*

- *A sequence* $\{g_n\}_{n\in\Lambda_p}$ *with n'th member being an isomorphism from the product principal SO(3) bundle over the* $|x| < 1$ *ball in* $\mathbb{R}^3$ *to* $\phi_n{}^*P$.
- *A pair* $(A_\diamond, \hat{a}_\diamond)$ *of* $L^2_{1;loc}$ *connection on the product principal* SU(2) *bundle over the*

$|x| < \kappa^{-1}$ *ball in* $\mathbb{R}^3$*, and an* $\mathfrak{su}(2)$*-valued,* $L^2_{2;loc}$ *1-form on this ball with* $|\hat{a}_\lozenge| = |\hat{a}_\lozenge|(p)$. *These are such that the sequence* $\{g_n{}^*A_{\lozenge n}\}_{n\in\Lambda_p}$ *converges weakly in the* $L^2_{1;loc}$ *topology on the* $|x| < 1$ *ball to* $A_\lozenge$ *and the sequence* $\{g_n{}^*\hat{a}_{\lozenge n}\}_{n\in\Lambda_\lozenge}$ *converges weakly in the* $L^2_{2;loc}$ *topology on the* $|x| < 1$ *ball to* $\hat{a}_\lozenge$.

The proof of Lemma 6.2 is in the Section 6b.

The lemma that follows asserts the parts of Proposition 6.1 that concern the points in M where $|\hat{a}_\lozenge|$ is greater than zero.

**Lemma 6.3**: *Fix a subsequence of* $\Lambda \subset \{1, 2, \ldots\}$ *so that* $\{(A_n, \hat{a}_n)\}_{n\in\Lambda}$ *is described by the six bullets in Proposition 2.2. Let* $|\hat{a}_\lozenge|$ *denote the limit function given by the second bullet of Proposition 2.2. Then the* $|\hat{a}_\lozenge| > 0$ *subset of* M *is open and* $|\hat{a}_\lozenge|$ *on this subset is Hölder continuous with exponent* $\frac{1}{4}$ .

***Proof of Lemma 6.3***: The proof has two parts.

*Part 1*: Fix $p \in M$ where $|\hat{a}_\lozenge|(p) > 0$ and let $\Lambda_p$ denote the corresponding subsequence from Lemma 6.2. Let $r_* > 0$ denote a lower bound for $\{r_{\lozenge n}\}_{n\in\Lambda_p}$ . The convergence of $\{g_n{}^*\hat{a}_{\lozenge n}\}_{n\in\Lambda_p}$ that is asserted by Lemma 6.2 with the fact that each $n \in \Lambda_p$ version of $r_{\lozenge n}$ is greater than $r_*$ imply that the sequence $\{|\hat{a}_n|\}_{n\in\Lambda_p}$ converges strongly in the exponent $\upsilon = \frac{1}{4}$ Hölder topology in the radius $\frac{1}{2} r_*$ ball in M centered at p. That this is so follows from a 3-dimensional Sobolev inequality that asserts in part that $L^2_2$ functions in are Hölder continuous for any Hölder exponent less than $\frac{1}{2}$ .

Since $|\hat{a}_\lozenge(p)| > 0$, so $|\hat{a}_n|(p) > \frac{1}{2} |\hat{a}_\lozenge|(p)$ for all $n \in \Lambda_p$ if n is sufficiently large. This last observation with the Hölder topology convergence implies the following: There exists $r_p \in (0, r_*)$ such that $|\hat{a}_n| > 0$ on the radius $r_p$ ball about p if $n \in \Lambda$ is sufficiently large. Let p´ denote a point in this radius $r_p$ ball. Then $|\hat{a}_\lozenge|(p') > 0$ because the second bullet of Proposition 2.2 has $|\hat{a}_\lozenge|(p') = \limsup_{n\to\infty} |\hat{a}_n|(p')$ and this lim sup is no smaller than $\lim_{n\in\Lambda_p} |\hat{a}_n|(p')$. This proves that the $|\hat{a}_\lozenge| > 0$ part of M is an open set.

*Part 2*: To see about the continuity and Hölder continuity assertions where $|\hat{a}_\lozenge|$ is positive, fix $p \in M$ with $|\hat{a}_\lozenge|(p) > 0$ and let $\Lambda_p$ again denote the subsequence from Lemma 6.2. Let $B_p \subset M$ denote the radius $r_p$ ball centered on p. For each $n \in \Lambda_p$, let $u_n$ denote the automorphism of P's restriction to $B_p$ given by $(\phi_n^{-1})^*(g_1^{-1}g_n)$. The convergence asserted by Lemma 6.2 for the sequence $\{g_n{}^*\hat{a}_{\lozenge n}\}_{n\in\Lambda_p}$ on the unit ball in $\mathbb{R}^3$ implies that the sequence $\{u_n{}^*A_n\}_{n\in\Lambda_p}$ converges weakly in the $L^2{}_{1;loc}$ topology on $B_p$ and that $\{u_n{}^*\hat{a}_n\}_{n\in\Lambda_p}$ converges weakly in the $L^2{}_{2;loc}$ topology on $B_p$. A standard Sobolev inequality implies that convergence occurs in the exponent $\frac{1}{4}$ Hölder topology on compact subsets of this ball. Let $\mathfrak{a}_p$ denote the limit section over $B_p$ of the bundle $(P\times_{SO(3)}\mathfrak{su}(2))\otimes T^*M$. Note that $|\mathfrak{a}_p| = \lim_{n\in\Lambda_p} |(\phi_n{}^{-1})^*\hat{a}_{\Diamond n}|$. It follows in particular from Proposition 2.2 that $|\hat{a}_{\Diamond}| \geq |\mathfrak{a}_p|$ on $B_p$ with equality at p and on a set of full measure, the measure being full because $\{|\hat{a}_n|\}_{n=1,2,\ldots}$ converges in the $L^2$ topology on $B_p$. Let p´ denote another point in $B_p$. The point p´ has a corresponding Hölder continuous $|\mathfrak{a}_{p'}|$ that is defined on a ball $B_{p'}$ about p´ and is such that $|\hat{a}_{\Diamond}| \geq |\mathfrak{a}_{p'}|$ with equality at p´ and on a set of full measure. Both $|\mathfrak{a}_p|$ and $|\mathfrak{a}_{p'}|$ are Hölder continuous functions defined on $B_p\cap B_{p'}$ that are equal on a set of full measure. This can happen only if they are equal at all points in $B_p \cap B_{p'}$. This last observation implies that $|\hat{a}_{\Diamond}| = |\mathfrak{a}_p|$ on $B_p$ and so $|\hat{a}_{\Diamond}|$ is Hölder continuous on $B_p$ with exponent $\frac{1}{4}$.

**b) Proof of Lemma 6.2**

Item c) of Proposition 2.2 has the following consequence: There exists $E > 1$ such that each $n \in \Lambda$ version of $(A_n, \hat{a}_n)$ is described by the $E$ and $r = r_n$ version of (3.2). With this understood, invoke Proposition 3.1 for each $n \in \Lambda$ version of $(A_n, \hat{a}_n)$. Proposition 3.1 supplies corresponding $(A_n, \hat{a}_n)$ versions of what it denotes by g and $\hat{a}_{A_{\Diamond}}$. These versions are denoted in what follows by $g_n$ and $\hat{a}_{A_{\Diamond n}}$.

It follows from what is said in Proposition 3.1 that $\{\hat{a}_{A_{\Diamond n}}\}_{n\in\Lambda}$ is bounded in the $L^2{}_1$ topology on the $|x| \leq 1$ ball. Proposition 3.1 implies that the $\{g_n{}^*\hat{a}_n\}_{n\in\Lambda}$ is also bounded in the $L^2{}_1$ on the $|x| \leq 1$ ball and, if $r < 1$, then it is also bounded in the $L^2{}_2$ topology on each $|x| \leq r$ ball in $\mathbb{R}^3$. It follows as a consequence that there is a subsequence in $\Lambda$, this to be denoted by $\Lambda_{\Diamond}$, such that $\{\hat{a}_{A_{\Diamond n}}\}_{n\in\Lambda_{\Diamond}}$ converges weakly in the $L^2{}_1$ toplogy on the $|x| \leq 1$ ball in $\mathbb{R}^3$ as does $\{g_n{}^*\hat{a}_{\Diamond n}\}_{n\in\Lambda_{\Diamond}}$. The latter sequence also converges weakly in the $L^2{}_{2;loc}$ topology on the $|x| \leq 1$ ball in $\mathbb{R}^3$.

Granted what is said in the preceding paragraph, the proof of Lemma 6.2 needs only a proof of the following assertion: There exists a subsequence $\Lambda_p \subset \Lambda_{\Diamond}$ with the property that $\{r_{\Diamond n}\}_{n\in\Lambda_p}$ is bounded away from zero. This assertion is proved by assuming it false so as to derive nonsense. The derivation of this nonsense has seven parts.

To set the notation, suppose that $n \in \Lambda_{\Diamond}$. The various parts of the proof use $h_n$ to denote the version of the function h in (5.2) that is defined using p and the pair $(A_n, \hat{a}_n)$. The corresponding versions of the functions H and $N$ are denoted by $H_n$ and $N_n$. The arguments that follow also use D to denote the assumed nonzero value of $|\hat{a}_\Diamond|(p)$. No generality is lost by assuming that $\frac{1}{2}D \le |\hat{a}_n|(p) \le 2D$ for all $n \in \Lambda$.

*Part 1*: Let p for the moment denote a given point in M. Fix $r \in (0, c_0^{-1})$ and let $\chi_{p,r}$ denote the function on M given by $\chi(r^{-1}\text{dist}(p,\cdot) - 1)$. This function equals 1 where $\text{dist}(p,\cdot) < r$ and it equals zero where the $\text{dist}(p,\cdot) > 2r$. Let $G_p$ denote the Green's function for the operator $d^\dagger d + 1$ on M with pole at p. Fix $n \in \{1, 2, \ldots\}$ and integrate the function $\chi_{p,r} G_p \langle \hat{a}_n \wedge * q_{A_n}(\hat{a}_n) \rangle$ over M, then integrate by parts so as to derive the following local version of (2.41):

$$\frac{1}{2}|\hat{a}_n|^2(p) + \int_M \chi_{p,r} G_p (|\nabla_A \hat{a}_n|^2 + r_n^{\ 2} |\hat{a}_n \wedge \hat{a}_n|^2) = \mathfrak{e}_{p0}$$

(6.1)

where $|\mathfrak{e}_{p0}|$ is bounded by

$$c_0 r^{1/2} + c_0 r^{-3} \int_{r \le \text{dist}(p,\cdot) \le 2r} |\hat{a}_n|^2 .$$

(6.2)

By way of an explanation, the bound on $|\mathfrak{e}_{p0}|$ follows from two facts, the first being that the $L^2$ norm of $q_{A_n}(\hat{a}_n)$ is bounded by $c_0$ and those of both $d_{A_n}\hat{a}_n$ and $d_{A_n} * \hat{a}_n$ are bounded by $c_0 r_n^{-1}$; these are the bounds asserted by Items c) and d) of Proposition 2.2. The third fact is that the $L^2$ norm of $G_p$ over $B_{2r}$ is less than $c_0 r^{1/2}$. Proposition 2.2 asserts in part that $|\hat{a}_\Diamond|(p) = \limsup_{n\in\Lambda} |\hat{a}_n|(p)$ and that $\{\hat{a}_n\}_{n\in\Lambda}$ converges strongly in the $L^2$ topology on M to $|\hat{a}_\Diamond|$. These facts with (6.1) and (6.2) imply that

$$\frac{1}{2}|\hat{a}_\Diamond|^2(p) \le c_0 r^{1/2} + c_0 r^{-3} \lim_{n\in\Lambda} \int_{r \le \text{dist}(p,\cdot) \le 2r} |\hat{a}_n|^2 .$$

(6.3)

Note that the sequence whose n'th term is the $n \in \Lambda$ version of the integral on the right hand side of (6.3) converges; and the limit is the integral of $|\hat{a}_\Diamond|^2$ over the indicated domain. This is so because $\{|\hat{a}_n|\}_{n\in\Lambda}$ converges strongly in the $L^2$ topology on M to $|\hat{a}_\Diamond|$.

*Part 2*: As explained momentarily, the following two assertions must hold:

- *If* $n \in \Lambda$ *is sufficiently large, then* $h_n(r) \ge c_0^{-1} D^2 r^2$ *for all* $r \in [\frac{1}{2} r_{\Diamond n}, 9r_{\Diamond n}]$.
- $\lim_{n\in\Lambda} \sup_{r\in[\frac{1}{2} r_{\Diamond n}, 9r_{\Diamond n}]} N_n(r) = 0$.

(6.4)

Given (6.4), nothing is lost by assuming that the $n \in \Lambda$ versions of $h_n$ and $N_n$ are such that $h_n(r) \geq c_0^{-2} D^2 r^2$ and $N_n(r) \leq 1$ for all $r \in [\frac{1}{2} r_{\Diamond n}, 9 r_{\Diamond n}]$. Moreover, the second bullet of (6.1) plus the assumption that $\lim_{n \in \Lambda} r_{\Diamond n} = 0$ implies the following: Given $\delta \in (0, \frac{1}{16})$, there exists an integer $N_\delta$ such that

$$\sup_{r \in [\frac{1}{2} r_{\Diamond n}, 9 r_{\Diamond n}]} N_n(r) < \delta \quad \textit{and} \quad r_{\Diamond n} < e^{-1/\delta} \quad \textit{when } n \geq N_\delta \tag{6.5}$$

The subsequent parts of the proof generate the required nonsense when $\delta < c_0^{-1}$.

To prove the top bullet in (6.4), use p's version of (6.3) to conclude that

$$D^2 \leq c_0 r^{1/2} + c_0 r^{-2} \sup_{s \in [r, 2r]} h_n(s) \tag{6.6}$$

when $r \in (0, c_0^{-1}]$. Fix $R > \frac{1}{2}$ for the moment. If n is sufficiently large, then $h_n$ is defined at $r = R r_{\Diamond n}$. Assuming this to be the case, then (6.6) demands a point $s \in [R r_{\Diamond n}, 2R r_{\Diamond n}]$ with $h_n(s) \geq c_0^{-1} D^2 R^2 r_{\Diamond n}^{\;2}$. Fix such a point and integrate the first bullet of the p and $(A_n, \hat{a}_n)$ version of Proposition 5.1 using the fact that $h(t) \leq c_0 t^2$ for any $t \in (0, c_0^{-1})$ to conclude that $h(R r_\Diamond) \geq c_0^{-1} D^2 R^2 r_{\Diamond n}^{\;2} - c_0 (r_n^{-2} R r_{\Diamond n} + r_n^{-1/2} R^2 r_{\Diamond n}^{\;2})$. Since $r_{\Diamond n} \geq c_0^{-1} r_n$ in any event, this last inequality leads directly to the top bullet (6.4) when n is large.

The proof of the second bullet of (6.4) assumes to the contrary that the lim sup in question is positive and derives nonsense. To start the derivation, fix $m > 16$ so that the lim sup in the second bullet of (6.4) is greater than $m^{-1}$. Let $\Lambda'$ denote a corresponding subsequence of $\Lambda$ with the following property: If $n \in \Lambda'$, then there exists $r_{1n} \in [\frac{1}{2} r_{\Diamond n}, 9 r_{\Diamond n}]$ with $N_n(r_{1n}) \geq m^{-1}$.

Fix $n \in \Lambda'$ and let $\kappa$ denote the constant that is supplied by the p and $(A_n, \hat{a}_n)$ version of Proposition 5.1. Note in particular that $\kappa$ is independent of n. The top bullet in (6.4) implies that $h_n(r) \geq r^{3-1/\kappa}$ for $r \in [\frac{1}{2} r_{\Diamond n}, 9 r_{\Diamond n}]$ when n is large. This being the case, there exists a maximal $r_1 \in [9 r_{\Diamond n}, c_0^{-1} r_0]$ such that $h_n(r) > r^{3-1/\kappa}$ for all $r \in [r_{\Diamond n}, r_1]$. Invoke the first bullet of Proposition 5.1 to conclude that

$$\tfrac{d}{dr} h_n \geq 2 r^{-1} (1 + m^{-1} - r^{-1/8} - c_E (r^{1/c_0} - r_*^{\,1/c_0})) h_n - c_E r^{-1/4} r^{5/4} \tag{6.7}$$

for all $r \in [r_{\Diamond n}, r_1]$. Integrate this equation to see that

$$h_n(r) \geq c_0^{-1} D^2 \left(\frac{r}{r_{\Diamond n}}\right)^{1/m} r^2. \tag{6.8}$$

for all $r \in [r_{\Diamond n}, r_1]$. It follows from this that $r_1$ must be equal to $c_0^{-1} r_0$.

Since $h_n(c_0^{-1}r_0) \le c_0 r_0^2$ in any event, the $r = c_0^{-1}r_0$ version of (6.8) can hold only if $r_{\diamond n} \ge c_0^{-1} r_0 D^{-2m}$. The latter conclusion constitutes the desired nonsense because it runs afoul of the assumption at the outset that $\limsup_{n \in \Lambda_\diamond} r_{\diamond n} = 0$.

*Part 3*: To set the notation that is used below, fix $n \in \Lambda$ and $p_{*n} \in M$. The notation uses $r_{*\diamond n}$ to denote the $p_{*n}$ version of $r_{\diamond n}$ as defined using the pair $(A_n, \hat{a}_n)$. The notation also has $h_{*n}$, $\hbar_{*n}$ and $\mathrm{N}_{*n}$ denoting the $p_{*n}$ versions of the functions $h$, $\hbar$ and $\mathrm{N}$. The function that measures the distance to $p_{*n}$ is denoted by $r_{*n}$.

Suppose that $n \in \Lambda$ and that $p_{*n}$ is a point in M with $\mathrm{dist}(p, p_{*n}) \le 3\, r_{\diamond n}$. As proved directly, the following must be true:

- $h_{*n}(r_{*n}) \ge r_{*n}^{2+1/16}$ *for all* $r_{*n} \in [\frac{1}{2} r_{*\diamond n}, 4 r_{\diamond n}]$.
- $\mathrm{N}_{*n}(r_{*n}) < c_0 D^{-2} \delta$ *for all* $r_{*n} \in [\frac{1}{2} r_{*\diamond n}, 4 r_{\diamond n}]$ *if* $n \ge N_\delta$.

(6.9)

If the top bullet in (6.9) were false for some $r_{*n}$ in the indicated range, then the fourth bullet in the $p_{*n}$ and $(A_n, \hat{a}_n)$ version of Proposition 5.1 would apply with $r_*$ being the relevant value of $r_{*n}$. In particular, this bullet would find $h_{*n} \le c_0 r_{*n}^{2+1/16}$ for all $r_{*n}$ between $4r_{\diamond n}$ and $c_0^{-1}$. This would imply that the integral of $|\hat{a}_n|^2$ over the spherical annulus centered at $p_{*n}$ where $4r_{\diamond n} \le \mathrm{dist}(p_{*n}, \cdot) \le 12 r_{\diamond n}$ is no greater than $c_0 r_{\diamond n}^{3+1/16}$. But this is not possible when $n \ge c_0^{-1}$ because the latter spherical annulus contains the spherical annulus centered at p where $7r_{\diamond n} \le \mathrm{dist}(p, \cdot) \le 9r_{\diamond n}$ and it follows from (6.4) that the integral of $|\hat{a}_n|^2$ over the $7r_{\diamond n} \le \mathrm{dist}(p, \cdot) \le 9r_{\diamond n}$ annulus is greater than $c_0^{-1} D^2 r_{\diamond n}^2$.

Granted that the top bullet is true, suppose for the sake of argument that the lower bullet is not true. The three steps that follow generate nonsense from this assumption. To set the stage for what is to come, fix for the moment $m > 1$ and suppose that there exists $n \in \Lambda$ with a corresponding $r_* \in [\frac{1}{2} r_{*\diamond n}, 4r_{\diamond n}]$ where $\mathrm{N}_{*n}(r_*) \ge 2m\delta$.

<u>Step 1</u>: If $\mathrm{N}_{*n}(r_*) \ge \frac{1}{16}$, then the argument in Part 2 of Section 5e can be repeated to see that (5.48) must hold. The latter conclusion is untenable when n is large because it implies that $h_n(r_{\diamond n}) < c_0 r_{\diamond n}^{2+1/16}$. This understood, assume that $m\delta \le \mathrm{N}_*(r_{*n}) < \frac{1}{16}$.

<u>Step 2</u>: There exists in any case $r > r_*$ such that $h_{*n}(s) > s^{2+3/4}$ for all $s \in [r_*, r]$. It then follows from the second bullet of Proposition 5.1 that

$$\mathrm{N}_{*n}(s) \ge m\delta - c_0 (s^{1/c_0} - r_*^{1/c_0}) \ \textit{for all}\ s \in [r_*, r]\ . \tag{6.10}$$

This inequality is untenable when n is large if $r = 4r_{\Diamond n}$ for the following reason: The ball of radius $4r_{\Diamond n}$ centered on $p_{*n}$ is contained in the ball of radius $7r_{\Diamond n}$ centered on p. This implies that the integral of $|\nabla \hat{a}_n|^2$ over the ball of radius $7r_{\Diamond n}$ is greater than $c_0^{-1} m \delta r_{\Diamond n}$ when n is large. But the latter integral can be no larger than $c_0 \mathrm{D}^{-2} \delta r_{\Diamond n}$, this being a consequence of (6.4). These two bounds can not hold simultaneously if $m > c_0 \mathrm{D}^{-2}$.

<u>Step</u> <u>3</u>: If $m = c_0 \mathrm{D}^{-2}$, then it follows from what was said in Step 2 that there exists $r_\ddagger \in [r_*, 4r_{\Diamond n}]$ such that $h_{*n}(r_\ddagger) \le r_\ddagger^{2+3/4}$. This understood, then the fourth bullet of Proposition 5.1 can be invoked using $r_\ddagger$ in lieu of $r_*$ to conclude that $h_{*n}(r_{*n}) \le c_0 r_{*n}^{2+1/16}$ when $r_{*n} \in [4r_{\Diamond n}, c_0^{-1}]$ and n is large. Repeat verbatim the argument for the first bullet in (6.9) to see that this is a nonsensical conclusion.

*Part 4*: Fix $n \ge N_\delta$ and $p_{*n} \in M$ with $\mathrm{dist}(p, p_{*n}) \le 3r_{\Diamond n}$. Introduce by way of notation $B_{*n}$ to denote the ball of radius $r_{*\Diamond n}$ centered on $p_{*n}$. Let $\phi_{*n}$ denote the $p_{*n}$ and $(A_n, \hat{a}_n)$ version of the map $\phi$ that is described at the outset of Section 3b, this to be viewed as a map from the radius 4 ball about the origin in $\mathbb{R}^3$ to the radius $4r_{*\Diamond n}$ ball centered on $p_{*n}$. Use this map to define the $p_{*n}$ version of the pair that is denoted in Section 3a by $(A_\Diamond, \hat{a}_\Diamond)$. The $p_{*n}$ version is denoted by $(A_{*\Diamond n}, \hat{a}_{*\Diamond n})$; it is a pair whose left hand member is connection on the $\phi_{*n}$ pull-back of P and whose right hand member is a section of the $\phi_{*n}$ pull-back of the bundle $(P \times_{SO(3)} \mathfrak{su}(2)) \otimes T^*M$. Denote by $\hat{a}_{**n}$ the $p_{*n}$ version of what is denoted in Section 3b by $\hat{a}_*$, this having the form $z_{*n}^{-1} \hat{a}_{*\Diamond n}$ with $z_{*n}$ denoting the $L^2$ norm of $\hat{a}_{*n}$ on $B_{*n}$. Thus, $\hat{a}_{**n}$ has $L^2$ norm 1 on the $|x| \le 1$ ball in $\mathbb{R}^3$. Given the definition of $\mathrm{N}$, the bound in the second bullet of (6.9) implies that

$$\int_{|x|\le 1} |\nabla_{A_{*\Diamond n}} \hat{a}_{**n}|^2 \le c_{\mathrm{D}} \delta, \tag{6.11}$$

with $c_{\mathrm{D}}$ denoting here and in what follows a number that is greater than 1 that depends only on $\mathrm{D}$. The value of $c_{\mathrm{D}}$ can be assumed to increase between successive appearances.

*Part 5*: Fix $n \ge N_\delta$ and introduce next $\tilde{a}_{**n}$ to denote $p_{*n}$ version of what is denoted by $\tilde{a}_*$ in Section 4a. The latter obeys the version of Proposition 4.1 that uses $(A_{*\Diamond n}, \hat{a}_{**n})$, $r_{*\Diamond n}$ and $r_{z*n}$ in lieu of $(A_{\Diamond n}, \hat{a}_{*n})$, $r_{\Diamond n}$ and $r_{zn}$. Given that $r_{z*n} \ge c_{\mathrm{D}}^{-1}$, $r_{*\Diamond n} \le c_0 r_{\Diamond n}$ and $r_{\Diamond n} \le e^{-1/\delta}$, the second bullet of Proposition 4.1 implies that $\|\tilde{a}_{**n}\|_2$ differs from 1 by at most $c_0 e^{-1/\delta}$. Meanwhile, the third bullet of Proposition 4.1 and (6.11) imply that

$$\int_{|x|\le 1} |\nabla_{A_{*\lozenge n}} \tilde{a}_{**n}|^2 \le c_D \delta .$$

(6.12)

With an appeal to Proposition 3.2 in mind, fix $\varepsilon \in (0, 1)$ and introduce $\kappa_\varepsilon$ to denote the $\varepsilon$ and $\mu = \frac{1}{2}$ and $E = c_0$ version of what is denoted by $\kappa_{E,\mu,\varepsilon}$ in Proposition 3.2. If it is the case that $\delta < c_D^{-1}\kappa_\varepsilon^{-1}$ and also $r_{\lozenge n} < c_D^{-1}\kappa_\varepsilon^{-1}$ then Proposition 3.2 and (3.8) say that the square of the $L^2$ norm of $F_{A_{*\lozenge n}}$ on the $|x| < \frac{3}{4}$ ball in $\mathbb{R}^3$ is less than $\varepsilon$.

*Part 6*: Fix an integer $n > N_\delta$. This part of the proof defines an iterative procedure that starts with p and generates from p a finite sequence of points. The k'th point in the sequence is denoted by $p_{*n,k}$. It is such that $\text{dist}(p, p_{*n,k}) < 3 r_{\lozenge n}$. The description that follows of the iteration step uses $p_{*n,0}$ to denote p. If some $k \ge 0$ version of $p_{*n,k}$ has been defined, the notation has $r_{*\lozenge n,k}$ denoting the $p_{*n,k}$ version of $r_{*\lozenge n}$.

THE ITERATION PROCEDURE: Fix $k \ge 0$ and suppose that that the point $p_{*n,k}$ has been defined. Let $p_{*n} \in M$ denote any given point with $\text{dist}(p_{*n,k}, p_{*n}) \le 2 r_{*\lozenge n,k}$. The assignment of $r_{*\lozenge n}$ to $p_{*n}$ defines a continuous function on the radius $2r_{*\lozenge n,k}$ ball about $p_{*n,k}$. If the minimum value is less than or equal to $\frac{1}{256} r_{*\lozenge n,k}$, take $p_{*n,k+1}$ to be a point with distance $2r_{*\lozenge n,k}$ or less from $p_{*n,k}$ where the minimum value is achieved. If the minimum value of this function is greater than $\frac{1}{256} r_{*\lozenge n,k}$, then there is no (k+1)'st point and $p_{*n,k}$ is the last point in the sequence. Note that the iteration must stop after a finite number of runs because each $p_{*n}$ version of $r_{*\lozenge n}$ is, in any event, greater than $c_0^{-1} r_n^{-1}$.

The sequence $\{p_{*n,k}\}_{k=1,2,\ldots}$ lies in the radius $(2 + \frac{1}{2}) r_{\lozenge n}$ ball centered on p. To understand why this is, keep in mind that

$$\text{dist}(p_{*n,k}, p_{*n,k+1}) \le \tfrac{1}{128} r_{*\lozenge n,k} .$$

(6.13)

Meanwhile, $r_{*\lozenge n,k} < \frac{1}{256} r_{*\lozenge n,k-1}$ and so $r_{*\lozenge n,k} \le (\frac{1}{256})^k r_{\lozenge n}$. Thus, $\text{dist}(p_{*n,k}, p_{*n,k+1}) \le 4(\frac{1}{256})^k r_{\lozenge n}$ and so $\text{dist}(p, p_{*\lozenge n,k+1}) \le (2 + 4(\sum_{m=1,2,\ldots k} (\frac{1}{256})^m)) r_{\lozenge n}$ which is no greater than $(2 + \frac{1}{32}) r_{\lozenge n}$.

*Part 7*: Let $p_{*n,k}$ now denote the final point in the iteration sequence $\{p, p_{*n,1}, \ldots\}$. It follows from the definition that the radius $2r_{*\lozenge n,k}$ ball centered on $p_{*n,k}$ is contained in the radius $3r_{\lozenge n}$ ball centered around p. If $p_{*n}$ is a given point in the radius $2r_{*\lozenge n,k}$ ball centered $p_{*n,k}$, then its corresponding $r_{*\lozenge n}$ is no less than $\frac{1}{256} r_{*\lozenge n,k}$. This understood, it follows that the radius $2r_{*\lozenge n,k}$ ball centered on $p_{*n,k}$ has a cover, $\mathfrak{U}$, with the following properties:

- $\mathfrak{U}$ *consist of at most* $c_0^{-1}$ *balls with centers in the radius* $2r_{*\lozenge n,k}$ *ball about* $p_{*n,k}$.

- *Let* $p_{*n}$ *denote a center point of a ball from* $\mathfrak{U}$. *The radius of its ball is* $\frac{1}{2}\, r_{*\lozenge n}$.

(6.14)

Let $p_{*n}$ denote the center point of a given ball from $\mathfrak{U}$ and let $B_{*n}$ denote the corresponding ball from $\mathfrak{U}$. Use what is said at the end of Part 5 to see that the square of the $L^2$ norm of $F_{A_n}$ over $B_{*n}$ is at most $\varepsilon^2 r_{*\lozenge n}{}^{-1}$ and thus at most $256\varepsilon^2 r_{*\lozenge n,k}{}^{-1}$. Since there are at most $c_0{}^{-1}$ balls in $\mathfrak{U}$, this bound implies in turn that

$$\int_{\mathrm{dist}(p_{*n,k},\,\cdot\,)\le 2r_{*\lozenge n,k}} |F_{A_n}|^2 \le c_D \varepsilon^2 r_{*\lozenge n,k}{}^{-1}\,.$$

(6.15)

The latter bound is nonsensical if $\varepsilon < c_D{}^{-1}$ because it runs afoul of the definition of $r_{*\lozenge n,k}$. This is the promised nonsense from the small $\delta$ versions of (6.5).

**c) $|\hat{a}_\lozenge|$ near its zero locus**

This subsection proves the assertions in Proposition 6.1 that concern $|\hat{a}_\lozenge|$ on its zero locus. The proof has two parts.

*Part 1*: The lemma that follows plays a central role in the subsequent arguments. To set the notation, suppose that $p \in M$. Given a positive integer n, the lemma uses $h_n$ to denote the version of the function h in (5.2) that is defined using p and the pair $(A_n, \hat{a}_n)$.

**Lemma 6.4**: *Fix* $p \in M$ *and* $c > 2$. *Suppose there exists a subsequence in Proposition 6.1's sequence* $\Lambda$*, this denoted by* $\Lambda_p$*, such that* $\lim_{n\in\Lambda_p} h_n(r) \le c\, r^{1/c} r^2$ *when* $r \le c^{-1}$. *Then* $|\hat{a}_\lozenge|(p) = 0$. *Moreover,* $|\hat{a}_\lozenge|$ *is continuous at* p *and there exists a constant* $\kappa > 1$ *depending only on* $c$ *such that* $|\hat{a}_\lozenge|(\cdot) \le \kappa\, \mathrm{dist}(p,\cdot)^{1/\kappa}$ *on the ball of radius* $\kappa^{-1}$ *centered at* p.

***Proof of Lemma 6.4***: Fix $r \in (0, c_0{}^{-1})$ and it follows from (6.3) that $|\hat{a}_\lozenge|^2(p)$ is bounded by $c_0(r^{1/2} + r^{1/c})$. It follows as a consequence that $|\hat{a}_\lozenge|(p) = 0$. Let q denote some other point in M with $\mathrm{dist}(p,q) < c_0{}^{-1}$. Take r to equal $4\,\mathrm{dist}(p,q)$ and use q's version of (6.3) to see that

$$|\hat{a}_\lozenge|^2(q) \le c_0 r^{1/2} + c_0 r^{-3} \lim_{n\in\Lambda} \int_{r/2\le \mathrm{dist}(p,\cdot)\le 4r} |\hat{a}_n|^2$$

(6.16)

The final arguments for Lemma 6.4 exploit (6.16). To this end, fix $n \in \Lambda$ for the moment and let $h_n(\cdot)$ denote the version of the function h in (5.2) that is defined by the point p and the pair $(A_n, \hat{a}_n)$. If $\text{dist}(p,q) < c_0^{-1}$, then $h_n$ is defined on the interval $[0, 4r]$. Granted this is the case, use the definition of $h_n$ to conclude that

$$\int_{B_{4r} - B_{r/2}} |\hat{a}_n|^2 \leq c_0\, r \sup_{s \in [r/2, 4r]} h_n(s)\ . \tag{6.17}$$

Let $c$ and $\Lambda_p$ be as described by Lemma 6.4. If 4r is less than $c^{-1}$, then the right hand side of (6.17) is no greater than $c_0\, c\, r^{1/c} r^3$ when $n \in \Lambda_p$ is large. This implies that the right hand side of (6.16) is no greater than $c_0 (r^{1/2} + c\, r^{1/c})$ because the limit of a convergent sequence is equal to the limit of any its subsequences. Thus $|\hat{a}_\Diamond|^2(q) \leq c_0 (r^{1/2} + c\, r^{1/c})$. This is the assertion made by Lemma 6.4 because $r = 4\,\text{dist}(p,q)$.

*Part 2*: The assertion in Proposition 6.1 to the effect that $|\hat{a}_\Diamond|$ is continuous across its zero locus follows from the assertion in the proposition to the effect that $|\hat{a}_\Diamond|$ is Hölder continuous at each of its zeros for a fixed Hölder exponent. The proof of the uniform Hölder bound along the zero section of $|\hat{a}_\Diamond|$ is given in five steps. These steps employ the following terminology: The *local Hölder property* is said to hold at a given point p if there exist numbers $\kappa > 1$ and $\rho > 0$ such that $|\hat{a}_\Diamond(q)| < \text{dist}(p,q)^{1/\kappa}$ for all $q \in M$ with $\text{dist}(p,q) < \rho$. The Hölder assertion in Proposition 6.1 follows with a proof that the each point in Z has the local Hölder property with $\kappa$ being independent of the given point.

By way of a heads-up, the following observation is used implicitly in what follows: Suppose that $p \in Z$ and that there exists $\kappa > 1$ and $\rho' > 0$ and $x_p > 1$ such that $|\hat{a}_\Diamond(q)| \leq x_p\, \text{dist}(p,q)^{2/\kappa}$ when $\text{dist}(p, q) \leq \rho'$. Then $|\hat{a}_\Diamond|(q) \leq \text{dist}(p,q)^{1/\kappa}$ when $\text{dist}(p,q)$ is less than the minimum of $\rho'$ and $x_p^{\kappa}$.

<u>Step 1</u>: Fix $p \in M$ with $|\hat{a}_\Diamond|(p) = 0$. Since $|\hat{a}_\Diamond|(p) = 0$, so $\lim_{n\to\infty} |\hat{a}_n|(p) = 0$ also, this a consequence of the second bullet in Proposition 2.2. There are now two cases to consider, the first being where p's version of the sequence $\{r_{\Diamond n}\}_{n\in\Lambda}$ has a subsequence that is bounded away from zero. If such is the case, let $r_* > 0$ denote a lower bound for this subsequence. The corresponding subsequence of $\{\hat{a}_n\}_{n\in\Lambda}$ is bounded in the $L^2_2$ topology on the radius $r_*$ ball in M centered at p, and so it has a subsequence that converges strongly in the exponent $\frac{1}{4}$ Hölder topology on this ball. Let $\Lambda_p \subset \Lambda$ denote the indexing set for the latter subsequence. The Hölder convergence of $\{\hat{a}_n\}_{n\in\Lambda_p}$ on the radius $r_*$ ball centered at p has the following consequence: Given $\varepsilon > 0$, there exists $N_\varepsilon$ such that if $n \in \Lambda_p$ and $n > N_\varepsilon$, then

$$|\hat{a}_n| \le \varepsilon + \mathrm{dist}(p,\cdot)^{1/4} \quad \textit{on the radius } r_* \textit{ ball centered at } p. \tag{6.18}$$

Fix $n \in \Lambda_p$ with $n > N_\varepsilon$. The bound in (6.18) implies that the p and $(A_n, \hat{a}_n)$ version of the function $h_n(\cdot)$ obeys $h_n(r) \le c_0(\varepsilon + r^{1/4})r^2$ when $r \in (0, r_*)$. Granted this last bound, invoke Lemma 6.4 to see that the local Hölder property holds at p with Hölder exponent $\frac{1}{4}$.

Step 2: Assume here and in the subsequent steps that $\{(A_n, \hat{a}_n)\}_{n\in\Lambda}$ is such that $\lim_{n\in\Lambda} |\hat{a}_n|(p) = 0$ and $\lim_{n\in\Lambda} r_{\Diamond n} = 0$. Granted this assumption, then at least one of the three cases in the subsequent list describes $\{(A_n, \hat{a}_n)\}_{n\in\Lambda}$. Step 3 contains the proof that the list is inclusive.

CASE 1: This case occurs if there is a subsequence $\Lambda_p \subset \Lambda$ with two properties, the first being the following: If $n \in \Lambda_p$, then there exists $r_{\ddagger n} \in [\frac{1}{2} r_{\Diamond n}, c_0^{-1}]$ which is such that $h_n(r_{\ddagger n}) \le r_{\ddagger n}^{2+1/16}$. The second property requires that $\lim_{n\in\Lambda_p} r_{\ddagger n} = 0$. Fix $n \in \Lambda_p$. Let $r_{1n} \in [r_{\ddagger n}, c_0^{-1}]$ denote the maximal value for r such that $h_n(s) \le s^{2+1/16}$ for all $s \in [r_{\ddagger n}, r_{1n}]$. It follows from the fourth bullet of Proposition 5.1 that $h_n(r) \le c_0 r^{2+1/16}$ for all $r \in [r_{1n}, c_0^{-1}]$; it follows from the definitions of $r_{1n}$ and $r_{\ddagger n}$ that $h_n(r) < r^{2+1/16}$ for all $r \in [r_{\ddagger n}, r_{1n}]$. Since $\lim_{n\in\Lambda_p} r_{\ddagger n} = 0$, the subsequence $\Lambda_p$ with any $c > c_0$ can be used as input to Lemma 6.4 to prove that the local Hölder property holds at p with Hölder exponent greater than $c_0^{-1}$.

CASE 2: This case occurs if there exists $\delta > 0$ and a subsequence $\Lambda' \subset \Lambda$ with the following property: If $n \in \Lambda'$, then $h_n(r) > r^{2+1/16}$ for all $r \in [\frac{1}{2} r_{\Diamond n}, 9r_{\Diamond n}]$ and there exists $r \in [\frac{1}{2} r_{\Diamond n}, 9r_{\Diamond n}]$ with $N_n(r) \ge \delta$. Let $r_{1n} \in (9r_{\Diamond n}, c_0^{-1}]$ denote the maximal r which is such that $h_n(s) \ge s^{2+1/16}$ for all $s \in [r_{\Diamond n}, r_{1n}]$.

Suppose first that $\liminf_{n\in\Lambda'} r_{1n} = 0$. Fix $n \in \Lambda'$. The fourth bullet of Proposition 5.1 implies that $h_n(r) \le c_0 r^{2+1/16}$ for all $r \in [r_{1n}, c_0^{-1}]$. Fix a subsequence $\Lambda_p \subset \Lambda'$ such that $\lim_{n\in\Lambda_p} r_{1n} = 0$. The fact that $\lim_{n\in\Lambda_p} r_{1n} = 0$ implies that $\Lambda_p$ and any $c > c_0$ version of Lemma 6.4 can again be used to prove that the local Hölder property holds at p.

Suppose on the other hand that there exists $r_0 < c_0^{-1}$ such that $\liminf_{n\in\Lambda'} r_{1n} > 2r_0$. Fix $n \in \Lambda'$ such that $r_{1n} > r_0$. Then $h_n(r) > r^{2+1/16}$ on $[r_{\Diamond n}, r_0]$. This being the case, the second bullet of Proposition 5.1 can be invoked to see that $N_n(r) \ge \frac{1}{2}\delta$ if $r \in [\frac{1}{2} r_{\Diamond n}, c_0^{-1} r_0]$. With this understood, invoke the first bullet of Proposition 5.1 to obtain the inequality

$$\frac{d}{dr} h_n \ge 2r^{-1}(1 + \tfrac{1}{2}\delta) h_n - c_0 \delta^{-1} r_n^{-1} r \tag{6.19}$$

where $r \in [\frac{1}{2} r_{\Diamond n}, c_0^{-1} r_0]$. Fix r in this range, and integrate (6.19) from r to $c_0^{-1} r_0$ and use the fact that $h_n(c_0^{-1} r_0) \le c_0 r_0^2$ to conclude that $h_n(r) \le c_0(r_0^{-1/(2\delta)} r^{2+1/(2\delta)} + r_n^{-1} r_0^2)$. Let $r_{\ddagger n}$ denote the number $r_n^{-1/(2+1/(2\delta))} r_0$ and let $c$ denote $c_0(\delta^{-1} + c_0 r_0^{-1/2\delta})$. Then the preceding bound on $h_n(r)$ implies that $h_n(r) \le c\, r^{2+1/c}$ when $r \in [r_{\ddagger n}, c_0^{-1} r_0]$. Noting that $\lim_{n\in\Lambda'} r_{\ddagger n} = 0$, Lemma 6.4 can be invoked using as input $\Lambda_p = \Lambda'$ and $c$ to prove that the local Hölder property assertion holds at p with Hölder exponent greater than $c_0^{-1}\delta$.

The statement of the third case reintroduces notation from Part 3 of Section 6b.

CASE 3: This case occurs when three conditions are met. The first condition requires that $h_n(r) > r^{2+1/16}$ for all $r \in [\frac{1}{2} r_{\lozenge n}, 9r_{\lozenge n}]$ when $n \in \Lambda$ is sufficiently large; and the second condition requires that $\lim_{n\in\Lambda} \sup_{r\in[\frac{1}{2}r_{\lozenge n}, 9r_{\lozenge n}]} N_n(r) < \delta$. The third condition requires there be a subsequence $\Lambda' \in \{1, 2, \ldots\}$ and an associated sequence $\{p_{*n}\}_{n\in\Lambda'} \subset M$ with the following properties.

- *Each* $n \in \Lambda'$ *version of* $p_{*n}$ *has distance less than* $3r_{\lozenge n}$ *from* p.
- *Either or both of the following statements are true.*
  i) *If* $n \in \Lambda'$*, then there exists* $r_{*\ddagger n} \in [\frac{1}{2} r_{*\lozenge n}, 9r_{*\lozenge n}]$ *such that* $h_{*n}(r_{*\ddagger n}) \le r_{*\ddagger n}^{2+1/16}$.
  ii) $\sup_{r\in[\frac{1}{2}r_{*\lozenge n}, 9r_{*\lozenge n}]} N_{*n}(r) \ge \delta$

(6.20)

Suppose there is a subsequence $\Lambda'' \subset \Lambda'$ such that Item i) in the second bullet of (6.20) holds for all $n \in \Lambda''$. Fix $n \in \Lambda''$ and let $h_{*n}$ denote the $p_{*n}$ version of h. But for cosmetic changes, the argument in Case 1 can be used with $p_{*n}$ replacing p to see that $h_{*n}(s) \le c_0 s^{2+1/16}$ for all $s \in [9r_{*\lozenge n}, c_0^{-1}]$. Keeping this in mind, use an integration by parts with the fact that $\|\nabla_{A_n} \hat{a}_n\|_2 \le c_0$ and $|\hat{a}_n| \le c_0$ to see that

$$h_n(s) \le h_{*n}(s + 4r_{\lozenge n}) + c_0 r_{\lozenge n}^{3/2},$$

(6.21)

when $s > 10 r_{\lozenge n}$. Fix $r > 0$ and (6.21) implies that $\lim_{n\in\Lambda'} h_n(r) \le c_0 r^{2+1/16}$. This being the case, then Lemma 6.4 can be invoked using as input $\Lambda_p = \Lambda''$ and any $c \ge c_0$ to prove that p has the local Hölder property.

Suppose next that Item i) of (6.20) is not true if $n \in \Lambda'$ is large. This understood, throw out the finite set of integer where Item i) is true and use $\Lambda'$ now to denote the remaining set. Each $n \in \Lambda'$ obeys the condition in Item ii) of (6.20). But for cosmetic changes, the argument in Case 2 can be used with each $n \in \Lambda'$ version of $p_{*n}$ replacing p to obtain the following data: A number $c > 1$ that depends only on $\delta$ and a sequence $\{r_{*\ddagger n}\}_{n\in\Lambda'} \subset (0, c^{-1})$ with limit zero and with the following additional property: If $n \in \Lambda'$, then $r_{*\ddagger n}$ is such that $h_{*n}(s) \le c\, s^{2+1/c}$ when $s \in [r_{*\ddagger n}, c^{-1}]$. Granted this data, use (6.21) to conclude that $\lim_{n\in\Lambda'} h_n(r) \le c_0\, c\, r^{2+1/c}$ for each $r \in (0, c^{-1})$. It follows from the latter bound that the sequence $\Lambda_p = \Lambda'$ and the given value of $c$ can be used as input to Lemma 6.4 to prove that p has the local Hölder property with Hölder exponent greater than $c_0^{-1}\delta$.

Step 3: Assume that $\{(A_n, \hat{a}_n)\}_{n\in\Lambda}$ is such that $\lim_{n\in\Lambda}|\hat{a}_n|(p) = 0$ and $\lim_{n\to\infty} r_{\Diamond n} = 0$. The paragraphs that follow proves that at least one of the three cases in Step 2 apply with $\delta$ being greater than $c_0^{-1}$. To this end, assume to the contrary that none of these cases for a given $\delta$. The existence of such a sequence is shown below to lead to nonsense when $\delta$ is smaller than $c_0^{-1}$.

After discarding a finite set of terms and then relabling the result as $\Lambda$, then any given $n \in \Lambda$ pair from the sequence $\{(A_n, \hat{a}_n)\}_{n\in\Lambda}$ must have the following properties:

- $h_n(r) \ge r^{2+1/16}$ *for all* $r \in [\frac{1}{2} r_{\Diamond n}, 4 r_{\Diamond n}]$ .
- $N_n(r) < \delta$ *for all* $r \in [\frac{1}{2} r_{\Diamond n}, 4 r_{\Diamond n}]$ .
- *Supposing that* $p_{*n}$ *has distance less than* $3 r_{\Diamond n}$ *from* p *and that* $r \in [\frac{1}{2} r_{*\Diamond n}, 9r_{*\Diamond n}]$, *then* $h_{*n}(r) \le r^{2+1/16}$ *and* $N_{*n}(r) < \delta$.

(6.22)

Indeed, the first bullet of (6.22) must be obeyed to avoid a CASE 1 label, the second bullet of (6.22) must be obeyed to avoid a CASE 2 label, and the third bullet of (6.22) must be obeyed to avoid a Case 3 label.

Step 4: Fix $n \in \Lambda$ and a point $p_{*n} \in M$ with $\text{dist}(p,p_{*n}) \le 3r_{\Diamond n}$. The constructions in Part 4 of Section 6b can be repeated to construct what is denoted there by $\hat{a}_{**n}$. The $L^2$ norm of $\hat{a}_{**n}$ on the $|x| \le 1$ ball in $\mathbb{R}^3$ is equal to 1. Moreover, the condition on $N_{*n}$ in the third bullet of (6.22) implies that

$$\int_{|x|\le 1} |\nabla_{A_{*\Diamond n}} \hat{a}_{**n}|^2 \le c_0 \delta .$$

(6.23)

With (6.23) understood, fix $\varepsilon \in (0, 1]$ and let $\kappa_\varepsilon$ denote the $E \le c_0$ and $\mu = \frac{1}{4}$ version of what is denoted by $\kappa_{E,\mu,\varepsilon}$ in Proposition 3.2 and (3.8). It follows from (6.23) that if $\delta \le c_0^{-1}\kappa_\varepsilon^{-1}$ and if n is large so that $r_{\Diamond n} \le c_0^{-1}\kappa_\varepsilon^{-1}$, then the square of the $L^2$ norm of $F_{A_{*\Diamond n}}$ on the $|x| \le \frac{3}{4}$ ball is bounded by $\varepsilon$.

Step 5: Having fixed $n \in \Lambda$, repeat the iteration procedure in Part 6 of Section 6b to construct an interation sequence $\{p, p_{*n,1}, \ldots, p_{*n,k}\}$. The final paragraph of Step 4 implies that what is said about $p_{*n,k}$ in Part 7 in Section 6b holds in this case also. In particular, (6.15) holds. The latter inequality is nonsensical of $\varepsilon < c_0^{-1}$ and n is large for the same reason it is nonsensical in Section 6b: It runs afoul of the definition of $r_{*\Diamond n,k}$.

## 7. The data Z, $\mathcal{I}$ and $\nu$

The forthcoming proposition is used to complete the proof of Theorem 1.1a. . By way of notation, the proposition denotes the zero locus of $|\hat{a}_{\diamond}|$ by Z. The following notation is also used: Fix $p \in M$ and $r > 0$. The proposition uses $B_r$ to denote the radius r ball centered at p and it uses $\partial B_r$ to denote the boundary of the closure of $B_r$. Keep in mind that the term *real line bundle* is used here to describe the associated line bundle to a principal $\mathbb{Z}/2\mathbb{Z}$ bundle.

**Proposition 7.1**: *The set* $Z \subset M$ *is a closed set and so* $M{-}Z$ *is open. There is a real line bundle over* $M{-}Z$*, this denoted by* $\mathcal{I}$*, and a harmonic section,* $\nu$*, of* $T^*(M{-}Z) \otimes \mathcal{I}$ *with the properties listed below.*

- $|\nu| = |\hat{a}_{\diamond}|$.
- $|\nabla\nu|$ *extends from* $M{-}Z$ *as an* $L^2$ *function on* $M{-}Z$.
- *For any point* p *in* M*, the function* $|\nabla\nu|\,\mathrm{dist}(\cdot, p)^{-1/2}$ *is an* $L^2$ *function on* $M{-}Z$ *and there is a* p*-independent bound on its* $L^2$ *norm.*
- *There exists* $\kappa \geq 1$ *with the following significance. Fix* $p \in M$ *to define functions* h *and* H *on* $(0, \kappa^{-1})$ *by the rules*

$$r \to h(r) = \int_{\partial B_r} |\nu|^2 \quad \textit{and} \quad r \to H(r) = \int_{B_r - Z} |\nabla\nu|^2 \ .$$

  a) *The function* h *is strictly positive on* $(0, \kappa^{-1})$.
  b) *Define the function* $\mathrm{N}$ *on* $(0, \kappa^{-1})$ *by the rule* $r \to \mathrm{N}(r) = r\,H(r)/h(r)$. *The function* h *is differentiable on* $[0, \kappa^{-1})$. *Moreover, its derivative on* $(0, \kappa^{-1})$ *can be written as* $\frac{d}{dr} h = 2r^{-1}(1 + \mathrm{N} + \mathfrak{e})h$ *with* $\mathfrak{e}$ *such that* $|\mathfrak{e}| \leq \kappa r^2$.
  c) *If* $s > r > 0$*, then* $\mathrm{N}(s) \geq e^{-\kappa(s^2 - r^2)}\,\mathrm{N}(r) - \kappa(s^2 - r^2)$.

With regards to the Item b) of the fourth bullet, keep in mind that $\mathrm{N}h$ is $rH$ and so $r^{-1}\mathrm{N}h$ has limit zero as r limits to 0. This is also the case for $r^{-1}h$ because $h \leq c_0 r^2$.

The function $\mathrm{N}$ plays the role here of the frequency function introduced by [Al] and used by [HHL] and [Han]. It is perhaps needless to say that Proposition 7.1's version of $\mathrm{N}$ is the analog for $\nu$ of the function in (5.3).

Granted that Z is closed and granted the second and third bullets of the proposition, the convention in what follows is to extend $|\nabla\nu|$ and any $p \in M$ version of $\mathrm{dist}(p, \cdot)^{-1}|\nabla\nu|$ as an $L^2$ functions on the whole of M by declaring them to be zero on Z. For example, this convention writes H(r) as $\int_{B_r} |\nabla\nu|^2$ .

Sections 7a and 7b contain the proof of Proposition 7.1. Section 7c contains a lemma that concern the behavior of the $r \to 0$ limit of Proposition 7.1's function N as the point p is varied in M. Section 7d proves that $|\nu|$ is uniformly Hölder continuous on M.

## a) The construction of $\mathcal{I}$ and $\nu$

Proposition 6.1 asserts in part that $|\hat{a}_\Diamond|$ is continuous, and this implies directly that its zero locus, Z, is a closed set. The subsequent lemma provides what is needed to define $\mathcal{I}$ and $\nu$. The notation is that used by Proposition 6.1

**Lemma 7.2**: *Let* $\Lambda \subset \{1, 2, \ldots\}$ *denote the subsequence from Proposition 6.1. There exists a subsequence* $\Lambda_\Diamond \subset \Lambda$ *such that the corresponding sequence* $\{(A_n, \hat{a}_n)\}_{n \in \Lambda_*}$ *has the properties listed below.*

- *The sequence of* $\{F_{A_n}\}_{n \in \Lambda_*}$ *has bounded $L^2$ norm on compact subsets of* M–Z *and* $\{\nabla_{A_n} \hat{a}_n\}_{n \in \Lambda_*}$ *has bounded $L^2$ norm on* M.
- *There exists a sequence* $\{h_n\}_{n \in \Lambda_*} \subset \mathrm{Aut}(P)$ *such that each* $n \in \Lambda_*$ *version of* $h_n{}^*A_n$ *can be written on* M–Z *as* $A_0 + \hat{a}_{h^*A_n}$ *with* $\{\hat{a}_{h^*A_n}\}_{n \in \Lambda_*}$ *having bounded $L^2_1$ norm on any given compact subset of* M–Z *and converging in the $L^2_{1;loc}$ topology on* M–Z. *Meanwhile* $\{h_n{}^*\hat{a}_n\}_{n \in \Lambda_*}$ *has bounded $L^2_2$ norm on any given compact set in* M–Z *and it converges weakly in the $L^2_{2;loc}$ topology on* M–Z.
- *Let* $(A_\Diamond, \mathfrak{a}_\Diamond)$ *denote the limit pair of $L^2_{1;loc}$ connection on* $P|_{M-Z}$ *and $L^2_{2;loc}$ section over* M–Z *of* $(P \times_{SO(3)} \mathfrak{su}(2)) \otimes T^*M$. *These are such that*
  a) $d_{A_\Diamond} \mathfrak{a}_\Diamond = 0$, $d_{A_\Diamond}{*}\mathfrak{a}_\Diamond = 0$ *and* $\mathfrak{a}_\Diamond \wedge \mathfrak{a}_\Diamond = 0$.
  b) $|\mathfrak{a}_\Diamond| = |\hat{a}_\Diamond|$.
  c) $|\nabla_{A_\Diamond} \mathfrak{a}_\Diamond| = \lim_{n \in \Lambda} |\nabla_{A_n} \hat{a}_n|$ *with the convergence being in the $L^2_{1;loc}$ topology on the set where* $|\hat{a}_\Diamond| > 0$.

This lemma is proved momentarily. Assume it to be true for the time being.

Parts 1 of what follows directly use Lemma 7.2 to define $\mathcal{I}$, and Part 2 of what follows defines $\nu$ and verifies the first bullet of Proposition 7.1 Part 3 of what follows uses Lemma 7.2 to prove the second and third bullets of Proposition 7.1.

*Part 1*: Let SM denote the unit sphere bundle TM. Since $|\hat{a}_\Diamond| < c_0$, Item b) of the third bullet in Lemma 7.2 finds $|\mathfrak{a}_\Diamond| < c_0$. This understood, define a quadratic map from the $SM/\{\pm 1\}$ to $[0, c_0)$ by the rule $v \to |\mathfrak{a}_\Diamond(v)|^2$. This map has a unique maximum at each point of M–Z, this being one consequence of the fact that $\mathfrak{a}_\Diamond \wedge \mathfrak{a}_\Diamond = 0$. The corresponding line in $TM|_{M-Z}$ defines a real line subbundle of TM over M–Z. This real line bundle subbundle is denoted by $\mathcal{I}_\ddagger$. It is associated in a canonical way to the principal $\mathbb{Z}/2\mathbb{Z}$ that is defined by the points where $\mathcal{I}_{\ddagger}$ intersects the unit sphere bundle in TM. The bundle $\mathcal{I}$ is the dual to $\mathcal{I}_{\ddagger}$.

*Part 2*: Let $\mathfrak{U}$ denote a countable, locally finite open cover of M–Z by balls, and let B denote a given ball from this cover. The fact that $\mathfrak{a}_\Diamond \wedge \mathfrak{a}_\Diamond = 0$ implies that $\mathfrak{a}_\Diamond$ can be written on B as $\sigma_B \nu_B$ with $\sigma_B$ being a Sobolev class $L^2_2$ map from B to the unit sphere in $\mathfrak{su}(2)$ and with $\nu_B$ being an $\mathbb{R}$ valued 1-form on B that annihilates the orthogonal complement in $TM|_B$ of the line subbundle $\mathcal{I}_{\ddagger}|_B$.

The equations $[\sigma_B, d_{A_\Diamond}\mathfrak{a}_\Diamond] = 0$ and $[\sigma_B, d_{A_\Diamond}{*}\mathfrak{a}_\Diamond] = 0$ are equivalent to the assertions $\nabla_{A_\Diamond}\sigma_B \wedge \nu_B = 0$ and $\nabla_{A_\Diamond}\sigma_B \wedge {*}\nu_B = 0$. Since $|\nu_B| = |\hat{a}_\Diamond| \neq 0$ on B, these together assert that $\nabla_{A_\Diamond}\sigma_B = 0$. Thus, $\sigma_B$ is $A_\Diamond$-covariantly constant. Meanwhile, the two equations $\langle \sigma_B\, d_{A_\Diamond}\mathfrak{a}_\Diamond \rangle = 0$ and $\langle \sigma_B\, d_{A_\Diamond}{*}\mathfrak{a}_\Diamond \rangle = 0$ say that $\nu_B$ is a harmonic 1-form on B. As a parenthetical remark, the fact that $\sigma_B$ is $A_\Diamond$ covariantly constant implies that $F_{A_\Diamond}$ can be written on B as $\sigma_B w_B$ with $w_B$ denoting a closed square integrable 2-form on B.

The following turns out to be a crucial observation about this writing of $\mathfrak{a}_\Diamond$ The definition of $\nu_B$ has a sign ambiguity because there are automorphisms of B × SU(2) that pull $\sigma_B$ back as $-\sigma_B$. This sign ambiguity disappears if and only if $\mathcal{I}_{\ddagger}$ is the product line bundle.

What follows gives a second view of this sign ambiguity. Fix a length 1 element, $\sigma \in \mathfrak{su}(2)$. There exists a Sobolev class $L^2_2$ automorphism of the bundle B × SU(2) that writes $A_\Diamond$ as $\theta_0 + \sigma A_B$ and writes $\mathfrak{a}_\Diamond$ on B as $\sigma\nu_B$ with $A_B$ denoting a 1-form on B of Sobolev class $L^2_1$. Now let $\tau \in \mathfrak{su}(2)$ denote a given element with $\langle \sigma\tau \rangle = 0$ and length $|\tau| = \pi$. View $e^\tau$ as an automorphism The latter pulls back $A_\Diamond$ as $\theta_0 - \sigma A_B$ and it pulls back $\mathfrak{a}_\Diamond$ as $-\sigma\nu_B$. The sign ambiguity is due to the fact that $\theta_0 - \sigma A_B$ and $-\sigma\nu_B$ can be written respectively as $\theta_0 + \sigma A_B'$ and $\sigma\nu_B'$ with $A_B' = -A_B$ and $\nu_B' = -\nu_B$.

Let B and B′ denote intersecting balls from $\mathfrak{U}$. Then $\nu_B$ can be written on B ∩ B′ as $z_{BB'}\nu_{B'}$ with $z_{BB'}$ being either 1 or -1. The collection $\{z_{BB'}\}_{B,B'\in\mathfrak{U}}$ defines the transition functions for the line bundle $\mathcal{I}$. Meanwhile, the collection $\{\nu_B\}_{B\in\mathfrak{U}}$ defines a smooth, harmonic section over M–Z of $T^*(M–Z) \otimes \mathcal{I}$ that vanishes on Z, this being Proposition 7.1's section $\nu$.

*Part 3*: The fact that $|\nu| = |\hat{a}_\Diamond|$ follows from Item b) of Lemma 7.2's third bullet because $|\nu| = |\mathfrak{a}_\Diamond|$. The fact that $|\nabla\nu|$ is square integrable follows from the first bullet of Lemma 7.2 because $|\nabla\nu| = |\nabla_{A_\Diamond}\mathfrak{a}_\Diamond|$. To elaborate, fix $\rho > 0$ and introduce by way of notation $Z_\rho \subset M$ to denote the set where $|\hat{a}_\Diamond| < \rho$. Lemma 7.2 implies directly that

$$\int_{M-Z_\rho} |\nabla \nu|^2 = \lim_{n\in\Lambda} \int_{M-Z_\rho} |\nabla_{A_n} \hat{a}_n|^2 .$$

(7.1)

With (7.1) in mind, define a function on $(0, c_0^{-1})$ by the rule that assigns a given number ρ in this integral the value of the left hand integral in (7.1). The first bullet of Lemma 7.2 asserts in part that this function is bounded on $(0, c_0^{-1})$ and since it is a decreasing function of ρ, so the dominated convergence theorem says that it has a unique $\rho \to 0$ limit. This limit is the integral of $|\nabla\nu|^2$.

Proposition 7.1's third bullet follows using an identical argument after invoking what is said in the sixth bullet of Proposition 2.2 about the sequence whose n'th term is the integral over M of the function $G_p|\nabla_{A_n} \hat{a}_n|^2$.

***Proof of Lemma 7.2***: The proof has three steps.

Step 1: Fix $\rho > 0$. This step proves the following assertion:

*The sequence* $\{\int_{M-Z_\rho} |F_{A_n}|^2\}_{n\in\Lambda}$ *is bounded.*

(7.2)

To see why this is, suppose to the contrary that (7.2) is false. Given a point p in the $M–Z_\rho$, and $n \in \Lambda$, let $r_{\Diamond n,p}$ denote p's version of the number $r_\Diamond$ that is defined in (3.6) by the connection $A_n$. If (7.2) is false, then it must be that $\liminf_{p\in M-Z_\rho} (\liminf_{n\in\Lambda} r_{\Diamond n,p}) = 0$. This understood, there exists a point $p \in M–Z_\rho$ a subsequence $\Lambda' \subset \Lambda$ and a sequence $\{p_n\}_{n\in\Lambda'} \subset M–Z_\rho$ that converges to p and is such that $\lim_{n\in\Lambda'} r_{\Diamond n,p_n} = 0$. Granted this, then it must be true that $\lim_{n\in\Lambda'} r_{\Diamond n,p} = 0$ also.

The sequence $\{(A_n, \hat{a}_n)\}_{n\in\Lambda'}$ can be used as input for Proposition 2.2. Let $|\hat{a}_\Diamond'|$ denote the $\Lambda'$ version of what is denoted in Proposition 2.2 by $|\hat{a}_\Diamond|$. Use $|\hat{a}_\Diamond|$ to denote the original version that is supplied by Λ. Proposition 6.1 asserts that $|\hat{a}_\Diamond'|$ is also a Hölder continuous function on the complement of its zero locus, and locally Hölder continuous on its zero locus. It then follows that $|\hat{a}_\Diamond'| = |\hat{a}_\Diamond|$ because both are Hölder continuous and because they define the same $L^2$ function. Since $p \in M–Z$, it follows that $|\hat{a}_\Diamond'|(p) > 0$ and so $\Lambda'$ has a subsequence, this denoted by $\Lambda_p$, with $\lim_{n\in\Lambda_p} |\hat{a}_n|(p) = |\hat{a}_\Diamond'|(p)$. Invoke Lemma 6.2 using the point p and $\Lambda_p$ to conclude that $\{r_{\Diamond n,p}\}_{n\in\Lambda_p}$ is bounded away from zero. But this is nonsense because $\Lambda_p \subset \Lambda'$ and $\Lambda'$ was chosen so $\lim_{n\in\Lambda'} r_{\Diamond n,p} = 0$.

Step 2: Granted (7.2), then the constructions of Uhlenbeck in [U] can be employed to obtain a subsequence $\Lambda' \subset \Lambda$ and a sequence $\{h_n\}_{n\in\Lambda'}$ of automorphisms of $P|_{M-Z}$ such that the sequence $\{h_n{}^*A_n\}_{n\in\Lambda'}$ has the properties asserted by the second bullet of Lemma 7.3. The $L^2{}_{1;loc}$ limit connection can be taken to be $A_\lozenge$.

Proposition 2.2 implies in part that the sequence $\{\nabla_{A_n}\hat{a}_n\}_{n\in\Lambda'}$ is bounded, and this with the first bullet of (3.11) implies that the sequence $\{h_n{}^*\hat{a}_n\}_{n\in\Lambda'}$ is bounded in the $L^2{}_1$ topology on compact subsets of M–Z. As the sequence $\{q_{A_n}(\hat{a}_n)\}_{n\in\Lambda}$ is also bounded, (7.2) with the same sort of integration by parts argument that is used to prove the fifth bullet of Proposition 4.1 proves that $\{\nabla_{A_n}(\nabla_{A_n}\hat{a}_n)\}_{n\in\Lambda'}$ has bounded $L^2$ norm on compact subsets of M–Z. It follows that (3.11) can be used again to see that $\{h_n{}^*\hat{a}_n\}_{n\in\Lambda'}$ has bounded $L^2{}_2$ norm on compact subsets of M–Z. This implies in particular that $\Lambda'$ has a subsequence, this being $\Lambda_*$, such that $\{h_n{}^*\hat{a}_n\}_{n\in\Lambda_*}$ converges weakly in the $L^2{}_{2;loc}$ topology on M–Z. Use $\mathfrak{a}_\lozenge$ to denote the limit.

<u>Step 3</u>: What is said in (7.2) has the following implication: Suppose that U is a compact set in M–Z. Then there exists $r_U > 0$ such that $\inf_{p\in U}\{\inf_{n\in\Lambda'}\{r_{\lozenge n,p}\}\} > r_U$. This being the case, the argument in Part 2 of the proof of Lemma 6.3 can be used to prove that $|\mathfrak{a}_\lozenge| = |\hat{a}_\lozenge|$. Item c) of Lemma 7.2 with what is said in Step 2 about convergence implies that $d_{A_\lozenge}\hat{a}_\lozenge = 0$ and $d_{A_\lozenge}{*\hat{a}_\lozenge} = 0$ and $\hat{a}_\lozenge \wedge \hat{a}_\lozenge = 0$.

**b) The fourth bullet of Proposition 7.1**

The eight parts of this subsection prove the fourth bullet of Proposition 7.1. The three items are proved in more or less reverse order.

*Part 1*: Fix $\rho > 0$ and let $Z_\rho \subset M$ again denote the set where $|\nu| < \rho$. The following observation is invoked repeatedly in subsequent arguments in this section and in later sections.

$$\lim_{\rho\to 0}\int_{Z_\rho}|\nabla\nu|^2 = 0 . \tag{7.3}$$

By way of an explanation, the integrand is a measurable function with support on $Z_\rho$–Z. Meanwhile, Z is a closed set and so the function of ρ given by the volume of $Z_\rho$–Z has limit zero as ρ limits to zero.

With $\rho > 0$ given, the subsequent arguments refer to the function $\chi_\rho = \chi(2-\rho^{-1}|\nu|)$, this being a function on M that equals 0 on $Z_\rho$ and 1 on $M-Z_{2\rho}$. It is introduced to avoid certain delicate issues with regard to derivatives of ν near Z.

Define $h_{(\rho)}$ and $H_{(\rho)}$ to be the functions

$$h_{(\rho)} = \int_{\partial B_r} \chi_\rho \, |\nu|^2 \quad \textit{and} \quad H_{(\rho)} = \int_{B_r} \chi_\rho \, |\nabla\nu|^2 \tag{7.4}$$

It follows from the definitions that $h_{(\rho)}(r) \le h(r) \le h_{(\rho)}(r) + c_0 r^2 \rho^2$. Meanwhile, $H_{(\rho)}(r)$ is no greater than $H(r)$, and (7.3) implies that $\lim_{\rho\to 0}|H(r) - H_{(\rho)}(r)| = 0$ with the limit being uniform in the following sense: Given $\varepsilon > 0$, there exists $\rho_\varepsilon > 0$ such that if $r \in [0, c_0^{-1}]$ and $\rho < \rho_\varepsilon$, then $|H(r) - H_\rho(r)| < \varepsilon$.

*Part 2*: The lemma stated and then proved below asserts in part that h is non zero on $(0, c_0^{-1}]$.

**Lemma 7.3**: *There exists* $\kappa > 1$ *with the following significance: Fix* $p \in M$ *so as to define the function* h. *If* $r \in (0, \kappa^{-1})$, *then* $\int_{B_r} |\nu|^2 \le (1-\kappa r^2)\, r\, h(r)$.

***Proof of Lemma 7.3***: Use Q to denote the symmetric section of the tensor bundle $\otimes_2 T^*M$ given by $Q = \nu\otimes\nu - \frac{1}{2}\mathfrak{m}|\nu|^2$. The fact that ν is closed and coclosed implies that Q is divergence free. This means that it has vanishing $L^2$ inner product with the covariant derivatives of 1-forms. Note that Q is an $L^2_1$ section of $\otimes_2 T^*M$ that is smooth on M–Z and continuous and locally Hölder continuous across Z. The fact that Q is smooth on M–Z follows from the fact that ν is harmonic on M–Z and thus smooth. The fact that Q is locally Hölder continuous across Z follows from what is said by Proposition 6.1.

Fix a Gaussian coordinate system centered at p, and use the coordinate differentials as a basis for T*M and the coordinate vector fields as a basis for TM. Integrate the inner product between $\nabla d(|x|^2)$ and Q over the ball $B_r$. Use the fact that Q is divergence free and that $|\nu|$ is continuous and zero on Z with an integrate by parts to identify the latter integral with a boundary term. The resulting identity can be written as

$$r \int_{\partial B_r} (|\nu_r|^2 - \tfrac{1}{2}|\nu|^2) + \tfrac{1}{2}\int_{B_r} |\nu|^2 = \mathfrak{d}(r)\,, \tag{7.5}$$

where $\nu_r$ denotes the inner product between $d|x|$ and ν and where $|\mathfrak{d}(r)| \le c_0 r^2 \int_{B_r} |\nu|^2$.

Granted this bound on $|\mathfrak{d}|$, the identity in (7.5) implies what is asserted by Lemma 7.3.

*Part 3*: Lemma 7.3 has the following consequence: If $r \in (0, c_0^{-1})$ and if $h(r) = 0$, then $h(s) = 0$ for all $s \le r$. This understood, let $D \in [0, c_0^{-1}]$ denote the maximum value of r where h(r) is zero. The function N is defined on (D, $c_0^{-1}$). The rest of this Part 3 and all of Parts 4-6 prove Item c) from Proposition 7.1's fourth bullet for s > r with r > D.

The proof of Item c) for when s and r are greater than D will appeal to the lemma that follows directly and the upcoming Lemma 7.5.

**Lemma 7.4**: *There exists* $\kappa \geq 1$ *with the following significance. Let* $\Lambda'$ *denote the subsequence from Lemma 7.2. Fix* $\rho > 0$. *If* $n \in \Lambda$ *is sufficiently large, then*

$$\int_{Z_\rho} (|\nabla_{A_n} \hat{a}_n|^2 + r_n^2 |\hat{a}_n \wedge \hat{a}_n|^2) < \kappa \rho^{1/\kappa} .$$

***Proof of Lemma 7.4***: Fix for the moment $s \in (0, c_0^{-1})$ and introduce the function $\chi_{n,s}$ on M that is given by $\chi(s^{-1}|\hat{a}_n| - 1)$. This function is equal to 1 where $|\hat{a}_n| < \frac{5}{4}$ s and it is equal to zero where $|\hat{a}_n| \geq \frac{7}{4}$ s. Take $f = \chi_{n,s}$ in the $(A_n, \hat{a}_n)$ version of (2.2) to obtain an equation with the form

$$- \int_M \chi'_{n,s} |\hat{a}_n| |d|\hat{a}_n||^2 + \int_M \chi_{n,s} (|\nabla_{A_n} \hat{a}_n|^2 + r_n^2 |\hat{a}_n \wedge \hat{a}_n|^2) = \mathfrak{e}_n(s) \tag{7.6}$$

where $|\mathfrak{e}_n(s)| \leq c_0 (r_n^{-1} + s^2)$.

Fix $\rho < c_0^{-1}$ and let $N_\rho$ denote the smallest integer with the following property: If $n \in \Lambda$ and $n \geq N_\rho$, then

$$r_n > \rho^{-2} \quad \textit{and} \quad \tfrac{7}{8} |\nu| \leq |\hat{a}_n| \leq \tfrac{9}{8} |\nu| \textit{ where } |\nu| > \tfrac{1}{8} \rho . \tag{7.7}$$

If $n \geq N_\rho$ and $s \geq \rho$, then (7.6) implies the inequality

$$\int_{Z_s} (|\nabla_{A_n} \hat{a}_n|^2 + r_n^2 |\hat{a}_n \wedge \hat{a}_n|^2) \leq c_0 \int_{Z_{2s} - Z_s} (|\nabla_{A_n} \hat{a}_n|^2 + r_n^2 |\hat{a}_n \wedge \hat{a}_n|^2) + c_0 s^2). \tag{7.8}$$

The integral on the left in (7.8) is the difference between the respective integrals of its integrand over $Z_{2s}$ and $Z_s$. Use this fact to rewrite (7.8) so as to read

$$\int_{Z_s} (|\nabla_{A_n} \hat{a}_n|^2 + r_n^2 |\hat{a}_n \wedge \hat{a}_n|^2) \leq (1 - c_0^{-1}) \int_{Z_{2s}} (|\nabla_{A_n} \hat{a}_n|^2 + r_n^2 |\hat{a}_n \wedge \hat{a}_n|^2) + c_0 s^2 . \tag{7.9}$$

Given $k \in \{0, 1, \ldots, \}$ but less than $\frac{1}{\ln 2} |\ln \rho| - c_0$ introduce $x_k$ to denote the $s = 2^{-k} c_0^{-1}$ version of the integral on the left hand side of (7.9). With this notation understood, then (7.9) asserts that

$$x_k \le (1 - c_0^{-1})\, x_{k-1} + c_0\, 2^{-2k} \ .$$

(7.10)

Iterating this finds $x_k \le (1 - c_0^{-1})^k (x_0 + c_0)$. This implies that

$$\int_{Z_\rho} (|\nabla_{A_n} \hat{a}_n|^2 + r_n^2 |\hat{a}_n \wedge \hat{a}_n|^2) \le c_0\, \rho^{k/c_0} \ .$$

(7.11)

Since (7.11) holds for all $n \ge N_\rho$, it leads directly to the claim made by Lemma 7.4.

The next lemma is also needed for the fourth bullet of the proposition:

**Lemma 7.5**: *Let* $\Lambda'$ *denote the subsequence from Lemma 7.2. Given* $\rho > 0$, *then there exists* $\kappa > 1$ *with the following significance: If* $n \in \Lambda'$, *then* $r_n^2 \int_{M-Z_\rho} |\hat{a}_n \wedge \hat{a}_n|^2 \le \kappa\, r_n^{-2}$.

***Proof of Lemma 7.5***: The integral of $|F_{A_n}|^2$ on $M–Z_\rho$ enjoys an n-independent upper bound and so this is also the case for the integral of $r_n^4 |\hat{a}_n \wedge \hat{a}_n|^2$ on $M–Z_\rho$.

*Part 4*: Given $n \in \Lambda'$, use $h_n$, $\hslash_n$ and $N_n$ to denote the version of the functions $\hslash$ and $N$ that are defined in (5.2) and (5.3) using the point p, $r = r_n$ and $(A, \hat{a}) = (A_n, \hat{a}_n)$. Suppose $\rho > 0$ has been specified. Fix $r > 0$ and let $h_{n(\rho)}(r)$ and $H_{n(\rho)}(r)$ denote the respective integrals

$$h_{n(\rho)}(r) = \int_{\partial B_r} \chi_\rho |\hat{a}_n|^2 \quad \textit{and} \quad H_{n(\rho)}(r) = \int_{B_r} \chi_\rho (|\nabla_{A_n} \hat{a}_n|^2 + r_n^{\,2} |\hat{a}_n \wedge \hat{a}_n|^2) \,.$$

(7.12)

These are such that

- $h_{n(\rho)} \le h_n$ *for all* $n \in \Lambda$ *and* $h_{n(\rho)} \ge h_n - c_0 r^2 \rho^2$ *if* $n \in \Lambda$ *is sufficiently large*.
- $H_{n(\rho)} \le H_n$ *for all* $n \in \Lambda$ *and* $H_{n(\rho)} \ge H_n - c_0\, \rho^{1/c_0}$ *if* $n \in \Lambda$ *is sufficiently large*.

(7.13)

By way of an explanation, the lower bound for $h_{n(\rho)}$ follows from the fact that $\{|\hat{a}_n|\}_{n\in\Lambda}$ converges to $|\nu|$, and that for $H_{n(\rho)}$ follows from Lemma 7.4.

Fix $r \in (D, c_0^{-1})$. It follows from what is said by (7.13) and Lemma 7.5 by taking limits first as $n \in \Lambda'$ gets ever larger and then as $\rho$ limits to zero that

$$N(r) = \lim_{n\in\Lambda'} N_n(r) \ , \tag{7.14}$$

and that this limit is uniform as r various on compact subsets of $(D, c_0^{-1})$.

*Part 5*: This part of the subsection invokes some of what is said in Sections 5c and 5d to derive the $r > D$ version of Item c) of the fourth bullet from a lemma that is proved in Part 6. The starting point is the $r = r_n$ and $(A, \hat{a}) = (A_n, \hat{a}_n)$ version of (5.30). Keep in mind that what is denoted in (5.30) by $\mathfrak{t}$ is discussed subsequent to (5.11). It is enough to know that $|\mathfrak{t}| < c_0 r$. What is denoted by $\mathfrak{z}$ is defined subsquent to (5.28). A bound for $|\mathfrak{z}|$ is supplied by (5.32), but the latter is not sufficient for the purposes at hand. The upcoming Lemma 7.6 says what is needed about $\mathfrak{z}$'s absolute value. There is also a term in (5.30) that is denoted by $\mathfrak{R}$. This term enters via the second bullet of Lemma 5.2. A bound for $|\mathfrak{R}|$ is supplied by (5.31), but a stronger bound is needed and Lemma 7.6 provides one.

**Lemma 7.6**: *There exists* $\kappa > 1$*, and given* $\varepsilon \in (0, 1]$*, there exists* $\kappa_\varepsilon > \kappa$ *with the following significance: Let* $\Lambda'$ *denote the subsequence from Lemma 7.2. If* $n \in \Lambda'$ *is greater than* $\kappa_\varepsilon$*, then the absolute values of the* $r = r_n$ *and* $(A, \hat{a}) = (A_n, \hat{a}_n)$ *versions of* $\mathfrak{z}$ *and* $\mathfrak{R}$ *at points* $r \in [\varepsilon, \kappa^{-1})$ *are such that* $|\mathfrak{z}| < \varepsilon + \kappa r h_n$ *and* $|\mathfrak{R}| < \varepsilon + \kappa (h_n + r H_n)$.

This lemma is proved in Part 6.

Accept Lemma 7.6 for the moment to complete the proof of the $s > D$ version of Item c) from Proposition 7.1's fourth bullet. To do this, fix first $r_* > D$ and suppose that r is greater than $r_*$. It then follows from Lemma 7.3 that $h(r) > m_*^{-1}$ with $m_* > 1$ a number that depends on $r_*$ but not on r. Fix ε positive but less than the smaller of $c_0^{-1} m_*^{-1}$ and $r_*$. If n is large, then $h_n(r)$ will differ by at most $\varepsilon^2$ from the integral of $|\nu|^2$ over $\partial B_r$, and so $h_n(r)$ will be greater than $c_0^{-1} m_*^{-1}$.

With the preceding in mind, use Lemma 7.6 with (5.30) to conclude that

$$\frac{d}{dr} N_n \geq -c_0\Big(r + \frac{r}{h_n}\varepsilon\Big)(1 + N_n) \quad \textit{when } \max(D, \varepsilon) < r \leq c_0^{-1} \textit{ and } n \textit{ is large}. \tag{7.15}$$

Since $h_n(r) \geq c_0^{-1} m_*^{-1}$, the inequality (7.15) leads to the bound

$$\frac{d}{dr} N_n \geq -c_0(r + m_* \varepsilon)(1 + N_n) \quad \textit{when } \max(D, \varepsilon) < r \leq c_0^{-1} \textit{ and } n \textit{ is large}. \tag{7.16}$$

Suppose now that ε is greater than $m_*^{-2}$. If both s and r are greater than $\max(D, \varepsilon)$ and less than $c_0^{-1}$, and if $s > r$, then integrating (7.17) finds that

$$N_n(s) > e^{-c_0(s^2-r^2+\varepsilon^{1/2})}\, N_n(r) - c_0(s^2-r^2+\varepsilon^{1/2}) \quad \textit{when } n \textit{ is large}. \tag{7.17}$$

Proposition 7.1's fourth bullet for $r > D$ follows directly from (7.14) and (7.17).

*Part 6*: This part of the subsection contains the proof of Lemma 7.6.

***Proof of Lemma 7.6***: The proof has four steps.

Step 1: What is said about $\mathfrak{z}$ subsequent to (5.28) finds

$$|\mathfrak{z}| \le c_E\Big(|H - h| + \int_{B_r} |\hat{a}|^2 + \Big|\int_{B_r} \big\langle \hat{a} \wedge * q_A(\hat{a})\big\rangle\Big|\Big) , \tag{7.18}$$

with $c_E > 1$ denoting here and in what follows, a number that depends only on the value of $E$ in (3.2) and whose value can increase between successive appearances. The absolute value of $H - h$ is no greater than $c_E\, r^{-1}$, this being the content of (5.12). Meanwhile, the integral that involves $q_A(\hat{a})$ is bounded by $c_E\, r^{-1/2}$, this being a consequence of (5.16). If the integer $n \in \Lambda'$ is large, then $r_n^{-1/2}$ will be less than $c_E^{-1}\varepsilon^2$ and in particular, the left most and right most terms in the $r = r_n$ and $(A, \hat{a}) = (A_n, \hat{a}_n)$ version of (7.16) will make a contribution to $|\mathfrak{z}|$ that is less than $\frac{1}{2}\varepsilon^2$. Meanwhile, the integral of $|\hat{a}_n|^2$ over any $r > \varepsilon$ version of $B_r$ differs from that $|\nu|^2$ by at most $c_E^{-1}\varepsilon^2$ if n is large, and Lemma 7.5 asserts that the latter integral is no greater than $r(1 - c_0 r^2)$ times that of $|\nu|^2$ over the boundary of $B_r$. This integral of $|\nu|^2$ will differ from $r(1 - c_0 r^2)\, h_n(r)$ by at most $c_E^{-1}\varepsilon^2$ when $r \ge \varepsilon$ if n is large. It follows as a consequence that $|\mathfrak{z}| \le \varepsilon + c_E\, r\, h_n$ when $r \ge \varepsilon$ and n is large.

Step 2: Reintroduce the notation from Step 1 of the proof of Lemma 5.2. The idnentities in (5.6)-(5.8) lead to a bound on $|\mathfrak{R}|$ that can be written as

$$|\mathfrak{R}| \le c_0 \Big|\int_{B_r} \big\langle (\nabla_{A,i}\hat{a})_k\, q_A(\hat{a})_k \big\rangle \partial_i \Delta \Big| + c_E\Big(h + rH + \int_{B_r} |\hat{a}|^2 + r^{-1/2}\Big) . \tag{7.19}$$

What is denoted by $\Delta$ designates the function $\operatorname{dist}(\cdot\,, p)^2$ and $\partial_i\Delta$ designates the directional derivative of the function $\Delta$ along the i'th basis vector. Remarks in Step 1 imply the following: If $r \in (\varepsilon, c_0^{-1})$ and n is large, then the $r = r_n$ and $(A, \hat{a}) = (A_n, \hat{a}_n)$ version of the right hand side of (7.19) is no greater than

$$c_0 \Big| \int_{B_r} \Big\langle (\nabla_{A,i}\hat{a})_k \, q_A(\hat{a})_k \Big\rangle \partial_i \Delta \Big| + c_E(h + rH) + c_E^{-1}\varepsilon^2 \quad .$$

(7.20)

Fix $\rho > 0$ and use the identity $1 = (1 - \chi_\rho) + \chi_\rho$ to break the integral in (7.20) into two parts, the first having the factor $(1 - \chi_\rho)$ in the integrand and the second having the factor $\chi_\rho$. The next step supplies bounds for the absolute values of these two integrals.

<u>Step</u> <u>3</u>: Use Lemma 7.4 with Item d) in Proposition 2.2 to see that the norm of

$$\int_{B_r} (1 - \chi_\rho) \Big\langle (\nabla_{A_n,i}\hat{a}_n)_k \, q_{A_n}(\hat{a}_n)_k \Big\rangle \partial_i \Delta$$

(7.21)

is no greater than $c_0 \, r \, \rho^{1/c_0}$ when n is sufficiently large.

The $\chi_\rho$ part of the integral in (7.20) is

$$\int_{B_r} \chi_\rho \Big\langle (\nabla_{A_n,i}\hat{a}_n)_k \, q_{A_n}(\hat{a}_n)_k \Big\rangle \partial_i \Delta \, .$$

(7.22)

To bound the absolute value of (7.22), fix $\mu \in (0, \frac{1}{64})$ for the moment and introduce by way of notation $\sigma_\mu$ to denote the function on $B_r$ given by $\chi(2 - \mu^{-1}(1 - r^{-1}\text{dist}(\cdot, p)))$. This function equals 0 where the distance to p is greater than $(1 - \mu) r$ and it equals 1 where the distance to p is less than $1 - 2\mu$. Insert the identity $1 = (1 - \sigma_\mu) + \sigma_\mu$ into the integrand in (7.22) to write it as a sum of two terms.

The absolute value of the term with the factor $1 - \sigma_\mu$ is no greater than

$$c_0 \, r \, \Big( \int_{1-2\mu<\text{dist}(\cdot,p)<1} \chi_\rho |\nabla_{A_n}\hat{a}_n|^2 \Big)^{1/2} ,$$

(7.23)

this being a consequence of Item d) of the second bullet of Proposition 2.2. Meanwhile, the integral in (7.23) is no greater than $c_0 r^{1/2} \mu^{1/4}$ times the $L^4$ norm of $\nabla_{A_n}\hat{a}_n$ over $B_r$. The latter norm has a $\rho$ dependent but n independent upper bound, this being a consequence of the first bullet of (3.11) and the second bullet of Lemma 7.2. This understood, it then follows that the contribution to the absolute value of the integral in (7.22) from the term with $1 - \sigma_\mu$ is no greater than $c_0 r^{3/2} \mu \, K_\rho$ with $K_\rho$ being a number that is determined by $\rho$ but independent of $\mu$ and n.

The contribution to the integral of (7.22) from the term with $\sigma_\mu$ is the integral over $B_r$ of $\sigma_\mu \chi_\rho \langle (\nabla_{A_n,i}\hat{a}_n)_j \, q_{A_n}(\hat{a}_n)_j \rangle \partial_i \Delta$. A bound for the absolute value of this integral is obtained by writing

$$q_{A_n}(\hat{a}_n) = *d_{A_n}(*d_{A_n}\hat{a}_n) - d_{A_n}(*d_{A_n}*\hat{a}_n) \ ,$$

(7.24)

and then integrating by parts to make an integrand which has terms that are linear in the components of $d_{A_n}\hat{a}_n$ and $d_{A_n}*\hat{a}_n$. Use Item c) of Proposition 2.2 and the $L^2_2$ bound in the second bullet of Lemma 7.2 to see from the resulting integral that the absolute value of the $\sigma_\mu$ contribution to (7.22) is bounded by $K_\rho\, \mu^{-1} r_n^{-1}$ with $K_\rho$ denoting again a number that is determined by ρ but is independent of both $\mu$ and n.

Step 4: By way of a summary, what is said in Steps 1-3 bound the absolute value of the explicit integral in (7.20) by

$$c_0\, \rho^{1/c_0} + K_\rho(\mu^{1/4} + \mu^{-1} r_n^{-1})$$

(7.25)

With ε given, first choose $\rho < \varepsilon^{c_0}$ so that the left most term in (7.25) is less than $\frac{1}{3}\varepsilon^2$. With ρ so chosen, choose a postive value of $\mu$, but sufficiently small so as to make the middle term in (7.25) less than $\frac{1}{3}\varepsilon^2$ also. With ρ and $\mu$ fixed, the right most term in (7.25) is less than $\frac{1}{3}\varepsilon^2$ when n is sufficiently large. Granted these bounds, then (7.20) leads directly to the bound asserted by Lemma 7.6.

*Part 7*: This part of the subsection derives a differential equation for h that is the same as that given in Item b) of Proposition 7.1 where r > D. To set the notation, let $\mathfrak{m}$ denote the metric inner product on T*M. Use $\partial_r$ to denote the derivative along the radial geodesics from p in $B_r$. Fix $\rho \in (0, c_0^{-1})$ and differentiate to obtain the identity

$$\frac{d}{dr} \int_{\partial B_r} \chi_\rho\, |\nu|^2 = 2r^{-1}(1+\mathfrak{e}_1) \int_{\partial B_r} \chi_\rho\, |\nu|^2 + 2 \int_{\partial B_r} \chi_\rho \mathfrak{m}(\nu,\partial_r\nu)$$

(7.26)

with $\mathfrak{e}_1$ being bounded by $c_0 r^2$. Integrate by parts in the second integral using the fact that $\nu$ is harmonic to write it as

$$\int_{\partial B_r} \chi_\rho \mathfrak{m}(\nu,\partial_r\nu) = \int_{B_r} \chi_\rho(\,|\nabla\nu|^2 + \mathrm{Ric}(\nu,\nu)) + \mathfrak{e}_2 \quad .$$

(7.27)

where $\mathfrak{e}_2$ is a term with absolute value obeying

$$|\mathfrak{e}_2| \le c_0\, \rho^{-1} \int_{B_r \cap Z_{2\rho}} |\nu|\, |\nabla\nu|^2$$

(7.28)

Given that $|\nu| < 2\rho$ on $Z_{2\rho}$, this bound for $|\mathfrak{e}_2|$ implies directly that

$$|\mathfrak{e}_2| \le c_0 \int_{Z_{2\rho}} |\nabla\nu|^2 \ . \tag{7.29}$$

As noted in (7.3), the $\rho \to 0$ limit in (7.29) is zero.

Fix $r > 0$ and $\varepsilon > 0$; then integrate (7.26) and use (7.27) to see that

$$\int_{\partial B_{r+\varepsilon}} \chi_\rho\, |\nu|^2 - \int_{\partial B_r} \chi_\rho\, |\nu|^2 \ = \int_r^{r+\varepsilon} \Big(s^{-1} \int_{\partial B_s} (1 + \mathfrak{e}_1)\, \chi_\rho |\nu|^2 + \int_{B_s} \chi_\rho (|\nabla\nu|^2 + \mathrm{Ric}(\nu,\nu)) \ + \ \mathfrak{e}_2 \Big)\, ds \ . \tag{7.30}$$

Take $\rho$ to zero on both sides of (7.30). What is said subsequent to (7.4) and what is said in (7.3) imply that the $\rho \to 0$ limit of (7.30) is the identity

$$h(r+\varepsilon) - h(r) = \int_r^{r+\varepsilon} \Big(s^{-1} \int_{\partial B_s} (1 + \mathfrak{e}_1) |\nu|^2 + \int_{B_s} (|\nabla\nu|^2 + \mathrm{Ric}(\nu,\nu)) \Big)\, ds \ . \tag{7.31}$$

Divide both sides of (7.31) by $\varepsilon$ and take the $\varepsilon \to 0$ limit. As the $\varepsilon \to 0$ limit of the right hand side exists, the result of taking the $\varepsilon \to 0$ limit is an identity for the derivative of h that reads

$$\tfrac{d}{dr} h = 2r^{-1}(1 + \mathfrak{e}_*)\, h + H + \int_{B_r} \mathrm{Ric}(\nu,\nu) \tag{7.32}$$

with $\mathfrak{e}_*$ being a function of r obeying $|\mathfrak{e}_*| \le c_0 r^2$.

The integral of $\mathrm{Ric}(\nu,\nu)$ that appears in (7.32) is bounded by $c_0$ times the integral of $|\nu|^2$ over $B_r$. This understood, use Lemma 7.3 to write (7.32) schematically as

$$\tfrac{d}{dr} h = 2r^{-1}(1 + \mathfrak{e})\, h + H \ , \tag{7.33}$$

where $\mathfrak{e}$ is a function of r that obeys $|\mathfrak{e}| \le c_0 r^2$. The equation in (7.33) is the equation in Item b) of Proposition 7.1's fourth bullet at values of $r > \mathrm{D}$, this because the definition of $\mathrm{N}$ can be invoked to write H as $H = r^{-1}\mathrm{N}\, h$.

*Part 8*: This step proves the assertion in Item a) of Proposition 7.1's fourth bullet to the effect that any given $p \in M$ version of h is positive on the whole of $(0, c_0^{-1})$. This is done by assuming that there exists $p \in M$ where the corresponding version of $\mathrm{D}$ is positive and deriving nonsense.

To set the stage for the derivation, note that if $\mathrm{D}$ is positive, then Z contains an open set. If such is the case, then there exists, for any $\varepsilon > 0$, a point $p \in Z$ with $\mathrm{D}$ positive but less than $\varepsilon$. In particular, there exists a point $p \in M$ with $\mathrm{D}$ positive, with h defined on $(0, r_0)$ with $r_0 > c_0^{-1}$ and with $h(\frac{1}{2} r_0) > 0$. Let $p \in M$ denote such a point.

Write H in (7.33) on the $r > \mathrm{D}$ part of $(0, r_0)$ as $r^{-1}\mathrm{N}h$. Fix $s > r > \mathrm{D}$ with both in $(0, r_0)$ and integrate this rewriting of (7.33) to see that

$$h(s) = (1+\mathfrak{d})\,(\tfrac{s}{r})^2 \exp\Big(2\int_r^s \tfrac{1}{t}\mathrm{N}(t)dt\Big)\, h(r) \ , \tag{7.34}$$

where $\mathfrak{d}(s)$ is such that $|\mathfrak{d}| \le c_0 s^2$ and $|\frac{d}{ds}\mathfrak{d}| \le c_0 s$. As the $r \to \mathrm{D}$ limit of h(r) is zero and as $\mathrm{D}$ is positive, the identity in (7.34) implies that the function $\mathrm{N}$ can not be a bounded function on $(\mathrm{D}, \frac{1}{2} r_0]$. This understood, fix for the moment $m > 1$ and some $r_m \in (\mathrm{D}, \frac{1}{2} r_0]$ with $\mathrm{N}(r_m) > m$. The $r > \mathrm{D}$ version of Item c) of the fourth bullet in Proposition 7.1 says that $\mathrm{N}(\frac{1}{2} r_0) \ge c_0^{-1} m$. Since $m$ can be as large as desired, this constitutes the desired nonsense. The reason it is nonsense is as follows: By assumption, $h(\frac{1}{2} r_0) > 0$. Since H is bounded and since $\mathrm{N} = rH/h$, so $\mathrm{N}(\frac{1}{2} r_0)$ is finite, and thus less than $m$ if $m$ is sufficiently large.

**c) The $r \to 0$ limit of $\mathrm{N}$**

Given a point $p \in M$, let $\mathrm{N}_{(p)}$ denote p's version of the function $\mathrm{N}$. The upcoming Lemma 7.7 concerns the behavior of the $r \to 0$ limit of $\mathrm{N}_p(r)$ as p varies in Z.

**Lemma 7.7**: *There exists* $\kappa \ge 1$ *with the following significance: Any given* $p \in M$ *version of* $\mathrm{N}_{(p)}$ *extends to* $[0, \kappa^{-1}]$ *as a continuous function. Moreover,*

- *If* $\{q_k\}_{k=1,2,\ldots} \subset M$ *converges to a given point* $p \in M$ *then* $\lim_{k\to\infty} \mathrm{N}_{(q_k)}(0) \le \mathrm{N}_{(p)}(0)$.
- *If* $p \in M{-}Z$, *then* $\mathrm{N}_{(p)}(0) = 0$*; and if* $p \in Z$, *then* $\mathrm{N}_{(p)}(0) > \kappa^{-1}$.

The remainder of this subsection contains the proof of this lemma.

***Proof of Lemma 7.7***: The proof has four parts.

*Part 1*: To see that a given $p \in M$ version of $\mathrm{N}_{(p)}$ is continuous on its original domain of definition, keep in mind that p's version of the function h is differentiable, this being an assertion of Proposition 7.1. This understood, it follows that $\mathrm{N}_{(p)}$ is continuous on $(0, c_0^{-1})$ if and only if p's version of the function H is continuous. The fact that H is continuous follows from what is said in Part 1 of Section 7b about H and (7.4)'s function $H_{(\rho)}$. To say more, fix $\varepsilon > 0$ and use what is said at the end of Part 1 of Section 7b to find $\rho_\varepsilon > 0$ such that $|H(\cdot) - H_{(\rho)}(\cdot)| < \varepsilon$ when $\rho < \rho_\varepsilon$. Fix such a value for $\rho$. Since $\nu$ is smooth on M–Z, the function $H_{(\rho)}$ is smooth. This understood fix $r \in (0, c_0^{-1}]$ and then fix $\Delta > 0$ so that $|H_{(\rho)}(r+\Delta) - H_{(\rho)}(r)| < \varepsilon$. It then follows that $|H(r+\Delta) - H(r)| < 3\varepsilon$ and so H is continuous.

The assertion that $N_{(p)}$ extends continuously to $[0, c_0^{-1}]$ follows from Item c) of the fourth bullet of Proposition 7.1 given the fact that $N_{(p)}$ is positive.

*Part 2*: The proof of the bulleted assertions requires a weak version of what is said by the lemma's second bullet, this being the following:

$$\textit{If } p \in M\text{–}Z, \textit{ then } N_{(p)}(0) = 0\textit{; and if } p \in Z, \textit{ then } N_{(p)}(0) > 0 \,. \tag{7.35}$$

The proof of (7.35) has two steps.

<u>Step 1</u>: If $p \in$ M–Z, then $|\nu|^2$ is greater than zero at p and smooth in a neighborhood of p. This being the case, p's version of the function h must have the form $4\pi r^2 |\nu|^2(p) + \mathfrak{e}$ with $\mathfrak{e}$ being a function of r with absolute value bounded by $c_0 r^3$. With this in mind, write $N_{(p)}$ near $r = 0$ as $N_{(p)}(0) + \mathfrak{o}$ with $\lim_{r\to 0} \mathfrak{o}(r) = 0$. Suppose for the sake of argument that $\varepsilon > 0$ and that $N_{(p)}(0)$ is greater than $\varepsilon$ so as to derive nonsense. To do this, fix $s > 0$ and use (7.34) to see that $h(r) \le c(s)\, r^{2+\varepsilon/2}$ when r is small.

<u>Step 2</u>: Suppose that $p \in Z$ and suppose for the sake of argument that $N_{(p)}(0) = 0$ so as to derive nonsense. To start, introduce by way of notation $c$ to denote the version of the constant $\kappa$ that is assigned to p by Proposition 6.1. It follows from Proposition 6.1 that p's version of the function h is such that $h(r) \le c\, r^{2+2/c}$ when $r \le c^{-1}$. Hold onto this bound for a moment.

Use the fact that $N_{(p)}$ is continuous on $[0, c_0^{-1}]$ and 0 at $r = 0$ to draw the following conclusion: Given $\varepsilon > 0$, there exists $r_\varepsilon > 0$ such that $N_{(p)}(r) < \varepsilon$ when $r < 2r_\varepsilon$. Granted this bound, use the $s = r_\varepsilon$ version of (7.34) to see that $h(r) \ge c_\varepsilon^{-1}\, r^{2+2\varepsilon}$ when $r < r_\varepsilon$ with $c_\varepsilon > 1$ being a number that depends on $\varepsilon$ but not on r.

The lower bound $h(r) \ge c_\varepsilon^{-1}\, r^{2+2\varepsilon}$ runs afoul of the upper bound $c\, r^{2+2/c}$ when $\varepsilon > 1/c$. The fact that the lower bound holds for all nonzero $\varepsilon$ constitutes the desired nonsense.

*Part 3*. The three steps that follow prove Lemma 7.7's assertion about the behavior of the function $p \to N_{(p)}(0)$. The notation used below has $h_{(p)}$ and $H_{(p)}$ denoting a given $p \in M$ version of the functions h and H.

<u>Step</u> <u>1</u>: Granted that $N_{(p)}(0) = 0$ when $p \in M{-}Z$, then what is said by Lemma 7.7 about the semi-continuity $N_{(\cdot)}(0)$ at a point $p \in M{-}Z$ follows directly from the fact that $M{-}Z$ is an open set.

<u>Step</u> <u>2</u>: Fix $p \in Z$ and $D > 0$. Let $q \in Z$ denote a point with distance $D$ or less from p. If $r \in (c_0 D, c_0^{-1})$, then the sphere of radius r centered at q is contained in the ball of radius $r + D$ centered at p and it contains the ball of radius $r - D$ centered at p. This understood, an application of the fundamental theorem of calculus finds

$$|h_{(q)}(r) - h_{(p)}(r)| \le c_0 \, (H_{(p)}(r + D))^{1/2} \, r \, D^{1/2} \, . \tag{7.36}$$

Use the third bullet of Proposition 7.1 to bound $H_{(p)}(r + D)$ by $c_0 \, (r + D)^{1/2}$ and use this bound with (7.35) to conclude that $|h_{(q)}(r) - h_{(p)}(r)| \le c_0 \, r^{3/2} D^{1/2}$ when $r \in (10D, c_0^{-1})$.

<u>Step</u> <u>3</u>: The fact that $N_{(p)}$ is continuous on $[0, c_0^{-1})$ has the following consequence: Fix $\varepsilon > 0$ and there exists $r_\varepsilon \in (0, \varepsilon)$ such that $N_{(p)}(0) - \varepsilon \le N_{(p)}(r) \le N_{(p)}(0) + \varepsilon$ if $r \in (0, 2r_\varepsilon]$. This being the case, then (7.34) implies that $h_{(p)}(r)$ for $r \in [0, 2r_\varepsilon]$ can be written as

$$h_{(p)}(r) = \mathcal{c}_\varepsilon \, r^{2(1 + N_{(p)}(0) + \mathfrak{e})} \, , \tag{7.37}$$

where $\mathcal{c}_\varepsilon > 0$ and where the function $\mathfrak{e}$ is such that $|\mathfrak{e}| \le c_0 \varepsilon$. If $D < c_0^{-1} r_\varepsilon$, then this equation for $h_{(p)}$ with (7.36) implies that the function $h_{(q)}(r)$ for $r \in (10D, 2r_\varepsilon]$ can be written as

$$h_{(q)}(r) = \mathcal{c}_\varepsilon \, r^{2(1 + N_{(p)}(0) + \mathfrak{e})} + \mathfrak{r} \tag{7.38}$$

where the function $\mathfrak{r}$ is such that $|\mathfrak{r}| \le c_0 \, r^{3/2} D^{1/2}$. If it is the case that $D < c_0^{-1} \, \mathcal{c}_\varepsilon^{\,2} \, r_\varepsilon^{\,c_0}$, then (7.38) holds for $r \in [r_\varepsilon, 2r_\varepsilon]$ with $\mathfrak{r} = 0$ but with a different version of $\mathfrak{e}$ such that $|\mathfrak{e}| \le c_0 \varepsilon$.

<u>Step</u> <u>4</u>: Suppose that $D$ obeys the bound $D \le c_0^{-1} \, \mathcal{c}_\varepsilon^{\,2} \, r_\varepsilon^{\,c_0}$. It follows from what is said at the end of Step 2 that if $s > r$ with both from $[r_\varepsilon, 2r_\varepsilon]$, then

$$h_{(q)}(r)/h_{(q)}(s) \ge \left(\tfrac{r}{s}\right)^{2(1 + N_{(p)}(0)) + c_0 \varepsilon} \, . \tag{7.39}$$

Meanwhile, the point q has its corresponding version of Proposition 7.1's fourth bullet. Item c) of the q version implies that $N_{(q)}(t) > N_{(q)}(0) - c_0 \, r_\varepsilon^{\,2}$ for $t \in [0, 2r_\varepsilon]$. This bound with q's version of (7.34) implies that

$$h_{(q)}(r)/h_{(q)}(s) \leq (\tfrac{r}{s})^{2(1+N_{(q)}(0)) - c_0\varepsilon} .$$

(7.40)

These upper and lower bounds are not compatible unless $N_{(q)}(0) \leq N_{(p)}(0) + c_0\varepsilon$. Since $\varepsilon$ can be made as small as desired by taking $D$ sufficiently small, this last inequality proves Lemma 7.7's claim about the semi-continuity of $N_{(\cdot)}(0)$ at p.

*Part 4*: This part proves that there is a positive lower bound for the value of the function $N_{(\cdot)}(0)$ on Z. The proof starts with the identity $|\nu| = |\hat{a}_\delta|$ from Proposition 7.1. Granted this identity, then Proposition 6.1 supplies $\kappa > 1$ such that if $p \in Z$ and q is a point in M that in a small radius ball centered at p, then $|\nu|(q) \leq \mathrm{dist}(p,q)^{1/\kappa}$. What follows is a consequence: If r is positive but small, then $h_{(p)}(r) \leq c_0\, x^2\, r^{2+2/\kappa}$. Meanwhile, Items b) and c) of the fourth bullet of Proposition 7.1 that $h_{(p)}(r) \geq x'\, r^{2+2N_{(p)}(0)}$ if r is small with $x'$ being independent of r if r is small. These two bounds are not compatible if $N_{(p)}(0) < \kappa$.

**d) Hölder continuity of $|\nu|$**

The three parts of this subsection prove that $|\nu|$ is Hölder continuous on M.

*Part 1*: Since $\nu$ is harmonic on M–Z it obeys the equation

$$\nabla^\dagger\nabla\nu + \mathrm{Ric}(\nu) = 0$$

(7.41)

with Ric denoting here the endomorphism of T*X that is defined using the metric and Ricci curvature. Supposing that $\delta \in (0,1]$, let $\beta_\delta$ denote the function $\chi(2 - \delta^{-1}\mathrm{dist}(\cdot, Z))$. This function is equal to 1 where the distance to Z is greater than $2\delta$ and it is equal to zero where the distance to Z is less than $\delta$. The equation in (7.41) leads to the bound

$$\|\beta_\delta \nabla^\dagger\nabla\nu\|_2^{\,2} \leq c_0 ;$$

(7.42)

and an integration by parts leads from (7.42) to an $L^2_2$ bound for $\nu$ of the form

$$\|\nabla\nabla(\beta_\delta\nu)\|_2^{\,2} \leq c_0\,\delta^{-4}\,\|\nabla\nu\|_2^{\,2}$$

(7.43)

This implies in particular that the $L^2_1$ norm of $\beta_\delta|\nabla\nu|$ is bounded by $c_0\delta^{-2}$. The latter bound leads turn a $c_0\delta^{-2}$ bound for the $L^6$ norm of $\beta_\delta|\nabla\nu|$ and hence for the $L^6$ norm of $\beta_\delta d|\upsilon|$. Granted this last $L^6$ bound, then a standard dimension 3 Sobolev inequality leads to the exponent $\frac{1}{4}$ Hölder norm bound

$$||\nu(p)|-|\nu(q)|| \le c_0\,\delta^{-2}\,\mathrm{dist}(p,q)^{1/4} \tag{7.44}$$

if both p and q have distance $2\delta$ or more from Z.

*Part 2*: As explained in Part 3, if p has distance less than $\delta$ from Z, then

$$|\nu(p)| \le c_0\,\delta^{\upsilon} \tag{7.45}$$

with $\upsilon$ positive and independent of both $\delta$ and p. Of course, the analogous inequality holds for $|\nu(q)|$ if q has distance less than $\delta$ from Z. To exploit (7.45), label p and q so that $\mathrm{dist}(p,Z) \le \mathrm{dist}(q,Z)$. If $\mathrm{dist}(p,q)^{1/16} \le \mathrm{dist}(p,Z)$, then it follows from (7.44), that

$$||\nu(p)|-|\nu(q)|| \le c_0|\mathrm{dist}(p,q)|^{1/8}. \tag{7.46}$$

Suppose on the other hand that $\mathrm{dist}(p,q)^{1/16} \ge \mathrm{dist}(p,Z)$. If $\mathrm{dist}(q,Z) > 10\,\mathrm{dist}(p,Z)$, then $\mathrm{dist}(q,Z) \le \frac{10}{9}\,\mathrm{dist}(p,q)$ and since

$$||\nu(p)|-|\nu(q)|| \le |\nu(p)|+|\nu(q)|, \tag{7.47}$$

it follows from (7.45) and its $|\nu(q)|$ analog that

$$||\nu(p)|-|\nu(q)|| \le c_0\,\mathrm{dist}(p,q)^{\upsilon}. \tag{7.48}$$

Meanwhile, if $\mathrm{dist}(p,q)^{1/16} \ge \mathrm{dist}(p,Z)$ and if $\mathrm{dist}(q,Z) \le 10\,\mathrm{dist}(p,Z)$, then it follows from (7.45), its $|\nu|(q)$ analog and from (7.47) that

$$||\nu(p)|-|\nu(q)|| \le c_0\,\mathrm{dist}(p,q)^{\upsilon/16}\ . \tag{7.49}$$

The claimed Hölder continuity of $|\nu|$ follows from (7.46), (7.48) and (7.49) with a proof of (7.45).

*Part 3*: This part proves what is asserted in (7.45). Supposing that p is a given point in M, let $\delta$ now denote the distance from p to Z. Fix a point $p' \in Z$ with distance at most $2\delta$ from p and let $\sigma_{p',\delta}$ denote the function $\chi(4 - \delta^{-1}\mathrm{dist}(\cdot,p'))$ This function is equal to 1 where the distance to $p'$ is greater than $4\delta$ and it is equal to 0 where the distance to $p'$ is less than $3\delta$. Keeping in mind that $|\nu| = |\hat{a}_{\Diamond}|$, subtract that $f = \sigma_{p',\delta}G_p$ version of the fourth bullet of Proposition 2.2 from the equation at the end of the sixth bullet of Proposition 2.2 to obtain an inequality that reads

$$|\nu|^2(p) \le c_0\, \delta^{-3} \int_{\mathrm{dist}(\cdot,p') \le 4\delta} |\nu|^2 \, . \tag{7.50}$$

This inequality uses the bound $|G_p| \le c_0\, \mathrm{dist}(\cdot\,, p)^{-1}$. It also uses the fact that $(1 - \sigma_{p',\delta})\, G_p$ is a non-negative function on M. Let $\kappa$ denote the constant from the second bullet of Lemma 7.7. It follows from Items b) and c) of Proposition 7.1 and from the second bullet of Lemma 7.7 that the right hand sice of (7.50) is bounded by $c_0\delta^{2\kappa}$ if $\delta \le c_0^{-1}$.